\documentclass[a4paper,11pt]{article}
\usepackage[latin1]{inputenc}
\usepackage[T1]{fontenc}
\usepackage[francais]{babel}
\usepackage{amsmath,amsbsy,amsfonts,amssymb}
\usepackage{fullpage}
\usepackage{enumerate}
\usepackage{graphicx}
\usepackage{titlesec}
\usepackage{textcomp}
\usepackage{epic}
\usepackage{eepic}
\usepackage{bbold}
\usepackage{soul}
\usepackage{pst-eucl}
\usepackage[all, ps]{xy}
\newtheorem{cor}{Corollaire}[subsection]
\newtheorem{lem}{Lemme}[subsection]
\newtheorem{duf}{Definition}[subsection]
\newtheorem{theom}{Théorème}
\newtheorem{theo}{Théorème}[subsection]

\newtheorem{prop}{Proposition}[subsection]
\title{La conjecture locale de Gross-Prasad pour les représentations tempérées des groupes unitaires}
\author{Raphaël Beuzart-Plessis}
\begin{document}

\maketitle

\bigskip

\begin{abstract}
Let $E/F$ be a quadratic extension of non-archimedean local fields of characteristic $0$ and let $G=U(n)$, $H=U(m)$ be unitary groups of hermitian spaces $V$ and $W$. Assume that $V$ contains $W$ and that the orthogonal complement of $W$ is a quasisplit hermitian space (i.e. whose unitary group is quasisplit over $F$). Let $\pi$ and $\sigma$ be smooth irreducible representations of $G(F)$ and $H(F)$ respectively. Then Gan, Gross and Prasad have defined a multiplicity $m(\pi,\sigma)$ which for $m=n-1$ is just the dimension of $Hom_{H(F)}(\pi,\sigma)$. For $\pi$ and $\sigma$ tempered, we state and prove an integral formula for this multiplicity. As a consequence, assuming some expected properties of tempered $L$-packets, we prove a part of the local Gross-Prasad conjecture for tempered representations of unitary groups. This article represents a straight continuation of recent papers of Waldspurger dealing with special orthogonal groups.
\end{abstract}

\textbf{Introduction}

\vspace{4mm}

 Soit $F$ un corps local non archimédien de caractéristique nulle et $E$ une extension quadratique de $F$. Soit $V$ un espace vectoriel de dimension finie sur $E$ muni d'une forme hermitienne symétrique non dégénérée $h$. Un tel couple sera appelé espace hermitien. Supposons donnée une décomposition orthogonale $V=W\oplus D\oplus Z$ où $D$ est une droite et $Z$ est une somme de plans hyperboliques i.e. admet une base $(v_i)_{i=\pm 1,\ldots,\pm r}$ telle que $h(v_i,v_j)=\delta_{i,-j}$ pour $i,j=\pm 1,\ldots,\pm r$. Notons $G$ et $H$ les groupes unitaires de $V$ et $W$ respectivement. Considérons le sous-groupe parabolique $P$ de $G$ des éléments qui conservent le drapeau de sous-espaces totalement isotropes
$$Ev_r\subset Ev_r\oplus Ev_{r-1}\subset\ldots\subset Ev_r\oplus\ldots\oplus Ev_1$$
Soit $M$ la composante de Levi de $P$ des éléments qui conservent les droites $Ev_i$ pour $i=\pm 1,\ldots,\pm r$. On a un isomorphisme naturel $M\simeq (R_{E/F}GL_1)^r\times G_0$ où $G_0$ est le groupe unitaire de $V_0=D\oplus W$ et $R_{E/F}$ désigne la restriction des scalaires à la Weil. Notons $U$ le radical unipotent de $P$. Soit $\xi$ un caractère de $U(F)$ invariant par conjugaison par $H(F)$ et qui est générique pour cette propriété (la définition exacte de $\xi$ est donnée dans la section 4). Pour $(\pi,E_\pi)$ et $(\sigma,E_\sigma)$ des représentations admissibles irréductibles de $G(F)$ et $H(F)$ respectivement on définit $Hom_{H,\xi}(\pi,\sigma)$ comme l'espace des applications linéaires $l:E_\pi\to E_\sigma$ qui vérifient

$$l(\pi(hu)e)=\xi(u)\sigma(h)l(e)$$

\noindent pour tous $e\in E_\pi$, $h\in H(F)$, $u\in U(F)$. On définit la multiplicité $m(\pi,\sigma)$ comme la dimension de cet espace, elle ne dépend que des représentations $\pi$ et $\sigma$ et pas des divers choix effectués. D'après [AGRS], on a $m(\pi,\sigma)\leqslant 1$. La conjecture de Gross-Prasad telle qu'énoncée dans [GGP] décrit exactement quand la multiplicité $m(\pi,\sigma)$ vaut 1 en terme des paramètres de Langlands de $\pi$ et $\sigma$. On démontre ici sous certaines hypothèses concernant les $L$-paquets des groupes unitaires une version affaiblie de la conjecture pour les représentations tempérées. Plus précisément, supposons $G$ et $H$ quasi-déployés et affectons d'un indice $i$ tout les objets définis jusque là: $G_i,H_i,V_i,W_i\ldots$. Il existe à isomorphisme près un unique espace hermitien $V_a$ de même dimension que $V_i$ mais qui ne lui est pas isomorphe. On définit de même l'espace hermitien $W_a$. On a encore une décomposition orthogonale $V_a=Z\oplus D\oplus W_a$. Admettons l'existence de $L$-paquets pour les groupes $G_i(F),H_i(F),G_a(F)$ et $H_a(F)$ ainsi que certaines propriétés de ces $L$-paquets (les hypothèses faites sont énoncées précisément en 18.1). Les paramètres de Langlands pour les groupes $G_i$ et $G_a$ sont les mêmes. Soit $\varphi$ un tel paramètre que l'on suppose tempéré et $\Pi_i$ resp. $\Pi_a$ le $L$-paquet de représentations de $G_i(F)$ resp. de $G_a(F)$ correspondant (on a éventuellement $\Pi_a=\emptyset$). Soit de même $\psi$ un paramètre de Langlands tempérés pour les groupes $H_i$ et $H_a$ et $\Sigma_i$ resp. $\Sigma_a$ le $L$-paquet de représentations de $H_i(F)$ resp. de $H_a(F)$ correspondant. Sous les hypothèses admises sur les $L$-paquets on montre le théorème suivant

\begin{theom}
Il existe un unique couple $(\pi,\sigma)\in (\Pi_i\times \Sigma_i)\bigsqcup (\Pi_a\times \Sigma_a)$ tel que $m(\pi,\sigma)=1$
\end{theom}

Le théorème 1 découle d'une formule intégrale pour la multiplicité $m(\pi,\sigma)$. L'idée de l'existence d'une telle formule intégrale ainsi que sa preuve s'inspirent directement des articles [W1] et [W2] de Waldspurger où est traité le cas des groupes spéciaux orthogonaux. Décrivons cette formule. On définit un ensemble $\underline{\mathcal{T}}$ de tores (en général non maximaux) de $H$. Ce sont les tores $T$ obtenus de la façon suivante: on a une décomposition orthogonale $W=W'\oplus W''$ telle que les groupes unitaires de $W''$ et $W''\oplus D\oplus Z$ soient quasi-déployés et $T$ est un tore maximal anisotrope du groupe unitaire de $W'$. Soit $\mathcal{T}$ un ensemble de représentants de $\underline{\mathcal{T}}$ pour l'action par conjugaison de $H(F)$. Soit $T\in \mathcal{T}$. On définit sur $T(F)$ trois fonctions déterminant $\Delta,D^H$ et $D^G$ qui importent peu et pour lesquels on renvoie au corps de l'article pour des définitions précises. Plus important, si $\pi$ et $\sigma$ sont des représentations admissibles irréductibles de $G(F)$ et $H(F)$ respectivement on définit deux fonctions $c_\pi$ et $c_\sigma$ sur $T(F)$. Ces deux fonctions sont définies à partir des caractères $\theta_\pi$ et $\theta_\sigma$ de $\pi$ et $\sigma$. Soit $t\in T(F)$ et soit $G_t$ le centralisateur connexe de $t$ dans $G$. D'après Harish-Chandra, au voisinage de $t$ le caractère $\theta_\pi$ est combinaison linéaire des transformées de Fourier des intégrales orbitales nilpotentes de $\mathfrak{g}_t(F)$ (où $\mathfrak{g}_t$ désigne l'algèbre de Lie de $G_t$). Plus précisément il existe des nombres complexes $c_{\pi,\mathcal{O}}(t)$ indexés par les orbites nilpotentes dans $\mathfrak{g}_t(F)$ et un voisinage $\omega$ de $0$ dans $\mathfrak{g}_t(F)$ de sorte que pour toute fonction $\varphi\in C_c^\infty(\omega)$ on ait

$$\displaystyle\int_{\mathfrak{g}_t(F)} \theta_\pi(texp(X))\hat{\varphi}(X) dX=\sum_\mathcal{O} c_{\pi,\mathcal{O}}(t) J_\mathcal{O}(\varphi)$$

\noindent où $\hat{\varphi}$ désigne la transformée de Fourier de $\varphi$ et $J_\mathcal{O}$ est l'intégrale orbitale sur $\mathcal{O}$. La somme porte sur l'ensemble des orbites nilpotentes de $\mathfrak{g}_t(F)$. Bien sûr les mesures et la transformée de Fourier doivent être définies précisément, on renvoie pour cela au corps de l'article. Soit $Nil_{reg}(\mathfrak{g}_t(F))$ l'ensemble des orbites nilpotentes régulières. On pose alors

$$\displaystyle c_\pi(t)=\frac{\sum_{\mathcal{O}\in Nil_{reg}(\mathfrak{g}_t(F))} c_{\pi,\mathcal{O}}(t)}{|Nil_{reg}(\mathfrak{g}_t(F))|}$$

\noindent On définit de même la fonction $c_\sigma$. Posons

$$m_{geom}(\pi,\sigma)=\displaystyle\sum_{T\in\mathcal{T}} |W(H,T)|^{-1} \lim\limits_{s\to 0^+} \int_{T(F)} c_{\sigma^\vee}(t)c_\pi(t)D^H(t)^{1/2}D^G(t)^{1/2}\Delta(t)^{s-1/2} dt$$

\noindent où $W(H,T)=Norm_{H(F)}(T)/Cent_{H(F)}(T)$ et $\sigma^\vee$ est la représentation contragrédiente de $\sigma$. L'expression ci-dessus est bien définie (les intégrales sont convergentes pour $Re(s)>0$ et admettent une limite lorsque $s$ tend vers $0$). Le résultat principal de cet article est le 17.1.1 dont voici un énoncé:

\begin{theom}
Si $\pi$ et $\sigma$ sont tempérées, on a
$$m(\pi,\sigma)=m_{geom}(\pi,\sigma)$$
\end{theom}

Décrivons brièvement comment l'on déduit le théorème 1 du théorème 2, l'idée étant issue de l'article [W1]. Sommons les multiplicités $m(\pi,\sigma)$ pour $(\pi,\sigma)$ parcourant $(\Pi_i\times \Sigma_i)\bigsqcup(\Pi_a\times \Sigma_a)$. Le théorème 2 permet d'exprimer la somme sur les $(\pi,\sigma)\in \Pi_i\times \Sigma_i$ comme une somme d'intégrales sur un ensemble de tores $\mathcal{T}_i$. De même le théorème 2 exprime la somme sur les $(\pi,\sigma)\in \Pi_a\times \Sigma_a$ comme une somme d'intégrales sur un ensemble de tores $\mathcal{T}_a$. On peut regrouper ces tores suivant leurs classes de conjugaison stable. On a une notion naturel de transfert entre les tores de $\mathcal{T}_a$ et les tores de $\mathcal{T}_i$.  Plus précisément on peut transférer une classe de conjugaison stable de $\mathcal{T}_a$ en une classe de conjugaison stable de $\mathcal{T}_i$. Le transfert est injectif et presque surjectif: toutes les classes sont atteintes sauf celle du tore $\{1\}\in\mathcal{T}_i$. On peut aussi transférer les différentes fonctions qui apparaissent dans les intégrales et ce transfert fait apparaître un signe. Il se trouve que ce signe vaut $-1$. Ainsi dans la somme totale que l'on considère toutes les intégrales disparaissent sauf une: celle correspondant au tore $\{1\}$. Un resultat de Rodier ([Ro]) permet alors de montrer que ce terme vaut 1. La démonstration du théorème 2 est plus longue et occupe la majeure partie de l'article (sections 5 à 17). Encore une fois la démonstration suit de très près celles de [W1] et [W2]. Une fonction $f\in C_c^\infty(G(F))$ est dite très cuspidale si pour tout sous-groupe parabolique propre $P=MU$ de $G$ on a

$$\forall m\in M(F), \; \displaystyle\int_{U(F)} f(mu) du=0$$

\noindent On définit une suite croissante exhaustive $(\Omega_N)_{N\geqslant 1}$ de compacts-ouverts de $H(F)U(F)\backslash G(F)$ et on note, pour $N\geqslant 1$, $\kappa_N$ la fonction caractéristique de $\Omega_N$. Posons

$$\displaystyle J_N(\theta_{\sigma^\vee},f)=\int_{H(F)U(F)\backslash G(F)} \kappa_N(g)\int_{H(F)}\theta_{\sigma^\vee}(h)\int_{U(F)} f(g^{-1}hug) \xi(u) dudhdg$$

C'est l'analogue de l'expression $I_N(\theta_{\sigma^\vee},f)$ de [W1] et [W2]. Cette expression admet une limite lorsque $N$ tend vers l'infini et comme dans la formule des traces locale d'Arthur, il y a deux façons de calculer cette limite: l'une qualifiée de géométrique, l'autre de spectrale. Le développement géométrique de la limite fait l'objet des sections 5 à 10, le développement spectral fait l'objet de la section 16. C'est de l'égalité entre ces deux expressions que l'on déduira le théorème 2. Décrivons à présent brièvement le contenu des différentes parties. \\

 La première partie rassemble un certain nombre de définitions, notations et résultats qui seront utilisés par la suite. On fixe notamment les choix de mesures qui apparaissent les plus naturels pour le développement géométrique et pour le développement spectral. On rappelle aussi quelques propriétés des fonctions cuspidales et très cuspidales démontrées dans [W1] et [W2]. \\
 
 La deuxième partie comprend des majorations sur le radical unipotent d'un sous-groupe parabolique dont on aura besoin dans la section 13. \\
 
 La troisième partie décrit brièvement certains objets associés aux groupes unitaires: sous-groupes paraboliques, sous-groupes compacts spéciaux, $R$-groupes et orbites nilpotentes régulières. \\
 
 La quatrième partie consiste principalement à définir l'expression $J_N(\theta,f)$. \\
 
 Dans la cinquième partie, se trouvent les définitions nécessaires à l'énoncé du développement géométrique (théorème 5.4.1). \\
 
 La démonstration du théorème 5.4.1 fait l'objet des sections 6 à 10. Le schéma de la preuve est le même que dans [W1]. Dans la section 6, on descend le problème à l'algèbre de Lie, dans la section 7 on transforme l'expression obtenue grâce à une transformée de Fourier. Une fois effectuée cette transformation, tout se passe beaucoup mieux et on montre dans la huitième partie que $J_N(\theta,f)$ admet une limite et on calcule cette limite. Le résultat principal de la section 9 permet de se ramener au cas où $f$ n'a pas d'éléments unipotents dans son support. Cela permet dans la section 10 une preuve par récurrence du théorème 5.4.1. \\
 
 Dans les sections 11 à 13, on montre les majorations qui seront nécessaires lors du développement spectral (section 16). Ce sont les analogues des majorations 4.3(2) à 4.3(8) de [W2]. Les démonstrations sont pourtant de nature différente. Dans la section 11 on montre une "`décomposition de Cartan relative"' dans le cas $r=0$ (inspirée des résultats de [BO], [DS] et [Sa]). Alliée aux majorations de la section 2, on obtient une démonstration plus directe des résultats. \\
 
 Dans la section 14, on montre que les éléments de $Hom_{H,\xi}(\pi,\sigma)$, pour $\pi$ et $\sigma$ tempérées, peuvent être calculés explicitement. Ce résultat est nécessaire pour faire apparaître la multiplicité $m(\pi,\sigma)$ dans le développement spectral. On en donne une démonstration différente de celle de [W2] basée sur l'existence d'une décomposition de Cartan relative ainsi que sur des éléments directement inspirés de la proposition 6.4.1 de [SV]. \\
 
 Dans la quinzième partie, on montre qu'en un certain sens la multiplicité $m(\pi,\sigma)$ est compatible à l'induction parabolique pour les représentations tempérées. \\
 
 C'est dans la seizième partie que l'on énonce et démontre le développement spectral. La démonstration suit de très près celle de [W2], à tel point que l'on s'est contenté en général d'indiquer les quelques modifications à apporter. \\
 
 On donne dans la section 17 la preuve du théorème 2. Comme on l'a dit, cela découle de l'égalité entre les développements géométrique et spectral. Néanmoins, il est nécessaire d'affaiblir la condition sur $f$ et de supposer qu'il s'agit seulement d'une fonction cuspidale. Suivant un procédé d'Arthur, on transforme l'égalité précédente en une égalité entre distributions invariantes. Le lemme 1.8.1 (qui n'est que le rappel du lemme 2.7 de [W2]) permet alors cet affaiblissement sur $f$. \\
 
 La dix-huitième et dernière partie est celle où l'on démontre le théorème 1. C'est dans cette partie que l'on fait des hypothèses sur les L-paquets (paragraphe 18.1).

\vspace{4mm}

\textbf{Remerciements}: Cet article est une partie de ma thèse sous la direction de Jean-Loup Waldspurger. Qu'il soit chaleureusement remercié pour m'avoir proposé ce sujet si riche, pour ses indications toujours justes et clairvoyantes et pour ses relectures attentives des différentes versions préliminaires du présent article.

\bigskip

\section{Groupes, mesures, notations}

\subsection{Groupes}
Soit $F$ un corps local non archimédien de caractéristique nulle. On note $\mathcal{O}_F$ son anneau d'entiers, $\mathfrak{p}_F$ l'idéal maximal de $\mathcal{O}_F$, $q_F$ le cardinal du corps résiduel et on fixe une uniformisante $\pi_F$. On désigne par $val_F$ la valuation normalisée par $val_F(\pi_F)=1$ et $|.|_F=q_F^{-val_F}$ la valeur absolue. On fixe $\psi$ un caractère additif de $F$ trivial sur $\mathcal{O}_F$ et non trivial sur $\mathfrak{p}_F^{-1}$. Si $X$ est un espace topologique totalement discontinu, nous noterons $C_c^\infty(X)$ l'espace des fonctions localement constantes à support compact de $X$ dans $\mathbb{C}$. Sauf mention explicite du contraire, tous les groupes et sous-groupes de cet article seront supposés définis sur $F$. \\

Soit $G$ un groupe réductif connexe défini sur $F$, $\mathfrak{g}$ son algèbre de Lie. On note $A_G$ le plus grand tore déployé central de $G$ et $X(G)$ le groupe des caractères algébriques de $G$ définis sur $F$. On note $G^1$ le sous-groupe des éléments $g\in G(F)$ vérifiant $|\chi(g)|_F=1$ pour tout $\chi\in X(G)$. On pose $\mathcal{A}_G=Hom(X(G),\mathbb{R})$ et $\mathcal{A}_G^*=X(G)\otimes \mathbb{R}$. On définit une application $H_G:G(F)\to \mathcal{A}_G$ par $<H_G(g),\chi>=log(|\chi(g)|_F)$ pour tout $\chi\in X(G)$. On note $\mathcal{A}_{G,F}$ et $\tilde{\mathcal{A}}_{G,F}$ les images par $H_G$ de $G(F)$ et $A_G(F)$ respectivement. Ce sont des réseaux dans $\mathcal{A}_G$. On définit $\mathcal{A}_{G,F}^\vee$ (resp. $\tilde{\mathcal{A}}_{G,F}^\vee$) comme l'ensemble des $\lambda\in\mathcal{A}_G^*$ vérifiant $\lambda(\zeta)\in 2\pi\mathbb{Z}$ pour tout $\zeta\in\mathcal{A}_{G,F}$ (resp. pour tout $\zeta\in\tilde{\mathcal{A}}_{G,F}$). On pose $i\mathcal{A}_{G,F}^*=i\mathcal{A}_G^*/i\mathcal{A}_{G,F}^\vee$. Soit $T$ un tore maximal de $G$, on pose $\delta(G)=dim(G)-dim(T)$, cela ne dépend pas du choix de $T$. \\

Pour $g\in G(F)$, on notera $X\mapsto gXg^{-1}$ l'action adjointe de $g$ sur $\mathfrak{g}$. On définit $G_{ss}$ et $\mathfrak{g}_{ss}$ comme les ensembles des éléments semi-simples dans $G$ et $\mathfrak{g}$ respectivement. On note $G_{reg},\mathfrak{g}_{reg}$ les ensembles des éléments réguliers semi-simples dans $G$ et $\mathfrak{g}$ respectivement. Pour $S$ une partie de $G$, on notera $Z_G(S)$ le centralisateur de $S$ dans $G$ et $G_S$ la composante connexe de l'élément neutre dans $Z_G(S)$. Si $f$ est une fonction sur $G(F)$ et $g\in G(F)$, on définit la fonction ${}^g f$ par ${}^g f(x)=f(g^{-1}xg)$ et on définit de même ${}^g f$ si $f$ est une fonction sur $\mathfrak{g}$. Pour $x\in G_{ss}(F)$ et $X\in \mathfrak{g}_{ss}(F)$, on pose
$$D^G(x)=|det(ad(x)_{\mathfrak{g}/\mathfrak{g}_x})|_F$$
$$D^G(X)=|det(ad(X)_{\mathfrak{g}/\mathfrak{g}_X})|_F$$

On appelle Levi de $G$ tout sous-groupe de $G$ qui est une composante de Levi d'un sous-groupe parabolique de $G$ tout deux définis sur $F$. Si $M$ est un Levi de $G$ on note $\mathcal{F}(M)$, $\mathcal{P}(M)$ et $\mathcal{L}(M)$ l'ensemble des sous-groupes paraboliques qui contiennent $M$ resp. de composante de Levi $M$ resp. des Levi qui contiennent $M$. Si $P$ est un sous-groupe parabolique de $G$ on notera $\overline{P}$ le parabolique opposé et $\delta_P$ le module usuel. Fixons $P_{min}=M_{min}U_{min}$ un sous-groupe parabolique minimal et $K$ un sous-groupe compact spécial en bonne position par rapport à $M_{min}$. On appellera sous-groupe parabolique standard (resp. antistandard) un sous-groupe parabolique qui contient $P_{min}$ (resp. qui contient $\overline{P}_{min}$) et sous-groupe parabolique semi-standard un sous-groupe parabolique qui contient $M_{min}$. Soit $P$ est un sous-groupe parabolique semi-standard, alors $P$ possède une unique composante de Levi qui contient $M_{min}$. Lorsque l'on notera $P=MU$ on sous-entendra que $M$ est ce Levi et que $U$ est le radical unipotent de $P$. On a alors $G(F)=M(F)U(F)K$ et pour $g\in G(F)$ on notera $g=m_P(g)u_P(g)k_P(g)$ une décomposition de $g$ pour laquelle $m_P(g)\in M(F),u_P(g)\in U(F)$ et $k_P(g)\in K$. L'application $G(F)\to \mathcal{A}_M$, $g\mapsto H_M(m_P(g))$ est bien définie, on la note $H_P$. Pour $P=MU$ un sous-groupe parabolique semi-standard on notera $W^M=Norm_{M(F)}(M_{min})/M_{min}(F)$ le groupe de Weyl de $M$ et $W(M)=Norm_{G(F)}(M)/M(F)$.

On désignera par $\mathcal{T}(G)$ un ensemble de représentants des classes de conjugaison par $G(F)$ dans l'ensemble des tores maximaux de $G$.

On notera $\Xi^G$ la fonction de Harish-Chandra sur $G(F)$ comme définie dans [W3]. Elle dépend du choix d'un sous-groupe parabolique minimal de $G$ et d'un sous-groupe compact spécial en bonne position. Changer de choix remplace $\Xi^G$ par une fonction équivalente. Puisque $\Xi^G$ ne sera utilisée que pour des questions de majorations, les choix n'importent pas. Soit $\rho: G\to GL(V)$ une représentation algébrique fidèle de $G$. Fixons une norme $||.||$ sur $V$. On en déduit une norme aussi notée $||.||$ sur l'espace des endomorphismes de $V$. Pour $g\in G(F)$ on pose $\sigma(g)=sup(1,log(||g||),log(||g^{-1}||))$. On a alors $\sigma(gg')\leqslant \sigma(g)+\sigma(g')\leqslant 2\sigma(g)\sigma(g')$ pour tous $g,g'\in G(F)$. Soit $b>0$ un réel, on note $\mathbf{1}_{\sigma>b}$ et $\mathbf{1}_{\sigma\leqslant b}$ les fonctions caractéristiques de l'ensemble des $g\in G(F)$ vérifiant $\sigma(g)>b$, $\sigma(g)\leqslant b$ respectivement. On adoptera la notation commode mais imprécise suivante: si $F_1$ et $F_2$ sont deux fonctions à valeurs réelles positives dépendant de paramètres $x_1,\ldots,x_N$, on notera

$$F_1(x_1,\ldots,x_N)<<F_2(x_1,\ldots,x_N)$$

\noindent pour signifier qu'il existe une constante $C>0$ telle que $F_1(x_1,\ldots,x_N)\leqslant CF_2(x_1,\ldots,x_N)$ pour toutes valeurs des paramètres $x_1,\ldots,x_N$. On dira alors que $F_1$ est essentiellement majorée par $F_2$. \\

Soit $\mathcal{S}(G(F))$ l'espace des fonctions de Schwartz-Harish-Chandra sur $G(F)$. C'est l'espace des fonctions $f:G(F)\to\mathbb{C}$ qui sont biinvariantes par un sous-groupe compact-ouvert de $G(F)$ et qui vérifient pour tout $D\in\mathbb{R}$ la majoration
$$|f(g)|<<\Xi(g)\sigma(g)^{-D}$$
pour tout $g\in G(F)$.
Soit $\pi$ une représentation admissible de $G(F)$. On notera alors $E_\pi$ un espace vectoriel sur lequel $\pi$ se réalise. Si $\pi$ est unitaire, on note $(,)$ une forme hermitienne définie positive $G(F)$-invariante. Nous dirons que $\pi$ est tempérée si elle est unitaire, de longueur finie et qu'il existe un réel $D$ tel que pour tous $e,e'\in E_\pi$ on ait
$$|(e',\pi(g)e)|<< \Xi(g) \sigma(g)^D$$
pour tout $g\in G(F)$. Si $\pi$ est tempérée et $G(F)$ est muni d'une mesure de Haar, l'action de $C_c^\infty(G(F))$ sur $E_\pi$ s'étend en une action de $\mathcal{S}(G(F))$. Pour tout $f\in \mathcal{S}(G(F))$ et tous $e,e'\in E_\pi$ on a alors
$$(\pi(f)e,e')=\displaystyle\int_{G(F)} f(g) (\pi(g)e,e') dg$$

On supposera jusqu'en 2.4 inclus fixé le groupe $G$ ainsi qu'un sous-groupe parabolique minimal $P_{min}=M_{min}U_{min}$ et un sous-groupe compact spécial $K$ de $G(F)$ en bonne position par rapport à $M_{min}$.

\subsection{Mesures}

La formule géométrique (théorème 5.4.1) ne dépend que du choix de mesures sur certains tores anisotropes et la formule spectrale (théorème 16.1.1) que du choix de mesures sur $i\mathcal{A}_{L,F}^*$ où $L$ est un sous-groupe de Levi semi-standard de $G$. Néanmoins, les preuves font intervenir des mesures sur d'autres groupes et les choix naturels varient. \\

Pour le développement géométrique (sections 5 à 10) on normalisera les mesures comme suit. Fixons une forme bilinéaire symétrique non dégénérée $G(F)$-invariante $<.,.>$ sur $\mathfrak{g}(F)$. Ceci permet de définir sur $\mathfrak{g}(F)$ une transformée de Fourier. Pour $f\in C_c^\infty(\mathfrak{g}(F))$, on pose
$$\hat{f}(X)=\displaystyle\int_{\mathfrak{g}} f(Y)\psi(<Y,X>) dY$$
où $dY$ est la mesure de Haar autoduale sur $\mathfrak{g}(F)$ c'est-à-dire telle que $\hat{\hat{f}}(X)=f(-X)$. De la même façon un sous-espace de $\mathfrak{g}(F)$ pour lequel la restriction de $<.,.>$ est non dégénérée est naturellement muni d'une mesure autoduale. On munit les sous-espaces de $\mathfrak{g}(F)$ pour lesquels la restriction de $<.,.>$ est dégénérée d'une mesure de Haar quelconque. Si $H$ est un sous-groupe algébrique de $G$ on munit $H(F)$ d'une mesure en relevant celle fixée sur $\mathfrak{h}(F)$ par l'exponentielle. Enfin on munira les sous-groupes compacts de $G(F)$ qui ne sont pas de façon évidente le groupe des $F$-points d'un sous-groupe algébrique de $G$ de la mesure de Haar de masse totale 1. Soit $T\subset G$ un tore, alors on a un deuxième choix de mesure $d_c t$ sur $T(F)$: si $T$ est déployé c'est la mesure qui donne au sous-groupe compact maximal de $T(F)$ la mesure $1$, dans le cas général $d_ct$ est compatible avec la mesure que l'on vient de définir sur $A_T(F)$ et $T(F)/A_T(F)$ est de mesure 1 pour $d_ct$. On note $\nu(T)$ le facteur tel que $d_c t=\nu(T)dt$. \\

On note $Nil(\mathfrak{g})$ l'ensemble des orbites nilpotentes de $\mathfrak{g}(F)$. Soit $\mathcal{O}$ une telle orbite. Pour $X\in\mathcal{O}$, la forme bilinéaire $(Y,Z)\mapsto <X,[Y,Z]>$ sur $\mathfrak{g}(F)$ se descend en une forme symplectique sur $\mathfrak{g}(F)/\mathfrak{g}_X(F)$, c'est-à-dire sur l'espace tangent à $\mathcal{O}$ au point $X$. Ainsi $\mathcal{O}$ est muni d'une structure de variété $F$ analytique symplectique et on en déduit une mesure "`autoduale"' sur $\mathcal{O}$. Cette mesure est invariante par conjugaison par $G(F)$.

\vspace{3mm}

La normalisation des mesures lors du développement spectral (section 16) sera la suivante. Rappelons que l'on a fixé $P_{min}=M_{min}U_{min}$ un sous-groupe parabolique minimal de $G$ et $K$ un sous-groupe compact spécial en bonne position par rapport à $A_{M_{min}}$. On choisit sur $K$ la mesure de Haar de masse totale 1 et sur $G(F)$ une mesure de Haar quelconque. Pour tout sous-groupe parabolique semi-standard $P=MU$ de $G$, on choisit sur $U(F)$ l'unique mesure de Haar vérifiant
$$\displaystyle\int_{U(F)} \delta_{\overline{P}}(m_{\overline{P}}(u)) du=1$$
Il existe alors une unique mesure de Haar $dm$ sur $M(F)$ qui vérifie
$$\displaystyle\int_{G(F)} f(g) dg=\int_{M(F)}\int_{U(F)}\int_K f(muk) dkdudm$$
pour tout $f\in C_c^\infty(G(F))$. Et la mesure $dm$ ne dépend pas du choix de $P\in\mathcal{P}(M)$. Les réseaux $i\mathcal{A}_{M,F}^\vee$ et $i\tilde{\mathcal{A}}_{M,F}^\vee$ sont munis des mesures de comptage. On munit $i\mathcal{A}_M^*$ de la mesure de Haar qui donne à $i\mathcal{A}_M^*/i\tilde{\mathcal{A}}_{M,F}^\vee$ la mesure $1$ et $i\mathcal{A}_{M,F}^*$ de la mesure quotient.

\vspace{2mm}

On utilisera ainsi la première normalisation des mesures pour toutes les définitions "géométriques" (section 1.4) et la deuxième normalisation pour toutes les définitions "spectrales" (sections 1.5 à 1.8).

\subsection{Bons voisinages}
Soit $x\in G(F)$ un élément semi-simple. On reprend la définition 3.1 de [W1] de la notion de bon voisinage de $\mathfrak{g}_x(F)$. C'est un voisinage de $0$ dans $\mathfrak{g}_x(F)$ qui est $Z_G(x)(F)$-invariant, compact modulo conjugaison sur lequel l'exponentielle est définie et qui vérifie les conditions 1-7 de la section 3.1 de [W1]. Soit $f$ une fonction sur $G$ et $\omega$ un bon voisinage de $\mathfrak{g}_x(F)$. On note $f_{x,\omega}$ la fonction sur $\mathfrak{g}_x(F)$ définie par

$$
f_{x,\omega}(X)= \left\{
    \begin{array}{ll}
        f(xexp(X)) & \mbox{si } X \in \omega \\
        0 & \mbox{sinon.}
    \end{array}
\right.
$$

\subsection{Intégrales orbitales invariantes et pondérées, quasi-caractères}
Pour $X\in\mathfrak{g}_{reg}(F)$, on définit l'intégrale orbitale

$$J_G(X,f)=D^G(X)^{1/2}\displaystyle\int_{G_X(F)\backslash G(F)} f(g^{-1} Xg) dg$$

\noindent pour tout $f\in C_c^\infty(\mathfrak{g}(F))$. Il existe une unique fonction $\hat{j}$ sur $\mathfrak{g}(F)\times \mathfrak{g}(F)$, localement intégrable, localement
constante sur $\mathfrak{g}_{reg}(F)\times \mathfrak{g}_{reg}(F)$, telle que, pour toute $f\in C_c^\infty(\mathfrak{g}(F))$ et tout $X\in \mathfrak{g}_{reg}(F)$,
on ait l'égalité :

$$\displaystyle J_G(X,\hat{f})=\int_{\mathfrak{g}(F)} f(Y) \hat{j}(X,Y) dY$$

Soient $M$ un Levi de $G$ et $X\in\mathfrak{m}(F)\cap\mathfrak{g}_{reg}(F)$. Pour $Y\in\mathfrak{g}_{reg}(F)$, fixons un ensemble de représentants $(Y_i)_{i=1,\ldots,r}$ des classes de conjugaison par $M(F)$ de l'ensemble des éléments de $\mathfrak{m}(F)$ qui sont conjugués à $Y$ par un élément de $G(F)$. On a alors l'égalité (cf la formule 2.6(5) de [W1])

$$\mbox{(1)} \;\;\; \displaystyle\hat{j}^G(X,Y)D^G(Y)^{1/2}=\sum_{i=1}^r \hat{j}^M(X,Y_i) D^M(Y_i)^{1/2}$$

Pour $\mathcal{O}\in Nil(\mathfrak{g})$, on définit l'intégrale orbitale nilpotente

$$J_\mathcal{O}(f)=\displaystyle\int_\mathcal{O} f(X) dX$$

\noindent pour tout $f\in C_c^\infty(\mathfrak{g}(F))$. On définit la transformée de Fourier de la distribution précédente par

$$\hat{J}_\mathcal{O}(f)=J_\mathcal{O}(\hat{f})$$

\noindent pour tout $f\in C_c^\infty(\mathfrak{g}(F))$, où intègre pou la mesure sur $\mathcal{O}$ fixée en 1.2. D'après Harish-Chandra, $\hat{J}_\mathcal{O}$ est représentable par une fonction $Y\mapsto \hat{j}(\mathcal{O},Y)$ localement intégrable et localement constante sur $\mathfrak{g}_{reg}(F)$. Pour $\lambda\in F^\times$, définissons $f^\lambda$ par $f^\lambda(X)=f(\lambda X)$. Pour tout $\lambda\in F^\times$, on a

$$J_\mathcal{O}(f^\lambda)=|\lambda|_F^{-dim(\mathcal{O})/2}J_{\lambda\mathcal{O}}(f)$$
On en déduit que pour tout $\lambda\in F^\times$ on a $\hat{j}(\mathcal{O},\lambda Y)=|\lambda|_F^{-dim(\mathcal{O})/2}\hat{j}(\lambda \mathcal{O},Y)$. Pour $\mathcal{O}\in Nil(\mathfrak{g})$, on note $\Gamma_\mathcal{O}$ le germe de Shalika associé à $\mathcal{O}$. Pour tout $f\in C_c^\infty(\mathfrak{g}(F))$ il existe donc un voisinage $\omega$ de $0$ tel qu'on ait pour tout $X\in \omega\cap \mathfrak{g}_{reg}(F)$

$$J_G(X,f)=\displaystyle\sum_{\mathcal{O}\in Nil(\mathfrak{g})} \Gamma_\mathcal{O}(X) J_\mathcal{O}(f)$$
Et de plus, $\Gamma_\mathcal{O}$ vérifie la condition habituelle d'homogénéité. \\

Soit $M$ un Levi semistandard de $G$. Pour tout $g\in G(F)$ la famille $(H_{P}(g))_{P\in \mathcal{P}(M)}$ est $(G,M)$-orthogonale positive. D'après Arthur ([A?]) on peut donc lui associer un nombre $v_M(g)$. Il faut pour cela avoir fixé une mesure sur $\mathcal{A}_M^G$, choix qui a été effectué en 1.2. La fonction $g\mapsto v_M(g)$ est invariante à gauche par $M(F)$ et à droite par $K$. Remarquons qu'elle ne dépend pas que de $M$ mais aussi du choix du sous-groupe compact spécial $K$ en bonne position par rapport à $M$. Pour $f\in C_c^\infty(G(F))$ et $x\in M(F)\cap G_{reg}(F)$, on définit l'intégrale orbitale pondérée de $f$ en $x$ par

$$\displaystyle J_M(x,f)=D^G(x)^{1/2}\int_{G_x(F)\backslash G(F)} f(g^{-1}xg) v_M(g) dg$$

\vspace{2mm}

Un quasi-caractère sur $G(F)$ est une fonction définie presque partout $\theta: G(F)\to \mathbb{C}$ invariante par conjugaison et tel que pour tout $x\in G_{ss}(F)$ il existe un bon voisinage $\omega$ de $\mathfrak{g}_x(F)$ et des coefficients $c_{\mathcal{O},\theta}(x)$ pour $\mathcal{O}\in Nil(\mathfrak{g}_x)$ vérifiant

\begin{center}
$\theta(xexp(X))=\displaystyle\sum_{\mathcal{O}\in Nil(\mathfrak{g}_x)} c_{\mathcal{O},\theta}(x) \hat{j}(\mathcal{O},X)$ pour presque tout $X\in\omega$
\end{center}

\subsection{Représentations, induites paraboliques, opérateurs d'entrelacements, caractères pondérés}

Toutes les représentations sont supposées lisses. On adopte les notations suivantes: $Irr(G)$, $Temp(G)$, $\Pi_2(G)$ et $\Pi_{ell}(G)$ désignent respectivement l'ensemble des classes d'isomorphismes de représentations irréductibles, irréductibles tempérées, irréductibles de la série discrète et irréductibles elliptiques de $G(F)$.
Si $P=MU$ est un sous-groupe parabolique de $G$ et $\tau$ est une représentation de $M(F)$ on notera $i^G_P(\tau)$ l'induite normalisée de $\tau$ à $G$. Elle se réalise sur l'espace $E^G_{P,\tau}$ des fonctions $\varphi: G(F)\to E_\tau$ localement constantes qui vérifient $\varphi(mug)=\delta_P(m)^{1/2}\tau(m)f(g)$ pour tous $m\in M(F),u\in U(F)$ et $g\in G(F)$, l'action de $G(F)$ se faisant par translation à droite. Pour $g\in G(F)$ et $f\in C_c^\infty(G(F))$, on note $i^G_P(\tau,g)$ et $i^G_P(\tau,f)$ les actions de $g$ et $f$ agissant sur $E_{P,\tau}^G$. \\

Pour $\pi$ une représentation de $G(F)$ et $\lambda\in \mathcal{A}_G^*\otimes_{\mathbb{R}} \mathbb{C}/i\mathcal{A}_{G,F}^\vee$, on note $\pi_\lambda$ la représentation suivante de $G(F)$: son espace est le même que celui de $\pi$ et $\pi_\lambda(g)=e^{<H_G(g),\lambda>}\pi(g)$ pour tout $g\in G(F)$. Soit $P=MU$ un sous-groupe parabolique semi-standard de $G$ et $\tau$ une représentation de $M(F)$. On note $\mathcal{K}^G_{P,\tau}$ l'espace des fonctions $\varphi:K\to E_\tau$ qui vérifient $\varphi(muk)=\tau(m)\varphi(k)$ pour tous $m\in M(F)\cap K, u\in U(F)\cap K$ et $k\in K$. Les représentations $i^G_P(\tau_\lambda)$ pour $\lambda\in \mathcal{A}_M^*\otimes \mathbb{C}$ se réalisent toutes sur l'espace $\mathcal{K}^G_{P,\tau}$ via l'isomorphisme qui à une fonction $\varphi\in E^G_{P,\tau_\lambda}$ associe sa restriction à $K$. Si $\tau$ est unitaire, on munit $\mathcal{K}^G_{P,\tau}$ d'un produit hermitien par
$$(e,e')=\displaystyle\int_K (e(k),e'(k)) dk$$
C'est un produit scalaire invariant pour la représentation $i^G_P(\tau_\lambda)$ pour tout $\lambda\in i\mathcal{A}_{M,F}^*$.\\
Pour $P,P'\in\mathcal{P}(M)$ et $\lambda\in \mathcal{A}_M^*\otimes_\mathbb{R} \mathbb{C}$, on peut définir des opérateurs d'entrelacements
$$J_{P'|P}(\tau_\lambda): E^G_{P,\tau_\lambda}\to E^G_{P',\tau_\lambda}$$
Quand la partie réelle de $\lambda$ est dans un certain cône, on a

$$(J_{P'|P}(\tau_\lambda)e)(g)=\displaystyle\int_{(U'(F)\cap U(F))\backslash U'(F)} e(u'g) du'$$

La fonction $\lambda\mapsto J_{P'|P}(\tau_\lambda)$ à valeurs dans les opérateurs $\mathcal{K}^G_{P,\tau}\to \mathcal{K}^G_{P',\tau}$ admet un prolongement méromorphe à tout $\mathcal{A}_M^*\otimes_\mathbb{R} \mathbb{C}$. Si $\tau$ est irréductible, pour tout $P\in \mathcal{P}(M)$ l'opérateur $J_{P|\overline{P}}(\tau_\lambda)J_{\overline{P}|P}(\tau_\lambda)$ est scalaire. On note $j(\tau_\lambda)$ ce scalaire, il ne dépend pas du choix de 
$P$. On peut normaliser les opérateurs d'entrelacements: il existe des facteurs $r_{P'|P}(\tau_\lambda)$ tels que les opérateurs $R_{P'|P}(\tau_\lambda)=r_{P'|P}(\tau_\lambda)^{-1} J_{P'|P}(\tau_\lambda)$ vérifient les conditions du théorème 2.1 de [A4]. Notamment
\begin{itemize}
\item Pour tous $P,P',P''\in\mathcal{P}(M)$, $R_{P''|P'}(\tau_\lambda)R_{P'|P}(\tau_\lambda)=R_{P''|P}(\tau_\lambda)$
\item Si $\tau$ est tempérée, les fonctions $\lambda\mapsto R_{P'|P}(\tau_\lambda)$ sont holomorphes en tout point de $i\mathcal{A}_{M,F}^*$ et pour $\lambda\in i\mathcal{A}_{M,F}^*$, l'adjoint de $R_{P'|P}(\tau_\lambda)$ est $R_{P|P'}(\tau_\lambda)$
\end{itemize}

La définition des opérateurs d'entrelacements normalisés s'étend au cas où $\tau$ est semi-simple. Soit $\tau$ une représentation tempérée et $M$ un Lévi semistandard de $G$. Fixons $P\in\mathcal{P}(M)$ et posons pour tout $P'\in\mathcal{P}(M)$ et tout $\lambda\in i\mathcal{A}_{M,F}^*$
$$\mathcal{R}_{P'}(\tau,\lambda)=R_{P|P'}(\tau)R_{P'|P}(\tau_\lambda)$$
C'est un endomorphisme de $\mathcal{K}^G_{P,\tau}$. la famille $(\mathcal{R}_{P'}(\tau,.))_{P'\in\mathcal{P}(M)}$ est une $(G,M)$-famille à valeur opérateurs ([A1] paragraphe 7). On peut donc suivant Arthur lui associer un opérateur que l'on note $\mathcal{R}_M(\tau)$. \\

Le caractère pondéré de la représentation $\tau$ est la distribution $f\mapsto J_M^G(\tau,f)$ qui à $f\in C_c^\infty(G(F))$ associe la trace de l'opérateur $i^G_P(\tau,f)\mathcal{R}_M(\tau)$ agissant sur $\mathcal{K}^G_{P,\tau}$. Cette distribution ne dépend en fait pas du sous-groupe parabolique $P$ que l'on a fixé. Dans le cas où $M=G$ on notera simplement $\theta_\tau(f)=J_G^G(\tau,f)$, c'est le caractère usuel de $\tau$.

\subsection{R-groupes}

On utilisera la théorie des $R$-groupes telle qu'elle est exposée dans la section 1.5 de [W2]. On reprend aussi les notations de cette référence : $Norm_{G(F)}(\tau)$, $W(\tau)$, $W'(\tau)$, $R(\tau)$, $R_P(w,\tau)$. Il est fait dans [W2] un certain nombre d'hypothèses qui permettent de simplifier la théorie. Nous verrons en 3.1 que ces hypothèses sont vérifiées dans le cas des groupes unitaires. Les hypothèses sont les suivantes: soit $M$ un Levi semistandard de $G$ et $\tau\in\Pi_2(M)$, alors
\begin{itemize}
\item $\tau$ se prolonge en une représentation $\tau^N$ de $Norm_{G(F)}(\tau)$.
\item L'homomorphisme naturel $K\cap Norm_{G(F)}(\tau)\to W(\tau)$ admet une section qui est un homomorphisme.
\item Le $R$-groupe $R(\tau)$ est abélien.
\end{itemize}

Rappelons comment l'on retrouve les représentations elliptiques de $G(F)$ à partir des $R$-groupes. Soit $M$ un Levi de $G$, $\tau\in \Pi_2(M)$ et $P\in\mathcal{P}(M)$. On a alors une décomposition en représentations irréductibles $\displaystyle i_P^G(\tau)=\bigoplus_{\zeta\in R(\tau)^\vee} i^G_P(\tau,\zeta)$, où $R(\tau)^\vee$ est le dual du groupe abélien $R(\tau)$ et $i^G_P(\tau,\zeta)$ est la sous-représentation de $i^G_P(\tau)$ où $R_P(w,\tau)$ agit comme $\zeta(w)$ pour tout $w\in R(\tau)$. Notons $W(M)_{reg}$ l'ensemble des éléments de $W(M)$ qui agissent sans point fixe sur $\mathcal{A}_M/\mathcal{A}_G$. Alors une représentation $i^G_P(\tau,\zeta)$ comme ci-dessus est elliptique si et seulement si $R(\tau)\cap W(M)_{reg}\neq \emptyset$ et si cette condition est vérifiée, on a $W'(\tau)=\{1\}$. On obtient ainsi toutes les représentations irréductibles elliptiques de $G(F)$.

\subsection{Formule de Plancherel-Harish-Chandra}

Pour tout $M\in\mathcal{L}(M_{min})$, on fixe un élément $P\in \mathcal{P}(L)$. Notons $\{\Pi_2(M)\}$ l'ensemble des orbites de $\Pi_2(M)$ pour l'action de $i\mathcal{A}_{M,F}^*$. Pour chaque orbite $\mathcal{O}$, fixons un élément $\tau$ de cette orbite. Notons $i\mathcal{A}_\mathcal{O}^\vee$ le stabilisateur de $\tau$ dans $i\mathcal{A}_M^*$. Pour tout $\lambda\in i\mathcal{A}_{M,F}^*$ on pose

$$m(\tau_\lambda)=j(\tau_\lambda)^{-1} d(\tau)$$
où $d(\tau)$ est le degré formel. Pour tout $f\in \mathcal{S}(G(F))$, on a alors l'égalité

$$\displaystyle f(g)=\sum_{M\in \mathcal{L}(M_{min})} |W^M||W^G|^{-1} \sum_{\mathcal{O}\in\{\Pi_2(M)\}} [i\mathcal{A}_\mathcal{O}^\vee:i\mathcal{A}_{M,F}^\vee]^{-1}$$
$$\displaystyle \int_{i\mathcal{A}_{M,F}^*} m(\tau_\lambda) Tr\left(i_P^G(\tau_\lambda,g^{-1}) i_P^G(\tau_\lambda,f)\right) d\lambda$$
pour tout $g\in G(F)$. C'est la formule de Plancherel-Harish-Chandra et c'est le théorème VIII.1.1 de [W3]. Soit $K_f$ un sous-groupe compact-ouvert de $G(F)$ tel que $f$ soit biinvariante par $K_f$. Seules interviennent de façon non nulle les orbites $\mathcal{O}$ pour lesquelles une représentation $i^G_P(\tau_\lambda)$ admet des invariants non nuls par $K_f$. Ces orbites
sont en nombre fini. \\

Fixons $P=MU\in\mathcal{F}(M_{min})$ et une représentation $\tau$ de $M(F)$ irréductible et de la série discrète. Notons $\pi_\lambda=i^G_P(\tau_\lambda)$ pour tout $\lambda\in i\mathcal{A}_{M,F}^*$. Soit $\varphi$ une fonction $C^\infty$ sur $i\mathcal{A}_{M,F}^*$ et $e,e'\in\mathcal{K}^G_{P,\tau}$. Posons

$$\displaystyle f_{e,e',\varphi}(g)=\int_{i\mathcal{A}_{M,F}^*} \varphi(\lambda) m(\tau_\lambda) (\pi_\lambda(g)e',e) d\lambda$$

\noindent pour tout $g\in G(F)$. Cette fonction appartient à $\mathcal{S}(G(F))$. Identifions tout élément de $W(M)$ à un représentant dans $K\cap Norm_{G(F)}(M)$. Notons ${\cal E}(\tau)$ l'ensemble des couples $(w,\mu)\in W(M)\times i{\cal A}_{M,F}^*$ tels que $w^{-1}\tau\simeq \tau_{\mu}$. Pour $(w,\mu)\in {\cal E}(\tau)$, fixons un automorphisme unitaire $\tau(w,\mu)$ de $E_\tau$ tel que 
$$\tau(w,\mu)\tau_{\mu}(m)=(w^{-1}\tau)(m)\tau(w,\mu)$$
pour tout $m\in M(F)$. Définissons l'homomorphisme $A(w,\mu):{\cal K}^G_{w^{-1}Pw,\tau}\to {\cal K}_{P,\tau}^G$ par 
$$(A(w,\mu)e)(g)=\tau(w,\mu)e(w^{-1}g).$$
Pour $\lambda\in i{\cal A}_{M,F}^*$, définissons l'endomorphisme $R(w,\mu,\lambda)$ de ${\cal K}^G_{P,\tau}$ par
$$R(w,\mu,\lambda)=A(w,\mu)R_{w^{-1}Pw\vert P}(\tau_{\lambda+\mu}).$$
Il vérifie la relation d'entrelacement 
$$R(w,\mu,\lambda)\pi_{\lambda+\mu}(g)=\pi_{w\lambda}(g)R(w,\mu,\lambda).$$
Soient $e_{0},e'_{0}\in {\cal K}^G_{P,\tau}$. Alors on a l'égalité
$$\mbox{(1)} \;\;\; \int_{G(F)}f_{e,e',\varphi}(g)(e'_{0},\pi_\lambda(g)e_{0})dg=\sum_{(w,\mu)\in {\cal E}(\tau)}\varphi(\mu+w^{-1}\lambda)(R(w,\mu,w^{-1}\lambda)e',e_{0})(e'_{0},R(w,\mu,w^{-1}\lambda)e).$$
Cf proposition VII.2 de [W3]. \\

On en déduit le point suivant \\
 
 (2) supposons le support de $\varphi$ vérifie la condition suivante
 
\begin{center}
 pour tout $(w,\mu)\in {\cal E}(\tau)$, $\mu\in Supp(\varphi)$ entraîne $\mu=0$.
\end{center}
alors on a l'égalité (avec les notations du paragraphe précédent)
$$\int_{G(F)}f_{e,e',\varphi}(g)(e'_{0},\pi_0(g)e_{0})dg=\vert W'(\tau)\vert \varphi(0)\sum_{w\in R(\tau)}(R_{P}(w,\tau)e',e_{0})(e'_{0},R_{P}(w,\tau)e).$$

\noindent cf [W2]1.6(1). On aura aussi besoin de la formule suivante

\vspace{2mm}

  (3) Soient $\varphi,\psi\in C_c^\infty\in C_c^\infty(i\mathcal{A}_{M,F}^*)$. Supposons que les supports de $\varphi$ et $\psi$ vérifie $(w^{-1} Supp(\psi)+\mu)\cap Supp(\varphi)=\emptyset$ pour tout $(w,\mu)\in \mathcal{E}(\tau)-\{(Id,0)\}$, alors on a l'égalité
$$\displaystyle \int_{G(F)} f_{e_0,e'_0,\varphi}(g_1gg_2)\overline{f_{e_1,e'_1,\psi}(g)} dg=\int_{i\mathcal{A}_{M,F}^*} m(\tau_\lambda) \varphi(\lambda)\overline{\psi}(\lambda) (\pi_\lambda(g_2)e'_0,e'_1)(\pi_\lambda(g_1)e_1,e_0) d\lambda$$

\noindent pour tous $e_0,e'_0,e_1,e'_1\in\mathcal{K}^G_{P,\tau}$ et pour tous $g_1,g_2\in G(F)$.
  
\vspace{2mm}

Après une permutation d'intégrales, l'égalité (3) découle de l'égalité (1) et de l'hypothèse faite sur $Supp(\varphi)$.

\subsection{Fonctions cuspidales et très cuspidales, quasi-caractères associés}

Une fonction $f\in C_c^\infty(G(F))$ est dite très cuspidale si elle vérifie
$$\displaystyle\int_{U(F)} f(mu)du=0$$
pour tout sous-groupe parabolique propre $P=MU$ de $G$ et pour tout $m\in M(F)$. On définit de façon analogue la notion de fonction très cuspidale sur l'algèbre de Lie. \\

Une fonction $f\in C_c^\infty(G(F))$ est dite cuspidale si pour tout sous-groupe de Levi propre $M$ de $G$ et pour tout $x\in M(F)\cap G_{reg}(F)$ on a

$$\displaystyle\int_{G_x(F)\backslash G(F)} f(g^{-1}xg) dg=0$$

Cette condition est équivalente à ce que $\theta_\pi(f)=0$ pour toute représentation $\pi$ de $G(F)$ qui est tempérée et proprement induite. Toute fonction très cuspidale est cuspidale. On a le résultat suivant

\begin{lem}
Soit $f\in C_c^\infty(G(F))$ une fonction cuspidale, il existe une fonction très cuspidale $f'\in C_c^\infty(G(F))$ telle que la propriété suivante soit vérifiée

\vspace{2mm}

Pour toute distribution $D$ sur $C_c^\infty(G(F))$ invariante par conjugaison on a $D(f')=0$ si et seulement si $D(f)=0$.
\end{lem}

\ul{Preuve}: C'est le lemme 2.7 de [W2] $\blacksquare$

Il est possible d'associer à une fonction très cuspidale $f$ un quasicaractère $\theta_f$ défini à partir des intégrales orbitales pondérées de $f$. Pour $x\in G(F)$ notons $M(x)$ le centralisateur de $A_{G_x}$ dans $G$. Soit $x\in G(F)$ et $y\in G(F)$ tel que $yM(x)y^{-1}$ soit semistandard. On pose alors

$$\theta_f(x)=(-1)^{a_{M(x)}-a_G} D^G(x)^{-1/2} \nu(G_x)^{-1} J_{yM(x)y^{-1}}(yxy^{-1},f)$$

\noindent cf lemme 5.2 de [W1] pour la preuve que cette définition ne dépend pas du choix de $y$ et les paragraphes suivants pour la preuve qu'il s'agit d'un quasicaractère. \\

Soit $f\in C_c^\infty(G(F))$. Waldspurger définit dans [W2] paragraphe 2.4 une fonction $I\theta_f$ sur $G_{reg}(F)$ qui est localement constante et invariante par conjugaison. Cette fonction est définie à partir des intégrales orbitales pondérées invariantes de Arthur. Pour tout $x\in G_{reg}(F)$ la distribution $f\mapsto I\theta_f(x)$ est invariante. D'après le lemme 2.5 de [W2], si $f$ est cuspidale, $I\theta_f$ est un quasi-caractère. On dispose donc si $f$ est très cuspidale de deux quasi-caractères construits à partir de $f$: $\theta_f$ et $I\theta_f$. Le lemme 2.6 de [W2] compare ces deux quasi-caractères. Avant de l'énoncer, il faut rappeler un résultat d'Arthur. Pour $Z\in \mathcal{A}_{G,F}$, notons $\mathbf{1}_{H_G=Z}$ la fonction caractéristique de l'ensemble des $x\in G(F)$ qui vérifie $H_G(x)=Z$. Soit $\mathcal{H}_{ac}(G(F))$ l'espace des fonctions $f:G(F)\to \mathbb{C}$ qui vérifient

\begin{itemize}
\item $f$ est biinvariante par un sous-groupe compact-ouvert de $G(F)$
\item pour tout $Z\in \mathcal{A}_{G,F}$, la fonction $f\mathbf{1}_{H_G=Z}$ est à support compact
\end{itemize}

Soient $L\in {\cal L}(M_{min})$ et $f\in C_{c}^{\infty}(G(F))$. Arthur montre qu'il existe une fonction $\phi_{L}(f)\in {\cal H}_{ac}(L(F))$ telle que, pour toute représentation $\pi\in Temp(L)$ et tout $Z\in {\cal A}_{L,F}$, on ait l'égalité
$$(1) \qquad \int_{i{\cal A}_{L,F}^*}J_{L}(\pi_{\lambda},f)exp(-\lambda(Z))d\lambda=mes(i\mathcal{A}_{L,F}^*)\theta_{\pi}(\phi_{L}(f){\bf 1}_{H_{L}=Z}).$$

On dira qu'une fonction $f\in \mathcal{H}_{ac}(G(F))$ est cuspidale si pour tout $Z\in\mathcal{A}_{G,F}$ la fonction $f\mathbf{1}_{H_G=Z}$ est cuspidale. Pour $f\in \mathcal{H}_{ac}(G(F))$, on définit la fonction $I\theta_f$ par 

\begin{center}
$I\theta_f(x)=I\theta_{f \mathbf{1}_{H_G=H_G(x)}}(x)$ pour tout $x\in G_{reg}(F)$
\end{center}
Si $f$ est cuspidale, la fonction $I\theta_f$ est toujours un quasi-caractère. Pour $L$ un sous-groupe de Levi de $G$ et $\theta$ un quasi-caractère de $L(F)$, on sait définir un quasi-caractère induit $Ind_L^G(\theta)$: c'est l'induction de $L$ à $G$ de la distribution invariante sur $L(F)$ définie par $\theta$ (cf le paragraphe 2.3 de [W2]). Le lemme 2.6 de [W2] est alors le suivant

\begin{lem}
Soit $f\in C_{c}^{\infty}(G(F))$ une fonction très cuspidale. Alors
  
  (i) pour tout $L\in {\cal L}(M_{min})$, la fonction $\phi_{L}(f)$ est cuspidale;
  
  (ii) on a l'égalité
  $$\theta_{f}=\sum_{L\in {\cal L}(M_{min})}\vert W^L\vert \vert W^G\vert ^{-1}(-1)^{a_{L}-a_{G}}Ind_{L}^G(I\theta^L_{\phi_{L}(f)}).$$

\end{lem}

\section{Majorations unipotentes}

Soit $G$ un groupe réductif sur $F$. On conserve les notations déjà fixées. En particulier, $P_{min}=M_{min}U_{min}$ est un sous-groupe parabolique minimal et $K$ est un sous-groupe compact spécial en bonne position par rapport à $M_{min}$. On notera $A_{min}$ la composante déployée du centre de $M_{min}$.

\subsection{ Une première majoration}

\begin{prop}
Soit $P=MU$ un sous-groupe parabolique de $G$. Pour tout entier $d$, il existe un entier $d'$ tel que
$$\displaystyle\int_{U(F)}\Xi^G(mu)^2 \sigma(mu)^d du << \Xi^M(m)^2\sigma(m)^{d'}$$
pour tout $m\in M(F)$
\end{prop}
\ul{Preuve}: On ne nuit pas à la généralité en supposant que $P$ est antistandard. On va d'abord montrer que pour tout $d\in\mathbb{N}$ il existe un entier $d'\in \mathbb{N}$ de sorte que
\begin{center}
(1) $\displaystyle\int_{U(F)} \Xi^G(mu)^2\sigma(mu)^d du<< \sigma(m)^{d'} \Xi^M(m)^2 \delta_{\overline{P}}(m)$ , pour tout $m\in M(F)$
\end{center}
Soit $M_{min}^{M,+}$ le sous-ensemble de $M_{min}(F)$ des éléments en position positive pour $P_{min}\cap M$. D'après Bruhat-Tits, il existe un sous-groupe compact-ouvert $K_M\subset M(F)$ tel que $M(F)=K_M M_{min}^{M,+} K_M$. Il suffit donc d'établir (1) pour $m\in M_{min}^{M,+}$. \\

D'après [W3] lemme II.4.4, il existe un entier $D\in \mathbb{N}$ tel que

\begin{center}
$\Xi^G(mu)<<\delta_{P_{min}}(m)^{1/2}\delta_{P_{min}}(m_{P_{min}}(u))^{1/2} \sigma(mu)^D$
\end{center}

\noindent pour tout $m\in M_{min}(F)$ et pour tout $u\in \overline{U}_{min}(F)$. Comme de plus $\sigma(mu)<<\sigma(m)\sigma(u)$, on en déduit la majoration

$$\displaystyle\int_{U(F)} \Xi^G(mu)^2 \sigma(mu)^d du<< \sigma(m)^{d+2D} \delta_{P_{min}}(m) \displaystyle\int_{U(F)} \delta_{P_{min}}(m_{P_{min}}(u)) \sigma(u)^{d+2D} du$$

\noindent pour tout $m\in M_{min}(F)$. Grâce au lemme II.4.1 de [W3], on vérifie aisément que l'intégrale $\displaystyle\int_{U(F)} \delta_{P_{min}}(m_{P_{min}}(u)) \sigma(u)^{d+D} du$ est convergente. D'après le lemme II.1.1 de [W3] on a $\sigma(m)^{d+D}\delta_{P_{min}}(m)<<\sigma(m)^{d+D} \Xi^M(m)^2 \delta_{\overline{P}}(m)$ pour tout $m\in M_{min}^{M,+}$. Cela établit (1). \\

Pour tout $g\in G(F)$, on a $\Xi^G(g^{-1})=\Xi^G(g)$ et $\sigma(g^{-1})=\sigma(g)$. On a donc aussi le majoration

$$\displaystyle\int_{U(F)} \Xi^G(um)^2\sigma(um)^d du<< \sigma(m)^{d'} \Xi^M(m)^2 \delta_{P}(m)$$

\noindent pour tout $m\in M(F)$. Après le changement de variable $u\mapsto m^{-1}um$, on obtient la majoration de l'énoncé $\blacksquare$

\subsection{Intégrales à paramètres}

Dans la suite on aura besoin du lemme suivant.

\begin{lem}
Soit $(X,\mathcal{A},\mu)$ un espace mesuré, $U\subset \mathbb{C}$ un ouvert et $f:X\times U\to \mathbb{C}$ une fonction vérifiant
\begin{enumerate}
\item Pour tout $z\in U$ la fonction $t\mapsto f(t,z)$ est mesurable.
\item Pour tout $t\in X$ la fonction $z\mapsto f(t,z)$ est holomorphe sur $U$.
\item Il existe un sous-ensemble discret $E$ de $U$ tel que pour tout compact $K\subset U-E$ il existe $g_K\in L^1(X)$ de sorte que pour tout $(t,z)\in X\times K$ on ait
$$|f(t,z)|\leqslant g_K(t)$$
\end{enumerate}
Alors la fonction $U\to \mathbb{C}$, $z\mapsto \displaystyle\int_{X} f(t,z) dt$ est partout définie par une intégrale absolument convergente et est holomorphe.
\end{lem}

\ul{Preuve}:  Par un théorème classique d'intégrale à paramètre la fonction $z\mapsto \displaystyle\int_X f(t,z) dt$ est définie par une intégrale absolument convergente et est holomorphe sur $U-E$. Soit $z_0\in E$ et $r>0$ de sorte que $B(z_0,r)\subset U$ et $B(z_0,r)\cap E=\{z_0\}$. On applique 3) à $K=S(z_0,r)=\{z; |z-z_0|=r\}$, ce qui fournit une fonction $g_K\in L^1(X)$. Par le principe du maximum, on a alors pour tout $(t,z)\in X\times B(z_0,r)$
$$|f(t,z)|\leqslant sup_{z\in K} |f(t,z)|\leqslant g_K(t)$$
Toujours par un théorème d'intégrale à paramètre, on en déduit que $z\mapsto \displaystyle\int_X f(t,z) dt$ est absolument convergente et holomorphe sur l'intérieur de $B(z_0,r)$. Le raisonnement étant valable pour tout $z_0\in E$, le résultat s'en suit. $\blacksquare$

\subsection{Fonctionnelles de Jacquet}

Introduisons l'ensemble $\Sigma$ des racines de $A_{min}$ dans $G$, $\Sigma^+$ l'ensemble des racines dans $P_{min}$ et $\Delta$ la base correspondante. Soit $P=MU$ un sous-groupe parabolique standard, on notera alors $\Sigma(P)$ l'ensemble des racines de $A_{min}$ dans $U$ et $\Delta(P)=\Delta\cap \Sigma(P)$. L'application $P\mapsto \Delta(P)$ est une bijection entre les sous-groupes paraboliques standards et les sous-ensembles de $\Delta$. On désignera par $A_M^+$ la chambre positive correspondant à $P$ c'est-à-dire l'ensemble des éléments $a$ de $A_M(F)$ vérifiant $|\alpha(a)|_F\leqslant 1$ pour tout $\alpha\in \Delta(P)$. Pour $\epsilon>0$, on pose

$$A_M(\epsilon)^+=\{a\in A_M(F); |\alpha(a)|_F<\epsilon \;\forall \alpha\in\Delta(P)\}$$

\noindent Soit $w_l\in W^G$ l'élément de plus grande longueur i.e. celui tel que $w_lP_{min}w_l^{-1}=\overline{P}_{min}$ où $\overline{P}_{min}$ est le sous-groupe parabolique opposé à $P_{min}$.
On désignera par $\overline{U}_{min}$ le radical unipotent de $\overline{P}_{min}$. Soit $\mathcal{D}(U_{min})$ le sous-groupe dérivé de $U_{min}$. On a un isomorphisme de groupes algébriques $U_{min}/\mathcal{D}(U_{min})\simeq \bigoplus_{\alpha\in \Delta} V_\alpha$, où pour tout $\alpha\in \Delta$, $V_\alpha$ est le sous-espace propre correspondant à $\alpha$ pour l'action par conjugaison de $A_{min}$. De plus, on a l'égalité $U_{min}(F)/\mathcal{D}(U_{min})(F)=(U_{min}/\mathcal{D}(U_{min}))(F)\simeq \bigoplus_{\alpha\in\Delta} V_\alpha(F)$. Fixons des formes linéaires non nulles $l_\alpha$ sur chacun des espaces vectoriels $V_\alpha(F)$. On définit un caractère de $U_{min}(F)/\mathcal{D}(U_{min})(F)$ via l'isomorphisme précédent par la formule

$$(x_\alpha)_{\alpha\in \Delta}\mapsto \psi(\sum_\Delta l_\alpha(x_\alpha))$$

On notera aussi $\psi$ ce caractère et on l'identifiera à un caractère de $U_{min}(F)$. On fait agir $F^\Delta$ sur $U_{min}(F)/\mathcal{D}(U_{min})(F)$ par $y.(x_\alpha)=(y_\alpha x_\alpha)$ où $y=(y_\alpha)_{\alpha\in\Delta}$ est un élément de $F^\Delta$. Soit $\mathcal{K}^G_{P_{min}}$ l'espace des fonctions $\phi: K\to \mathbb{C}$ qui sont localement constantes et invariantes à gauche par $K\cap P_{min}(F)$, c'est un espace vectoriel sur lequel se réalisent toutes les représentations $\pi_s=i^G_{P_{min}}(\delta_{P_{min}}^s)$ pour $s\in \mathbb{C}$ via l'isomorphisme donné par la restriction des fonctions à $K$. Soit $(,)$ la forme bilinéaire non dégénérée sur $\mathcal{K}^G_{P_{min}}$ définie par

$$(e_1,e_2)=\displaystyle\int_K e_1(k)e_2(k) dk$$

Soit $e\in\mathcal{K}^G_{P_{min}}$ la fonction constante égale à 1 et pour tout $s$, soit $e_s$ la fonction sur $G(F)$ déduite via l'isomorphisme précédent entre $i^G_{P_{min}}(\delta_{P_{min}}^s)$ et $\mathcal{K}^G_{P_{min}}$. Pour $y\in F^\Delta$, $g\in G(F)$ et $s\in \mathbb{C}$, on pose

$$\displaystyle\Phi_{y,s}(g)=\int_{U_{min}(F)}e_s(w_lug)\psi(y.u) du$$

L'intégrale est absolument convergente pour $Re(s)>0$. Pour $c\geqslant 0$ un réel, notons $U_{min}(F)_c$ le sous-groupe de $U_{min}(F)$ constitué des éléments $u$ tels que $val_F(l_\alpha(u))\geqslant -c$ pour tout $\alpha\in \Delta$. Par transformée de Fourier inverse, on a

$$\displaystyle\int_{\mathcal{O}_F^\Delta} \Phi_{y,s}(1) dy=\int_{U_{min}(F)_0} e_s(w_lu) du$$

Pour $a\in A_{min}^+$ on note $y(a)$ le $\Delta$-uplet $(\alpha(a))_{\alpha\in\Delta}$ d'éléments de $(\mathcal{O}_F\backslash \{0\})^\Delta$. L'application $a\mapsto y(a)$ est un isomorphisme de $A_{min}^+/(Z_G(F)\cap A_{min}(F))$ sur un sous-monoïde ouvert $\mathcal{M}$ de $(\mathcal{O}_F\backslash \{0\})^\Delta$. Il existe un nombre fini d'éléments $y_i$ tels que $(\mathcal{O}_F\backslash \{0\})^\Delta=\bigsqcup_i y_i \mathcal{M}$. On a par conséquent

$$\displaystyle\int_{\mathcal{O}_F^\Delta} \Phi_{y,s}(1) dy=\sum_i \int_{y_i \mathcal{M}}  \Phi_{y,s}(1) dy$$

Le changement de variable $A_{min}^+/(A_{min}(F)\cap Z_G(F))\to y_i \mathcal{M}$, $a\mapsto y_iy(a)$ donne
$$\displaystyle\int_{y_i \mathcal{M}}  \Phi_{y,s}(1) dy=|y_i|_F\int_{A_{min}^+/(Z_G(F)\cap A_{min}(F))} \Phi_{y_iy(a),s}(1) \prod_{\Delta} |\alpha(a)|_F da$$

Par le changement de variable $u\mapsto aua^{-1}$ on a $\Phi_{y_iy(a),s}(1)=\delta_{P_{min}}(a)^{s-1/2}\Phi_{y_i,s}(a)$. On en déduit

$$\mbox{(1)} \;\;\; \int_{U_{min}(F)_0} e_s(w_lu) du=\sum_i |y_i|_F\int_{A_{min}^+/(A_{min}(F)\cap Z_G(F))} \delta_{P_{min}}(a)^{s-1/2} \Phi_{y_i,s}(a) \prod_{\Delta} |\alpha(a)|_F da$$

\begin{lem}
Soient $e_1,e_2\in\mathcal{K}^G_{P_{min}}$, $P=MU$ un sous-groupe parabolique standard de $G$. Il existe $\epsilon>0$ et un sous-ensemble discret $E\subset \mathbb{C}$ tels que, pour tout $\gamma\in A_{min}^+$, il existe des fonctions holomorphes $\varphi_w^\gamma:\mathbb{C}-E\to\mathbb{C}$ indexées par $w\in W^G/W^M$ telles que, pour tout $s\in\mathbb{C}-E$ et tout $a_M\in A_M(\epsilon)^+$, on ait
$$(\pi_s(a_M\gamma)e_1,e_2)=\displaystyle\sum_{w\in W^G/W^M} \delta_{P_{min}}(a_M)^{1/2}\delta_{P_{min}}(w.a_M)^s \varphi_w^\gamma(s)$$
\end{lem}

\ul{Preuve}: Soit $\mathcal{X}$ le tore complexe des caractères non ramifiés de $M_{min}(F)$ et \\
\noindent $B=\mathbb{C}[M_{min}(F)/M_{min}^1]$ la $\mathbb{C}$-algèbre des fonctions régulières sur $\mathcal{X}$. Soit $\chi_{nr}:M_{min}(F)\to B^\times$ le caractère non ramifié générique qui à $m_0$ associe la fonction régulière $\chi\mapsto \chi(m_0)$. On pose $(\pi_B,V_B)=i_{P_{min}}^G(\chi_{nr})$, c'est une $(G,B)$-représentation lisse et admissible au sens de [BDKV]. Son dual lisse est la représentation induite $(\pi_B^\vee,V_B^\vee)=i^G_{P_{min}}(\chi_{nr}^{-1})$. Pour tout sous-groupe parabolique $Q$ de $G$, on notera $(\pi_{B,Q},(V_B)_Q)$ le module de Jacquet de cette représentation relativement à $Q$ et $j_Q: V_B\to (V_B)_Q$ la projection naturelle. Un théorème de Casselman affirme l'existence d'une dualité $(V_B)_P\times (V_B^\vee)_{\overline{P}}\to B$ qui est $B$-bilinéaire et $G$-équivariante vérifiant la condition suivante:

\begin{center}
Pour $f\in V_B$ et $f^\vee\in V_B^\vee$, il existe $\epsilon>0$ tel que pour tout $a\in A_{min}^+$ vérifiant $|\alpha(a)|_F<\epsilon$ pour tout $\alpha\in \Delta(P)$ on ait
$$<\pi_B(a)f,f^\vee>=\delta_P(a)^{1/2} <\pi_{B,P}(a)j_P(f), j_{\overline{P}}(f^\vee)>$$
\end{center}

La représentation $\pi_{B,P}$ admet une filtration indexée par $w\in W^G/W^M$ pour laquelle les quotients successifs admettent pour caractère central les $w^{-1}\chi_{nr}$. Puisque ces caractères sont tous différents, la filtration se scinde en une somme directe après extension des scalaires à $Frac(B)$. On a donc une décomposition $j_P(f)=\displaystyle\sum_{w\in W^G/W^M} j_P(f)_w$ où $j_P(f)_w\in Frac(B)\otimes_B (V_B)_P$ pour tout $w$ et $A_M(F)$ agit par $w^{-1}\chi_{nr}$ sur $j_P(f)_w$. On obtient pour $a=a_M\gamma$ avec $a_M\in A_M(\epsilon)^+$ et $\gamma\in A_{min}^+$, 

$$<\pi_B(a)f,f^\vee>=\delta_P(a_M\gamma)^{1/2}\displaystyle\sum_{w\in W^G/W^M} \chi_{nr}(w.a_M) <\pi_{B,P}(\gamma)j_P(f)_w,j_{\overline{P}}(f^\vee)>$$

Par restriction des fonctions à $K$, les espaces $V_B$ et $V_B^\vee$ sont tout deux isomorphes à l'espace $\mathcal{K}^G_{P_{min}}\otimes B$ des fonctions $\varphi: K\to B$ localement constantes et invariantes à gauche par $P_{min}(F)\cap K$. Les éléments $e_1$ et $e_2$ peuvent être vus comme des éléments de cet espace (puisque $\mathbb{C}\subset B$). On notera encore par $e_1$ (resp. $e_2$) l'image réciproque de $e_1$ dans $V_B$ (resp. $V_B^\vee$). En appliquant ce qui précède à $f=e_1$ et $f^\vee=e_2$ on obtient le résultat par spécialisation de $\chi_{nr}$ à $\delta_{P_{min}}^s$. $\blacksquare$

\vspace{4mm}

Il existe des entiers naturels strictement positifs $k_\alpha$ pour $\alpha\in\Delta$ tels que $\delta_{P_{min}}(a)=\displaystyle\prod_{\alpha\in\Delta} |\alpha(a)|_F^{k_\alpha}$ pour tout $a\in A_{min}(F)$. Soit $k=sup(k_\alpha)$, on a alors pour tout $a\in A_{min}^+$
$$\displaystyle\prod_{\alpha\in\Delta} |\alpha(a)|_F\leqslant \delta_{P_{min}}(a)^{1/k}$$
Soient $e_1,e_2\in\mathcal{K}^G_{P_{min}}$ et $s$ un nombre complexe, on pose
$$F_{e_1,e_2}(s)=\displaystyle\int_{A_{min}^+/A_G(F)}(\pi_s(a)e_1,e_2)\delta_{P_{min}}(a)^{s-1/2} \prod_{\alpha\in \Delta} |\alpha(a)|_F da$$

\begin{lem}
Pour $Re(s)>-\frac{1}{2k}$ l'intégrale définissant $F_{e_1,e_2}(s)$ est absolument convergente et la fonction $s\mapsto F_{e_1,e_2}(s)$ est holomorphe sur le demiplan $Re(s)>-\frac{1}{2k}$.
\end{lem}

\ul{Preuve}: Soit $K_1$ un sous-groupe ouvert de $A_{min}^1$ qui laisse stable $e_1$. Pour $1>\epsilon>0$, il existe des sous-ensembles finis $\Gamma_P^\epsilon\subset A_{min}^+$ pour $P$ parcourant les sous-groupes paraboliques standards tels que
$$A_{min}^+=\bigsqcup_{P_{min}\subset P=MU}\bigsqcup_{\gamma\in \Gamma_P^\epsilon} A_M(\epsilon)^+\gamma K_1$$
Il suffit donc de montrer que pour $P=MU$ un sous-groupe parabolique standard et $\gamma\in \Gamma_P^\epsilon$, l'intégrale à un paramètre
$$\displaystyle\int_{A_M(\epsilon)^+/A_G(F)}(\pi_s(a_M\gamma)e_1,e_2)\delta_{P_{min}}(a_M)^{s-1/2} \prod_{\alpha\in \Delta} |\alpha(a_M)|_F da_M$$
est convergente pour $Re(s)>-\frac{1}{2k}$ et qu'elle définit une fonction holomorphe sur ce demi-plan. Choisissons $\epsilon$ de sorte que la conclusion du lemme 2.3.1 soit vérifiée pour tout les sous-groupes paraboliques standards. Le lemme 2.3.1 fournit des fonctions holomorphes $\phi^\gamma_w$ telles que pour $s$ en dehors d'un  ensemble discret et tout $a_M\in A_M(\epsilon)^+$,

$$ (\pi_s(a_M\gamma)e_1,e_2)\delta_{P_{min}}(a_M)^{s-1/2}=\displaystyle\sum_{w\in W^G/W^M} \delta_{P_{min}}(a_M(w.a_M))^s \phi_w^\gamma(s)$$
On va appliquer le lemme 2.2.1 à notre intégrale à paramètre. Les points 1 et 2 sont aisément vérifiés. Soit $C$ un sous-ensemble compact de $\{Re(s)>-\frac{1}{2k}\}$ sur lequel les fonctions holomorphes $\phi_w^\gamma$ sont bien définies. Soit $r_-=inf_{s\in C} Re(s)>-\frac{1}{2k}$ et $r_+=sup_{s\in C} Re(s)>-\frac{1}{2k}$. Pour tous $w\in W^G/W^M$, $a_M\in A_M^+$ et pour tout $s\in C$ on a

$$\mbox{(2)}\;\;\; |\delta_{P_{min}}(a_M(w.a_M))^s|\leqslant \delta_{P_{min}}(a_M(w.a_M))^{r_+}+\delta_{P_{min}}(a_M(w.a_M))^{r_-}$$

Il existe une constante $c>0$ telle que pour tout $w\in W^G/W^M$

$$sup_{s\in C} |\phi_w^\gamma(s)|\leqslant c $$

En utilisant la majoration (2), on obtient que la fonction

$$s\mapsto (\pi_s(a_M)e_1,e_2)\delta_{P_{min}}(a_M)^{s-1/2} \prod_{\alpha\in \Delta} |\alpha(a_M)|_F$$
est majorée sur $C$ par

$$c\displaystyle\sum_{w\in W^G/W^M} (\delta_{P_{min}}(a_M(w.a_M))^{r_+}+\delta_{P_{min}}(a_M(w.a_M))^{r_-})\delta_{P_{min}}(a_M)^{1/k}$$

Pour tous $w\in W^G$ et $a\in A_{min}^+$ on a $\delta_{P_{min}}(a)^2\leqslant \delta_{P_{min}}(a(w.a))\leqslant 1$. Puisque la fonction $a_M\mapsto \delta_{P_{min}}(a_M)^r$ est intégrable sur $A_M(\epsilon)^+/A_G(F)$ pour $r$ réel strictement positif, la fonction précédente est bien intégrable sur $A_M(\epsilon)^+/A_G(F)$. On peut donc appliquer le lemme 2.2.1 $\blacksquare$

\begin{lem}
Soit $y\in F^\Delta$, la fonction 
$$s\mapsto \displaystyle\int_{A_{min}^+/Z_G(F)\cap A_{min}(F)} \delta_{P_{min}}(a)^{s-1/2}\Phi_{y,s}(a)\prod_{\alpha\in \Delta} |\alpha(a)|_F da$$
bien définie et holomorphe pour $Re(s)>0$ admet un prolongement holomorphe à $Re(s)>-\frac{1}{2k}$.
\end{lem}

\ul{Preuve}: Pour $Re(s)>0$, la forme linéaire $\Omega_s: i_{P_{min}}^G(\delta_{P_{min}}^s)\to \mathbb{C}$, $\varphi \mapsto \displaystyle\int_{U_{min}(F)} \varphi(w_l u)\psi(y.u) du$ est une fonctionnelle de Whittaker pour le caractère générique $\xi: u\mapsto \psi(y.u^{-1})$. On a $\Phi_{y,s}(a)=\Omega_s(\pi_s(a)e_s)$. Il existe un sous-groupe compact-ouvert $K_1$ de $G(F)$ tel que pour tout $a\in A_{min}^+$ on a $K_1\subset Ker(\xi)aKa^{-1}$. Soit $\mathcal{B}$ un base de $(\mathcal{K}^G_{P_{min}})^{K_1}$ et $\mathcal{B}^{*}=\{e'^*; e'\in\mathcal{B}\}$ la base duale pour $(,)$. On a alors
$$\Omega_s(\pi_s(a)e_s)=\displaystyle\sum_{e'\in\mathcal{B}} (\pi_s(a)e,e'^*)\Omega_s(e')$$
Il est bien connu que $\Omega_s(e')$ qui est la fonctionnelle de Jacquet admet un prolongement holomorphe au plan complexe (cf par exemple [CS]). Le résultat du lemme est alors une conséquence du lemme 2.3.2 $\blacksquare$

\subsection{Majorations de mesures}

Pour $c\geqslant 0$ un réel, on pose $\overline{U}_{min}(F)_c=w_lU_{min}(F)_c w_l^{-1}$. \\
\noindent Soient $A_n=\{\overline{u}\in\overline{U}_{min}(F)_0; \delta_{P_{min}}(m_{P_{min}}(\overline{u}))\geqslant q^{-n}\}$ et $a_n=mes(A_n)$.
\begin{lem}
Il existe $\epsilon>0$ tel que
$$a_n<<q^{n(1/2-\epsilon)}$$
\end{lem}

\ul{Preuve}: Pour $Re(s)>0$, on a
$$\displaystyle\int_{U_{min}(F)_0} e_s(w_lu) du=\sum_{n\geqslant 0} (a_n-a_{n-1}) q^{-n/2-ns}$$
D'après le lemme 2.3.3 et 2.3(1), cette fonction admet un prolongement holomorphe à $Re(s)>-\frac{1}{2k}$. Par conséquent le rayon de convergence de la série entière
$\displaystyle\sum_n (a_n-a_{n-1}) z^n$ est strictement plus grand que $q^{-1/2}$ d'où le résultat $\blacksquare$

\vspace{2mm}

Soit $\overline{Q}=M\overline{U}$ un sous-groupe parabolique antistandard. Pour $c\geqslant 0$ on pose $\overline{U}(F)_c=\overline{U}(F)\cap\overline{U}_{min}(F)_c$. On définit $A_{n,c}(\overline{U})=\{\overline{u}\in\overline{U}(F)_c;\; \delta_{P_{min}}(m_{P_{min}}(\overline{u}))\geqslant q^{-n}\}$ pour tout $n\in \mathbb{N}$ et pour tout réel $c\geqslant 0$.

\begin{prop}
Il existe deux réels $\epsilon>0$ et $\alpha>0$ tels que

$$mes(A_{n,c}(\overline{U}))<<q^{n(1/2-\epsilon)+\alpha c}$$
pour tout $n\in\mathbb{N}$ et pour tout $c\geqslant 0$.
\end{prop}

\ul{Preuve}: On a les faits suivants

\vspace{2mm}

\noindent(1) Il existe un réel $\alpha_1>0$ tel que pour tout $c\geqslant 0$, il existe $a_c\in A_M(F)$ qui vérifie \\
    (i) $\sigma(a_c)\leqslant \alpha_1 (c+1)$ \\
    (ii) $\overline{U}(F)_c\subset a_c\overline{U}(F)_0a_c^{-1}$
    
\vspace{2mm}

\noindent(2) Il existe un réel $\alpha_2>0$ tel que \\
     (i) Pour tout $g\in G(F)$, $\delta_{P_{min}}(m_{P_{min}}(g))\leqslant exp(\alpha_2 \sigma(g))$; \\
     (ii) Pour tout $a\in A_M(F)$, $\delta_{\overline{Q}}(a)\leqslant exp(\alpha_2 \sigma(a))$.
  
\vspace{2mm}

Soient $c\geqslant 0$ et $a_c\in A_M(F)$ qui vérifie (1)(i) et (ii). On a alors 

$$m_{P_{min}}(a_cua_c^{-1})=a_cm_{P_{min}}(u)m_{P_{min}}(k_{P_{min}}(u)a_c^{-1})$$
pour tout $u\in \overline{U}(F)_0$. On peut toujours supposer que la fonction $\sigma$ est invariante à gauche par $K$. On a alors

$$\delta_{P_{min}}(m_{P_{min}}(a_cua_c^{-1}))\leqslant \delta_{P_{min}}(m_{P_{min}}(u))exp(2\alpha_1\alpha_2(c+1))$$
pour tout $u\in \overline{U}(F)_0$. Soit $k_c=E\left(2\alpha_1\alpha_2(c+1)/log(q)\right)$ où $E(.)$ désigne la partie entière. On déduit de ce qui précède que

$$A_{n,c}(\overline{U})\subset a_cA_{n+k_c,0}(\overline{U})a_c^{-1}$$
pour tout $n\in\mathbb{N}$. D'où

$$\mbox{(3)} \;\;\; mes(A_{n,c}(\overline{U}))\leqslant \delta_{\overline{Q}}(a_c)mes\left(A_{n+k_c,0}(\overline{U})\right)\leqslant exp(\alpha_1\alpha_2(1+c)) mes\left(A_{n+k_c,0}(\overline{U})\right)$$

\noindent pour tout $n\in\mathbb{N}$ et pour tout $c\geqslant 0$. \\

 Soit $K_{\overline{U}_{min}}=K\cap\overline{U}_{min}(F)$, on a alors $A_{n,0}(\overline{U})K_{\overline{U}_{min}}\subset A_n$ pour tout $n\in\mathbb{N}$. On en déduit que
$$mes(A_{n,0}(\overline{U}))mes((\overline{U}(F)\cap K_{\overline{U}_{min}})\backslash K_{\overline{U}_{min}})\leqslant a_n$$
pour tout $n\in\mathbb{N}$. Puisque $(\overline{U}(F)\cap K_{\overline{U}_{min}})\backslash K_{\overline{U}_{min}}$ est compact donc de mesure finie, le lemme 2.4.1 et l'inégalité (3) permettent d'obtenir le résultat de la proposition $\blacksquare$

\begin{cor}
Il existe deux réels $\epsilon>0$ et $\alpha>0$ tels que

$$mes\{u\in \overline{U}(F)_c;\; q^{-n-1}<\Xi^M(m_Q(u))\delta_Q(m_Q(u))^{1/2}\leqslant q^{-n}\} <<q^{n(1-\epsilon)+\alpha c}$$
pour tout $n\in\mathbb{N}$ et pour tout $c\geqslant 0$.
\end{cor}

\ul{Preuve}: D'après le lemme II.3.2 de [W3], il existe un réel $D$ tel que

$$\Xi^M(m_Q(u))\delta_Q(m_Q(u))^{1/2}<<\delta_{P_{min}}(m_{P_{min}}(u))^{1/2} \sigma(u)^D$$
pour tout $u\in \overline{U}(F)$. D'après le lemme II.3.4 de [W3] on a aussi

$$\sigma(u)<<1-log \left(\delta_{P_{min}}(m_{P_{min}}(u))\right)$$
pour tout $u\in \overline{U}(F)$. Le résultat est alors une conséquence facile de la proposition précédente $\blacksquare$

\section{Les groupes unitaires}

\subsection{Généralités}

Dorénavant on fixe une extension quadratique $E$ de $F$. Les notations $\mathcal{O}_E$, $k_E$, $q_E$, $val_E$... seront les analogues de celles définies pour $F$. On désignera par $N$ et $Tr$ respectivement la norme et la trace de l'extension $E/F$. On définit un caractère additif $\psi_E$ de $E$ par $\psi_E(x)=\psi(Tr(x))$. On notera $x\mapsto \overline{x}$ le $F$-automorphisme non trivial de $E$ et $\chi_E$ le caractère quadratique de $F^\times$ donné par 

$$\chi_E:F^\times\to F^\times/N(E^\times)=\{\pm 1\}$$

Si $G$ est un groupe algébrique définie sur $E$, on notera $R_{E/F}G$ le groupe algébrique défini sur $F$ obtenu par restriction des scalaires à la Weil de $E$ à $F$.

\vspace{2mm} 
 
 Un espace hermitien sera pour nous un couple $(V,h)$ constitué d'un espace vectoriel de dimension finie $V$ sur $E$ et d'une forme hermitienne $h:V\times V\to E$ non dégénérée (avec la convention que $h$ est linéaire en la première variable). Quand la forme hermitienne sur $V$ est évidente on dira juste que $V$ est un espace hermitien (c'est le cas par exemple pour un sous-espace d'un espace hermitien sur lequel la restriction de la forme hermitienne est non dégénérée).  
Soit $(V,h)$ un espace hermitien et soit $G$ son groupe unitaire. On appelle système hyperbolique de $V$ toute famille $(v_i)_{i=\pm 1,\ldots,\pm n}$ d'éléments de $V$ qui vérifient $h(v_i,v_j)=\delta_{i,-j}$ pour tous $i,j\in\{\pm 1,\ldots,\pm n\}$. Si $V$ admet un système hyperbolique qui est aussi une base, on dira que $V$ est hyperbolique. Soit $V_{hyp}$ un sous-espace hermitien hyperbolique maximal de $V$ et $V_{an}$ son orthogonal. Alors la forme hermitienne sur $V_{an}$ est anisotrope. Posons $d_{an}(V)=dim(V_{an})$. Notons $d(V)$ la dimension de $V$ (comme $E$-espace vectoriel). On a toujours $d(V)\equiv d_{an}(V) \; [2]$ et $d_{an}(V)\leqslant 2$. Le groupe $G$ est quasi-déployé si et seulement si $d_{an}(V)\leqslant 1$.

Pour $v',v''\in V$, on note $c(v',v'')$ l'endomorphisme de $V$ défini par

$$c(v',v'')(v)=h(v,v')v''-h(v,v'')v'$$

\noindent pour tout $v\in V$. C'est un élément de $\mathfrak{g}(F)$ et les $c(v',v'')$ pour $v',v''\in V$ engendrent $\mathfrak{g}(F)$ comme $F$-espace vectoriel.

On définit une fonction $\Delta$ sur $G_{ss}(F)$ de la façon suivante. Pour $x\in G_{ss}(F)$, notons $V''(x)$ le noyau de $x-1$ dans $V$ et posons

$$\Delta(x)=|N(det(1-x)_{|V/V''(x)})|_F$$

\noindent Soit $W$ un sous-espace non dégénérée de $V$ et notons $H$ son groupe unitaire. On peut considérer $H$ comme un sous-groupe de $G$ en laissant agir ses éléments trivialement sur l'orthogonal de $W$.

\begin{lem}
Posons $k=d(V)-d(W)$. On a alors

$$D^G(x)=D^H(x)\Delta(x)^k$$

\noindent pour tout $x\in H_{ss}(F)$.
\end{lem}

\ul{Preuve}: Posons $W'=W^\perp$, $H'$ le groupe unitaire de $W'$ et 

$$\mathfrak{g}(W,W')=\{X\in\mathfrak{g};\; X.W\subset W',\; X.W'\subset W\}$$

\noindent On a alors deux décompositions $\mathfrak{g}=\mathfrak{h}'\oplus\mathfrak{h}\oplus \mathfrak{g}(W,W')$ et $\mathfrak{g}_x=\mathfrak{h}_x\oplus \mathfrak{g}(W,W')_x$ où $\mathfrak{g}(W,W')_x$ désigne le commutant de $x$ dans $\mathfrak{g}(W,W')$. Par définition, on a donc

$$D^G(x)=D^H(x)|det(1-ad(x))_{|\mathfrak{g}(W,W')/\mathfrak{g}(W,W')_x}|_F$$

L'application $\mathfrak{g}(W,W')\to Hom_E(W',W)$, $X\mapsto X_{|W'}$ est un isomorphisme (où $X_{|W'}$ désigne la restriction de $X$ à $W'$). Via cet isomorphisme $ad(x)$ agit sur $\mathfrak{g}(W,W')$ par $X\mapsto x\circ X$. On a donc un isomorphisme de $F[x]$-modules $\mathfrak{g}(W,W')\simeq W^{\oplus k}$ où l'action de $x$ sur chaque facteur est l'action naturelle de $x$ sur $W$. On en déduit facilement que

$$|det(1-ad(x))_{|\mathfrak{g}(W,W')/\mathfrak{g}(W,W')_x}|_F=\Delta(x)^k$$

$\blacksquare$

\subsection{Sous-groupes compacts spéciaux, paraboliques, Levi et R-groupes}

Soit $(v_i)_{i=\pm 1,\ldots,\pm r}$ un système hyperbolique maximal de $V$. Soient $k_1,\ldots,k_s\geqslant 1$ des entiers vérifiant $k_1+\ldots+k_s\leqslant r$. Notons alors $Z_i$ (resp. $Z_{-i}$), pour $i=1,\ldots,s$, le sous-espace de $V$ engendré par les $v_j$ (resp. par les $v_{-j}$) pour $j=(k_1+\ldots+k_{i-1}+1),\ldots,(k_1+\ldots+k_i)$. Soit $P$ le stabilisateur dans $G$ du drapeau

$$Z_1\subset Z_1\oplus Z_2\subset\ldots\subset Z_1\oplus\ldots\oplus Z_s$$

C'est un sous-groupe parabolique de $G$. Le sous-groupe $M$ de $P$ des éléments qui préservent $Z_{\pm i}$, $i=\pm 1,\ldots,\pm s$ est une composante de Levi de $P$. 

Tous les couples $(P,M)$ formés d'un sous-groupe parabolique et d'une composante de Levi de $P$ arrivent de cette façon, c'est-à-dire qu'on peut trouver un système hyperbolique maximal et des entiers $k_1,\ldots,k_s$ tels que $P$ et $M$ soient déterminés comme ci-dessus. \\

 Supposons à nouveau fixés notre système hyperbolique maximal et nos entiers $k_1,\ldots,k_s$ cela détermine un sous-groupe parabolique $P$ et un Levi $M$. Soit $\tilde{V}$ l'orthogonal de $Z_1\oplus Z_{-1}\oplus\ldots\oplus Z_s\oplus Z_{-s}$ et notons $\tilde{G}$ son groupe unitaire. On a alors un isomorphisme naturel

$$M\simeq R_{E/F} GL(Z_1)\times\ldots\times R_{E/F} GL(Z_s)\times \tilde{G}$$

Le groupe de Weyl $W(M)$ s'identifie alors naturellement à un sous-groupe de $\mathcal{S}_s\ltimes (\mathbb{Z}/2\mathbb{Z})^s$ où $\mathcal{S}_s$ agit sur $(\mathbb{Z}/2\mathbb{Z})^s$ par permutation sur les entrées. Soit $\tau$ une représentation irréductible de la série discrète de $M(F)$. Alors $R(\tau)$ est un sous-groupe de $(\mathbb{Z}/2\mathbb{Z})^s$ et $R(\tau)_{reg}\neq \emptyset$ si et seulement si $R(\tau)$ contient l'élément $t=(-1,\ldots,-1)$. Dans ce cas, $R(\tau)=(\mathbb{Z}/2\mathbb{Z})^s$, $\tau$ se prolonge en une représentation de $Norm_{G(F)}(\tau)$, on a $R(\tau)_{reg}=\{t\}$ et $|det(t-1_{|\mathcal{A}_M})|=2^{a_M}$. \\

Soit $\pi$ une représentation irréductible et elliptique de $G(F)$. On peut trouver $M$, $\tau$ comme ci-dessus, et $\zeta\in R(\tau)^{\vee}$ de sorte que $\pi=Ind_{P}^G(\tau,\zeta)$, où $P$ est un élément de ${\cal P}(M)$. Puisque la classe de conjugaison du couple $(M,\tau)$ est bien déterminée, on peut poser $t(\pi)=2^{a_{M}}$. \\

Plus généralement, soit $M$ comme précédemment et soit $L$ un sous-groupe de Levi de $M$. Alors $L$ admet une décomposition

$$L=L_1\times\ldots\times L_s\times\tilde{L}$$

\noindent où pour $i=1,\ldots,s$, $L_i$ est un sous-groupe de Levi de $R_{E/F}GL(Z_i)$ et $\tilde{L}$ est un sous-groupe de Levi de $\tilde{G}$. Soit $\tau\in\Pi_2(L)$, on a alors une décomposition en produit tensoriel $\tau\simeq \tau_1\otimes\ldots\otimes \tau_s\otimes \tilde{\tau}$ où $\tau_i\in \Pi_2(L_i)$ pour $i=1,\ldots,s$ et $\tilde{\tau}\in\Pi_2(\tilde{L})$. La théorie du R-groupe est triviale pour les $R_{E/F}GL(Z_i)$ donc $R^M(\tau)$ s'identifie à $R^{\tilde{G}}(\tilde{\tau})$. Il s'ensuit que $R^M(\tau)_{reg}$ est vide sauf si $L_i=R_{E/F}GL(Z_i)$ pour $i=1,\ldots,s$ et $R^{\tilde{G}}(\tilde{\tau})_{reg}\neq\emptyset$. Dans ce cas là $R^M(\tau)_{reg}$ est réduit à un élément $t$ et on a $|det(t-1)_{|\mathcal{A}_M/\mathcal{A}_L}|=|R^M(\tau)|=2^{a_M-a_L}$. \\

Il existe un ensemble de $\mathcal{O}_E$-réseaux de $V$ qualifiés de spéciaux tel que pour tout réseau spécial $R$ de $V$ le stabilisateur de $R$ dans $G(F)$ soit un sous-groupe compact spécial. Reprenons les constructions précédentes dans le cas où $k_1=\ldots=k_r=1$. Alors $P$ est un sous-groupe parabolique minimal. On peut choisir de prendre $P_{min}=P$ et $M_{min}=M$.  Il existe des vecteurs $v'_i$ non nuls proportionnels à $v_i$ pour $i=\pm 1,\ldots,\pm r$  vérifiant $h(v'_i,v'_{-i})=h(v'_{i'},v'_{-i'})$ pour tout $i,i'=1,\ldots,r$ et un réseau spécial $\tilde{R}\subset \tilde{V}$ tel que le réseau $R=R_Z\oplus \tilde{R}$, où $R_Z=\mathcal{O}_Ev'_1\oplus\mathcal{O}_Ev'_{-1}\oplus\ldots\oplus \mathcal{O}_Ev'_r\oplus\mathcal{O}_E v'_{-r}$, soit spécial. Le stabilisateur $K$ de $R$ dans $G(F)$ qui est un sous-groupe compact spécial est alors en bonne position par rapport à $M_{min}$. Il est utile de remarquer (c'est une des hypothèses du paragraphe 1.6) qu'alors le morphisme naturel $Norm_{G(F)}(M_{min})\cap K\to W^G$ admet une section. Notons $W^K$ le sous-groupe des éléments $k\in K$ qui agissent par permutation sur les $(v_i)_{i=\pm 1,\ldots,\pm r}$ et qui agissent comme l'identité sur $\tilde{V}$. Alors la restriction du morphisme précédent à $W^K$ est un isomorphisme, ce qui montre l'existence d'une section.

\subsection{Orbites nilpotentes régulières}

Supposons $G$ quasidéployé et soit $(v_{\pm i})_{i=1,\ldots,r}$ un système hyperbolique maximal. On a alors $dim(V)=2r$ ou $2r+1$. Si $dim(V)=2r+1$ alors $\mathfrak{g}(F)$ ne contient qu'une seule orbite nilpotente régulière. Supposons que $dim(V)=2r$. Considérons $P$ et $M$ le sous-groupe parabolique et le Levi associé aux entiers $k_1=\ldots=k_r=1$. Alors $P$ est un sous-groupe parabolique minimal. Notons $U$ son radical unipotent et $\mathfrak{u}$ son algèbre de Lie. Soit $\eta\in Ker(Tr)-\{0\}$ un élément de $E$ de trace nulle qui est non nul. Soit $N_\eta$ l'élément de $\mathfrak{u}(F)$ définit par

\begin{center}
$N_\eta v_r=0, N_\eta v_i=v_{i+1}$ pour $i=r-1,\ldots,1,-2,\ldots,-r$ et $N_\eta v_{-1}=\eta v_1$
\end{center}

Alors $N_\eta$ est un élément régulier qui engendre une orbite nilpotente $\mathcal{O}_\eta$ qui ne dépend que de l'image de $\eta$ dans $(Ker(Tr)-\{0\})/N(E^\times)$. L'application $\eta\mapsto \mathcal{O}_\eta$ est une bijection de $(Ker(Tr)-\{0\})/N(E^\times)$ sur $Nil(\mathfrak{g})(F)_{reg}$, en particulier cet ensemble est de cardinal 2 et la multiplication par n'importe quel élément $\lambda\in F^\times-N(E^\times)$ permute ces deux orbites.

\section{Position du problème}

 Soit $(V,h)$ un espace hermitien et $G$ le groupe unitaire associé. On se donne une décomposition orthogonale de $V$ de la forme $V=Z\oplus D\oplus W$ où $D$ est une droite non dégénérée, dont on fixera un générateur $v_0$, et $Z$ possède une base hyperbolique $(v_i)_{i=\pm1,\ldots,\pm r}$: on a donc $h(v_i,v_j)=\delta_{i,-j}$ pour tous $i,j\in\{\pm 1,\ldots, \pm r\}$. On notera $Z_+$ resp. $Z_-$ les sous-espaces de $V$ engendrés par les $v_i$, $i=1,\ldots,r$ resp. les $v_{-i}$, $i=1,\ldots,r$ et $V_0=D\oplus W$. Soient $G$, $G_0$ et $H$ les groupes unitaires respectifs de $V$, $V_0$ et $W$. On identifie $H$ (resp. $G_0$) avec le sous-groupe de $G$ des éléments qui agissent comme l'identité sur $Z\oplus D$ (resp. $Z$). Soit $P$ le sous-groupe parabolique de $G$ qui conserve le drapeau de sous-espaces isotropes

$$Ev_r\subset Ev_r+Ev_{r-1}\subset\ldots\subset Ev_r+\ldots+Ev_1$$

Soit $A$ le tore des éléments de $G$ qui conservent chaque droite $Ev_i$, $i=\pm 1,\ldots,\pm r$ et qui agissent comme l'identité sur $D\oplus W$, $U$ le radical unipotent de $P$. Posons $M=AG_0$, c'est une composante de Levi de $P$. Pour $a\in A(F)$ on note $a_i$ la valeur propre de $a$ agissant sur $v_i$. Soient $\xi_0,\ldots,\xi_{r-1}$ des éléments de $F^{\times}$, la formule suivante défini alors un caractère de $U(F)$ invariant par conjugaison par $H(F)$

$$\xi(u)=\psi_E(\displaystyle\sum_{i=0}^{r-1} \xi_i h(uv_i,v_{-i-1}))$$

Pour $\pi\in Irr(G)$ et $\sigma\in Irr(H)$ on définit $Hom_{H,\xi}(\pi,\sigma)$ comme l'espace des homomorphismes $l:E_\pi\to E_\sigma$ qui vérifient

$$l(\pi(hu)e)=\xi(u)\sigma(h)l(e)$$

\noindent pour tous $h\in H(F),u\in U(F)$ et $e\in E_\pi$. On note $m(\pi,\sigma)$ la dimension de cet espace. D'après [AGRS] et [GGP], cet espace est de dimension au plus 1. \\
Soit $R$ un $\mathcal{O}_E$-réseau spécial de $V$ en bonne position par rapport à un parabolique minimal inclus dans $P$ qui se décompose sous la forme $R=R_0\oplus R_Z$ où $R_0$ et $R_Z$ sont des $\mathcal{O}_E$-réseaux de $V_0$ et $Z$ respectivement. On a alors $G(F)=A(F)G_0(F)U(F)K$. Pour $N\geqslant 1$, on définit une fonction $\kappa_N$ sur $G(F)$ ainsi: c'est la fonction caractéristique de l'ensemble des $g\in G(F)$ vérifiant

\begin{itemize}
\item $g^{-1}(v_r)\cap \pi_E^{-N}R\neq \emptyset$
\item $g^{-1}(v_{r-1}+Ev_r)\cap \pi_E^{-2N}R\neq \emptyset$
\item $\ldots$
\item $g^{-1}(v_0+Ev_1+\ldots+Ev_r)\cap \pi_E^{-2N}R\neq \emptyset$
\item $g^{-1}(v_{-1}+W+Ev_0+Ev_1+\ldots+Ev_r)\cap \pi_E^{-2N}R\neq \emptyset$
\item $\ldots$
\item $g^{-1}(v_{-r}+Ev_{-r+1}+\ldots+Ev_{-1}+W+Ev_0+\ldots+Ev_r)\cap \pi_E^{-2N}R\neq \emptyset$.
\end{itemize}

\vspace{3mm}

La fonction $\kappa_N$ est invariante à gauche par $H(F)U(F)$ et à droite par $K$. Elle est de plus à support compact dans $H(F)U(F)\backslash G(F)$. Soient $\theta$ un quasi-caractère de $H(F)$ et $f\in C_c^{\infty}(G(F))$ une fonction très cuspidale. Pour $g\in G(F)$, on définit alors la fonction ${}^g\!f^\xi$ sur $H(F)$ de la façon suivante:

$${}^g\!f^\xi(h)=\displaystyle\int_{U(F)} f(g^{-1}hug)\xi(u) du$$

\noindent Cette fonction est localement constante à support compact, on peut donc définir

$$J(\theta,f,g)=\displaystyle\int_{H(F)} \theta(h){}^g\!f^\xi(h) dh$$
qui est une fonction localement constante en $g$. Pour $N\geqslant 1$, on pose alors:

$$J_N(\theta,f)=\displaystyle\int_{H(F)U(F)\backslash G(F)} J(\theta,f,g) \kappa_N(g) dg$$

La raison de l'expression précédente est qu'elle apparaît naturellement lorsque l'on essaye de calculer une multiplicité $m(\pi,\sigma)$ pour $\pi\in Irr(G)$ et $\sigma\in Irr(H)$ cuspidales (prendre pour $f$ un coefficient de $\pi$ et $\theta=\theta_{\sigma^\vee}$, cf [W1] proposition 13.1 pour le cas des groupes spéciaux orthogonaux). Le but de cet article est de démontrer que $J_N(\theta,f)$ admet une limite lorsque $N$ tend vers l'infini et de donner deux expressions pour cette limite. L'une qualifiée de géométrique est énoncée dans la section 5 et démontrée dans les sections 6 à 10. L'autre spectrale est énoncée et démontrée dans la section 16. De l'égalité entre ces deux développements, l'on déduira dans la section 17 la formule intégrale pour la multiplicité $m(\pi,\sigma)$ ($\pi\in Temp(G)$, $\sigma\in Temp(H)$) annoncée dans l'introduction.

\section{Le développement géométrique: définitions et énoncé}

\subsection{Un ensemble de tores}

 On définit $\underline{\mathcal{T}}$ comme l'ensemble des sous-tores $T$ de $H$ pour lesquels il existe une décomposition orthogonale $W=W'\oplus W''$ telle que, en notant $H'$, $H''$ et $G''$ les groupes unitaires de $W'$, $W''$ et $V''=Z\oplus D\oplus W''$ respectivement, les conditions suivantes soient vérifiées:
\begin{itemize}
\item $T$ est un tore maximal totalement anisotrope de $H'$.
\item $G''$ et $H''$ sont quasidéployés sur $F$.
\end{itemize}
On fixe $\mathcal{T}$ un ensemble de représentants des classes de conjugaison par $H(F)$ dans $\underline{\mathcal{T}}$.

\vspace{3mm}

 Soit $T\in\underline{\mathcal{T}}$ et $W',W''$ comme précédemment. On note $T_{\natural}$ l'ouvert de Zariski des éléments dont la restriction à $W'$ n'a pas $1$ pour valeur propre et a toutes ses valeurs propres de multiplicité 1. Soient $\theta$ et $\theta'$ deux quasicaractères sur $H(F)$ et $G(F)$ respectivement. On définit deux fonctions $c_\theta$ et $c_{\theta'}$ sur $T_{\natural}(F)$ de la façon suivante: pour $t\in T_{\natural}(F)$, on a $G_t(F)=G''(F)\times T(F)$ et $H_t(F)=H''(F)\times T(F)$. Par conséquent on a des identifications naturelles $Nil(\mathfrak{g}''(F))= Nil(\mathfrak{g}_t(F))$ et $Nil(\mathfrak{h}''(F))= Nil(\mathfrak{h}_t(F))$. On pose alors

$$c_{\theta'}(t)= \displaystyle\sum_{\mathcal{O}\in Nil_{reg}(\mathfrak{g}''(F))} \frac{c_{\theta',\mathcal{O}}(t)}{|Nil_{reg}(\mathfrak{g}''(F))|}$$

et

$$c_{\theta}(t)= \displaystyle\sum_{\mathcal{O}\in Nil_{reg}(\mathfrak{h}''(F))} \frac{c_{\theta,\mathcal{O}}(t)}{|Nil_{reg}(\mathfrak{h}''(F))|} $$

Si $f\in C_c^\infty(G(F))$ est une fonction très cuspidale, on posera $c_f=c_{\theta_f}$.

\subsection{Un critère de convergence}

 Soit $T\in\underline{\mathcal{T}}$. On note $W'_T$ et $W_T''$ les espaces notés $W'$ et $W''$ précédemment. On note $H'_T$ le groupe unitaire de $W'_T$. Pour $t\in T(F)$, on définit $E''(t)=Ker(t-1)_{|W'_T}$ et $E'(t)$ son supplémentaire orthogonal dans $W'_T$. On notera respectivement $J'(t)$ et $J''(t)$ les groupes unitaires de $E'(t)$ et de $E''(t)$ et $\mathfrak{z}_t$ désignera le centre de l'algèbre de Lie $\mathfrak{j}'(t)_t$. Remarquons que l'on a toujours $\mathfrak{z}_t\subset \mathfrak{t}$. \\
 
\begin{lem}
Soit $t\in T(F)$. On a

\begin{enumerate}[(i)]
\item $dim(\mathfrak{t})-dim(\mathfrak{z}_t)\geqslant dim(E''(t))$ avec égalité si et seulement si $J'(t)_t$ est un tore;

\item Il existe un bon voisinage $\omega\subset \mathfrak{t}(F)$ qui vérifie les deux conditions suivantes
\begin{itemize}
\item Pour tous $X\in \omega$ et $\lambda\in F^\times$ tels que $\lambda X\in\omega$ on a $\mathfrak{z}_{texp(X)}=\mathfrak{z}_{texp(\lambda X)}$
\item Pour tout $X\in \omega\backslash\mathfrak{z}_t$ on a $\mathfrak{z}_t\subsetneq\mathfrak{z}_{texp(X)}$
\end{itemize}
\end{enumerate}
\end{lem}

\ul{Preuve}: (i) $T$ est un sous-tore maximal de $H'_{T,t}=J''(t)\times J'(t)_t$. On a donc une décomposition $T=T''\times T'$ avec $T''$ et $T'$ sous-tores maximaux de $J''(t)$ et $J'(t)_t$ respectivement. En particulier, on a $dim(T'')=dim(E''(t))$ et $\mathfrak{z}_t\subset \mathfrak{t}'$. On en déduit que

\[\begin{aligned}
\displaystyle dim(\mathfrak{t})-dim(\mathfrak{z}_t) & =dim(T'')+dim(\mathfrak{t}')-dim(\mathfrak{z}_t) \\
 & \geqslant dim(T'')=dim(E''(t))
\end{aligned}\]

\noindent avec égalité si et seulement si $\mathfrak{t}'=\mathfrak{z}_t$ c'est-à-dire si et seulement si $J'(t)_t$ est un tore.

 (ii) Il existe un bon voisinage $\omega\subset \mathfrak{t}(F)$ tel que pour tout $X\in \omega$ on ait $E''(texp(X))=Ker(X_{|E''(t)})$ et $J'(texp(X))_{texp(X)}=J'(texp(X))_{t,X}$. Soit $\omega$ un tel voisinage. Pour tous  $X\in \omega$ et $\lambda\in F^\times$ tels que $\lambda X\in\omega$ on a alors $J'(texp(X))_{texp(X)}=J'(texp(\lambda X))_{texp(\lambda X)}$ d'où $\mathfrak{z}_{texp(X)}=\mathfrak{z}_{texp(\lambda X)}$. \\
Soit $X\in \omega$ et $\tilde{G}$ le groupe unitaire de l'orthogonal de $E'(t)$ dans $E'(texp(X))$. On a alors $J'(texp(X))_t=\tilde{G}\times J'(t)_t$ et $J'(texp(X))_{texp(X)}=\tilde{G}_X\times J'(t)_{t,X}$. De plus, $Z(J'(t)_t)^0\subset J'(t)_{t,X}$ d'où $\mathfrak{z}_t\subset\mathfrak{z}_{texp(X)}$. Puisque $X$ appartient au centre de l'algèbre de Lie $\mathfrak{j}'(texp(X))_{t,X}$, on a $X\in \mathfrak{z}_{texp(X)}$ donc $\mathfrak{z}_t\neq\mathfrak{z}_{texp(X)}$ si $X\notin \mathfrak{z}_t$ $\blacksquare$

\vspace{4mm}

 Pour $X$ un $F$-espace vectoriel et $i\in\mathbb{R}$, on note $C_i(X)$ l'espace des fonctions mesurables $\varphi: X\to \mathbb{C}$ vérifiant $\varphi(\lambda x)=|\lambda|_F^i\varphi(x)$ pour tout $x\in X$ et tout $\lambda\in F^{\times 2}$. On note $C_{\geqslant i}(X)$ l'espace des fonctions à valeurs complexes sur $X$ engendré par les $C_j(X)$ pour $j\geqslant i$. Soit $\chi:F^{\times}\to \{\pm 1\}$ un caractère quadratique, pour $i\in \mathbb{R}$ on notera $C_{i,\chi}(X)$ l'espace des fonctions mesurables $\varphi:X\to \mathbb{C}$ telle que $\varphi(\lambda x)=\chi(\lambda)|\lambda|_F^i\varphi(x)$ pour tout $x\in X$ et tout $\lambda\in F^{\times}$. Enfin, $C_{\geqslant i,\chi}(X)$ désignera l'espace de fonctions sur $X$ engendré par $C_{i,\chi}(X)$ et les $C_j(X)$ pour $j>i$. \\
 
 Soit $\delta: T(F)\to \mathbb{R}$ une fonction, on définit alors $C_{\geqslant \delta}(T)$ (resp. $C_{\geqslant \delta,\chi}(T)$) comme l'espace des fonctions $f$ définies presque partout sur $T(F)$ et telle que pour tout $t\in T(F)$, il existe un bon voisinage $\omega$ de $\mathfrak{t}(F)$ et une fonction $\varphi\in C_{\geqslant \delta(t)}(\mathfrak{t}(F)/\mathfrak{z}_t(F))$ (resp. $\varphi\in C_{\geqslant \delta(t),\chi}(\mathfrak{t}(F)/\mathfrak{z}_t(F))$) vérifiant

\begin{center}
$\forall X\in\omega$, $f_{t,\omega}(X)=\varphi(\overline{X})$ p.p. 
\end{center}

\noindent où $\overline{X}$ désigne la projection de $X$ sur $\mathfrak{t}(F)/\mathfrak{z}_t(F)$. Rappelons que l'on a défini en 3.1 une fonction $\Delta$ sur $G_{ss}(F)$.

\vspace{3mm}

\begin{lem}
Posons 

$$\displaystyle \delta(t)=min\big(\frac{dim(\mathfrak{z}_t)-dim(\mathfrak{t})-dim(E''(t))}{2},-1\big)$$

\noindent Soit $f\in C_{\geqslant\delta, \chi_E}(T)$, alors pour tout $s\in\mathbb{C}$ tel que $Re(s)>0$ l'intégrale $\displaystyle\int_{T(F)} f(t)\Delta(t)^s dt$ converge absolument et de plus la limite 

$$\lim\limits_{s\to 0^+} \displaystyle\int_{T(F)} f(t)\Delta(t)^s dt$$

\noindent existe.
\end{lem}

\ul{Preuve}: Il suffit de montrer que pour tout $t\in T(F)$, il existe un voisinage $\Omega$ de $t$ tel que $\displaystyle\int_{\Omega} f(x)\Delta(x)^s dx$ converge absolument pour $Re(s)>0$ et $\lim\limits_{s\to 0} \displaystyle\int_{\Omega} f(x)\Delta(x)^s dt$ existe. (*)\\

Pour $t\in T(F)$ et pour $X$ dans un bon voisinage assez petit de $0$ dans $\mathfrak{t}(F)$, on a $\Delta(t exp(X))=\Delta(t) |N(det(X_{|E''(t)}))|_F$ presque partout. La fonction $\Delta_t: X\mapsto |N(det(X_{|E''(t)}))|_F$ est invariante par $\mathfrak{z}_t(F)\subset \mathfrak{h}'_{T,t}(F)$ et presque partout homogène de degré $2dim(E''(t))$. \\

On pose $T(F)_n=\{t\in T(F): dim(\mathfrak{t})-dim(\mathfrak{z}_t)\leqslant n\}$ et on va montrer par récurrence sur $n$ que (*) est vraie pour $t\in T(F)_n$. Pour $t\in T(F)_0$, la fonction $f$ est localement constante au voisinage de $t$ (car $\mathfrak{z}_t=\mathfrak{t}$), donc le résultat est vérifié. Supposons le résultat vérifié pour $n-1$ et soit $t\in T(F)_n$. On peut supposer que $t$ n'est pas dans $T(F)_0$, on a alors $\delta(t)=\big(dim(\mathfrak{z}_t)-dim(\mathfrak{t})-dim(E''(t))\big)/2$. Choisissons un bon voisinage $\omega$ de $0$ dans $\mathfrak{t}(F)$ et $\varphi\in C_{\geqslant \delta(t),\chi_E}(\mathfrak{t}(F)/\mathfrak{z}_t(F))$ tels que $f(texp(X))=\varphi(\overline{X})$ p.p. pour $X\in \omega$. On peut décomposer $\varphi$ sous la forme 

$$\varphi=\displaystyle\sum_{i\geqslant \delta(t)} \varphi_i$$

\noindent où $\varphi_{\delta(t)}\in C_{\delta(t),\chi_E}(\mathfrak{t}(F)/\mathfrak{z}_t(F))$ et $\varphi_i \in C_i(\mathfrak{t}(F)/\mathfrak{z}_t(F))$ pour $i>\delta(t)$. Quitte à restreindre $\omega$ on peut supposer que $\omega$ vérifie les conditions du lemme 5.2.1, que l'exponentielle sur $\omega$ préserve les mesures et que $\omega$ admet une décomposition $\omega=\omega_{\mathfrak{z}}\times\omega'$ où $\omega_{\mathfrak{z}}\subset \mathfrak{z}_t$ est un voisinage de 0 et $\omega'$ est un réseau d'un supplémentaire de $\mathfrak{z}_t(F)$ dans $\mathfrak{t}(F)$. Pour $i\geqslant \delta(t)$ soit $f_i:T(F)\to \mathbb{C}$ la fonction définie par $f_i(texp(X))=\varphi_i(X)$ pour $X\in\omega$ et $f_i(t')=0$ si $t'\notin texp(\omega)$. On a alors

\begin{center}
(1)  $f_i\in C_{\geqslant\delta,\chi_E}(T)$
\end{center}  

 En effet, pour tout $\lambda\in \mathcal{O}_F^{\times 2}$, la fonction $g_\lambda$ définie par $g_\lambda(texp(X))=f(texp(\lambda X))$ pour $X\in\omega$ et $g_\lambda(t')=0$ si $t'\notin texp(\omega)$, est dans $C_{\geqslant\delta,\chi}(T)$ car pour tout $X\in \omega$ on a $\mathfrak{z}_{t exp(X)}=\mathfrak{z}_{t exp(\lambda X)}$. Comme $f_i$ est combinaison linéaire des fonctions $g_\lambda$ on a bien (1). \\
 
 Pour $X\in\omega\backslash \mathfrak{z}_t(F)$ on a $\mathfrak{z}_t\subsetneq\mathfrak{z}_{texp(X)}$ donc $texp(X)\in T(F)_{n-1}$. D'après l'hypothèse de récurrence, les $f_i$ vérifient (*) pour $t'=texp(X)$. Puisque l'exponentielle préserve les mesures sur $\omega$, pour tout compact $\omega''\subset\omega\backslash \mathfrak{z}_t(F)$ et pour $\mathfrak{R}(s)>0$ l'intégrale

$$\displaystyle\int_{\omega''} |\varphi_i(X)| \Delta_t(X)^s dX$$

\noindent converge et admet une limite lorsque $s$ tend vers 0. Choisissons une base du réseau $\omega'$ et soit $\Omega'$ l'ensemble des éléments de ce réseau dont la décomposition dans cette base s'écrit avec des coefficients dont la valuation minimum est 0 ou 1. Pour montrer la convergence de $\displaystyle\int_{texp(\omega)} f(t') \Delta(t')^sdt'$ pour $\mathfrak{R}(s)>0$, il suffit de montrer pour tout $i\geqslant\delta(t)$ la convergence de:

\[\begin{aligned}
\mbox{(2)} \;\;\; \displaystyle\int_{\omega} |\varphi_i(X)| \Delta_t(X)^s dX=\sum_{k\geqslant 0} \displaystyle\int_{\omega_{\mathfrak{z}}\times(\pi_F^{2k}\Omega')} |\varphi_i(\overline{X})| \Delta_t(X)^s dX &  \\
=vol(\omega_{\mathfrak{z}})\left (\sum_{k\geqslant 0} |\pi_F|_F^{2k\left(i+2sdim(E''(t))+dim(\mathfrak{t})-dim(\mathfrak{z}_t)\right)} \right) & \int_{\Omega'} |\varphi_i(\overline{X})| \Delta_t(X)^s dX
\end{aligned}\]

On a 

$$i+dim(\mathfrak{t})-dim(\mathfrak{z}_t)\geqslant \delta(t)+dim(\mathfrak{t})-dim(\mathfrak{z}_t)=\big(dim(\mathfrak{t})-dim(\mathfrak{z}_t)-dim(E''(t))\big)/2$$

\noindent et d'après le (i) du lemme 5.2.1, on a aussi $dim(\mathfrak{t})-dim(\mathfrak{z}_t)\geqslant dim(E''(t))$. Par conséquent

$$\displaystyle Re\big(i+2sdim(E''(t))+dim(\mathfrak{t})-dim(\mathfrak{z}_t)\big)\geqslant 2Re(s) dim(E''(t))\geqslant 0$$

\noindent On ne peut avoir égalité dans l'inégalité précédente que si $dim(\mathfrak{t})-dim(\mathfrak{z}_t)= dim(E''(t))=0$, c'est-à-dire si $t\in T(F)_0$, cas que l'on a exclu. On en déduit que l'intégrale (2) est égale à une série convergente. Pour $i>\delta(t)$ le calcul précédent nous montre même que l'intégrale 

$$\displaystyle\int_{\omega'} \varphi_i(X) \Delta_t(X)^s dX$$

\noindent admet une limite lorsque $s\to 0$. Pour montrer l'existence de la limite en $0$, il reste à montrer que pour $i=\delta(t)$ l'intégrale 

$$\displaystyle\int_{\omega'} \varphi_i(X) \Delta_t(X)^s dX$$

\noindent admet une limite lorsque $s$ tend vers $0$. Si $E/F$ est ramifiée alors on peut trouver $\lambda\in \mathcal{O}_F^{\times}$ tel que $\chi_E(\lambda)=-1$ et le changement de variable $X\mapsto \lambda X$ montre que cette intégrale est toujours nulle. Reste le cas où $E/F$ est non ramifiée. Définissons $\Omega''$ comme le sous-ensemble des éléments de $\Omega'$ dont les coefficients dans la base déjà fixée ont une valuation minimale nulle. On a alors:

\[\begin{aligned}
\displaystyle\int_{\omega'} \varphi_i(X) \Delta_t(X)^s dX=\displaystyle\sum_{k\geqslant 0} \int_{\pi_F^{k}\Omega''} \varphi_i(X) \Delta_t(X)^s dX \\
=\left (\sum_{k\geqslant 0} (-1)^k|\pi_F|_F^{k\big(\delta(t)+2sdim(E''(t))+dim(\mathfrak{t})-dim(\mathfrak{z}_t)\big)}\right) \int_{\Omega'} \varphi_i(X & ) \Delta_t(X)^s dX \\
=\big(1+|\pi_F|_F^{\delta(t)+2sdim(E''(t))+dim(\mathfrak{t})-dim(\mathfrak{z}_t)}\big)^{-1} \int_{\Omega'} \varphi_i(X) \Delta_t(X)^s
\end{aligned}\]

\noindent et cette expression admet une limite lorsque $s\to0$. $\blacksquare$

\subsection{Définition de $J_{geom}(\theta,f)$}

Soient $\theta$ et $\theta'$ des quasicaractères sur $H(F)$ et $G(F)$ respectivement.

\begin{lem}
Soit $T\in\underline{\mathcal{T}}$ alors la fonction $t\mapsto c_{\theta'}(t)c_\theta(t)D^H(t)^{1/2}D^G(t)^{1/2}\Delta(t)^{-1/2}$ est dans $C_{\geqslant \delta,\chi_E}(T)$.
\end{lem}

\vspace{4mm}

Soit $f\in C_c^\infty(G(F))$ une fonction très cuspidale. On peut donc d'après le lemme 5.2.2 définir la quantité suivante:

$$J_{geom}(\theta,f)=\displaystyle\sum_{T\in\mathcal{T}}\; |W(H,T)|^{-1} \nu(T) \lim\limits_{s\to 0^+} \displaystyle\int_{T(F)} c_\theta(t) c_f(t) D^H(t)^{1/2}D^G(t)^{1/2} \Delta(t)^{s-1/2}  dt$$

\vspace{3mm}

\ul{Preuve}: On reprend les notations du paragraphe 5.1. A $T$ sont associés des sous-espaces $W'$ et $W''$ de $W$ et des groupes unitaires $H'$, $H''$ et $G''$. Pour $t\in T(F)$, on note $E''(t)$ le noyau de $t-1$ dans $W'$ et $E'(t)$ le supplémentaire orthogonal de $E''(t)$ dans $W'$. On désignera par $J'(t)$ et $J''(t)$ les groupes unitaires de $E'(t)$ et $E''(t)$ respectivement et par $\tilde{G}(t)$ le groupe unitaire de $Z\oplus D\oplus W''\oplus E''(t)$. On a alors $G_t=\tilde{G}(t) J'(t)_t$ et $\mathfrak{g}_t=\tilde{\mathfrak{g}}(t)\oplus \mathfrak{j}'(t)_t$.

\vspace{3mm}

On définit les fonctions suivantes sur $T(F)$:

$$\delta_0,\delta_G,\delta_H:T(F)\to \mathbb{Z}$$

$$\delta_0(t)=\frac{\big(\delta(H_t)-\delta(H'')+\delta(G_t)-\delta(G'')\big)}{2}-dim(E''(t))$$

\vspace{2mm}

$$\delta_G(t)=\left\{
    \begin{array}{ll}
        \frac{1}{2}(\delta(G'')-\delta(G_t)+2) & \mbox{si } dim(E''(t))\geqslant 2 \mbox{ et } dim(W''\oplus E''(t))\equiv 0 [2]\\
        \frac{1}{2}(\delta(G'')-\delta(G_t)) & \mbox{sinon.}
    \end{array}
\right.$$

$$\delta_H(t)=\left\{
    \begin{array}{ll}
        \frac{1}{2}(\delta(H'')-\delta(H_t)+2) & \mbox{si } dim(E''(t))\geqslant 2 \mbox{ et } dim(W''\oplus E''(t))\equiv 1 [2]\\
        \frac{1}{2}(\delta(H'')-\delta(H_t)) & \mbox{sinon.}
    \end{array}
\right.$$

On va montrer dans un premier temps que $c_{\theta'}\in C_{\geqslant\delta_G}(T)$, $c_\theta\in C_{\geqslant \delta_H}(T)$ et $\big(D^HD^G\Delta^{-1}\big)^{1/2}\in C_{\geqslant\delta_0,\chi_0}(T)$ où $\chi_0$ est le caractère trivial de $F^{\times}$.

\vspace{2mm}

 Soit $t\in T(F)$ alors pour $X$ dans un bon voisinage assez petit de 0 dans $\mathfrak{t}(F)$, on a $D^H(texp(X))=D^H(t)D^{H_t}(X)$, $D^G(texp(X))=D^G(t)D^{G_t}(X)$ et \\
\noindent $\Delta(texp(X))=\Delta(t)|N(det(X''_{E''(t)/Ker(X'')}))|_F$ où $X''$ est la restriction de $X$ à $E''(t)$. La fonction $X\mapsto |N(det(X''_{E''(t)/Ker(X'')}))|_F$ est homogène de degré $2dim(E''(t))$ en dehors d'un fermé de mesure nulle et elle est clairement invariante par $\mathfrak{z}_t(F)$. Pour $X$ en dehors d'un ensemble de mesure nulle on a $H_{t,X}=TH''$ et $G_{t,X}=TG''$. Sur cet ensemble $D^{H_t}$ et $D^{G_t}$ sont homogènes de degrés $dim(H_t)-dim(H_{t,X})=\delta(H_t)-\delta(H'')$ et $dim(G_t)-dim(G_{t,X})=\delta(G_t)-\delta(G'')$ respectivement. Elles sont aussi invariantes par $\mathfrak{z}_t(F)$ car $\mathfrak{z}_t$ est central dans $\mathfrak{h}_t$ et $\mathfrak{g}_t$. On en déduit que $\big(D^HD^G\Delta^{-1}\big)^{1/2}\in C_{\geqslant\delta_0,\chi_0}(T)$.

\vspace{3mm}

 Soit $t\in T(F)$, il existe un bon voisinage $\omega\subset\mathfrak{g}_t(F)$ de $0$ de sorte que l'on ait
 
$$\displaystyle\theta'(texp(X))=\sum_{\mathcal{O}\in Nil(\mathfrak{g}_t)} c_{\theta',\mathcal{O}}(t)\hat{j}^{G_t}(\mathcal{O},X)$$

\noindent pour presque tout $X\in\omega$. Par linéarité on peut supposer que $\theta'(texp(X))=\hat{j}^{G_t}(\mathcal{O},X)$ pour presque tout $X\in\omega$. Cette fonction est invariante par translation par les éléments du centre de $\mathfrak{g}_t(F)$ donc aussi par $\mathfrak{z}_t(F)$. Soit $X\in \omega\cap\mathfrak{t}(F)$ tel que $texp(X)\in T_{\natural}(F)$. Pour presque tout $Y$ dans un voisinage de 0 dans $\mathfrak{g}_{t,X}(F)$ on a un développement 

$$\hat{j}^{G_t}(\mathcal{O},X+Y)=\displaystyle\sum_{\mathcal{O}'\in Nil(\mathfrak{g''}(F))} c_{\theta',\mathcal{O}'}(texp(X)) \hat{j}^{G''}(\mathcal{O}',Y)$$

\noindent La fonction $X\mapsto c_{\theta'}(texp(X))$ est une combinaison linéaire des coefficients $c_{\theta',\mathcal{O}'}(texp(X))$ pour $\mathcal{O}'\in Nil_{reg}(\mathfrak{g}''(F))$. Des propriétés d'homogénéité de $\hat{j}(\mathcal{O},.)$ et $\hat{j}(\mathcal{O}',.)$ on déduit que la fonction $X\in\omega \mapsto c_{\theta'}(texp(X))$ est dans $C_{\geqslant \frac{\delta(G'')-dim(\mathcal{O})}{2}}(\mathfrak{t}(F)/\mathfrak{z}_t(F))$. On a $dim(\mathcal{O})\leqslant \delta(G_t)$. Cela conclut si $dim(E''(t))\leqslant 1$ ou $dim(W''\oplus E''(t))$ est impaire. Si $dim(E''(t))\geqslant 2$ et $dim(W''\oplus E''(t))$ est paire il suffit de montrer que pour $\mathcal{O}_{reg}\in Nil_{reg}(\mathfrak{g}_t(F))$ et $X\in \omega\cap\mathfrak{t}(F)$ en position générale la fonction $Y\mapsto\hat{j}(\mathcal{O}_{reg},X+Y)$ est nulle dans un voisinage de 0 dans $\mathfrak{g}''(F)$. Dans le cas considéré $\tilde{G}(t)$ est un groupe unitaire de dimension impaire, donc $\tilde{\mathfrak{g}}(t)(F)$ admet une unique orbite nilpotente régulière $\tilde{\mathcal{O}}$. On a une factorisation $\hat{j}(\mathcal{O}_{reg},.)=\hat{j}(\tilde{\mathcal{O}},.)\tau(.)$ relativement à la décomposition $\mathfrak{g}_t(F)=\tilde{\mathfrak{g}}(t)(F)\oplus \mathfrak{j}'(t)_t(F)$, où $\tau$ est un quasi-caractère de $\mathfrak{j}'(t)_t(F)$. Or, $\tilde{\mathcal{O}}$ est induite à partir de l'orbite $\{0\}$ d'une sous-algèbre de Borel et donc est à support dans l'ensemble des éléments qui appartiennent à une sous-algèbre de Borel. Soit $X\in\omega\cap \mathfrak{t}(F)$ que l'on décompose en $X=X''+X'$ avec $X'\in\mathfrak{j}'(t)_t(F)$ et $X''\in\mathfrak{j}''(t)(F)\cap\mathfrak{t}(F)\subset \tilde{\mathfrak{g}}(t)(F)$. Le groupe $J''(t)\cap T$ est un tore anisotrope de dimension plus grande que $2$. Donc pour $X$ dans un ouvert de Zariski, $X''$ possède un voisinage qui ne contient aucun élément dans une sous-algèbre de Borel de $\tilde{\mathfrak{g}}(t)(F)$. Par conséquent $\hat{j}(\tilde{\mathcal{O}},.)$ s'annule au voisinage de $X''$. On en déduit le résultat pour $c_{\theta'}$. On procède de même pour $c_\theta$.

\vspace{3mm}

Un calcul indolore nous donne pour $t\in T(F)$, 
$$\delta_G(t)+\delta_H(t)+\delta_0(t)=\left\{
    \begin{array}{ll}
        -dim(E''(t))+1 & \mbox{si } dim(E''(t))\geqslant 2 \\
        -dim(E''(t)) & \mbox{sinon.}
    \end{array}
\right.
$$
 Si $dim(E''(t))\geqslant 2$, on a $\delta(t)=\big(dim(\mathfrak{z}_t)-dim(\mathfrak{t})-dim(E''(t))\big)/2$ et d'après le (i) du lemme 5.2.1, $\delta(t)\leqslant -dim(E''(t))$. Si $dim(E''(t))=0$, on a $\delta(t)\leqslant -1$. Enfin, si $dim(E''(t))=1$, toujours d'après le (i) du lemme 5.2.1, on a $\delta(t)\leqslant -dim(E''(t))$ avec égalité seulement si $J'(t)_t$ est un tore. On en déduit dans tout les cas que $\delta_G(t)+\delta_H(t)+\delta_0(t)\geqslant \delta(t)$ et que l'on ne peut avoir égalité que si $J'(t)_t$ est un tore et $dim(E''(t))=1$. Pour obtenir le lemme il nous reste à établir le fait suivant: pour $t\in T(F)$ tel que $J'(t)_t$ est un tore et $dim(E''(t))=1$ alors il existe un bon voisinage $\omega_T$ de $0$ dans $\mathfrak{t}(F)$ tel que:

\vspace{2mm}

\begin{itemize}
\item Si $dim(W'')$ est paire alors la fonction $X\in\omega_T\mapsto c_{\theta'}(texp(X))$ est la restriction à $\omega_T$ d'un élément de $C_{\geqslant \delta_G(t),\chi_E}(\mathfrak{t}(F)/\mathfrak{z}_t(F))$, et la fonction $X\in\omega_T\mapsto c_\theta(texp(X))$ est la restrcition à $\omega_T$ d'un élément de $C_{\geqslant \delta_H(t),\chi_0}(\mathfrak{t}(F)/\mathfrak{z}_t(F))$.

\item Si $dim(W'')$ est impaire alors la fonction $X\in\omega_T\mapsto c_{\theta'}(texp(X))$ est la restriction à $\omega_T$ d'un élément de $C_{\geqslant \delta_G(t),\chi_0}(\mathfrak{t}(F)/\mathfrak{z}_t(F))$, et la fonction $X\in\omega_T\mapsto c_\theta(texp(X))$ est la restriction à $\omega_T$ d'un élément de $C_{\geqslant \delta_H(t),\chi_E}(\mathfrak{t}(F)/\mathfrak{z}_t(F))$.
\end{itemize}

\vspace{2mm}

 On démontre ceci pour $c_{\theta'}$ la preuve étant la même pour $c_\theta$. Puisque $J'(t)_t$ est un tore on a une identification $Nil(\mathfrak{g}_t(F))=Nil(\tilde{\mathfrak{g}}(t)(F))$. D'après ce qui précède il n'y a rien à dire si $Nil_{reg}(\tilde{\mathfrak{g}}(t)(F))=\emptyset$ on suppose donc que $Nil_{reg}(\tilde{\mathfrak{g}}(t)(F))\neq\emptyset$.
 
\vspace{2mm}
 
  Dans le cas $dim(W'')$ impaire, $\tilde{\mathfrak{g}}(t)(F)$ possède une unique orbite nilpotente régulière $\tilde{\mathcal{O}}$ et $\mathfrak{g}''(F)$ possède deux orbites nilpotentes régulières $\mathcal{O}^{+}$ et $\mathcal{O}^{-}$. Soit $\omega$ un bon voisinage de $0$ dans $\mathfrak{g}_t(F)$ assez petit. Par linéarité, on peut supposer que $\theta'_{t,\omega}=\hat{j}(\tilde{\mathcal{O}},.)$ sur $\omega$. On a pour $X\in\omega\cap\mathfrak{t}(F)$ en position générique un développement au voisinage de $0$ dans $\mathfrak{g}''(F)$ de la forme $\hat{j}(\tilde{\mathcal{O}},X+Y)=\displaystyle\sum_{\mathcal{O}\in\ Nil(\mathfrak{g}''(F))} c_{\hat{j}(\tilde{\mathcal{O}},.),\mathcal{O}}(X)\hat{j}(\mathcal{O},Y)$. Par définition, on a
  
$$c_{\theta'}(texp(X))=\frac{c_{\hat{j}(\tilde{\mathcal{O}},.),\mathcal{O}^{+}}(X)+c_{\hat{j}(\tilde{\mathcal{O}},.),\mathcal{O}^{-}}(X)}{2}$$

\noindent Pour $\lambda\in F^{\times}$ et presque tout $X\in \mathfrak{t}(F), \; Y\in\mathfrak{g}''(F)$, on a $\hat{j}(\tilde{\mathcal{O}},\lambda X+\lambda Y)=|\lambda|_F^{-\frac{\delta(G_t)}{2}}\hat{j}(\tilde{\mathcal{O}},X+Y)$ et $\hat{j}(\mathcal{O}^\pm,\lambda Y)=|\lambda|_F^{-\frac{\delta(G'')}{2}}\hat{j}(\lambda\mathcal{O}^\pm,Y)$ et la multiplication par $\lambda$ conserve $\{\mathcal{O}^{+},\mathcal{O}^{-}\}$: on en déduit le résultat.
  
\vspace{2mm}

 Dans le cas $dim(W'')$ paire, $\tilde{\mathfrak{g}}(t)(F)$ possède deux orbites nilpotentes régulières $\tilde{\mathcal{O}}^{+}$ et $\tilde{\mathcal{O}}^{-}$ alors que $\mathfrak{g}''(F)$ ne possède qu'une orbite nilpotente régulières $\mathcal{O}''$. Comme auparavant, on peut fixer un bon voisinage $\omega$ de $0$ dans $\mathfrak{g}_t(F)$ assez petit et supposer que $\theta'_{t,\omega}=\hat{j}(\tilde{\mathcal{O}}^{+},.)$ presque partout sur $\omega$. On a pour $X\in\omega\cap\mathfrak{t}(F)$ en position générique un développement au voisinage de $0$ dans $\mathfrak{g}''(F)$ de la forme 

$$\hat{j}(\tilde{\mathcal{O}}^{+},X+Y)=\displaystyle\sum_{\mathcal{O}\in\ Nil(\mathfrak{g}''(F))} c_{\hat{j}(\tilde{\mathcal{O}}^{+},.),\mathcal{O}}(X)\hat{j}(\mathcal{O},Y)$$

\noindent Par définition, on a

$$c_{\theta'}(texp(X))=c_{\hat{j}(\tilde{\mathcal{O}}^{+},.),\mathcal{O}''}(X)$$

La fonction $\hat{j}(\tilde{\mathcal{O}}^{+},.)+\hat{j}(\tilde{\mathcal{O}}^{-},.)$ est à support dans l'ensemble des éléments de $\tilde{\mathfrak{g}}(t)(F)$ qui appartiennent à une sous-algèbre de Borel. Pour $X\in\mathfrak{t}(F)$ en position générale, la projection $X''$ de $X$ sur $\tilde{\mathfrak{g}}(t)(F)$ suivant $\mathfrak{j}'(t)_t(F)$ est non nulle et appartient à un tore anisotrope de dimension 1, donc n'appartient à aucune sous-algèbre de Borel de $\tilde{\mathfrak{g}}(t)(F)$ (car $\tilde{G}(t)$ est un groupe unitaire de dimension paire). Par conséquent pour presque tout $X\in \mathfrak{t}(F)$ la fonction $\hat{j}(\mathcal{O}^+,.)+\hat{j}(\mathcal{O}^-,.)$ s'annule au voisinage de $X$. Pour $\lambda\in F^{\times}$ et presque tout $X\in\mathfrak{t}(F),\; Y\in\mathfrak{g}''(F)$ on a donc

$$\hat{j}(\tilde{\mathcal{O}}^{+},\lambda X+\lambda Y) =|\lambda|_F^{-\frac{\delta(G_t)}{2}}\hat{j}(\lambda\tilde{\mathcal{O}}^{+},X+Y) =|\lambda|_F^{-\frac{\delta(G_t)}{2}}\chi_E(\lambda)\hat{j}(\tilde{\mathcal{O}}^{+},X+Y)$$

\noindent d'après ce qui précède et $\hat{j}(\mathcal{O}'',\lambda Y)=|\lambda|_F^{-\frac{\delta(G'')}{2}}\hat{j}(\mathcal{O}'',Y)$. On en déduit alors ce que l'on voulait  $\blacksquare$

\vspace{4mm}

\subsection{Enoncé du développement géométrique}

C'est le théorème suivant. Sa preuve occupera les sections 6 à 10.

\begin{theo}
Soient $\theta$ un quasicaractère de $H(F)$ et $f\in C_c^\infty (G(F))$ une fonction très cuspidale on a alors:
\begin{center}
$\lim\limits_{N\to\infty} J_N(\theta,f)=J_{geom}(\theta,f)$
\end{center}
\end{theo}

\vspace{4mm}

\subsection{Enoncé du théorème pour les algèbres de Lie}

 Le caractère $\xi$ de $U(F)$ se descend à l'algèbre de Lie via l'exponentielle, on en déduit un caractère encore noté $\xi$ de $\mathfrak{u}(F)$. Pour $\theta$ un quasicaractère de $\mathfrak{h}(F)$ et $f\in C_c^\infty(\mathfrak{g}(F))$ une fonction très cuspidale on définit alors les fonctions et quantités suivantes:

\begin{itemize}
\item Pour $g\in G(F)$ et $X\in\mathfrak{h}(F)$ on pose ${}^g\!f^\xi (X)=\displaystyle\int_{\mathfrak{u}(F)} f(g^{-1}(X+N)g)\xi(N) dN$.
\item Pour $g\in G(F)$ on pose $J(\theta,f,g)=\displaystyle\int_{\mathfrak{h}(F)} \theta(X) {}^gf^\xi(X) dX$.
\item Pour $N\geqslant 1$ on pose $J_N(\theta,f)=\displaystyle\int_{H(F)U(F)\backslash G(F)} \kappa_N(g) J(\theta,f,g) dg$.
\end{itemize}

\vspace{3mm}

Pour $X\in \mathfrak{g}_{ss}(F)$ on pose $\Delta(X)=|N(det(X_{|W/W''(X)}))|_F$ où $W''(X)$ est le noyau de $X$.

\vspace{3mm}

Pour $T\in\mathcal{T}$ et $X\in\mathfrak{t}(F)_\natural$ on définit

\begin{center}
$c_f(X)=\displaystyle\sum_{\mathcal{O}\in Nil_{reg}(\mathfrak{g}''(F))} \frac{c_{\theta_f,\mathcal{O}}(X)}{|Nil_{reg}(\mathfrak{g}''(F))|}$ et $c_\theta(X)=\displaystyle\sum_{\mathcal{O}\in Nil_{reg}(\mathfrak{h}''(F))} \frac{c_{\theta,\mathcal{O}}(X)}{|Nil_{reg}(\mathfrak{h}''(F))|}$.
\end{center}

\vspace{3mm}

De façon similaire à ce que l'on a fait sur le groupe on peut définir:

$$\mbox{(1)}\;\; J_{geom}(\theta,f)=\displaystyle\sum_{T\in\mathcal{T}} |W(H,T)|^{-1} \nu(T) \lim\limits_{s\to 0^+}\displaystyle\int_{\mathfrak{t}(F)} c_\theta(X)c_f(X)D^H(X)^{1/2}D^G(X)^{1/2}\Delta(X)^{s-1/2} dX$$

\vspace{3mm}

L'analogue du théorème 5.4.1 pour les algèbres de Lie est alors:
\begin{theo}
Soient $\theta$ un quasicaractère de $\mathfrak{h}(F)$ et $f\in C_c^\infty(\mathfrak{g}(F))$ une fonction très cuspidale on a
$$\lim\limits_{N\to\infty} J_N(\theta,f)=J_{geom}(\theta,f)$$
\end{theo}

Les théorèmes 5.4.1 et 5.5.1 seront démontrés dans la section 10.

\section{Descente à l'algèbre de Lie}

Fixons $x\in H_{ss}(F)$. Par la suite, on notera respectivement $W''$, $V_0''$ et $V''$ le noyau de $x-1$ dans respectivement $W$, $V_0$ et $V$, on notera $W'$ l'orthogonal de $W''$ dans $W$ (on a alors $W=W''\oplus W'$) et on notera $G'=H'$, $G''$, $H''$ et $G_0''$ les groupes unitaires de respectivement $W'$, $V''$, $W''$ et $V_0''$ (on a alors $G_x=G'_xG''$ et $H_x=H'_xH''$). Soit $\omega'\subset \mathfrak{g}'_x(F)$ et $\omega''\subset \mathfrak{g}''(F)$ deux bons voisinages de $0$ on pose alors $\omega=\omega'\times \omega''\subset \mathfrak{g}_x(F)$ c'est un bon voisinage de $0$. Soit $\Omega=(xexp(\omega))^G$. On supposera jusqu'à la section 8 que l'on a $Supp(f)\subset \Omega$.

\vspace{2mm}

\subsection{Localisation de $J_N(\theta,f)$}

\vspace{2mm}

On définit pour $g\in G(F)$ la fonction suivante sur $\mathfrak{h}_x(F)$
$${}^g\! f^\xi_{x,\omega}(X)=\displaystyle\int_{\mathfrak{u}_x(F)} {}^g\!f_{x,\omega}(X+N) \xi(N) dN$$
On définit aussi
$$J_{x,\omega}(\theta,f,g)=\displaystyle\int_{\mathfrak{h}_x(F)} \theta_{x,\omega}(X) {}^g\!f^\xi_{x,\omega}(X) dX$$
et
$$J_{x,\omega,N}(\theta,f)=\displaystyle\int_{H_x(F)U_x(F)\backslash G(F)} \kappa_N(g) J_{x,\omega}(\theta,f,g) dg$$
Posons $C(x)=D^G(x)^{1/2}D^H(x)^{1/2}\Delta(x)^{-1/2}$.

\begin{lem}
On a l'égalité
$$J_N(\theta,f)=C(x)J_{x,\omega,N}(\theta,f)$$
\end{lem}

\ul{Preuve}:
Il suffit de reprendre la preuve du lemme 8.2 de [W1]. Rappelons pour la commodité du lecteur comment l'on procède. \\

Par la formule d'intégration de Weyl, on a:
$$J(\theta,f,g)=\displaystyle\sum_{T\in\mathcal{T}(H)} |W(H,T)|^{-1} \displaystyle\int_{T(F)} \theta(t) D^H(t) \displaystyle\int_{T(F)\backslash H(F)}{}^{hg}\!f^\xi(t) dh dt$$
Pour $T$ et $T'$ deux tores maximaux de $H$ on note $W(T,T')=\{w\in H(F): wTw^{-1}=T'\}/T(F)$. On a alors:

\vspace{2mm}

1)Soient $T\in \mathcal{T}(H)$ et $t\in T(F)\cap H_{reg}(F)$ tel que $\displaystyle\int_{T(F)\backslash H(F)}{}^{hg}\!f^\xi(t) dh\neq 0$ alors

$$\displaystyle t\in\bigcup_{T_1\in \mathcal{T}(H_x)}\bigcup_{w\in W(T_1,T)} w(xexp(\mathfrak{t}_1(F)\cap \omega))w^{-1}$$

\vspace{2mm}

 En effet, il existe alors $u\in U(F)$ tel que $tu\in(xexp(\omega))^G$. La partie semisimple de $tu$ étant conjuguée à $t$ on a aussi $t\in (xexp(\omega))^G$. Soient donc $X\in \omega$ et $y\in G(F)$ tels que $yty^{-1}=xexp(X)$. On a alors $y(Z\oplus D)\subset Ker(xexp(X)-1) \subset Ker(x-1)=V''$ et il existe d'après le théorème de Witt un élément $y'\in G''(F)\subset G_x(F)$ tel que $y'y(Z\oplus D)=Z\oplus D$. Quitte à remplacer $y$ par $y'y$ et $X$ par $y'Xy'^{-1}$, on peut donc supposer que $y(Z\oplus D)=Z\oplus D$. Alors $yty^{-1}$ agit trivialement sur $Z\oplus D$ donc $xexp(X)$ aussi et $X\in\mathfrak{h}_x(F)$. Quitte à conjuguer encore $X$ par un élément de $H_x(F)$ on peut certainement supposer qu'il existe $T_1\in\mathcal{T}(H_x)$ tel que $X\in \mathfrak{t}_1(F)\cap \omega$. Soit $h$ la restriction de $y^{-1}$ à $W$ on a alors $t=hxexp(X)h^{-1}$. La conjugaison par $h$ envoie donc $Z_H(xexp(X))^0$ sur $Z_H(t)^0$. Puisque $t$ est régulier dans $H$ on a $Z_H(t)^0=T$ et $Z_H(xexp(X))^0=T_1$.

\vspace{2mm}

2) Pour $T\in \mathcal{T}(H)$, $T_1,T_2\in \mathcal{T}(H_x)$ et $w_1\in W(T_1,T)$, $w_2\in W(T_2,T)$ les deux ensembles $w_1(xexp(\mathfrak{t}_1(F)\cap\omega))w_1^{-1}$ et $w_2(xexp(\mathfrak{t}_2(F)\cap\omega))w_2^{-1}$ sont disjoints ou égaux et s'ils sont égaux on a $T_1=T_2$.

\vspace{2mm}

Soient $y_1,y_2\in H(F)$ qui relèvent $w_1$ et $w_2$ on pose $y=y_2^{-1}y_1$. Si les deux ensembles en question ne sont pas disjoints, on a $y(xexp(\omega))y^{-1}\cap xexp(\omega)\neq \emptyset$ donc $y\in Z_H(x)(F)=H_x(F)$ et la conjugaison par $y$ envoie $T_1$ sur $T_2$ et $\omega$ sur $\omega$.

\vspace{2mm}

3)Soient $T\in\mathcal{T}(H)$, $T_1\in\mathcal{T}(H_x)$ et $w_1\in W(T_1,T)$. Le nombre d'éléments $w_2\in W(T_1,T)$ tels que $w_2(xexp(\mathfrak{t}_1(F)\cap\omega))w_2^{-1}=w_1(xexp(\mathfrak{t}_1(F)\cap\omega))w_1^{-1}$ est $|W(H_x,T_1)|$.

\vspace{2mm}

En effet d'après ce qu'on vient de voir l'ensemble de ces éléments est en bijection avec $\{y\in H_x(F): yT_1y^{-1}=T_1\}/T_1(F)$ qui est précisément $W(H_x,T_1)$.

\vspace{2mm}

4)Soit $T_1\in\mathcal{T}(H_x)$, alors il existe un unique tore $T\in \mathcal{T}(H)$ tel que $W(T_1,T)\neq \emptyset$ et on a alors $|W(T_1,T)|=|W(H,T)|$.

\vspace{2mm}

L'existence et l'unicité de $T$ sont évidentes puisque $\mathcal{T}(H)$ est précisément un ensemble de représentants des classes de conjugaison de tores maximaux de $H$. Soit donc $T$ l'unique élément de $\mathcal{T}(H)$ tel que $W(T_1,T)\neq\emptyset$ et fixons $w_1\in W(T_1,T)$. On a alors une bijection entre $W(T_1,T)$ et $W(H,T)$ donnée par $w\mapsto ww_1^{-1}$.

\vspace{2mm}

On déduit des points 1) à 4) ci-dessus que l'on a

\[\begin{aligned}
J(\theta,f,g) & =\displaystyle\sum_{T_1\in\mathcal{T}(H_x)}\sum_{T\in\mathcal{T}(H)}\sum_{w_1\in W(T_1,T)} |W(H,T)|^{-1} |W(H_x,T_1)|^{-1} \\
  & \int_{\mathfrak{t}_1(F)\cap\omega} D^H(w_1(xexp(X))w_1^{-1})^{1/2} \theta(w_1(xexp(X))w_1^{-1}) J_H(w_1(xexp(X))w_1^{-1},{}^g\!f^\xi) dX \\
  & =\sum_{T_1\in\mathcal{T}(H_x)} |W(H_x,T_1)|^{-1} \int_{\mathfrak{t}_1(F)\cap\omega} D^H(xexp(X))^{1/2} \theta(xexp(X)) J_H(xexp(X),{}^g\!f^\xi) dX
\end{aligned}\]

Pour $X\in\omega\cap\mathfrak{t}_1(F)\cap\mathfrak{h}_{reg}(F)$ on a $D^H(xexp(X))=D^H(x)D^{H_x}(X)$ et

$$\displaystyle J_H(xexp(X),{}^g\!f^\xi)= D^H(x)^{1/2}\int_{H_x(F)\backslash H(F)} J_{H_x}(xexp(X), {}^{hg}\!f^\xi) dh$$

D'après la formule de Weyl sur $H_x$, on a donc:

$$J_N(\theta,f)=D^H(x)\displaystyle\int_{H_x(F)U(F)\backslash G(F)} \kappa_N(g) \int_{\mathfrak{h}_x(F)} \varphi_g(X) dX dg$$

où

$$
\varphi_g(X)= \left\{
    \begin{array}{ll}
        \theta_{x,\omega}(X) {}^g\!f^\xi(xexp(X)) & \mbox{si } X \in \omega \\
        0 & \mbox{sinon.}
    \end{array}
\right.
$$

Soit $\mathfrak{u}^x(F)=Im(ad(x)-1_{|\mathfrak{u}(F)})$, on a alors $\mathfrak{u}(F)=\mathfrak{u}_x(F)\oplus\mathfrak{u}^x(F)$. On a fixé une mesure sur $\mathfrak{u}(F)$ et on suppose que cette mesure est produit d'une mesure sur $\mathfrak{u}_x(F)$ et d'une mesure sur $\mathfrak{u}^x(F)$. Soit $U^x(F)=exp(\mathfrak{u}^x(F))$ qu'on munit de la mesure image par l'exponentielle. Pour $X\in\omega\cap\mathfrak{h}_{x,reg}(F)$ et $g\in G(F)$, on a
$${}^g\! f^\xi(xexp(X))=\int_{U_x(F)\backslash U(F)}\int_{U_x(F)} {}^g\!f(xexp(X)uv) \xi(uv) du dv$$
Pour $u\in U_x(F)$ l'application $v\mapsto \overline{(xexp(X)u)^{-1}v^{-1}xexp(X)uv}$ est un isomorphisme de $U^x(F)$ sur $U_x(F)\backslash U(F)$ (la barre désignant la classe d'un élément dans $U_x(F)\backslash U(F)$). Son jacobien est $|det(1-ad(xexp(X)u)_{|\mathfrak{u}^x(F)})|_F=|det(1-ad(xexp(X))_{|\mathfrak{u}^x(F)})|_F$ et on a \\
\noindent $|det(1-ad(xexp(X))_{|\mathfrak{u}^x(F)})|_F=|det(1-ad(x)_{|\mathfrak{u}^x(F)})|_F$ d'après la propriété 3.1(7) de [W1] sur les bons voisinages. On a un isomorphisme de $F$-espace vectoriel

$$W'\otimes_F (Fv_1\oplus\ldots\oplus Fv_r)\to u^x(F)$$

$$w'\otimes z\to c(z,w')$$

On en déduit que le jacobien de notre application est $\Delta(x)^r$. L'application

$$\mathfrak{u}_x(F)\to U_x(F)$$

$$N\mapsto exp(-X)exp(X+N)$$

\noindent est une bijection qui préserve les mesures. Puisque $\xi(exp(-X)exp(X+N))=\xi(N)$, on a

\[\begin{aligned}
{}^g\!f^\xi(xexp(X)) & =\Delta(x)^r\displaystyle\int_{U^x(F)}\int_{U_x(F)} {}^g\!f(v^{-1}xexp(X)uv) \xi(u) du dv \\
 & =\Delta(x)^r \int_{U^x(F)}\int_{\mathfrak{u}_x(F)} {}^{vg}\!f(xexp(X+N)) \xi(N) dN dv
 & = \Delta(x)^r \int_{U^x(F)} {}^{vg}\!f_{x,\omega}^\xi(X) dv
\end{aligned}\]

\noindent La flèche naturelle $U^x(F)\to U_x(F)\backslash U(F)$ est une bijection qui préserve les mesures on en déduit

$$J_N(\theta,f)=D^H(x)\Delta(x)^r\displaystyle\int_{U_x(F)H_x(F)\backslash G(F)} \kappa_N(g) \int_{\mathfrak{h}_x(F)} \theta_{x,\omega}(X) {}^g\!f^\xi_{x,\omega}(X) dX dg$$

\noindent On vérifie facilement, grâce au lemme 3.1.1, que $C(x)=D^H(x)\Delta(x)^r$, ce qui permet de conclure $\blacksquare$

\vspace{4mm}

\subsection{Localisation de $J_{geom}(\theta,f)$}

Pour $T\in \underline{\mathcal{T}}$ on notera maintenant $W'_T$ et $W''_T$ les sous-espaces que l'on avait noté $W'$ et $W''$ dans le paragraphe 5.1. Remarquons que l'ensemble des tores $\underline{\mathcal{T}}$ ne dépend que des espaces hermitiens $W$ et $D$. Quand on voudra préciser par rapport à quels espaces l'ensemble $\underline{\mathcal{T}}$ est défini on le notera plutôt $\underline{\mathcal{T}}(W,D)$. Soit $\underline{\mathcal{T}}_x$ l'ensemble des tores $T\in\underline{\mathcal{T}}$ qui vérifient $T\subset H_x$ et $W'\subset W'_T$ c'est équivalent à dire que $T=T'T''$ où $T'$ est un tore maximal anisotrope de $H'_x$ et $T''\in\underline{\mathcal{T}}(W'',D)$. On notera $\mathcal{T}_x$ un ensemble de représentants des classes de conjugaison par $H_x(F)$ dans $\underline{\mathcal{T}}_x$. Pour $T\in \underline{\mathcal{T}}_x$, on définit les fonctions $c_{\theta,x,\omega}$ et $c_{f,x,\omega}$ sur $\mathfrak{t}_\natural(F)$ de la façon suivante

$$c_{\theta,x,\omega}(X)= \left\{
    \begin{array}{ll}
        c_{\theta}(xexp(X)) & \mbox{si } X \in \omega \\
        0 & \mbox{sinon.}
    \end{array}
\right.$$

et

$$c_{f,x,\omega}(X)= \left\{
    \begin{array}{ll}
        c_f(xexp(X)) & \mbox{si } X \in \omega \\
        0 & \mbox{sinon.}
    \end{array}
\right.$$

On note $\Delta''$ la fonction définie sur $\mathfrak{h}_{x,ss}(F)$ par

$$\Delta''(X)=|N(det(X_{|W''/W''(X)}))|_F$$

où $W''(X)$ est le noyau de $X$ dans $W''$. On pose alors

\[\begin{aligned}
\displaystyle\mbox{(1)} \;\;\; & J_{geom,x,\omega}(\theta,f)= \\
 & \sum_{T\in \mathcal{T}_x} |W(H_x,T)|^{-1} \nu(T) \lim\limits_{s\to 0^+}\int_{\mathfrak{t}(F)} c_{\theta,x,\omega}(X)c_{f,x,\omega}(X)D^{H_x}(X)^{1/2}D^{G_x}(X)^{1/2}\Delta''(X)^{s-1/2}dX
\end{aligned}\]

\begin{lem}
On a l'égalité
$$J_{geom}(\theta,f)=C(x)J_{geom,x,\omega}(\theta,f)$$
\end{lem}

\vspace{3mm}

\ul{Preuve}: Il suffit encore une fois de reprendre la preuve du lemme 8.3 de [W1]. On a les propriétés suivantes

\vspace{2mm}

1) Pour $T\in \mathcal{T}$ et $t\in T_\natural(F)$, si $c_f(t)\neq 0$ on a
$$t\in\displaystyle\bigcup_{T_1\in \mathcal{T}_x}\bigcup_{w_1\in W(T_1,T)} w_1(xexp(\mathfrak{t}_1(F)\cap\omega))w_1^{-1}$$

\vspace{2mm}

En effet on a alors $t\in Supp(\theta_f)\subset \Omega$. Par conséquent il existe $X\in\omega$ et $y\in G(F)$ tel que $yty^{-1}=xexp(X)$. Encore une fois quitte à multiplier $y$ par un élément de $G''(F)$ on peut supposer que $y\in H(F)$. Alors $yTy^{-1}\in\underline{\mathcal{T}}$ et on a $W'\subset Ker(xexp(X)-1)^\perp =Ker(yty^{-1}-1)^\perp=W'_{yTy^{-1}}$ et $yTy^{-1}\subset Z_H(xexp(X))^0\subset  H_x$ donc $yTy^{-1}\in \underline{\mathcal{T}}_x$. Quitte à multiplier $y$ par un élément de $H_x(F)$ on a donc $yTy^{-1}\in \mathcal{T}_x$.

\vspace{2mm}

2) Soient $T_1,T_2\in \mathcal{T}_x$, $T\in \mathcal{T}$ et $w_1\in W(T_1,T)$, $w_2\in W(T_2,T)$, alors les ensembles $w_1(xexp(\mathfrak{t}_1(F)\cap\omega))w_1^{-1}$ et $w_2(xexp(\mathfrak{t}_2(F)\cap\omega))w_2^{-1}$ sont disjoints ou confondus et si ils sont confondus on a $T_1=T_2$.

\vspace{2mm}

3) Soient $T\in\mathcal{T}$, $T_1\in \mathcal{T}_x$ et $w_1\in W(T_1,T)$ alors le nombre d'éléments $w_2\in W(T_1,T)$ tels que $w_1(xexp(\mathfrak{t}_1(F)\cap\omega))w_1^{-1}=w_2(xexp(\mathfrak{t}_1(F)\cap\omega))w_2^{-1}$ est $W(H_x,T_1)$.

\vspace{2mm}

4) Soit $T_1\in \mathcal{T}_x$ alors il existe un unique $T\in \mathcal{T}$ tel que $W(T_1,T)\neq \emptyset$. On a alors $|W(T_1,T)|=|W(H,T)|$ et $\nu(T)=\nu(T_1)$.

\vspace{2mm}

Les points 2), 3) et 4) se démontrent exactement de la même manière que les points correspondants dans la démonstration du lemme 6.1.1. On déduit des trois points précédents que pour tout $s\in\mathbb{C}$ vérifiant $Re(s)>0$, on a

\[\begin{aligned}
\displaystyle\sum_{T\in\mathcal{T}} |W(H,T)|^{-1} \nu(T) \int_{\mathfrak{t}(F)} c_f(t)c_{\theta}(t)D^H(t)^{1/2}D^G(t)^{1/2}\Delta(t)^{s-1/2} dt &  \\
=\sum_{T_1\in\mathcal{T}_x}\sum_{T\in\mathcal{T}}\sum_{w_1\in W(T_1,T)} |W(H,T)|^{-1} |W(H_x,T_1)|^{-1} \nu(T_1) & \\
\int_{\mathfrak{t}_1(F)\cap\omega}  c_f(w_1xexp(X)w_1^{-1}) c_{\theta}(w_1xexp(X)w_1^{-1}) D^H(w_1xex & p(X)w_1^{-1})^{1/2} \\
D^{G}(w_1xexp(X) w_1^{-1} )^{1/2} \Delta(w_1xexp(X)w_1^{-1})^{s-1/2}  dX & 
\end{aligned}\]

Les fonctions $c_f$, $c_\theta$, $D^H$ et $\Delta$ sont invariantes par conjugaison par $H(F)$. On peut donc écrire

\[\begin{aligned}
\displaystyle\sum_{T\in\mathcal{T}} |W(H,T)|^{-1} \nu(T) \int_{\mathfrak{t}(F)} c_{\theta}(t)c_f(t & )D^H(t)^{1/2}D^G(t)^{1/2}\Delta(t)^{s-1/2} dt  \\
=\sum_{T_1\in\mathcal{T}_x} |W(H_x,T_1)|^{-1} \nu(T_1) \int_{\mathfrak{t}_1(F)\cap\omega} & c_f(xexp(X)) c_{\theta}(xexp(X)) \\
 & D^H(xexp(X))^{1/2} D^G(xexp(X))^{1/2} \Delta(xexp(X))^{s-1/2} dX
\end{aligned}\]

On a $D^H(xexp(X))=D^H(x)D^{H_x}(X)$, $D^G(xexp(X))=D^G(x)D^{G_x}(X)$ et $\Delta(xexp(X))=\Delta(x)\Delta''(X)$, par conséquent

\[\begin{aligned}
\displaystyle\sum_{T\in\mathcal{T}} |W(H,T)|^{-1} \nu(T) \int_{\mathfrak{t}(F)} c_{\theta}(t)c_f(t)D^H(t)^{1/2} D^G(t)^{1/2} & \Delta(t)^{s-1/2} dt \\
=C(x)\sum_{T_1\in\mathcal{T}_x} |W(H_x,T_1)|^{-1} \nu(T_1) \Delta(x)^s & \\
\int_{\mathfrak{t}_1(F)\cap\omega} c_{f,x,\omega}(X) c_{\theta,x,\omega}(X) D^{H_x}(X)^{1/2} D^{G_x}(X)^{1/2} & \Delta''(X)^{s-1/2} dX
\end{aligned}\]

Lorsque $s$ tend vers $0$, on obtient le résultat voulu. $\blacksquare$

\section{Utilisation de la transformée de Fourier}

Posons $U''=U\cap G''$. Soient $\varphi \in C_c^{\infty}(\mathfrak{g}''(F))$ et $\theta''$ un quasi-caractère de $\mathfrak{h}''(F)$. On pose

$$J_{\kappa''}(\theta'',\varphi)=\displaystyle\int_{U''(F)H''(F)\backslash G''(F)}\kappa''(g) J(\theta'',\varphi,g)dg$$

\noindent où 

$$J(\theta'',\varphi,g)= \displaystyle\int_{\mathfrak{h}''(F)} \theta''(X) \displaystyle\int_{\mathfrak{u}''(F)} {}^g\varphi(X+N) \xi(N) dN dX$$

\noindent On suppose de plus qu'il existe un élément $S\in\mathfrak{h}''_{reg}(F)$ de noyau nul dans $W''$ tel que l'on ait $\theta''(X)=\hat{j}^{H''}(S,X)$ pour tout $X\in \mathfrak{h}''(F)\cap Supp(\varphi)^{G''}$. \\
On va exprimer $J_{\kappa''}(\theta'',\varphi)$ en fonction de la transformée de Fourier de $\varphi$.

\vspace{4mm}

Fixons un élément $\eta\in E$ non nul de trace nulle. Posons $\nu_0=h(v_0,v_0)$. Soit $\Xi$ l'unique élément de $\overline{\mathfrak{u}''}(F)$ tel que pour tout $N\in \mathfrak{u}''(F)$ on ait $\xi(N)=\psi(<\Xi,N>)$. On a alors $\Xi(W'')=0$, \\
$\Xi v_{i+1}=\xi_iv_i$ pour $i\in\{0,\ldots,r-1\}$, $\Xi v_0=-\nu_0\xi_0v_{-1}$ et $\Xi v_{-i}=-\xi_iv_{-i-1}$ pour $i\in \{1,\ldots,r-1\}$. Soit $\Sigma$ l'orthogonal de $\mathfrak{u}''(F)\oplus \mathfrak{h}''(F)$ dans $\mathfrak{g}''(F)$. On a 

$$\Sigma=\mathfrak{a}(F)\oplus\mathfrak{u}''(F)\oplus \Lambda_0$$

\noindent où $\Lambda_0=Fc(v_0,\eta v_0)\oplus \{c(v_0,w)|w\in W''\}$.

\begin{lem}
Pour tout $\varphi\in C_c^{\infty}(\mathfrak{g}''(F))$ et tout $Y\in \mathfrak{h}''(F)$, on a l'égalité:
\begin{center}
$\widehat{\varphi^{\xi}}(Y)=\displaystyle\int_\Sigma \hat{\varphi}(\Xi+Y+X) dX$.
\end{center}
\end{lem}

\ul{Preuve}: C'est exactement la même que celle du lemme 9.2 de [W1] $\blacksquare$

\subsection{Etude  des classes de conjugaison dans $\Xi+S+\Sigma$}

On définit $\Lambda=\Lambda_0\oplus \Lambda_{\mathfrak{u}''}$ où $\Lambda_{\mathfrak{u}''}$ est le sous-$F$-espace vectoriel de $\mathfrak{u}''(F)$ engendré par les $c(v_i,v_{i+1})$ pour $i\in \{0,\ldots,r-1\}$ et les $c(v_i,\eta v_i)$ pour $i\in \{1,\ldots,r\}$. On a alors le:

\begin{lem}
$\Xi+S+\Sigma$ est stable par conjugaison par $U''$ et la conjugaison par $U''$ définit un isomorphisme de variétés algébriques:
\begin{center}
$U''\times \left(\Xi+S+\Lambda\right) \to \Xi+S+\Sigma$
\end{center}
\end{lem}

\ul{Preuve}: Soient $X\in \Sigma$, $u\in U''$, $U\in\mathfrak{u}''$ et $H\in\mathfrak{h}''$. Puisque $X$ est orthogonal à $\mathfrak{u}''\oplus\mathfrak{h}''$, $S$ est orthogonal à $\mathfrak{u}''$ et $\Xi$ est orthogonal à $\mathfrak{h}''$, on a les égalités

\[\begin{aligned}
<u(\Xi+S+X)u^{-1},U+H> & =<\Xi+S+X,u^{-1}(U+H)u> \\
& =<\Xi,u^{-1}Uu>+<\Xi,u^{-1}Hu-H>+<S,H>
\end{aligned}\]
Pour tout $N\in\mathfrak{u}''$, on a $<\Xi,N>=\displaystyle\sum_{j=0}^{r-1} \xi_jh(Nv_j,v_{-j-1})$. On en déduit facilement que \\
\noindent $<\Xi,u^{-1}Uu>=<\Xi,U>$ et $<\Xi,u^{-1}Hu-H>=0$. On a montré que 
$$<u(\Xi+S+X)u^{-1}-\Xi-S,U+H>=0$$ 
pour tout $H\in\mathfrak{h}''$ et tout $U\in\mathfrak{u}''$. Par conséquent, $u(\Xi+S+X)u^{-1}-\Xi-S$ est un élément de l'orthogonal de $\mathfrak{u}''\oplus\mathfrak{h}''$ c'est-à-dire de $\Sigma$. C'est la première partie du lemme.

\vspace{2mm}

 On définit les sous-groupes suivants de $U''$
\begin{itemize}
\item $U_1=U''$
\item $U_2=\{u\in U_1; u_{|Z_+}=Id\}$
\item $U_3=\{u\in U_2; u(v_0)=v_0\}$
\item $U_4=\{u\in U_3; u_{|W''}=Id\}$
\item $U_5=\{1\}$
\end{itemize}
On notera $\mathfrak{u}_i$ l'algèbre de Lie de $U_i$ pour $i=1,\ldots,5$. On définit aussi les sous-espaces suivants de $\Sigma$
\begin{itemize}
\item $\Sigma_1=\Sigma$
\item $\Sigma_2=\Lambda_0\oplus \mathfrak{u}_2\oplus \{c(v_{-1},v); v\in Z_+\}$
\item $\Sigma_3=\Lambda_0\oplus \mathfrak{u}_2$
\item $\Sigma_4=\Lambda_0\oplus \mathfrak{u}_4\oplus \{c(v_0,v); v\in Z_+\}$
\item $\Sigma_5=\Lambda$
\end{itemize}

\begin{center}
(1) $\Sigma_1,\Sigma_2$ et $\Sigma_3$ sont stables par conjugaison par $U_1$. L'ensemble $\Sigma_4$ est stable par conjugaison par $U_2$ et $\Sigma_5$ est stable par conjugaison par $U_4$
\end{center}

Pour $\Sigma_1$ c'est direct. Pour $\Sigma_3$, soit $P_2$ le sous-groupe parabolique de $G''$ des éléments qui stabilisent $Z_+$.Alors $\mathfrak{u}_2$ est l'algèbre de Lie du radical unipotent de $P_2$. Puisque $P''\subset P_2$, $\mathfrak{u}_2$ est stable par conjugaison par $P''$ donc a fortiori par $U_1$. De plus pour $u\in U_1$ et $X_{\Lambda_0}\in\Lambda_0$ on a $uX_{\Lambda_0}u^{-1}-X_{\Lambda_0}\in \mathfrak{u}_1$ et $(uX_{\Lambda_0}u^{-1}-X_{\Lambda_0})_{|Z_+}=0$. On en déduit que $uX_{\Lambda_0}u^{-1}\in \Sigma_3$ puis que $\Sigma_3$ est stable par conjugaison par $U_1$.  Pour $\Sigma_2$, soit $P_\sharp$ le sous-groupe parabolique de $G$ des éléments qui stabilisent $Ev_r\oplus\ldots\oplus Ev_2$ et $\mathfrak{u}_\sharp$ l'algèbre de Lie du radical unipotent de $P_\sharp$. On remarque que $\mathfrak{u}_2\oplus \{c(v_{-1},v); v\in Z_+\}=\mathfrak{u}_\sharp\oplus Ec(v_{-1},v_1)\subset \mathfrak{p}_\sharp$. Puisque $P''\subset P_\sharp$, la conjugaison par $U_1$ envoie $Ec(v_{-1},v_1)$ dans $Ec(v_{-1},v_1)\oplus\mathfrak{u}_\sharp$ et laisse stable $\mathfrak{u}_\sharp$. On a déjà vu que la conjugaison par $U_1$ envoie $\Lambda_0$ dans $\Sigma_3$. D'où le résultat pour $\Sigma_2$. \\
Pour $\Sigma_4$ et $\Sigma_5$, on remarque que $U_4$ est le centre de $U_2$ et on utilise la formule

\begin{center}
$gc(v,v')g^{-1}=c(gv,gv')$ pour tous $g\in G''$, $v,v'\in V''$
\end{center}

On a les caractérisations suivantes de $\Sigma_{i+1}$ dans $\Sigma_i$ pour $i=1,\ldots,3$

\begin{itemize}
\item $\Sigma_2$ est l'ensemble des $X\in\Sigma_1$ tels que $Xv_i=0$ pour $i=2,\ldots,r$
\item $\Sigma_3$ est l'ensemble des $X\in\Sigma_2$ tels que $Xv_1=0$
\item $\Sigma_4$ est l'ensemble des $X\in\Sigma_3$ tels que $Xw=0$ pour tout $w\in W''$
\end{itemize}

\begin{center}
(2) Pour $i=1,\ldots,4$ on a $u_i(\Xi+S)u_i^{-1}-(\Xi+S)\in \Sigma_i$ pour tout $u_i\in U_i$
\end{center}

En effet, pour $i=1$ c'est la première partie du lemme. Pour $i=2,3,4$, il suffit d'utiliser les caractérisations précédentes de $\Sigma_i$ dans $\Sigma_{i-1}$. \\
D'après (1) et (2) $\Xi+S+\Sigma_i$ est stable par conjugaison par $U_i$ pour $i=1,\ldots,5$. Pour $i=1,\ldots,4$, on peut considérer le quotient $U_i\times_{U_{i+1}} (\Xi+S+\Sigma_{i+1})$ de $U_i\times (\Xi+S+\Sigma_{i+1})$ pour la relation d'équivalence définie par $(u_iv,X)\sim (u_i,vXv^{-1})$ pour tous $u_i\in U_i$, $v\in U_{i+1}$, $X\in\Xi+S+\Sigma_{i+1}$. Le lemme découlera alors directement de

\begin{center}
(3) Pour $i=1,\ldots,4$, la conjugaison de $U_i$ sur $\Xi+S+\Sigma_{i+1}$ se descend en un isomorphisme entre $U_i\times_{U_{i+1}} (\Xi+S+\Sigma_{i+1})$ et $\Xi+S+\Sigma_i$
\end{center}

Pour $i=4$, $U_4$ centralise $\Sigma_5$ et on a un isomorphisme
$$\mathfrak{u}_4\to U_4$$
$$Y\mapsto I+Y$$
Il suffit donc de montrer que l'application
$$\mathfrak{u}_4\times \Sigma_5\to \Sigma_4$$
$$(Y,X)\mapsto (I+Y)(\Xi+S)(I-Y)-(\Xi+S)+X=Y\Xi-\Xi Y+X$$
est un isomorphisme. Il s'agit d'une application linéaire entre deux espaces vectoriels dont on vérifie facilement qu'ils ont même dimension. Il suffit donc de vérifier que c'est une application injective. Soit $Y\in \mathfrak{u}_4$ tel que $Y\Xi-\Xi Y\in \Sigma_5$ et montrons par récurrence descendante que $Yv_{-i}=0$ pour $i=1,\ldots,r$, ce qui établira l'injectivité de l'application. \\

Puisque $Y\Xi-\Xi Y\in \Sigma_5$ on a $-\Xi Yv_{-r}=(Y\Xi-\Xi Y)v_{-r}\in F\eta v_r+Fv_{r-1}$. Or on a 
$$h(\Xi Yv_{-r},v_{-r})=-h(Yv_{-r},\Xi v_{-r})=0$$
et 
$$h(\Xi Yv_{-r},v_{-r+1})=-h(Yv_{-r},\Xi v_{-r+1})=\xi_{r-1}h(Yv_{-r},v_{-r})\in F\eta$$

D'où $\Xi Yv_{-r}=0$ et puisque $Y(Z_-)\subset Z_+$ et que $\Xi$ est injective sur $Z_+$ on a $Y v_{-r}=0$. \\

Soit $1\leqslant i\leqslant r-1$. Supposons que $Yv_{-j}=0$ pour $j=i+1,\ldots,r$. \\
On a alors $(Y\Xi-\Xi Y)v_{-i}=\xi_i Yv_{-i-1}-\Xi Yv_{-i}=-\Xi Yv_{-i}$. Comme d'autre part on a aussi $Y\Xi-\Xi Y\in \Sigma_5$, $\Xi Yv_{-i}$ appartient à $Fv_{-i+1}+Fv_{-i-1}+F\eta v_{-i}$. Pour avoir $Yv_{-i}=0$, il suffit donc de vérifier que $h(\Xi Yv_{-i},v_{-i})=h(\Xi Yv_{-i},v_{-i-1})=h(\Xi Yv_{-i},v_{-i+1})=0$. On a 
$$h(\Xi Yv_{-i},v_{-i})=h(v_{-i},Y\Xi v_{-i})=-\xi_i h(v_{-i},Yv_{-i-1})=0$$, 
$$h(\Xi Yv_{-i},v_{-i-1})=h(v_{-i},Y\Xi v_{-i-1})=0$$
car $\Xi v_{-i-1}\in Ev_{-i-2}$ et 
$$h(\Xi Yv_{-i},v_{-i+1})=-h(Yv_{-i},\Xi v_{-i+1})\in F\eta$$
puisque $\Xi v_{-i+1}\in Fv_{-i}$. \\
On en déduit $h(\Xi Yv_{-i},v_{-i})=h(\Xi Yv_{-i},v_{-i-1})=h(\Xi Yv_{-i},v_{-i+1})=0$. \\

Pour $i=3$, $U_3$ normalise $\Sigma_4$ d'après (1) et on a un isomorphisme
$$Hom_E(W'',Z_+)\times U_4\to U_3$$
$$(Y,u_4)\to exp(Y-Y^*)u_4$$
Où pour $l$ un endomorphisme linéaire de $V''$ on a noté $l^*\in End_F(V'')$ l'application duale (on a identifié $V''$ à son propre dual via $h$). Il suffit donc de montrer que le morphisme suivant est un isomorphisme
$$Hom_E(W'',Z_+)\to \Sigma_3/\Sigma_4$$
$$Y\mapsto exp(Y-Y^*)(\Xi+S)exp(Y^*-Y)-(\Xi+S)$$
Soit $\underline{\Xi}$ l'endomorphisme de $Z_+$ qui vaut $\Xi$ sur $v_2,\ldots,v_r$ et qui envoie $v_1$ sur $0$. $\Sigma_3/\Sigma_4$ est aussi isomorphe à $Hom_E(W'',Z_+)$ par $X\mapsto p_{Z_+}\circ X_{|W''}$ où $p_{Z_+}$ est la projection sur $Z_+$ parallèlement à $Z_-\oplus D\oplus W''$. Via cet isomorphisme, l'application précédente devient
$$Hom_E(W'',Z_+)\to Hom_E(W'',Z_+)$$
$$Y\mapsto \underline{\Xi}Y-YS$$
Les endomorphismes linéaires de $Hom_E(W'',Z_+)$, $Y\mapsto YS$ et $Y\mapsto \underline{\Xi}Y$ commutent, le premier est semi-simple inversible et le second est nilpotent. On en déduit que l'application précédente est un isomorphisme. \\

Pour $i=2$, $U_2$ normalise $\Sigma_3$ et on a un isomorphisme
$$Hom_E(D,Z_+)\times U_3\to U_2$$
$$(Y,u_3)\mapsto exp(Y-Y^*) u_3$$
Comme pour $i=3$, il suffit donc de montrer que $Hom_E(D,Z_+)\to \Sigma_2/\Sigma_3$, $Y\mapsto exp(Y-Y^*)(\Xi+S)exp(Y^*-Y)$ est un isomorphisme. L'application $\Lambda_2/\Lambda_3\to Z_+$, $X\mapsto Xv_1$ est un isomorphisme. Via cet isomorphisme l'application précédente devient
$$Hom_E(D,Z_+)\to Z_+$$
$$Y\mapsto \xi_0Yv_0$$
qui est clairement un isomorphisme. \\

Pour $i=1$, soit $M_2$ la composante de Lévi de $P_2$ des éléments qui stabilisent $Z_+$ et $Z_-$ et posons $U_B=U_1\cap M_2$. Ce groupe s'identifie au radical unipotent du sous-groupe de Borel $B$ de $GL_E(Z_+)$ qui conserve le drapeau
$$Ev_r\subset \ldots\subset Ev_r\oplus\ldots\oplus Ev_1$$

L'application produit de $U_B\times U_2$ sur $U_1$ est un isomorphisme. Il suffit donc de montrer que
$$U_B\times (\Xi+S+\Sigma_2)\to \Xi+S+\Sigma_1$$
$$(u,X)\mapsto uXu^{-1}$$
est un isomorphisme. Soit $\overline{\Xi}\in End_E(Z_+)$ l'élément qui envoie $v_i$ sur $\Xi v_i$ pour $i=2,\ldots,r$ et $v_1$ sur $0$. Notons $\mathfrak{r}$ l'espace vectoriel engendré par les $c(v_{-1},v)$ pour $v\in Z_+$. On a alors $\Sigma_2=\mathfrak{r}\oplus \Sigma_3$ et $\Sigma_1=\mathfrak{b}\oplus \Sigma_3$. Puisque $\Sigma_3$ est invariant par conjugaison par $U_B$ et $S$ commute à $U_B$, on est ramené à prouver que l'application suivante est un isomorphisme
$$U_B\times(\overline{\Xi}+\mathfrak{r})\to \overline{\Xi}+\mathfrak{b}$$
$$(u,X)\mapsto uXu^{-1}$$
Tout se passe maintenant dans $GL_E(Z_+)$ et l'assertion est bien connue $\blacksquare$ \\

Donnons nous un élément $X\in \Lambda$, on peut l'écrire sous la forme
\begin{center}
(1) $X=c(v_0,w)+\displaystyle\sum_{i=0,\ldots,r-1} \lambda_ic(v_i,v_{i+1})+\displaystyle\sum_{i=0,\ldots,r} \mu_i c(v_i,\eta v_i)$
\end{center}
où $w$ est un vecteur de $W''$ et les $\lambda_i$, $\mu_i$ sont des éléments de $F$. Soit $D$ l'endomorphisme $E$-linéaire de $E[T]$ qui envoie $T^{i+1}$ sur $T^i$ et $1$ sur 0. Notons $R_S(T)$ le polynôme caractéristique de $S$ agissant sur $W''$. On a alors:

\begin{lem}
Le polynôme caractéristique de $\Xi+S+X$ agissant sur $V''$ est:
\begin{center}
$P_{\Xi+S+X}(T)=T^{2r}\displaystyle\sum \limits_{i=0}^{d_{W''}-1} \nu_0 h(S^i w,w) D^{i+1}(R_S(T))+R_S(T)\left[T^{2r+1}+\displaystyle\sum \limits_{i=0}^{r-1} (-1)^{i+1}2\nu_0 T^{2r-1-2i}\lambda_i\xi_i\left (\displaystyle\prod_{j=0}^{i-1} \xi_j^2\right) + \displaystyle\sum \limits_{i=0}^r(-1)^{i+1}T^{2r-2i}\nu_0\eta\mu_i \left(\displaystyle\prod_{j=0}^{i-1}\xi_j^2\right)\right]$
\end{center}
\end{lem}

\ul{Preuve}: Fastidieuse mais directe, elle est laissée au lecteur $\blacksquare$

\begin{lem}
Les polynômes $T^{2r} D^{i}(R_S(T))$, $i=1,\ldots,d_{W''}$, $R_S(T)T^j$, $j=0,\ldots,2r+1$ forment un base de l'espace des polynômes sur $E$ de degré $\leqslant 2r+d_{W''}+1$
\end{lem}

\ul{Preuve}: Notons $r(P)$ le reste de la division euclidienne de $P$ par $T^{2r}$. Il suffit de montrer que les familles $$(T^{2r+1}R_S(T),T^{2r}R_S(T),T^{2r}D(R_S(T)),\ldots,T^{2r}D^{d_{W''}}(R_S(T)))$$
 
et
 
$$(r(R_S(T)),\ldots,r(T^{2r-1}R_S(T)))$$

 sont des familles libres. Dans la première famille les degrés sont strictement décroissant et dans la deuxième l'indice du premier coefficient non nul est strictement croissant. D'où le résultat $\blacksquare$

\vspace{2mm}

Soit $\Lambda^S$ l'ensemble des éléments $X$ tels que $\Xi+S+X\in\mathfrak{g}''_{reg}(F)$, $\lambda_i\neq 0$, $\mu_i\neq 0$ et la famille $(w,Sw,\ldots,S^{d_{W''}-1}w)$ engendre $W''$. On déduit des deux lemmes précédents que l'application $\Lambda^S\to \{P\in E[T]| (-1)^{d''}\overline{P}(-T)=P(T)\} $ qui à $X$ associe le polynôme caractéristique de $\Xi+S+X$ est partout submersive. On notera $\Xi+S+\Sigma^S$ l'ensemble des éléments de $\Xi+S+\Sigma$ qui sont conjugués par $U''$ à un élément de $\Xi+S+\Lambda^S$. On a alors:

\begin{lem}
Le groupe $H_S''U''$ agit librement par conjugaison sur $\Xi+S+\Sigma^S$ et deux éléments de ce dernier sont conjugués par un élément de $G''$ si et seulement si ils le sont par un élément de $H_S''U''$.
\end{lem}

\vspace{2mm}

\ul{Preuve}: La conjugaison par $H_S''$ laisse stable $\Xi$, $S$ et $\Lambda^S$ donc d'après le lemme 7.1.1 il suffit de montrer que l'action de $H_S''$ sur $\Lambda^S$ est libre pour obtenir la première assertion du lemme. Pour $X=c(v_0,w)+\displaystyle\sum_{i=0,\ldots,r-1} \lambda_ic(v_i,v_{i+1})+\displaystyle\sum_{i=0,\ldots,r} \mu_i c(v_i,\eta v_i)\in \Lambda^S$ et $h\in H_S''$ on a $hXh^{-1}=c(v_0,h(w))+\displaystyle\sum_{i=0,\ldots,r-1} \lambda_ic(v_i,v_{i+1})+\displaystyle\sum_{i=0,\ldots,r} \mu_i c(v_i,\eta v_i)$. Par conséquent $hXh^{-1}=X$ si et seulement si $hw=w$ mais les conditions $h\in H_S''$ et $(w,Sw,\ldots,S^{d_{W''}-1}w)$ engendre $W''$ montre que cette condition est équivalente à $h=1$. \\

Soient $X,X'\in \Xi+S+\Sigma^S$ deux éléments conjugués par $G''$ et montrons qu'ils le sont par $H''_SU''$. D'après le lemme 7.1.1, on peut supposer que $X$ et $X'$ sont dans $\Xi+S+\Lambda^S$ et montrer qu'ils sont conjugués par un élément de $H''_S$. On peut alors écrire

$$X=\Xi+S+c(v_0,w)+\displaystyle\sum_{i=0,\ldots,r-1} \lambda_ic(v_i,v_{i+1})+\displaystyle\sum_{i=0,\ldots,r} \mu_i c(v_i,\eta v_i)$$
et
$$X'=\Xi+S+c(v_0,w')+\displaystyle\sum_{i=0,\ldots,r-1} \lambda_i'c(v_i,v_{i+1})+\displaystyle\sum_{i=0,\ldots,r} \mu_i' c(v_i,\eta v_i)$$

Puisque $X$ et $X'$ sont conjugués sous $G''$ ils ont même polynôme caractéristique. D'après le lemme 7.1.2, on a donc $\lambda_i=\lambda_i'$ et $\mu_i=\mu_i'$ pour tout $i$ et $h(S^jw,w)=h(S^jw',w')$ pour $0\leqslant j\leqslant d_{W''}-1$. On a alors

$$h(S^jw,S^iw)=(-1)^ih(S^{j-i}w,w)=(-1)^ih(S^{j-i}w',w')=h(S^jw',S^iw')$$
pour $0\leqslant i\leqslant j\leqslant d_{W''}-1$. Par conséquent il existe un unique élément $h\in H''$ qui envoie la base $(w,Sw,\ldots,S^{d_{W''}-1}w)$ sur la base $(w',Sw',\ldots,S^{d_{W''}-1}w')$ et on vérifie facilement que $h\in H''_S$ et $hXh^{-1}=X'$ $\blacksquare$

\vspace{5mm}

\subsection{Un calcul de Jacobien}

Pour $T\in\mathcal{T}(G'')$, on notera $\mathfrak{t}(F)^S$ le sous-ensemble de $\mathfrak{t}(F)$ des éléments conjugués par $G''(F)$ à un élément de $\Xi+S+\Sigma^S$ c'est un sous-ensemble ouvert. On a alors une application analytique bijective
\begin{center}
$\Xi+S+\Sigma^S/H_S''(F)U''(F)\to \displaystyle\bigsqcup_{T\in\mathcal{T}(G'')} \mathfrak{t}(F)^S/W(G'',T)$
\end{center}
Puisque l'application qui à un élément de $\Xi+S+\Sigma^S$ associe son polynôme caractéristique est partout submersive, l'application précédente est elle aussi partout submersive. On en déduit une mesure sur l'espace d'arrivée en transférant celle fixée sur l'espace de départ, si on note $d_\Sigma Y$ cette nouvelle mesure, on a alors:

\begin{lem}
$d_\Sigma Y=D^{H''}(S)^{-\frac{1}{2}}D^{G''}(Y)^{\frac{1}{2}}dY$
\end{lem}

\vspace{4mm}

\ul{Preuve}:  Plaçons nous sur $\overline{F}$. Soit $T\in\mathcal{T}(G'')$ et notons $\mathfrak{t}^S$ l'ensemble des éléments de $\mathfrak{t}$ qui sont conjugués à un éléments de $\Xi+S+\Sigma^S$. C'est un ouvert algébrique de $\mathfrak{t}^S$. On considère aussi $\Sigma^S$ comme un ouvert algébrique dans $\Sigma$. Soit $\overline{W}(G'',T)$ le groupe de Weyl de $T$ sur $\overline{F}$ alors l'application

\begin{center}
(1) $\Xi+S+\Sigma^S/H''_SU''\to \mathfrak{t}^S/\overline{W}(G'',T)$
\end{center}

\noindent est bien définie et est un isomorphisme de variétés algébriques. De plus, l'application sur les $F$-points

\begin{center}
(2) $\mathfrak{t}(F)^S/W(G'',T)\to \mathfrak{t}^S/\overline{W}(G'',T)(F)$
\end{center}

\noindent est un isomorphisme local de variétés analytiques. La forme bilinéaire non dégénérée $<.,.>$ sur $\mathfrak{g}$ permet d'obtenir sur chaque sous-espace non dégénéré de $\mathfrak{g}''$ une forme différentielle de degré maximale invariante bien définie au signe près. Sur les autres sous-espaces de $\mathfrak{g}''$, on choisit une forme différentielle de degré maximale invariante quelconque. Cela permet de fixer une telle forme différentielle sur $\mathfrak{u}''$ mais ce choix n'importe pas. On relève les formes différentielles invariantes sur les algèbres de Lie aux groupes. Alors l'application (2) préserve localement les mesures. On peut donc se contenter de calculer le jacobien algébrique de (1). Comme on vient de l'expliquer, ce jacobien n'est bien défini qu'au signe près, mais cela ne change pas sa valeur absolue. \\

Soit $X\in\mathfrak{g}''$ un élément régulier semisimple. Fixons une sous-algèbre de Borel $\mathfrak{b}_X$ qui contient $X$ et notons $\mathfrak{u}_X$ son radical nilpotent. Alors $det(ad(X)_{|\mathfrak{u}_X})$ ne dépend du choix de $\mathfrak{b}_X$ qu'au signe près. Comme les calculs qui suivent sont tous au signe près, on prend la liberté de noter $d^{G''}(X)$ ce déterminant sans référence au choix d'une sous-algèbre de Borel. On définit de même $d^{G_0''}$ et $d^{H''}$. Remarquons que l'on a $|d^{G''}(X)|_F=D^{G''}(X)^{1/2}$. \\

Supposons d'abord $r=0$. On a alors $U''=\{1\}$ et $\Xi=0$. Sur $\overline{F}$, $G_0''$ s'identifie au groupe linéaire de $V_0''\otimes_E \overline{F}$. Posons $d=d_{W''}$ et soit $w_1,\ldots,w_d$ une base de $W''\otimes_E \overline{F}$ constituée de vecteurs propres pour $S$. Dans la base $v_0,w_1,\ldots,w_d$, $\Sigma$ est l'ensemble des matrices de la forme

$$\begin{pmatrix}
\lambda_0 & \lambda_1 & \ldots & \lambda_d \\
\mu_1 & 0 & \ldots & 0 \\
\vdots & \vdots & & \vdots \\
\mu_d & 0 & \ldots & 0
\end{pmatrix}$$

Un simple calcul montre qu'avec ces notations la forme différentielle autoduale sur $\Sigma$ est au signe près
$$d\lambda_0 \wedge d\lambda_1\ldots \wedge d\lambda_d \wedge d\mu_1\ldots \wedge d\mu_d$$
$\Sigma^S$ est l'ensemble des éléments de $\Sigma$ pour lesquels les $\lambda_i$ et $\mu_j$ sont tous non nuls et $H''_S$ est le tore des éléments qui fixent les droites $\overline{F} w_i$ pour $i=1,\ldots,d$. Pour $h\in H''_S$ on note $h_i$ la valeur propre de $h$ agissant sur $w_i$. La mesure autoduale sur $H''_S$ est alors au signe près $\displaystyle\frac{dh_1\wedge\ldots \wedge dh_d}{h_1\ldots h_d}$. L'action par conjugaison de $H''_S$ sur $\Sigma$ est donnée par la formule

$$h\begin{pmatrix}\lambda_0 & \lambda_1 & \ldots & \lambda_d \\ \mu_1 & 0 & \ldots & 0 \\ \vdots & \vdots & & \vdots \\ \mu_d & 0 & \ldots & 0 \end{pmatrix}h^{-1}=\begin{pmatrix}
\lambda_0 & h_1^{-1}\lambda_1 & \ldots & h_d^{-1}\lambda_d \\ h_1\mu_1 & 0 & \ldots & 0 \\ \vdots & \vdots & & \vdots \\ h_d\mu_d & 0 & \ldots & 0
\end{pmatrix}$$

$\Sigma^S/H''_S$ s'identifie donc à $F\times \overline{F}^{\times,d}$ par l'application 

$$(\lambda_0,\ldots,\lambda_d)\mapsto \begin{pmatrix}
\lambda_0 & \lambda_1 & \ldots & \lambda_d \\
1 & 0 & \ldots & 0 \\
\vdots & \vdots & & \vdots \\
1 & 0 & \ldots & 0
\end{pmatrix}$$

La forme différentielle autoduale sur $\Sigma^S/H''_S$ étant alors $d\lambda_0\wedge \ldots \wedge d\lambda_d$ (toujours au signe près). Pour $i=1,\ldots,d$ soit $\alpha_i$ la valeur propre de $S$ agissant sur $w_i$. Le polynôme caractéristique de

$$S+\begin{pmatrix}
\lambda_0 & \lambda_1 & \ldots & \lambda_d \\
1 & 0 & \ldots & 0 \\
\vdots & \vdots & & \vdots \\
1 & 0 & \ldots & 0
\end{pmatrix}$$

est d'après le lemme 7.1.2

$$(T-\lambda_0)R_S(T)+\displaystyle\sum_{i=0}^{d-1} (-1)^{i+1} (\alpha_1^i\lambda_1+\ldots+\alpha_d^i\lambda_d) D^{i+1}(R_S(T)) $$

On munit l'espace des polynômes unitaires de degré $d+1$ de la forme différentielle $\delta=dc_0\wedge dc_1\wedge\ldots\wedge dc_{d-1}$ où $c_i$ est le coefficient de degré $i$ du polynôme. Alors le jacobien de l'application qui à un élément de $S+\Sigma^S/H''_S$ associe son polynôme caractéristique est au signe près le déterminant de la matrice suivante

$$\begin{pmatrix} 1 & \ldots & 1 \\
\alpha_1 & \ldots & \alpha_d \\
\vdots &  & \vdots \\
\alpha_1^{d-1} & \ldots & \alpha_d^{d-1}
\end{pmatrix}$$

C'est $d^{H''}(S)$ . Soient $e_1,\ldots,e_{d+1}$ une base commune de diagonalisation des éléments de $\mathfrak{t}$. Pour $Y\in \mathfrak{t}$ on note $Y_i$ la valeur propre de $X$ agissant sur $e_i$. Alors la forme différentielle autoduale sur $\mathfrak{t}$ est au signe près $\delta_\mathfrak{t}=dY_1\wedge \ldots\wedge dY_{d+1}$. On calcule aisément le jacobien de l'application qui à $Y\in\mathfrak{t}$ associe son polynôme caractéristique (c'est essentiellement le même calcul de déterminant que précédemment). Ce jacobien vaut $d^{G''}(Y)$. Le jacobien de l'application $S+\Sigma^S/H''_S\to \mathfrak{t}/\overline{W}(G'',T)$ est donc $d^{G''}(Y)d^{H''}(S)^{-1}$. \\

Traitons maintenant le cas général. L'application

$$U''\times (S+(\Lambda_0\oplus\mathfrak{a}))\to \Xi+S+\Sigma$$
$$(u,X)\mapsto u(\Xi+X)u^{-1}$$

est un isomorphisme local sur l'ouvert dense des $(u,X)$ où $X$ est régulier dans $G''$. On peut facilement calculer son jacobien. L'application étant équivariante pour l'action de $U''$ évidente sur chacun des deux espaces et les formes différentielles étant invariantes par action de $U''$, on peut se contenter de calculer le jacobien de l'application en $(1,X)$. La différentielle de l'application en $(1,X)$ est

$$\mathfrak{u}''\oplus \Lambda_0\oplus\mathfrak{a}\to \Sigma$$
$$(N,X_{\Lambda_0}+A)\mapsto [N,\Xi]+[N,X]+X_{\Lambda_0}+A$$

Dans une base de $\mathfrak{u}''$ bien choisie l'application $N\mapsto [N,X]$ est diagonale et l'application qui à $N$ associe la projection sur $\mathfrak{u}''$ de $[N,\Xi]$ est triangulaire supérieure stricte. Le jacobien est donc le déterminant de l'application $N\mapsto [N,X]$. Au signe près c'est $d^{G''}(X)d^{G''_0}(X)^{-1}$. Sur un ouvert dense l'application

\begin{center}
(3) $\mathfrak{a}\times (S+\Lambda_0/H_S'')\to \left(\Xi+S+\Sigma\right)/H_S''U''$
\end{center}

est donc un isomorphisme local de jacobien $d^{G''}(d^{G_0''})^{-1}$. On peut toujours supposer quitte à conjuguer que $\mathfrak{t}$ est de la forme $\mathfrak{a}\times\mathfrak{t}_0$ où $\mathfrak{t}_0$ est un tore maximal de $\mathfrak{g}''_0$. Puisque pour $X\in S+(\Lambda_0\oplus\mathfrak{a})$ régulier $X+\Xi$ est conjugué à $X$, la composée de l'application (3) et de l'application

$$\left(\Xi+S+\Sigma\right)/H_S''U''\to \mathfrak{t}/\overline{W}(G'',T)$$

est alors le produit de l'application identité sur $\mathfrak{a}$ et de l'application

$$S+\Lambda_0/H_S''\to \mathfrak{t}_0/\overline{W}(G''_0,T_0)$$

On a déjà calculé le jacobien de cette application: au signe près c'est $d^{H''}(S)(d^{G_0''})^{-1}$. On en déduit que le jacobien de (1) est au signe près $(d^{G''})^{-1}d^{G_0''}d^{H''}(S)(d^{G_0''})^{-1}=d^{H''}(S)(d^{G''})^{-1}$ $\blacksquare$

\subsection{Choix de sections localement analytiques}

Soit $T\in \mathcal{T}(G'')$. Puisque l'application $(\Xi+S+\Lambda^S)\to\displaystyle\bigsqcup_{T\in\mathcal{T}(G'')} \mathfrak{t}(F)^S/W(G'',T)$ est partout submersive, on peut trouver une section localement analytique $\mathfrak{t}(F)^S/W(G'',T)\to \Xi+S+\Lambda^S$. En composant avec la projection naturelle $\mathfrak{t}(F)^S\to \mathfrak{t}(F)^S/W(G'',T)$ on obtient une application localement analytique $\mathfrak{t}(F)^S\to\Xi+S+\Lambda^S$, $X\mapsto X_\Lambda$. Montrons que l'on peut choisir l'application $X\mapsto X_\Lambda$ de sorte que pour tout compact $\omega_T\subset \mathfrak{t}(F)$ elle envoie $\omega_T\cap\mathfrak{t}(F)^S$ dans un sous-ensemble compact de $\Xi+S+\Lambda$.

\vspace{2mm}

\ul{Preuve}: Soient $\lambda_i$, $\mu_i$ et $w\in W''$ les coordonnées de $X_\Lambda$ comme introduites en 7.1(1). Puisque les $\lambda_i$ et les $\mu_i$ sont des fonctions linéaires en les coefficients du polynôme caractéristique ces coordonnées seront bornées pour $X\in\omega_T$. Il suffit donc de montrer qu'on peut choisir l'application $X\mapsto X_\Lambda$ de sorte que la coordonnée $w$ reste bornée pour $X\in \omega_T$. Soit $W''^S$ l'ensemble des éléments $w\in W''$ tels que $(w,Sw,\ldots,S^{d_{W''}-1}w)$ engendre $W''$. D'après les lemmes 7.1.2 et 7.1.4 il suffit de montrer que l'application

$$W''^S\to E^{d_{W''}}$$
$$w\mapsto (h(w,w),\ldots,h(w,S^{d_{W''}-1}w))$$

admet une section localement analytique sur son image qui est bornée sur tout sous-ensemble compact de $E^{d_{W''}}$. D'après [W4]1.3, on peut décomposer $W''$ en somme directe de sous-espaces non dégénérés stables par $S$ qui sont isomorphes à des espaces hermitiens de la forme $F_i=F_i^\sharp\otimes_F E$ où $F_i^\sharp$ est une extension finie de $F$, la forme hermitienne étant donnée par $(x,y)\mapsto Tr_{F_i/E}(c_iy\tau_E(x))$ où $\tau_E$ est l'automorphisme de $F_i$ définit par $\tau_E(x\otimes y)=x\otimes \overline{y}$, $c_i\in F_i^\sharp$ et tel que l'action de $S$ sur chacun de ces sous-espaces soit la multiplication par un élément $a_i\in F_i$ vérifiant $\tau_E(a_i)=-a_i$ et $E[a_i]=F_i$. On se ramène alors facilement au cas où le couple $(W'',S)$ est de cette sorte pour une certaine extension $L^\sharp$ de $F$ et des éléments $a,c\in L=L^\sharp\otimes_F E$. L'application $W''^S\to E^{d_{W''}}$ est alors $x\mapsto (Tr(ca^kx\tau_E(x)))_{k=0,\ldots,d_{W''}-1}$. C'est l'application qui donne les coordonnées de $x\tau_E(x)$ dans la base duale de $(c,ca,\ldots,ca^{d_{W''}-1})$. On est donc ramené à prouver que l'application $N_{L/L^\sharp}:L^\times\to (L^\sharp)^{\times}$ admet une section analytique sur son image qui est bornée sur tout compact. Si l'extension $L/L^\sharp$ est une extension de corps c'est vrai pour toute section continue. Si $L=L^\sharp\times L^\sharp$ il suffit de prendre la section $x\mapsto (x,1)$. $\blacksquare$

\vspace{2mm}

Pour tout $T\in\mathcal{T}(G'')$ on peut trouver une application localement analytique $\mathfrak{t}(F)^S\to G''(F)$, $X\mapsto \gamma_{X}$ telle que pour $ X\in\mathfrak{t}(F)$, $X_{\Lambda}=\gamma_{X}X\gamma_{X}^{-1}\in \Xi+S+\Lambda^S$. On suppose que les fonctions analytiques $X\mapsto X_{\Lambda}$ et $X\mapsto \gamma_X$ vérifient les deux conditions suivantes

\begin{enumerate}
\item  Pour tout sous-ensemble compact $\omega_T\subset \mathfrak{t}(F)$, la fonction $X\mapsto X_\Lambda$ pour $X\in \mathfrak{t}(F)^S\cap \omega_T$ prend ses valeurs dans un sous-ensemble compact de $\Xi+S+\Lambda$.
\item  Pour tout sous-ensemble compact $\omega_T\subset \mathfrak{t}(F)$, il existe une constante $c>0$ telle que pour tout $X\in\mathfrak{t}(F)^S\cap \omega_T$ on ait $\sigma(\gamma_X)\leqslant c(1+|log D^{G''}(X)|)$.
\end{enumerate}

Le premier point est loisible d'après ce que l'on vient de voir et le deuxième découle du premier d'après le lemme 4.2 de [A3] que voici

\begin{lem}[4.2 de \text{[A3]}]
Soit $\Gamma\subset G''(F)$ et $\omega_T\subset \mathfrak{t}(F)$ des sous-ensembles compacts. Il existe $c>0$ de sorte que pour tout $X\in \omega_T\cap \mathfrak{g}''_{reg}(F)$ et pour tout $g\in G''(F)$ tels que $gXg^{-1}\in \Gamma$ on ait

$$inf\{\sigma(gt); \; t\in T(F)\}\leqslant c(1+|log\; D^{G''}(X)|)$$

\end{lem}

On a alors le lemme suivant qui nous servira par la suite:

\begin{lem}
\begin{itemize}
\item Pour tout $X\in\mathfrak{t}^S(F)$ la famille $(v_r,X_\Lambda v_r,\ldots,X_\Lambda^{d''-1} v_r)$ est une base de $V''$. On note $R_{X_\Lambda,v_r}$ le $\mathcal{O}_E$-réseau engendré par cette base.
\item Soit $R''$ un $\mathcal{O}_E$-réseau de $V''$. Il existe une fonction polynomiale non nulle $Q_S$ sur $\mathfrak{t}(F)$ et pour tout sous-ensemble compact $\omega_T\subset\mathfrak{t}(F)$ une constante $\alpha\in F^\times$ telles que

\begin{center}
$\forall X\in\mathfrak{t}(F)^S\cap\omega_T$, $\alpha Q_S(X)R''\subset R_{X''_\Lambda,v_r}$
\end{center}

\end{itemize}
\end{lem}

\ul{Preuve}: Soit $X\in\mathfrak{t}(F)^S$. On peut écrire $X_\Lambda$ sous la forme

$$X_\Lambda=c(v_0,w)+\displaystyle\sum_{i=0,\ldots,r-1} \lambda_ic(v_i,v_{i+1})+\displaystyle\sum_{i=0,\ldots,r} \mu_i c(v_i,\eta v_i)$$

D'après le lemme 7.1.2 les fonctions $X\mapsto \lambda_i$, $X\mapsto\mu_i$ et $X\mapsto h(S^iw,S^jw)$ sont polynomiales. Soit $R_{X_\Lambda,v_r}$ le sous-$\mathcal{O}_E$-module de $V''$ engendré par $v_r,X_\Lambda v_r,\ldots,X_\Lambda^{d''-1} v_r$ et soit $Q_S$ la fonction polynomiale sur $\mathfrak{t}(F)$ définie par

\begin{center}
$X\mapsto det\left((h(S^iw,S^jw))_{0\leqslant i\leqslant d_{W''}-1,0\leqslant j\leqslant d_{W''}-1}\right)$
\end{center}

(Remarquons que cette fonction à priori seulement définie sur $\mathfrak{t}(F)^S$ est en fait définie sur $\mathfrak{t}(F)$ d'après le lemme 7.1.2 et qu'elle est polynomiale toujours d'après le même lemme). Soit $\omega_T\subset \mathfrak{t}(F)$ un sous-ensemble compact. On va montrer le fait suivant dont découlent les deux points du lemme:

\begin{center}
(1) Soit $R_{W''}=\mathcal{O}_Ee_1\oplus\ldots \mathcal{O}_E e_{d_{W''}}$ un $\mathcal{O}_E$-réseau de $W''$. Il existe une constante $\alpha\in F^{\times}$ telle que \\
$\forall X\in\mathfrak{t}(F)^S\cap\omega_T$, $\forall i=0,\ldots,r$, $\alpha Q_S(X)v_{\pm i}\in R_{X_\Lambda,v_r}$ et $\alpha Q_S(X)R_{W''}\subset R_{X_\Lambda,v_r}$
\end{center}

Ramarquons, d'après l'hypothèse faite sur la fonction $X\mapsto X_\Lambda$, que la fonction $X\in\mathfrak{t}(F)^S\cap\omega_T\mapsto w$ prend ses valeurs dans un sous-ensemble compact de $W''$.On montre d'abord

\begin{center}
(2) Il existe $\beta\in F^\times$ telle que pour tout $X\in \mathfrak{t}(F)^S\cap \omega_T$  avec $X_\Lambda=c(v_0,w)+\displaystyle\sum_{i=0,\ldots,r-1} \lambda_ic(v_i,v_{i+1})+\displaystyle\sum_{i=0,\ldots,r} \mu_i c(v_i,\eta v_i)$, on ait

$$\displaystyle\left\{
    \begin{array}{ll}
        \beta v_i\in R_{X_\Lambda,v_r} \mbox{pour } i=0,\ldots,r \\
        \beta (S^{i-1}w+(-1)^i\prod_{j=0}^{i-1} \xi_j v_{-i})\in R_{X_\Lambda,v_r} \mbox{pour } i=1,\ldots,r \\
        \beta S^iw\in R_{X_\Lambda,v_r} \mbox{pour } i=0,\ldots,d_{W''}-1
    \end{array}
\right.
$$
\end{center}

On a déjà en effet $v_r,\xi_{r-1}v_{r-1}=X_\Lambda v_r,\ldots,\displaystyle\prod_{i=0}^{r-1} \xi_iv_0=X_\Lambda^rv_r\in R_{X_\Lambda,v_r}$ donc il existe une constante $\beta_1\in F^{\times}$ telle que pour tout $X\in \mathfrak{t}(F)^S\cap\omega_T$ on ait

\begin{center}
$\beta_1 v_i\in R_{X_\Lambda,v_r}$ pour $i=0,\ldots,r$
\end{center}

Puisque $X_\Lambda$ varie dans un sous-ensemble compact de $\Xi+S+\Lambda$, on peut certainement trouver $\gamma\in F^{\times}$ tel que $\gamma P_{X_\Lambda}$ soit un polynôme à coefficients entiers pour tout $X\in\mathfrak{t}(F)^S\cap\omega_T$. Le réseau $R_{X_\Lambda,v_r}$ est alors stable par $v\mapsto \gamma X_\Lambda v$.

\vspace{2mm}

On montre maintenant par récurrence que pour tout $i=1,\ldots,r$, il existe une constante $\alpha_i\in F^{\times}$ telle que pour tout $X\in\mathfrak{t}(F)^S\cap\omega_T$,  $\alpha_i(S^{i-1}w+(-1)^i\displaystyle\prod_{j=0}^{i-1}\xi_j v_{-i})\in R_{X_\Lambda,v_r}$. En effet pour $i=1$, on a

$$\gamma X_\Lambda(\beta v_0)=\gamma\beta\nu_0(w+\lambda_0v_1+2\mu_0\eta v_0-\xi_0v_{-1})\in R_{X_\Lambda,v_r}$$
et $\mu_0$, $\lambda_0$ restent dans des compacts d'où le résultat pour $i=1$. Si on a le résultat pour $i\leqslant r-1$ alors on a

$$\alpha_i\gamma X_\Lambda(S^{i-1}w+(-1)^i\displaystyle\prod_{j=0}^{i-1}\xi_j v_{-i})=$$

$$\alpha_i\gamma \left(S^iw+(-1)^i\displaystyle\prod_{j=0}^{i-1}\xi_j( -\xi_i v_{-(i+1)}+2\eta \mu_i v_i+\lambda_i v_{i+1})-h(w,S^{i-1}w)v_0 \right)\in R_{X_\Lambda,v_r}$$

Comme $\mu_i$, $\lambda_i$ et $h(w,S^{i-1}w)$ restent dans un certain compact on en déduit le résultat pour $i+1$.

\vspace{2mm}

On montre de la même façon que pour $i=r,\ldots,r+d_{W''}-1$, il existe une constante $\alpha_i$ telle que $\alpha_i S^iw\in R_{X_\Lambda,v_r}$ pour tout $X\in\mathfrak{t}(F)^S\cap\omega_T$. De plus, il existe une constante $\alpha_S\in F^{\times}$ (ne dépendant que de $S$) telle que

$$\alpha_S(\mathcal{O}_Ew+\mathcal{O}_ESw+\ldots+\mathcal{O}_ES^{d_{W''}-1}w)\subset \mathcal{O}_ES^rw+\mathcal{O}_ES^{r+1}w+\ldots+\mathcal{O}_ES^{r+d_{W''}-1}w$$
pour tout $w\in W''$ (car $S$ agit sans noyau sur $W''$). Cela prouve (2). Pour prouver (1) il suffit donc de montrer l'existence d'une constante $\alpha\in F^{\times}$ telle que pour tout $X\in\mathfrak{t}(F)^S\cap \omega_T$ on ait $\alpha Q(X)R_{W''}\subset \mathcal{O}_Ew+\mathcal{O}_ESw+\ldots+\mathcal{O}_ES^{d_{W''}-1}w$. Les coordonnées de $e_i$ dans la base $w,Sw,\ldots,S^{d_{W''}-1}w$ sont données par le vecteur colonne

\begin{center}
$M(X)^{-1}\begin{pmatrix} h(w,e_i) \\ h(Sw,e_i) \\ \vdots \\ h(S^{d_{W''}-1}w,e_i) \end{pmatrix}$
\end{center}

\noindent où $M(X)=(h(S^iw,S^jw))_{0\leqslant i\leqslant d_{W''}-1,0\leqslant j\leqslant d_{W''}-1}$. On en déduit le résultat puisque pour $X\in\mathfrak{t}(F)^S\cap\omega_T$, les coefficients de la matrice $det(M(X))M(X)^{-1}$ et les $h(S^kw,e_i)$ pour $k=0,\ldots,d_{W''}-1$ sont bornés. $\blacksquare$

\vspace{4mm}

\subsection{Calcul de $J_{\kappa''}(\theta'',\varphi)$}

Les résultats de la section précédente nous permettent maintenant d'évaluer $J_{\kappa''}(\theta'',\varphi)$ en fonction de la transformée de Fourier de $\varphi$. On a d'après l'hypothèse faite sur $\theta''$ et le lemme 7.0.2:

$$J(\theta'',\varphi,g)=D^{H''}(S)^{1/2}\displaystyle\int_{H_S''(F)\backslash H''(F)} \widehat{{}^g\!\varphi^\xi}(h^{-1}Sh)dh$$

$$=D^{H''}(S)^{1/2}\displaystyle\int_{H_S''(F)\backslash H''(F)} \int_{\Sigma} {}^g\!\widehat{\varphi}(\Xi+h^{-1}Sh+X) dXdh$$

La conjugaison par $H''(F)$ laisse stable $\Xi$ et $\Sigma$ et ne change pas la mesure sur ce dernier donc:

$$J(\theta'',\varphi,g)=D^{H''}(S)^{1/2}\displaystyle\int_{H_S''(F)\backslash H''(F)} \int_{\Sigma} {}^{hg}\!\widehat{\varphi}(\Xi+S+X) dXdh$$

On en déduit que:

$$J_{\kappa''}(\theta'',\varphi)=D^{H''}(S)^{1/2}\displaystyle\int_{H_S''(F)U''(F)\backslash G''(F)} \kappa''(g) \int_\Sigma {}^g\!\widehat{\varphi}(\Xi+S+X) dX dg$$

Le lemme 7.2.1, nous permet alors d'écrire:

\[\begin{aligned}
J_{\kappa''} & (\theta'',\varphi) \\
 & =\displaystyle\sum_{T\in\mathcal{T}(G'')} |W(G'',T)|^{-1}\int_{H_S''(F)U''(F)\backslash G''(F)} \kappa''(g) \int_{H_S''(F)U''(F)} \int_{\mathfrak{t}(F)^S} {}^{\gamma_Yyg}\!\widehat{\varphi}(Y)D^{G''}(Y)^{1/2}dYdy \\
 & =\displaystyle\sum_{T\in\mathcal{T}(G'')} \nu(A_T)^{-1} |W(G'',T)|^{-1}\times \int_{\mathfrak{t}(F)^S}\int_{A_T(F)\backslash G''(F)} \hat{\varphi}(g^{-1}Yg)\kappa_Y''(g)dg D^{G''}(Y)^{1/2}dY
\end{aligned}\]

où on a posé pour $T\in\mathcal{T}(G'')$ et $Y\in\mathfrak{t}(F)^S$ et $g\in G''(F)$

$$\kappa_Y''(g)=\nu(A_T)\displaystyle\int_{A_T(F)} \kappa''(\gamma_Y^{-1}ag) da$$

\section{Calcul de la limite de $J_{x,\omega,N}(\theta,f)$}

\subsection{Première transformation de $J_{x,\omega,N}(\theta,f)$}

\vspace{3mm}

On garde les notations du § 6. En reprenant le raisonnement de [W1] section 10.1, on montre qu'il existe un ensemble fini $\mathcal{S}$ d'éléments de $\mathfrak{h}''_{reg}(F)$ et une famille $(\hat{j}_S)_{S\in\mathcal{S}}$ de quasicaractères de $\mathfrak{h}'_x(F)$ indexée par $\mathcal{S}$ tels que pour tout $S\in\mathcal{S}$, le noyau de $S$ agissant sur $W''$ soit nul et

$$\forall X=X'+X''\in\mathfrak{h}_{x,reg}(F), \theta_{x,\omega}(X)=\displaystyle\sum_{S\in\mathcal{S}} \hat{j}_S(X')\hat{j}^{H''}(S,X'')$$

Pour $g\in G(F)$, on note $X=X'+X''\mapsto {}^gf^\sharp_{x,\omega}(X'+X'')$ la transformée de Fourier partielle de ${}^g f_{x,\omega}$ par rapport à la deuxième variable. C'est-à-dire que l'on a

$${}^g f^\sharp_{x,\omega}(X'+X'')=\displaystyle\int_{\mathfrak{h}''(F)} {}^g f_{x,\omega}(X'+Y'') \psi(<Y'',X''>) dY''$$

Les mêmes manipulations calculatoires que dans la section 10.1 de [W1] permettent d'obtenir l'expression suivante de $J_{x,\omega,N}$

$$\mbox{(1)} \;\;\;J_{x,\omega,N}(\theta,f)=\displaystyle\sum_{S\in\mathcal{S}} \displaystyle\sum_{T=T'T''\in\mathcal{T}(G_x)} \nu(A_{T''})^{-1} |W(G_x,T)|^{-1}$$ \\
$$\times \displaystyle\int_{\mathfrak{t}'(F)\times\mathfrak{t}''^S(F)} \hat{j}_S(X')D^{G'_x}(X')D^{G''}(X'')^{1/2}$$ \\
$$\times \displaystyle\int_{T'(F)A_{T''}(F)\backslash G(F)} {}^g\!f_{x,\omega}^{\sharp}(X'+X'')\kappa_{N,X''}(g) dg dX'dX''$$

où on a posé

$$\kappa_{N,X''}(g)=\nu(A_{T''})\displaystyle\int_{A_{T''}(F)} \kappa_N(\gamma_{X''}^{-1}ag) da$$

Pour $S\in\mathcal{S}$ et $T=T'T''\in\mathcal{T}(G'_x)\times\mathcal{T}(G'')$ soit $\mathcal{Q}_{S,T}$ l'ensemble des trois fonctions polynomiales suivantes sur $\mathfrak{t}(F)$:

\begin{itemize}
\item La fonction $X=X'+X''\in \mathfrak{t}=\mathfrak{t}'\oplus\mathfrak{t}''\mapsto det(ad(X')_{|\mathfrak{g}'_x/\mathfrak{t}'})$.
\item La fonction $X=X'+X''\in \mathfrak{t}=\mathfrak{t}'\oplus\mathfrak{t}''\mapsto det(ad(X'')_{|\mathfrak{g}''/\mathfrak{t}''})$.
\item La fonction polynomiale $Q_S$ sur $\mathfrak{t}''$ fournie par le lemme 7.3.2.
\end{itemize}

\vspace{2mm}

On note $\mathfrak{t}(F)[>\epsilon]$ (resp. $\mathfrak{t}(F)[\leqslant\epsilon]$) l'ouvert des éléments $X\in\mathfrak{t}(F)$ qui vérifient les inégalités $|Q(X)|_F>\epsilon$ (resp. $|Q(X)|_F\leqslant\epsilon$) pour $Q\in \mathcal{Q}_{S,T}$. Notons $J_{x,\omega,N,> N^{-b}}(\theta,f)$ (resp. $J_{x,\omega,N,\leqslant N^{-b}}(\theta,f)$) l'expression définie par la même formule que (1) où au lieu d'intégrer sur $\mathfrak{t}'(F)\times\mathfrak{t}''(F)^S$ on intègre sur $\mathfrak{t}'(F)\times\mathfrak{t}''(F)^S\cap \mathfrak{t}(F)[>N^{-b}]$ (resp. $\mathfrak{t}'(F)\times\mathfrak{t}''(F)^S\cap \mathfrak{t}(F)[\leqslant N^{-b}]$).  On justifie les manipulations ayant permis d'arriver à l'expression ci-dessus par la convergence absolue des intégrales qui y apparaissent. Notons $|J|_{x,\omega,N}(\theta,f)$ et $|J|_{x,\omega,N,\leqslant N^{-b}}(\theta,f)$ les mêmes expressions que précédemment où on a remplacé les fonctions intégrées par leurs valeurs absolues. On a alors comme dans [W1] les majorations suivantes:

\begin{itemize}
\item Il existe $k\in\mathbb{N}$ et $c>0$ tels que pour tout $N\geqslant 1$, $|J|_{x,\omega,N}(\theta,f)\leqslant cN^k$.
\item Il existe un entier $b\in\mathbb{N}$ et une constante $c>0$ tels que pour tout $N\geqslant 1$, \\
$|J|_{x,\omega,N,\leqslant N^{-b}}(\theta,f)\leqslant c N^{-1}$.
\end{itemize}

La démonstration est identique à celle présentée dans [W1], une fois la majoration suivante établie:

\begin{lem}
Soit $T=T'T''\in\mathcal{T}(G_x)$ et $\omega_{T''}$ un sous-ensemble compact de $\mathfrak{t}''(F)$. Soit $Q_S$ la fonction polynomiale sur $\mathfrak{t}''(F)$ obtenue à partir du lemme 7.3.2 pour le réseau $R''=R\cap V''$. \\
Alors il existe $k\in\mathbb{N}$ et $c>0$ tels que

$$\kappa_{N,X''}(g)\leqslant cN^k\sigma(g)^k(1+|log |Q_S(X'')|_F|)^k(1+|log D^{G''}(X'')|)^k$$
pour tout $g\in G(F)$, pour tout $X''\in\mathfrak{t}''(F)^S\cap \omega_{T''}$ et pour tout $N\geqslant 1$.
\end{lem}

\ul{Preuve}: Notons $\alpha\in F^\times$, la constante fournie par le lemme 7.3.2 pour le sous-ensemble compact $\omega_{T''}\subset \mathfrak{t}''(F)$. On adopte la notation suivante: si $x$ est un réel quelconque on pose $\pi_E^x=\pi_E^{E(x)}$. Il existe une constante $C>0$ telle que pour tout $g\in G(F)$, on ait:
$$ \pi_{E}^{C\sigma(g)}R\subset g(R)\subset \pi_{E}^{-C\sigma(g)}R$$
D'après le paragraphe 7.3, l'application $X''\in\mathfrak{t}''(F)^S\cap\omega_{T''}\mapsto X''_\Lambda\in \Xi+S+\Lambda$ a son image incluse dans un sous-ensemble compact de $\Xi+S+\Lambda$. Par conséquent, on peut trouver une constante $C'>0$ telle que pour tout $X''\in\mathfrak{t}''(F)^S\cap\omega_{T''}$ on ait

$${X_{\Lambda}''}^{k}R\subset \pi_E^{-C'}R$$
pour $k=0,\ldots,d_{W''}-1$. \\

Par définition $\kappa_{N,X''}(g)$ est le volume de l'ensemble des éléments $a\in A_{T''}(F)$ vérifiant $\kappa_N(\gamma_{X''}^{-1}ag)=1$. On a les implications:

\[\begin{aligned}
\kappa_N(\gamma_{X''}^{-1}ag)=1 & \Rightarrow g^{-1}a^{-1}\gamma_{X''}(v_r)\in\pi_{E}^{-N}R \\
 & \Rightarrow \gamma_{X''}^{-1}a^{-1}\gamma_{X''}(v_r)\in \pi_E^{-N} \gamma_{X''}^{-1}g(R) \\
 & \Rightarrow \gamma_{X''}^{-1}a^{-1}\gamma_{X''}(v_r)\in \pi_E^{-N-C\sigma(g)-C\sigma(\gamma_{X''})}R \\
 & \Rightarrow \gamma_{X''}^{-1}a^{-1}\gamma_{X''}(R_{X''_\Lambda,v_r})\subset \pi_E^{-N-C\sigma(g)-C\sigma(\gamma_{X''})-C'} R \; \text{(Car $\gamma_{X''}^{-1}a^{-1}\gamma_{X''}$ commute avec $X''_\Lambda$)} \\
 & \Rightarrow \alpha Q_S(X'')\pi_E^{C\sigma(\gamma_{X''})}a^{-1}(R'')\subset \pi_E^{-N-C\sigma(g)-C\sigma(\gamma_{X''})-C'} \gamma_{X''}(R) \\
 & \Rightarrow a^{-1}(R'')\subset \alpha^{-1} Q_S(X'')^{-1}\pi_E^{-N-C\sigma(g)-3C\sigma(\gamma_{X''})-C'}R
\end{aligned}\]

Puisque pour $a\in A_{T''}(F)$ on a $a^{-1}(R'')\subset W''$, on déduit des implications précédentes l'existence d'une constante $C_0>0$ telle que pour tout $X''\in \mathfrak{t}(F)^S\cap\omega_T''$ on ait

$$\kappa_N(\gamma_{X''}^{-1}ag)=1 \Rightarrow a^{-1}(R'')\subset Q_S(X'')^{-1}\pi_E^{-N-C_0\sigma(g)-C_0\sigma(\gamma_{X''})-C_0}R''$$

On a donc $\kappa_{N,X''}(g)\leqslant vol(a\in A_{T''}(F): a^{-1}(R'')\subset Q_S(X'')^{-1}\pi_E^{-N-C_0\sigma(g)-C_0\sigma(\gamma_{X''})-C_0}R'')$. Il existe une contante une constante $C_1>0$ et un entier $k\in\mathbb{N}$ tels que pour tout $\beta\in E^{\times}$ on ait

$$vol(\{a\in A_{T''}(F): a^{-1}(R'')\subset \beta R''\})\leqslant C_1 max(1,-val_E(\beta))^k$$

En utilisant l'inégalité $\sigma(\gamma_{X''})<< (1+|log(D^{G''}(X''))|)$ pour tout $X''\in\mathfrak{t}(F)^S\cap\omega_{T''}$, on en déduit l'inégalité voulue. $\blacksquare$

\vspace{4mm}

On fixe dorénavant un tore $T=T'T''\in\mathcal{T}(G_x)$ et on note $M_{\natural}$ le centralisateur de $A_{T''}$ dans $G$. C'est un Lévi de $G$.

\begin{lem}
On a $A_{M_{\natural}}=A_{T''}$.
\end{lem}

\ul{Preuve}: L'inclusion $A_{T''}\subset A_{M_{\natural}}$ est évidente. Puisque $T''H'\subset M_{\natural}$, $A_{M_{\natural}}$ est inclus dans le centre déployé de $T''H'$ qui est $A_{T''}$ car $A_{H'}=\{1\}$. $\blacksquare$

\subsection{Changement de la fonction de troncature}

Pour $g\in G(F)$ on pose $\sigma_T(g)=inf\{\sigma(tg); t\in T(F)\}$.

\begin{lem}
Il existe un réel $C_1>0$ et un sous-ensemble compact $\omega_T\subset \mathfrak{t}(F)$ tels que pour tout $g\in G(F)$, pour tout entier $N\geqslant 2$ et pour tout $X\in \mathfrak{t}'(F)\times \mathfrak{t}''(F)^S\cap\mathfrak{t}(F)[>N^{-b}]$ la non nullité de ${}^g\! f^\sharp_{x,\omega}(X)\neq 0$ entraîne $X\in \omega_T$ et $\sigma_T(g)\leqslant C_1 log(N)$.
\end{lem}

\ul{Preuve}: On sait qu'il existe un sous-ensemble compact-ouvert $\Gamma\subset G(F)$ tel que ${}^g\!f^\sharp_{x,\omega}(X)=0$ si $g\notin G_x(F)\Gamma$. De plus, il existe un sous-groupe ouvert-compact $K_f\subset G(F)$ tel que pour tout $g\in G(F)$ et tout $k\in K_f$, on ait ${}^{gk}\!f_{x,\omega}={}^g\!f_{x,\omega}$ et donc ${}^{gk}\!f^\sharp_{x,\omega}={}^g\!f^\sharp_{x,\omega}$. Quitte à agrandir $\Gamma$, on peut supposer que $K_f\Gamma=\Gamma$. Soit $\Gamma_f$ un ensemble de représentants des classes à droite de $\Gamma$ modulo $K_f$, c'est un ensemble fini. Posons $\Omega=\displaystyle \bigcup_{\gamma\in \Gamma_f} Supp({}^\gamma\!f^\sharp_{x,\omega})$ et $\omega_T=Cl(\Omega^{G_x}\cap \mathfrak{t}(F))$ où $Cl$ désigne l'adhérence. Ce sont des sous-ensembles compacts de $\mathfrak{g}_x(F)$ et $\mathfrak{t}(F)$ respectivement. Il existe deux réels $C,C'>0$  telles que

\begin{itemize}
\item $\forall g\in G_x(F)$, $\forall X\in\omega_T\cap\mathfrak{g}_{x,reg}(F)$, $g^{-1}Xg\in \Omega\Rightarrow \sigma_T(g)\leqslant C(1+|log D^{G_x}(X)|)$ (c'est le lemme 7.3.1)
\item $\forall X\in \mathfrak{t}'(F)\times\mathfrak{t}''(F)^S\cap \omega_T[>N^{-b}]$, $|log D^{G_x}(X)|\leqslant C'log(N)$
\end{itemize}

Pour $X$, $g$ et $N$ comme dans l'énoncé on a $X\in\omega_T$ d'après ce qui précède et on peut écrire $g=g_x\gamma$ avec $g_x\in G_x(F)$ et $\gamma\in \Gamma_f$. On a alors $g_x^{-1}Xg_x\in \Omega$ et on conclut avec les majorations ci-dessus $\blacksquare$

\vspace{2mm}

Soit $\mathcal{Y}=(Y_{P_{\natural}})_{P_{\natural}\in\mathcal{P}(M_{\natural})}$ une famille d'éléments de $\mathcal{A}_{M_{\natural}}$, $(G,M_{\natural})$-orthogonale et positive. On reprend les notations de [W1] qui reprennent celles d'Arthur: pour $Q=LU_Q\in \mathcal{F}(M_{\natural})$, on note $\sigma_{M_{\natural}}^Q(.,\mathcal{Y})$ resp. $\tau_Q$ la fonction caractéristique dans $\mathcal{A}_{M_{\natural}}$ de la somme de $\mathcal{A}_L$ et de l'enveloppe convexe des $(Y_{P_{\natural}})_{P_{\natural}\subset Q}$ resp. la fonction caractéristique de $\mathcal{A}^L_{M_{\natural}}+\mathcal{A}_Q^+$. On a alors

\begin{center}
$\displaystyle\sum_{Q\in\mathcal{F}(M_{\natural})} \sigma_{M_{\natural}}^Q(\zeta,\mathcal{Y})\tau_Q(\zeta-Y_Q)=1$ pour tout $\zeta\in \mathcal{A}_{M_{\natural}}$
\end{center}
où pour $Q=LU_Q\in \mathcal{F}(M_\natural)$, $Y_Q$ désigne la projection de $Y_{P_\natural}$ sur $\mathcal{A}_L$ pour n'importe quel parabolique $P_\natural\in\mathcal{P}(M_\natural)$ vérifiant $P_\natural\subset Q$ (cela ne dépend pas du choix de $P_\natural$). \\

On fixe $M_{min}$ un Lévi minimal de $G$ inclus dans $M_\natural$, un sous-groupe compact spécial $K$ en bonne position par rapport à $M_{min}$, et un sous-groupe parabolique minimal $P_{min}\in\mathcal{P}(M_{min})$. Soit $Y_{P_{min}}\in \mathcal{A}_{P_{min}}^{+}$. Pour tout $P\in\mathcal{P}(M_{min})$ il existe un unique $w\in W^G$ tel que $wP_{min}=P$ on pose alors $Y_P=wY_{P_{min}}\in\mathcal{A}_P^{+}$. La famille $(Y_P)_{P\in\mathcal{P}(M_{min})}$ est alors $(G,M_{min})$-orthogonale positive. Elle définit donc une $(G,M_\natural)$-famille orthogonale positive $(Y_{P_\natural})_{P_\natural\in\mathcal{P}(M_\natural)}$. Pour tout $g\in G(F)$ on notera $\mathcal{Y}(g)$ la $(G,M_\natural)$-famille orthogonale définie par $\mathcal{Y}(g)_Q=Y_Q-H_{\overline{Q}}(g)$. Posons

$$\tilde{v}(Y_{P_{min}},g)=\nu(A_{T''})\displaystyle\int_{A_{T''}(F)} \sigma^G_{M_\natural}(H_{M_\natural}(a),\mathcal{Y}(g)) da$$

Il existe une constante $C_2>0$ telle que pour $g\in G(F)$ vérifiant $\sigma(g)\leqslant C_2 inf_{\alpha\in \Delta_{min}} \alpha(Y_{P_{min}})$ la famille $(\mathcal{Y}(g)_{P_\natural})_{P_\natural\in\mathcal{P}(M_\natural)}$ vérifie $\mathcal{Y}(g)_{P_\natural}\in\mathcal{A}_{P_\natural}^+$ pour tout $P_\natural\in\mathcal{P}(M_\natural)$ donc est orthogonale positive.

\begin{lem}
Il existe deux constantes $C_3,C_4>0$ et un entier naturel $N_0$ telle que pour tout entier $N\geqslant N_0$ et pour tout $Y_{P_{min}}\in \mathcal{A}_{P_{min}}^{+}$ vérifiant

\begin{center}
\hspace{48mm} $C_3 log(N) \leqslant inf_{\alpha\in \Delta_{min}} \alpha(Y_{P_{min}})$ \hspace{41mm} (A)
\end{center}

\begin{center}
\hspace{48mm} $sup_{\alpha\in \Delta_{min}} \alpha(Y_{P_{min}})\leqslant C_4 N$ \hspace{48mm} (B)
\end{center}

on ait l'égalité

$$\displaystyle\int_{T'(F)A_{T''}(F)\backslash G(F)} {}^g\!f^{\sharp}_{x,\omega}(X'+X'')\kappa_{N,X''}(g)dg=\int_{T'(F)A_{T''}(F)\backslash G(F)} {}^g\!f^{\sharp}_{x,\omega}(X'+X'') \tilde{v}(Y_{P_{min}},g) dg$$
pour tout $X=X'+X''\in\mathfrak{t}'(F)\times \mathfrak{t}''(F)^S\cap \mathfrak{t}(F)[>N^{-b}]$.
\end{lem}

\ul{Preuve}: D'après le lemme 8.2.1 il n'y a rien à démontrer si $X\notin \omega_T$ on supposera donc dans la suite que $X\in \omega_T$. On fixe des constantes $C_3,C_4>0$, un entier naturel $N_0$, un élément $Y_{P_{min}}\in \mathcal{A}_{P_{min}}^+$ vérifiant les inégalités (A) et (B) et un entier $N\geqslant N_0$ et on va montrer que l'égalité de l'énoncé est vérifiée si $C_3$ et $N_0$ sont assez grands et $C_4$ assez petit. On peut déjà supposer $C_3\geqslant C_1$. Alors d'après le lemme 8.2.1, pour tout entier $N\geqslant 2$, pour tout $g\in G(F)$ et pour tout $X\in \mathfrak{t}'(F)\times \mathfrak{t}''(F)^S\cap \mathfrak{t}(F)[>N^{-b}]$, si le terme sous l'intégrale est non nul, la $(G,M_\natural)$ famille $(\mathcal{Y}(g)_{P_\natural})_{P_\natural\in\mathcal{P}(M_\natural)}$ est orthogonale positive. Si tel est le cas, pour tout $a\in A_{T''}(F)$ on a $\displaystyle\sum_{Q\in\mathcal{F}(M_\natural)} \sigma_{M_\natural}^Q(H_{M_\natural}(a),\mathcal{Y}(g))\tau_Q(H_{M_\natural}(a)-Y(g)_Q)=1$ et on peut donc écrire

$$\kappa_{N,X''}(g)=\displaystyle\sum_{Q\in\mathcal{F}(M_\natural)}\kappa_{N,X''}(Q,g)$$

\noindent où $\kappa_{N,X''}(Q,g)=\nu(A_{T''})\displaystyle\int_{A_{T''}(F)} \sigma_{M_\natural}^Q(H_{M_\natural}(a),\mathcal{Y}(g))\tau_Q(H_{M_\natural}(a)-Y(g)_Q) \kappa_N(\gamma_{X''}^{-1}ag) da$. Les fonctions $\kappa_{N,X''}(Q,.)$ sont invariantes à gauche par $T'(F)A_{T''}(F)$. On a donc une décomposition
$$\displaystyle\int_{T'(F)A_{T''}(F)\backslash G(F)} {}^g\!f^{\sharp}_{x,\omega}(X'+X'')\kappa_{N,X''}(g)dg=\sum_{Q\in\mathcal{F}(M_\natural)} I(Q,X)$$
où $I(Q,X)=\displaystyle\int_{T'(F)A_{T''}(F)\backslash G(F)} {}^g\!f^{\sharp}_{x,\omega}(X)\kappa_{N,X''}(Q,g) dg$. On va montrer les deux faits suivants dont découlent la proposition:

\begin{enumerate}
\item Si $C_4$ est assez petit et $N_0$ est assez grand, on a

$$I(G,X)=\displaystyle\int_{T'(F)A_{T''}(F)\backslash G(F)} {}^g\!f^{\sharp}_{x,\omega}(X) \tilde{v}(Y_{P_{min}},g) dg$$
pour tout $X\in \mathfrak{t}'(F)\times \mathfrak{t}''(F)^S\cap \omega_T[>N^{-b}]$.

\item Si $C_3$ et $N_0$ sont assez grand, on a $I(Q,X)=0$ pour tout $Q\in\mathcal{F}(M_{\natural})$ différent de $G$ et pour tout $X\in \mathfrak{t}'(F)\times \mathfrak{t}''(F)^S\cap \omega_T[>N^{-b}]$.
\end{enumerate}

Preuve de (1): D'après le lemme 8.2.1, il suffit d'avoir l'implication

$$\sigma_{M_\natural}^G(H_{M_\natural}(a),\mathcal{Y}(g))=1\Rightarrow \kappa_{N}(\gamma_{X''}^{-1}ag)=1$$
pour tout $a\in A_{T''}(F)$, pour tout $g\in G(F)$ vérifiant $\sigma(g)\leqslant C_1 log(N)$ et pour tout $X\in\mathfrak{t}'(F)\times \mathfrak{t}''(F)^S\cap \omega_T[>N^{-b}]$. Pour $a\in A_{T''}(F)$, $\sigma_{M_\natural}^G(H_{M_\natural}(a),\mathcal{Y}(g))=1$ entraîne $\sigma(a)<<(sup_{\alpha\in\Delta_{min}} \alpha(Y_{P_{min}}) +\sigma(g))$. Il existe une constante $C>0$ telle que $\sigma(g)<CN$ entraîne $\kappa_N(g)=1$ pour tout $g\in G(F)$. Comme de plus $\sigma(\gamma_{X''})<<log(N)$ on en déduit le point (1).

\vspace{2mm}

Preuve de (2): Soit $Q=LU_Q\in \mathcal{F}(M_\natural)$ différent de $G$ et $\overline{Q}=LU_{\overline{Q}}$ le parabolique opposé. On a

\[\begin{aligned}
I(Q,X) & =\displaystyle\int_{T'(F)A_{T''}(F)\backslash G(F)} {}^g\!f^{\sharp}_{x,\omega}(X)\kappa_{N,X''}(Q,g) dg \\
 & =\int_{T'(F)A_{T''}(F)\backslash L_Q(F)} \int_K \int_{U_{\overline{Q}}(F)}{}^{\overline{u}lk}\!f^\sharp_{x,\omega}(X) \kappa_{N,X''}(Q,\overline{u}lk) d\overline{u}dkdl
\end{aligned}\]

Si on a $\kappa_{N,X''}(Q,\overline{u}lk)=\kappa_{N,X''}(Q,lk)$ pour $\overline{u}\in U_{\overline{Q}}(F)$, $l\in L_{Q}(F)$ et $k\in K$ vérifiant ${}^{\overline{u}lk}\!f^\sharp_{x,\omega}(X)\neq 0$ alors l'intégrale intérieure sera nulle car $f$ est très cuspidale (c'est le lemme 5.5 de [W1]) et on aura bien $I(Q,X)=0$. D'après le lemme 8.2.1 il existe $C_5>0$ telle que ${}^{\overline{u}lk}\!f^\sharp_{x,\omega}(X)\neq 0$ entraîne $\sigma(\overline{u}lk),\sigma(lk),\sigma(\overline{u})\leqslant C_5 log(N)$.  Puisque l'on a 

$$\kappa_{N,X''}(Q,g)=\nu(A_{T''})\displaystyle\int_{A_{T''}(F)} \sigma_{M_\natural}^Q(H_{M_\natural}(a),\mathcal{Y}(g))\tau_Q(H_{M_\natural}(a)-Y(g)_Q) \kappa_N(\gamma_{X''}^{-1}ag) da$$
et que les fonctions $\sigma_{M_\natural}^Q(.,\mathcal{Y}(g))$ et $\tau_Q(.-Y(g)_Q)$ sont invariantes par translation à gauche de $g$ par $U_{\overline{Q}}(F)$ (car elle ne dépendent que de $H_{\overline{P}_\natural}(g)$ pour $P_\natural\subset Q$), il suffit donc d'avoir

\begin{center}
$\forall X\in\mathfrak{t}'(F)\times \mathfrak{t}''(F)^S\cap \omega_T[>N^{-b}]$, $\forall a\in A_{T''}(F)$, $\forall g\in G(F)$, $\forall \overline{u}\in U_{\overline{Q}}(F)$ tels que $\sigma_{M_\natural}^Q(H_{M_\natural}(a),\mathcal{Y}(g))\tau_Q(H_{M_\natural}(a)-Y(g)_Q)=1$ et $\sigma(\overline{u}g),\sigma(g),\sigma(\overline{u})\leqslant C_5 log(N)$ on a $\kappa_{N,X''}(\gamma_{X''}^{-1}a\overline{u}g)=\kappa_{N,X''}(\gamma_{X''}^{-1}ag)$
\end{center}

Soit $g\in G(F)$ et $a\in A_{T''}(F)$ tel que $\sigma(g)\leqslant C_5 log(N)$ et $\sigma_{M_\natural}^Q(H_{M_\natural}(a),\mathcal{Y}(g))\tau_Q(H_{M_\natural}(a)-Y(g)_Q)=1$. La fonction $\sigma_{M_\natural}^Q(.,\mathcal{Y}(g))\tau_Q(.-Y(g)_Q)$ est la fonction caractéristique de la somme de $\mathcal{A}_Q^+$ et de l'enveloppe convexe des $Y(g)_{P_\natural}$ pour $P_\natural \subset Q$. On a donc pour tout $\alpha\in\Sigma_Q^+$ l'inégalité $\alpha(H_{M_\natural}(a))\geqslant \displaystyle inf_{P_\natural\subset Q} \alpha(Y_{P_\natural}-H_{\overline{P}_\natural}(g))$. Il existe une constante $C'>0$ telle que $inf_{P_\natural\subset Q} \alpha(Y_{P_\natural}-H_{\overline{P}_\natural}(g))\geqslant C'inf_{\beta\in \Delta_{min}} \beta(Y_{P_{min}}) -C'\sigma(g)$. Par conséquent $\alpha(H_{M_\natural}(a))\geqslant C'(C_3-C_5)log(N)$. Le lemme suivant nous permet de conclure à l'existence d'une constante $C_3$ et d'un entier $N_2$ qui conviennent.

\begin{lem}
Soit $C_5>0$ une constante alors il existe $C_6>0$ et un entier $N_1$ tels que la propriété suivante soit vérifiée: pour tous $a\in A_{T''}(F)$, $g\in G(F)$, $\overline{u}\in U_{\overline{Q}}(F)$, $X\in \mathfrak{t}'(F)\times\mathfrak{t}''(F)^S\cap \omega_T[S;>N^{-b}]$ et pour tout entier $N\geqslant N_1$ vérifiant $\sigma(g),\sigma(\overline{u}g),\sigma(\overline{u})\leqslant C_5 log(N)$ et $\alpha(H_{M_\natural}(a))\geqslant C_6 log(N)$ pour tout $\alpha\in \Sigma_Q^+$ on a $\kappa_N(\gamma_{X''}^{-1}a\overline{u}g)=\kappa_N(\gamma_{X''}^{-1}ag)$.
\end{lem}

\ul{Preuve du Lemme 8.2.3}: Montrons d'abord qu'il existe un entier $N_1$ tel que pour tout entier $N\geqslant N_1$ et pour tous $a\in A_{T''}(F)$, $g\in G(F)$ et $X\in \mathfrak{t}'(F)\times \mathfrak{t}''(F)^S\cap\omega_T[S;>N^{-b}]$ avec $\sigma(g)\leqslant C_5 log(N)$ alors $\kappa_N(\gamma_{X''}^{-1}ag)=1$ si et seulement si $g^{-1}a^{-1}\gamma_{X''}(v_r)\in \pi_E^{-N} R$. \\
La condition est de toute façon nécessaire. Fixons $N,a,g$ et $X$ vérifiant les conditions précédentes et supposons que $g^{-1}a^{-1}\gamma_{X''}(v_r)\in \pi_E^{-N} R$. On va montrer que si $N_1$ est assez grand on a $\kappa_N(\gamma_{X''}^{-1}ag)=1$. On considère comme dans la preuve du lemme 8.1.1 une constante $C>0$ telle que pour tout $g'\in G(F)$ on ait $\pi_E^{C\sigma(g')}R\subset g'(R)\subset \pi_E^{-C\sigma(g')}R$ et une contante $C'>0$ telle que ${X_{\Lambda}''}^{k}(R)\subset \pi_E^{-C'}R$ pour tout $k=0,\ldots,d_{W''}-1$ et pour tout $X\in\mathfrak{t}'(F)\times \mathfrak{t}''(F)^S\cap \omega_T$. Le même raisonnement que dans la démonstration du lemme 8.1.1 montre que $g^{-1}a^{-1}\gamma_{X''}(R)\subset Q_S(X'')^{-1}\pi_E^{-N-2C\sigma(g)-2C\sigma(\gamma_{X''})-C'}R$. De plus, on a $|log |Q_S(X'')||<< log(N)$  et $\sigma(\gamma_{X''})<< log(N)$. Donc pour $N_1$ assez grand $g^{-1}a^{-1}\gamma_{X''}(v_r)\in \pi_E^{-N} R$ entraîne $g^{-1}a^{-1}\gamma_{X''}(R)\subset \pi_E^{-N-\sqrt{N}} R$ ce qui permet de conclure.

\vspace{2mm}

 Soit $N\geqslant N_1$ un entier et soit $a\in A_{T''}(F)$. On peut supposer qu'il existe $g\in G(F)$ et $X\in \mathfrak{t}'(F)\times \mathfrak{t}''(F)^S\cap\omega_T[S;>N^{-b}]$ tels que $\sigma(g)\leqslant C_5log(N)$ et $\kappa_N(\gamma_{X''}^{-1}ag)=1$ sinon il n'y a rien à dire. Notons $\mathcal{A}_N$ l'ensemble des $a\in A_{T''}(F)$ qui vérifient cette condition d'existence. En adaptant un tout petit peu le raisonnement précédent on montre qu'il existe une constante $C''>0$ tel que $a^{-1}(R)\subset \pi_E^{-N-C''log(N)}R$ pour tout $a\in \mathcal{A}_N$. On a alors pour tout $v\in V-\{0\}$, $val_R(a^{-1}v)-val_R(v)\geqslant -N-C''log(N)$. On peut fixer la constante $C''$ telle qu'on ait aussi $val_R(g^{-1}v)-val_R(v)\geqslant -C'' log(N)$ et $val_R(\gamma_{X''}v)-val_R(v)\geqslant -C'' log(N)$ pour tout $g\in G(F)$ vérifiant $\sigma(g)\leqslant C_5log(N)$, pour tout $X\in \mathfrak{t}'(F)\times\mathfrak{t}''(F)^S\cap \omega_T[S;>N^{-b}]$ et pour tout $v\in V-\{0\}$. Enfin, on peut trouver une constante $C_6$ de sorte que si $\alpha(H_{M_\natural}(a))\geqslant C_6 log(N)$ pour tout $\alpha\in \Sigma_Q^+$ alors pour tout $\overline{u}\in U_{\overline{Q}}(F)$ tel que $\sigma(\overline{u})\leqslant C_5 log(N)$ on a $val_R(a\overline{u}a^{-1}v-v)\geqslant 3C'' log(N)+val_R(v)$ pour tout $v\in V-\{0\}$. On a alors pour $a$, $g$, $\overline{u}$, $X$ comme dans l'énoncé
\[\begin{aligned}
val_R(g^{-1}\overline{u}^{-1}a^{-1}\gamma_{X''} v_r-g^{-1}a^{-1}\gamma_{X''}v_r) & \geqslant val_R(\overline{u}^{-1}a^{-1}\gamma_{X''} v_r-a^{-1}\gamma_{X''}v_r)-C''log(N) \\
 & =val_R(a^{-1}a\overline{u}^{-1}a^{-1}\gamma_{X''} v_r-a^{-1}\gamma_{X''}v_r)-C''log(N) \\
 & \geqslant val_R(a\overline{u}^{-1}a^{-1}\gamma_{X''} v_r-\gamma_{X''}v_r)-N-2C''log(N) \\
 & \geqslant val_R(\gamma_{X''}v_r)-N+C''log(N) \\
 & \geqslant -N
\end{aligned}\]
On obtient le résultat voulu d'après le premier point. $\blacksquare$

\vspace{4mm}

On déduit du lemme 8.2.2 la proposition suivante.

\begin{prop}
Il existe un entier naturel $N_1$ telle que pour tout $N\geqslant N_1$ et pour tout $Y_{P_{min}}\in \mathcal{A}_{P_{min}}^{+}$ vérifiant

$$C_3 log(N) \leqslant inf_{\alpha\in \Delta_{min}} \alpha(Y_{P_{min}})$$

on ait l'égalité

$$\displaystyle\int_{T'(F)A_{T''}(F)\backslash G(F)} {}^g\!f^{\sharp}_{x,\omega}(X'+X'')\kappa_{N,X''}(g)dg=\int_{T'(F)A_{T''}(F)\backslash G(F)} {}^g\!f^{\sharp}_{x,\omega}(X'+X'') \tilde{v}(Y_{P_{min}},g) dg$$
pour tout $X=X'+X''\in\mathfrak{t}'(F)\times \mathfrak{t}''(F)^S\cap \mathfrak{t}(F)[>N^{-b}]$.
\end{prop}

\ul{Preuve}:
En effet, il existe un entier $N_1$ tel que pour $N\geqslant N_1$ les deux ensembles

$$\{Y_{P_{min}}\in \mathcal{A}_{P_{min}}^{+}:C_3 log(N) \leqslant inf_{\alpha\in \Delta_{min}} \alpha(Y_{P_{min}}) \; \text{et} \; sup_{\alpha\in \Delta_{min}} \alpha(Y_{P_{min}})\leqslant C_4 N\}$$

et

$$\{Y_{P_{min}}\in \mathcal{A}_{P_{min}}^{+}:C_3 log(N+1) \leqslant inf_{\alpha\in \Delta_{min}} \alpha(Y_{P_{min}}) \; \text{et} \; sup_{\alpha\in \Delta_{min}} \alpha(Y_{P_{min}})\leqslant C_4 (N+1)\}$$

soient d'intersection non vide. Alors pour $X=X'+X''\in\mathfrak{t}'(F)\times \mathfrak{t}''(F)^S\cap \mathfrak{t}(F)[>N^{-b}]\subset \mathfrak{t}'(F)\times \mathfrak{t}''(F)^S\cap \mathfrak{t}(F)[>(N+1)^{-b}]$, on a d'après le lemme 8.2.2

$$\displaystyle\int_{T'(F)A_{T''}(F)\backslash G(F)} {}^g\!f^{\sharp}_{x,\omega}(X'+X'')\kappa_{N,X''}(g)dg=\displaystyle\int_{T'(F)A_{T''}(F)\backslash G(F)} {}^g\!f^{\sharp}_{x,\omega}(X'+X'')\kappa_{N+1,X''}(g)dg$$

D'où par une récurrence immédiate

$$\displaystyle\int_{T'(F)A_{T''}(F)\backslash G(F)} {}^g\!f^{\sharp}_{x,\omega}(X'+X'')\kappa_{N,X''}(g)dg=\displaystyle\int_{T'(F)A_{T''}(F)\backslash G(F)} {}^g\!f^{\sharp}_{x,\omega}(X'+X'')\kappa_{N',X''}(g)dg$$
pour tout $N'\geqslant N$. Soit $Y_{P_{min}}\in \mathcal{A}_{P_{min}}^{+}$ tel que $C_3 log(N) \leqslant inf_{\alpha\in \Delta_{min}} \alpha(Y_{P_{min}})$, on alors peut trouver un entier $N'\geqslant N$ tel que $Y_{P_{min}}$ vérifie les inégalités (A) et (B) du lemme 8.2.2 pour $N'$. Cela permet de conclure $\blacksquare$

\vspace{4mm}

\subsection{Le résultat final}

Notons $\theta_{f,x,\omega}$ la fonction sur $\mathfrak{g}_x(F)$ définie presque partout par

$$
\theta_{f,x,\omega}(X)= \left\{
    \begin{array}{ll}
        \theta_f(xexp(X)) & \mbox{si } X \in \omega \\
        0 & \mbox{sinon.}
    \end{array}
\right.
$$

On définit une fonction $\theta_{f,x,\omega}^\sharp$ sur $\mathfrak{g}_{x,reg}(F)$ par

$$\displaystyle \theta_{f,x,\omega}^\sharp(X)=(-1)^{a_{M(X)}}\nu(G_{x,X})^{-1} \int_{G_{x,X}(F)\backslash G(F)} {}^g f_{x,\omega}^\sharp(X) v_{M(X)}(g) dg$$

\noindent où pour $X\in\mathfrak{g}_{x,reg}(F)$, on a noté $M(X)$ le commutant de $A_{G_{x,X}}$ dans $G$ (c'est un Levi de $G$). Alors les fonctions $\theta_{f,x,\omega}$ et $\theta_{f,x,\omega}^\sharp$ sont localement constantes sur $\mathfrak{g}_{x,reg}(F)$, localement intégrables sur $\mathfrak{g}_x(F)$ et $\theta_{f,x,\omega}^\sharp$ est la transformée de Fourier partielle de $\theta_{f,x,\omega}$ par rapport à la deuxième variable. Plus précisément, cela signifie que l'on a l'égalité

$$\displaystyle \int_{\mathfrak{g}_x(F)} \theta_{f,x,\omega}^\sharp(X)\varphi(X)dX=\int_{\mathfrak{g}_x(F)} \theta_{f,x,\omega}(X)\varphi^\sharp(X) dX$$

\noindent pour tout $\varphi\in C_c^\infty(\mathfrak{g}_x(F))$ et où

$$\displaystyle \varphi^\sharp(X)=\int_{\mathfrak{g}''(F)} \varphi(X'+Y'') \psi(Y'',X''>) dY''$$

\noindent pour tout $X=X'+X''\in\mathfrak{g}_x(F)=\mathfrak{g}'_x(F)\oplus\mathfrak{g}''(F)$ (cf proposition 5.8 de [W1]).

\begin{lem}
Pour $N\geqslant N_1$ et $X=X'+X''\in \mathfrak{t}'(F)\times \mathfrak{t}''(F)^S\cap \mathfrak{t}(F)[>N^{-b}]$, on a les égalités

$$\displaystyle\int_{T'(F)A_{T''}(F)\backslash G(F)} {}^g\!f^{\sharp}_{x,\omega}(X)\kappa_{N,X''}(g)dg=0$$
si $A_{T'}\neq \{1\}$ et

$$\displaystyle\int_{T'(F)A_{T''}(F)\backslash G(F)} {}^g\!f^{\sharp}_{x,\omega}(X)\kappa_{N,X''}(g)dg=\nu(T')\nu(A_{T''})\theta^\sharp_{f,x,\omega}(X)$$
si $A_{T'}=\{1\}$
\end{lem}

\ul{Preuve}: C'est exactement la même que celle de la proposition 10.9 de [W1] $\blacksquare$

\vspace{4mm}

Si $A_{G'_x}=\{1\}$ on pose

\[\begin{aligned}
\mbox{(1)} \;\;K_{x,\omega}(\theta,f)=\displaystyle\sum_{S\in\mathcal{S}}\sum_{T=T'T''\in \mathcal{T}_{ell}(G'_x)\times \mathcal{T}(G'')} \nu(T') |W(G_x,T)|^{-1} \\ \int_{\mathfrak{t}'(F)\times\mathfrak{t}''(F)^S} \hat{j}_S(X') D^{G'_x}(X') D^{G''}(X'')^{1/2} \theta^\sharp_{f,x,\omega}(X) dX
\end{aligned}\]
et si $A_{G'_x}\neq \{1\}$, on pose

$$K_{x,\omega}(\theta,f)=0$$
\begin{prop}
L'expression (1) est absolument convergente et on a l'égalité
$$\lim\limits_{N\to\infty} J_{N,x,\omega}(\theta,f)=K_{x,\omega}(\theta,f)$$
\end{prop}

\ul{Preuve}: C'est exactement la même que celle de la proposition 10.10 de [W1] $\blacksquare$

\section{Cas des supports nilpotents}

Dans cette section on se propose de démontrer la proposition suivante:

\begin{prop}
Supposons que pour tout quasicaractère $\theta$ de $\mathfrak{h}(F)$ et toute fonction très cuspidale $f\in C_c^\infty(\mathfrak{g}(F))$ qui ne contient pas d'élément nilpotent dans son support on ait \\
\noindent $\lim\limits_{N\to \infty} J_N(\theta,f)=J_{geom}(\theta,f)$. Alors on a $\lim\limits_{N\to \infty} J_N(\theta,f)=J_{geom}(\theta,f)$ pour tout quasicaractère $\theta$ de $\mathfrak{h}(F)$ et toute fonction très cuspidale $f\in C_c^\infty(\mathfrak{g}(F))$.
\end{prop}

On fait donc dans la suite de la section l'hypothèse suivante

\begin{center}
\textbf{(Hyp)}: pour tout quasicaractère $\theta$ de $\mathfrak{h}(F)$ et toute fonction très cuspidale $f\in C_c^\infty(\mathfrak{g}(F))$ qui ne contient pas d'élément nilpotent dans son support on a 
$$\lim\limits_{N\to \infty} J_N(\theta,f)=J_{geom}(\theta,f)$$
\end{center}

\subsection{Calcul de $\lim\limits_{N\to\infty} J_N(\theta,f)$}

\vspace{2mm}

Soient $\theta$ un quasi-caractère de $\mathfrak{h}(F)$ et $f\in C_c^\infty(\mathfrak{g}(F))$ une fonction très cuspidale. Comme dans [W1] p.1265 on peut fixer un ensemble fini $\mathcal{S}$ d'éléments de $\mathfrak{h}_{reg}(F)$ dont tout les éléments ont un noyau nul dans $W$ et des nombres complexes $c_S$ pour $S\in\mathcal{S}$ tels que
$$\forall X\in\mathfrak{h}(F)\cap Supp(f)^G, \; \theta(X)=\displaystyle\sum_{S\in\mathcal{S}} c_S \hat{j}^H(S,X)$$
On pose alors

$$\mbox{(1)} \;\;\; K(\theta,f)=\displaystyle\sum_{S\in\mathcal{S}} \sum_{T\in\mathcal{T}(G)} c_S |W(G,T)|^{-1} \int_{\mathfrak{t}(F)^S} D^G(X)^{1/2} \hat{\theta}_f(X) dX$$
où l'ensemble $\mathfrak{t}(F)^S$ est défini comme dans la section 7. Puisque la fonction $\hat{\theta}_f$ est à support compact modulo conjugaison, les intégrales intervenant dans la définition de $K(\theta,f)$ sont à support compact, comme de plus $\hat{\theta}_f$ est un quasi-caractère (cf lemme 6.1 de [W1]) on en déduit d'après Harish-Chandra que les intégrales convergent.

\begin{lem}
On a $\lim\limits_{N\to\infty} J_N(\theta,f)=K(\theta,f)$
\end{lem}

\ul{Preuve}: Soit $\omega\subset \mathfrak{g}(F)$ un bon voisinage de $0$ et supposons que $Supp(f)\subset \omega$. Posons alors $\theta_\omega=\theta \bf{1}_{\omega}$. Par l'exponentielle, on relève les fonctions $f$ et $\theta_\omega$ en des fonctions sur respectivement $G(F)$ et $H(F)$ que l'on note $\bf{f}$ et $\Theta_\omega$. On a alors l'égalité $J_N(\theta,f)=J_N(\Theta_\omega,\bf{f})$. En appliquant la proposition 8.3.1 à $\bf{f}$, $\Theta_\omega$ et $x=1$ on obtient alors le résultat voulu.

\vspace{2mm}

Dans le cas général on peut toujours trouver $\lambda\in F^{\times 2}$ tel que $Supp(f)\subset\lambda\omega$. On pose alors $f'(X)=f(\lambda X)$, $\xi'(N)=\xi(\lambda N)$ et $\theta'(X)=\theta(\lambda X)$. Notons d'un indice $\xi'$ les expressions où on a remplacé $\xi$ par $\xi'$. Par les changements de variables $N\mapsto \lambda N$ et $X\mapsto \lambda X$ sur $\mathfrak{u}(F)$ et $\mathfrak{h}(F)$ respectivement, on a $J_{N,\xi'}(\theta',f')=|\lambda|_F^{-dim(U)-dim(H)}J_{N,\xi}(\theta,f)$. Grâce aux propriétés d'homogénéité des fonctions $\hat{j}^H$, on a

$$\displaystyle\theta'(Y)=\sum_{S\in\mathcal{S}} c_S \hat{j}^H(S,\lambda Y)=\sum_{S\in\mathcal{S}} |\lambda|_F^{-\delta(H)/2}c_S \hat{j}^H(\lambda S,Y)$$
pour tout $Y\in Supp(f')^G\cap\mathfrak{h}(F)$. On peut donc prendre $\mathcal{S'}=\lambda \mathcal{S}$ et les constantes $c_{\lambda S}=|\lambda|_F^{-\delta(H)/2}c_S$ pour définir $K_{\xi'}(\theta',f')$. Pour $T\in\mathcal{T}(G)$, on vérifie que l'ensemble $\mathfrak{t}(F)^{\lambda S}$ défini par $\xi'$ et $\lambda S$ est égal à $\lambda \mathfrak{t}(F)^S$. De plus, $\hat{\theta}_{f'}(X)=|\lambda|_F^{-dim(G)}\hat{\theta}_f(\lambda^{-1} X)$ et $D^G(\lambda X)=|\lambda|_F^{-\delta(G)}D^G(X)$. Par changement de variable on obtient 

$$K_{\xi'}(\theta',f')=|\lambda|_F^{-dim(G)+dim(T)+\delta(G)/2-\delta(H)/2}K_\xi(\theta,f)$$

On a $-dim(G)+dim(T)+\delta(G)/2=-\delta(G)/2$. On laisse au lecteur le soin de vérifier que $\delta(G)/2+\delta(H)/2=dim(U)+dim(H)$. On s'est alors ramené au premier cas. $\blacksquare$

\subsection{Une première approximation}

On pose $E(\theta,f)=K(\theta,f)-J_{geom}(\theta,f)$ pour $\theta$ un quasi-caractère de $\mathfrak{h}(F)$ et $f\in C_c^\infty(\mathfrak{g}(F))$ une fonction très cuspidale, le but étant de montrer que $E(.,.)$ est identiquement nul. Le lemme suivant donne une première approximation de ce terme d'erreur.

\begin{lem}
La forme bilinéaire $(\theta,f)\mapsto E(\theta,f)$ est proportionnelle à la forme bilinéaire 

$$(\theta,f)\mapsto \displaystyle\sum_{\mathcal{O}^H\in Nil_{reg}(\mathfrak{h}(F))}c_{\theta,\mathcal{O}^H}(0) \sum_{\mathcal{O}\in Nil_{reg}(\mathfrak{g}(F))} c_{\theta_f,\mathcal{O}}(0)$$
\end{lem}

\ul{Preuve}: Montrons dans un premier temps que $E$ est combinaison linéaire des formes bilinéaires $(\theta,f)\mapsto c_{\theta,\mathcal{O}^H}(0)c_{\theta_f,\mathcal{O}}(0)$ pour $\mathcal{O}^H\in Nil(\mathfrak{h}(F))$ et $\mathcal{O}\in Nil(\mathfrak{g}(F))$. \\

 Si $c_{\theta,\mathcal{O}^H}(0)=0$ pour tout orbite $\mathcal{O}^H\in Nil(\mathfrak{h}(F))$ alors il existe un $G$-domaine $\omega\subset \mathfrak{g}(F)$ compact modulo conjugaison contenant 0 tel que $\theta_{|\omega\cap\mathfrak{h}(F)}=0$. On pose alors $f'=f-f_{|\omega}$. On a $E(\theta,f)=E(\theta,f')+E(\theta,f_{|\omega})$. Puisque $f'$ ne contient aucun élément nilpotent dans son support on a $E(\theta,f')=0$ (c'est l'hypothèse (Hyp)). En revenant aux définitions on voit que $J_{geom}(\theta,f_{|\omega})=0$ et $J_N(\theta,f_{|\omega})=0$ pour tout $N\geqslant 1$ donc $K(\theta,f_{|\omega})=0$ et par conséquent $E(\theta,f)=0$. \\

 Si $c_{\theta_f,\mathcal{O}}(0)=0$ pour tout $\mathcal{O}\in Nil(\mathfrak{g}(F))$ on peut trouver un $G$-domaine compact modulo conjugaison tel que $\theta_{f}(X)=0$ pour tout $X\in\omega$. On pose alors $f'=f-f_{|\omega}$. On a $E(\theta,f)=E(\theta,f')+E(\theta,f_{|\omega})$. On a encore $E(\theta,f')=0$ car le support de $f'$ ne contient aucun élément nilpotent. Comme de plus on a $\theta_{f_{|\omega}}=\theta_f \bf{1}_\omega=0$ on vérifie immédiatement que $K(\theta,f_{|\omega})=J_{geom}(\theta,f_{|\omega})=0$. Ainsi on a bien $E(\theta,f)=0$.

\vspace{2mm}

 Montrons maintenant que $E$ est combinaison linéaire des formes bilinéaires \\
\noindent $(\theta,f)\mapsto c_{\theta,\mathcal{O}^H}(0)c_{\theta_f,\mathcal{O}}(0)$ pour $\mathcal{O}^H\in Nil_{reg}(\mathfrak{h}(F))$ et $\mathcal{O}\in Nil_{reg}(\mathfrak{g}(F))$.
Pour $\lambda\in F^{\times 2}$ on pose $\theta^\lambda(X)=\theta(\lambda X)$ et $f^\lambda (X)=f(\lambda X)$. On a alors $\theta_{f^\lambda}=\theta_f^\lambda$. Par les propriétés d'homogénéité des germes de Shalika, on a donc $c_{\theta^\lambda,\mathcal{O}^H}(0)c_{\theta_{f^\lambda},\mathcal{O}}(0)=|\lambda|_F^{-[dim(\mathcal{O}^H)+dim(\mathcal{O})]/2} c_{\theta,\mathcal{O}^H}(0)c_{\theta_f,\mathcal{O}}(0)$. \\ 
Puisque $dim(\mathcal{O}^H)+dim(\mathcal{O})\leqslant \delta(H)+\delta(G)$ avec égalité si et seulement si $\mathcal{O}^H$ et $\mathcal{O}$ sont des orbites nilpotentes régulières, il suffit maintenant de montrer que pour tout $\lambda\in F^{\times 2}$ on a $E(\theta^\lambda, f^\lambda)=|\lambda|_F^{-[\delta(H)+\delta(G)]/2} E(\theta,f)$. On a vérifié lors de la démonstration du lemme 9.1.1 une telle propriété d'homogénéité pour $K(\theta,f)$ à condition de changer $\xi$ en $\xi^\lambda$. Il est clair que $J_{geom}(\theta,f)$ ne dépend pas du choix de $\xi$. De plus, par les propriétés d'homogénéité des fonctions $\hat{j}(\mathcal{O},.)$ on vérifie que

$$J_{geom}(\theta^\lambda,f^\lambda)=|\lambda|_F^{-[\delta(H)+\delta(G)]/2}J_{geom}(\theta,f)$$

 L'expression $K(\theta,f)$ ne dépend de $\xi$ que par la définition des ensembles $\mathfrak{t}(F)^S$ qui dépendent du choix de $\Xi$ mais changer $\xi$ en $\xi'$ revient à changer $\Xi$ en un élément $\Xi'$ qui lui est conjugué par $A(F)$ on vérifie alors que les ensembles $\mathfrak{t}(F)^S$ ne dépendent en fait pas du choix de $\xi$. Cela suffit pour obtenir la propriété d'homogénéité voulue.

\vspace{2mm}

 Il existe donc des nombres complexes $c_{\mathcal{O}^H,\mathcal{O}}$ tels que pour tout quasi-caractère $\theta$ de $\mathfrak{h}(F)$ et toute fonction très cuspidale $f\in C_c^\infty(\mathfrak{g}(F))$ on ait

$$E(\theta,f)=\displaystyle\sum_{\mathcal{O},\mathcal{O}^H} c_{\mathcal{O}^H,\mathcal{O}} c_{\theta,\mathcal{O}^H}(0) c_{\theta_f,\mathcal{O}}(0)$$

où la somme porte sur $\mathcal{O}^H\in Nil_{reg}(\mathfrak{h}(F))$ et $\mathcal{O}\in Nil_{reg}(\mathfrak{g}(F))$.\\

Reste à vérifier que les coefficients $c_{\mathcal{O}^H,\mathcal{O}}$ sont tous égaux. Si $G$ ou $H$ n'est pas quasidéployé il n'y a rien à montrer. Si $G$ et $H$ sont quasidéployés alors l'un des deux a une seule orbite nilpotente régulière et l'autre en a deux qui sont permutées par n'importe quel élément $\lambda\in F^\times -N_{E/F}(E^\times)$. Pour $\lambda\in F^\times$ quelconque l'égalité $E(\theta^\lambda,f^\lambda)=|\lambda|_F^{-dim(U)-dim(H)} E(\theta,f)$ reste valable et on a  
$$c_{\theta^\lambda,\mathcal{O}^H}(0)c_{\theta_{f^\lambda},\mathcal{O}}(0)=|\lambda|_F^{-[dim(\mathcal{O}^H)+dim(\mathcal{O})]/2} c_{\theta,\lambda\mathcal{O}^H}(0)c_{\theta_f,\lambda\mathcal{O}}(0)$$
En choisissant $\lambda\in F^\times-N_{E/F}(E^\times)$ on obtient l'égalité des deux coefficients recherchée. $\blacksquare$

\subsection{Calcul de germes de Shalika}

D'après le lemme précédent, la proposition 9.2.1 est établie si $G$ ou $H$ n'est pas quasi-déployé. On suppose donc dans la suite que $G$ et $H$ sont quasi-déployés.

\begin{lem}
Soit $B$ un sous-groupe de Borel de $G$ et $T_{qd}$ un tore maximal de $B$ tous définis sur $F$. Soit $X_{qd}\in \mathfrak{t}_{qd}(F)\cap \mathfrak{g}_{reg}(F)$ et $\mathcal{O}\in Nil(\mathfrak{g}(F))$ on a alors

$$
\Gamma_\mathcal{O}(X_{qd})= \left\{
    \begin{array}{ll}
        0 & \mbox{si } \mathcal{O} \mbox{ n'est pas régulière} \\
        1 & \mbox{si } \mathcal{O} \mbox{ est régulière}
    \end{array}
\right.
$$

\end{lem}

\ul{Preuve}: Le tore $T_{qd}$ est un Lévi de $G$ et la distribution $f\mapsto J_G(X_{qd},f)$ est induite à partir de la distribution $f\mapsto f(X_{qd})$. Donc $\Gamma_\mathcal{O}(X_{qd})$ est non nul si et seulement si $\mathcal{O}$ intervient dans l'orbite induite de l'obite $\{0\}$ de $\mathfrak{t}_{qd}(F)$. Cette condition équivaut à ce que $\mathcal{O}$ soit régulière. On en déduit la première égalité. Puisque $\Gamma_{\{0\}}^{T_{qd}}(X_{qd})=1$ on en déduit aussi la deuxième égalité. $\blacksquare$

\vspace{3mm}

\subsection{Preuve de la proposition 1}

D'après le lemme 9.2.1 il existe un nombre complexe $c$ tel que pour tout quasi-caractère $\theta$ de $\mathfrak{h}(F)$ et toute fonction très cuspidale $f\in C_c^\infty(\mathfrak{g}(F))$ on ait

$$E(\theta,f)=c\displaystyle\sum_{\mathcal{O}^H\in Nil_{reg}(\mathfrak{h}(F))}c_{\theta,\mathcal{O}^H} \sum_{\mathcal{O}\in Nil_{reg}(\mathfrak{g}(F))} c_{\theta_f,\mathcal{O}}$$

On doit donc montrer que $c=0$. \\

Soit $T_{qd}\in \mathcal{T}(G)$ resp. $T_{qd}^H\in\mathcal{T}(H)$ un tore maximal d'un sous-groupe de Borel défini sur $F$ de $G$ resp. $H$. Soit $X_{qd}\in\mathfrak{t}_{qd}(F)\cap\mathfrak{g}_{reg}(F)$ et $X_{qd}^H\in\mathfrak{t}_{qd}^H(F)\cap \mathfrak{h}_{reg}(F)$. D'après le résultat 6.3(3) de [W1] on peut trouver un voisinage $\omega_{X_{qd}}$ de $X_{qd}$ dans $\mathfrak{t}_{qd}(F)$ aussi petit que l'on veut et qui ne contient que des éléments réguliers et une fonction très cuspidale $f=f[X_{qd}]\in C_c^\infty(\mathfrak{g}(F))$ vérifiant les propriétés suivantes:

\begin{enumerate}
\item Pour $T\in \mathcal{T}(G)$ si $T\neq T_{qd}$, la restriction de $\hat{\theta}_{f}$ à $\mathfrak{t}(F)$ est nulle.
\item Pour toute fonction localement intégrable $\varphi$ sur $\mathfrak{t}_{qd}(F)$, invariante par $W(G,T_{qd})$ on a
$$|W(G,T)|^{-1}\displaystyle\int_{\mathfrak{t}_{qd}(F)} \varphi(X) D^G(X)^{1/2} \hat{\theta}_{f}(X) dX=vol(\omega_{X_{qd}})^{-1} \int_{\omega_{X_{qd}}} \varphi(X) dX$$
\item Pour $Y\in\mathfrak{g}_{reg}(F)$,
$$\theta_{f}(Y)=vol(\omega_{X_{qd}})^{-1}\displaystyle\int_{\omega_{X_{qd}}} \hat{j}^G(X,Y) dX$$
\end{enumerate}

La dernière égalité au voisinage de 0 devient

$$\theta_{f[X_{qd}]}(Y)=\hat{j}^G(X_{qd},Y)$$

Par conséquent $c_{\theta_f,\mathcal{O}}=\Gamma_{\mathcal{O}}(X_{qd})$.\\
Soit $\theta$ un quasi-caractère de $\mathfrak{h}(F)$ tel que $\theta(X)=\hat{j}^H(X_{qd}^H,X)$ pour $X\in Supp(f)^G\cap\mathfrak{h}(F)$. D'après le lemme 9.3.1 pour tout $\mathcal{O}\in Nil_{reg}(\mathfrak{g}(F))$ on a $c_{\theta_f,\mathcal{O}}=\Gamma_{\mathcal{O}}(X_{qd})=1$ et pour tout $\mathcal{O}^H\in Nil_{reg}(\mathfrak{h}(F))$ on a $c_{\theta,\mathcal{O}^H}=\Gamma_{\mathcal{O}^H}(X_{qd}^H)=1$. On a donc $E(\theta,f)=2c$.

\vspace{2mm}

Pour $T\in \mathcal{T}$, $T\neq\{1\}$ et $Y\in \mathfrak{t}(F)$ en position générale il existe selon la parité de $dim(G)$:

\begin{itemize}
\item un voisinage de $Y$ dans $\mathfrak{h}(F)$ qui ne rencontre aucune sous-algèbre de Borel de $\mathfrak{h}(F)$ si $dim(G)$ est impaire.
\item un voisinage de $Y$ dans $\mathfrak{g}(F)$ qui ne rencontre aucune sous-algèbre de Borel de $\mathfrak{g}(F)$ si $dim(G)$ est paire.
\end{itemize}

\noindent Dans les deux cas on en déduit que $c_{\theta}(Y)c_f(Y)=0$. Donc la contibution du tore $T$ dans $J_{geom}(\theta,f)$ est nulle. La contribution du tore $\{1\}$ étant

$$\displaystyle\frac{1}{2}\sum_{\mathcal{O}^H\in Nil_{reg}(\mathfrak{h}(F))}c_{\theta,\mathcal{O}^H} \sum_{\mathcal{O}\in Nil_{reg}(\mathfrak{g}(F))} c_{\theta_f,\mathcal{O}}$$
on a $J_{geom}(\theta,f)=1$.

\vspace{2mm}

D'après la propriété (2) de $f$ on a $K(\theta,f)=vol(\omega_{X_{qd}})^{-1} vol(\omega_{X_{qd}}\cap\mathfrak{t}_{qd}(F)^{X_{qd}^H})$. Il suffit donc d'avoir $vol(\omega_{X_{qd}}\cap\mathfrak{t}_{qd}(F)^{X_{qd}^H})=vol(\omega_{X_{qd}})$ pour en déduire que $c=0$. Puisqu'on peut choisir $\omega_{X_{qd}}$ aussi petit que l'on veut, cela revient à montrer qu'on peut trouver $X_{qd}$ et $X_{qd}^H$ comme précédemment tels que $X_{qd}\in \mathfrak{t}_{qd}(F)^{X_{qd}^H}$ (car $\mathfrak{t}_{qd}(F)^{X_{qd}^H}$ est un ouvert de $\mathfrak{t}_{qd}(F)$). Puisque $\Xi+X_{qd}^H+\Sigma^{X_{qd}^H}$ est dense dans $\Xi+X_{qd}^H+\Sigma$, il suffit encore de montrer que $\left(\Xi+X_{qd}^H+\Sigma\right)\cap (\mathfrak{t}_{qd}(F)\cap \mathfrak{g}_{reg}(F))^G\neq \emptyset$. On distingue alors deux cas:

\begin{itemize}
\item Si $dim(H)$ est paire, alors il existe $A\in\mathfrak{a}(F)$ et $\mu_0\in F$ tels que $\Xi+X_{qd}^H+\mu_0c(v_0,\eta v_0)+A\in (\mathfrak{t}_{qd}(F)\cap \mathfrak{g}_{reg}(F))^G$.
\item Si $dim(H)$ est impaire, alors on peut trouver une base $((w_{\pm i})_{i=1,\ldots,p},w_{p+1})$ de $W$ formée de vecteurs propres pour $X_{qd}^H$ vérifiant $h(w_i,w_j)=\delta_{i,-j}$ si $i,j\in\{\pm 1,\ldots,\pm p\}$, $h(w_{p+1},w_i)=0$ si $i\in\{\pm 1,\ldots,\pm p\}$, $h(w_{p+1},w_{p+1})=\nu_1$ et $X_{qd}^H w_{p+1}=s_{p+1}w_{p+1}$. Puisque $G$ est quasi-déployé $G_0$ est aussi quasi-déployé. On peut donc supposer $\nu_1=-\nu_0$. On fixe $\mu_0\in F$ tel que $s_{p+1}=2\nu_0 \mu_0\eta$. Pour $\lambda\in F$ on vérifie facilement que $X_{qd}^H+\mu_0c(v_0+\eta v_0)+c(v_0,\lambda w_{p+1})$ admet $v_0+w_{p+1}$ et $v_0-w_{p+1}$ pour vecteurs propres avec pour valeurs propres respectives $s_{p+1}+\nu_0 \lambda$ et $s_{p+1}-\nu_0 \lambda$. On en déduit qu'on peut choisir $\lambda\in F$ et $A\in \mathfrak{a}(F)$ tels que l'élément $X=X_{qd}^H+A+\mu_0 c(v_0,\eta v_0)+c(v_0,\lambda w_{p+1})$ admet pour vecteurs propres $v_{\pm i}$ $i=1,\ldots,r$, $w_{\pm i}$ $i=1,\ldots,p$ et $w_{p+1}\pm v_0$ avec des valeurs propres distinctes. Alors on a bien $\Xi+X\in (\mathfrak{t}_{qd}(F)\cap\mathfrak{g}_{reg}(F))^G$.
\end{itemize}

$\blacksquare$

\section{Preuve des théorèmes 5.4.1 et 5.5.1}
On démontre les théorèmes 5.4.1 et 5.5.1 par récurrence sur la dimension de $V$. Supposons les résultats établis pour tout les couples $(V',W')$ avec $dim(V')<dim(V)$

\begin{lem}
Soient $\theta$ un quasicaractère de $H(F)$ et $f\in C_c^\infty(G(F))$ une fonction très cuspidale qui ne contient aucun élément unipotent dans son support. On a alors
$$\lim\limits_{N\to\infty} J_N(\theta,f)=J_{geom}(\theta,f)$$
\end{lem}

\ul{Preuve}: Pour tout $x\in G_{ss}(F)$ fixons un bon voisinage $\omega_x$ de $0$ dans $\mathfrak{g}_x(F)$. Etant donné l'hypothèse sur $f$, on peut supposer que le support de $f$ ne rencontre pas $exp(\omega_1)$. Par un procédé de partition de l'unité, on peut supposer qu'il existe $x\in G_{ss}(F)$ tel que $f$ soit à support dans $(xexp(\omega_x))^G$. Si ${}^gf^\xi$ est non nul pour un certain $g\in G(F)$, il existe $h\in H(F)$ et $u\in U(F)$ tels que $hu\in (xexp(\omega_x))^G$. La partie semisimple de $hu$ étant conjuguée à $h$, on a aussi $h\in (xexp(\omega_x))^G$. Donc il existe $X\in\omega_x$ et $g\in G(F)$ tels que $h=g^{-1}xexp(X)g$. On a alors $g(D\oplus Z)\subset Ker(xexp(X)-1)\subset Ker(x-1)$ donc $g^{-1}xg\in H(F)$. On en déduit que si $x$ n'est conjugué à aucun élément de $H(F)$ on a $J_N(\theta,f)=0$ pour tout $N\geqslant 1$. On vérifie aussi aisément que si $x$ n'est conjugué à aucun élément de $H(F)$ alors $\theta_f$ est nul au voisinage de $H(F)$ donc $J_{geom}(\theta,f)=0$. On peut donc supposer que $x$ est conjugué à un élément de $H(F)$. Quitte à conjuguer $x$ dès le départ, on peut aussi bien supposer que $x\in H(F)$. Notons $\omega=\omega_x$ et reprenons les notations du début de la section 6. Quitte à changer les choix des bons voisinages $\omega_x$, on peut supposer que $\omega$ se décompose en un produit $\omega'\times\omega''$ avec $\omega'\subset \mathfrak{g}'_x(F)$ et $\omega''\subset \mathfrak{g}''(F)$ des bons voisinages de $0$. D'après les lemmes 6.1.1 et 6.2.1, on est ramené à prouver l'égalité

$$\lim\limits_{N\to\infty} J_{x,\omega,N}(\theta,f)=J_{geom,x,\omega}(\theta,f)$$

\noindent Comme au début de la section 8, on peut trouver une famille $\mathcal{S}$ d'éléments de $\mathfrak{h}''_{reg}(F)$ et une famille $(\hat{j}_S)_{S\in\mathcal{S}}$ de quasicaractères de $\mathfrak{h}'_x(F)$ tels que pour tout $X=X'+X''\in\mathfrak{h}_{x,reg}(F)$ on ait
$$\theta_{x,\omega}(X)=\displaystyle\sum_{S\in\mathcal{S}} \hat{j}_S(X')\hat{j}^{H''}(S,X'')$$
Notons $\theta''_S=\hat{j}^{H''}(S,.)$. On peut de la même façon trouver deux familles finies $(\theta'_{f,b})_{b\in B}$ et $(\theta''_{f,b})_{b\in B}$ de quasicaractères de $\mathfrak{g}'_x(F)$ et $\mathfrak{g}''(F)$ respectivement telles que pour tout $X=X'+X''\in\mathfrak{g}_x(F)$

$$\theta_{f,x,\omega}(X)=\displaystyle\sum_{b\in B} \theta'_{f,b}(X')\theta''_{f,b}(X'')$$
D'après la proposition 6.4 de [W1], pour tout $b\in B$, on peut trouver une fonction très cuspidale $f''_b\in C_c^\infty(\mathfrak{g}''(F))$ telle que $\theta''_{f,b}=\theta_{f''_b}$. En comparant 5.5(1) et 6.2(1) on a

$$J_{geom,x,\omega}(\theta,f)=\displaystyle\sum_{b\in B,S\in \mathcal{S}} J'(S,b)J_{geom}(\theta''_S,f''_b)$$
où

$$J'(S,b)=\displaystyle\sum_{T'\in\mathcal{T}_{ell}(H'_x)} |W(H'_x,T')|^{-1} |\nu(T')| \int_{\mathfrak{t}'(F)} \hat{j}_S(X')\theta'_{f,b}(X') D^{H'_x}(X') dX'$$
Dans la section 8, on a montré que $J_{x,\omega,N}(\theta,f)$ admet une limite et on a calculé cette limite: c'est $K_{x,\omega}(\theta,f)$. Si $A_{G_x}\neq \{1\}$ alors $K_{x,\omega}(\theta,f)=0$ et l'ensemble de tores $\underline{\mathcal{T}}_x$ est vide donc $J_{geom,x,\omega}(\theta,f)=0$. Sinon, en comparant les formules 8.3(1) et 9.1(1) on a

$$K_{x,\omega}(\theta,f)=\displaystyle\sum_{b\in B, S\in\mathcal{S}} J'(S,b) K(\theta''_S,f''_b)$$
D'après l'hypothèse de récurrence appliquée au couple $(V'',W'')$ on a $K(\theta''_S,f''_b)=J_{geom}(\theta''_S,f''_b)$ d'où le résultat $\blacksquare$

\subsection{Preuve de 5.5.1}

D'après la section 9, on peut se restreindre aux fonctions très cuspidales $f\in C_c^\infty(\mathfrak{g}(F))$ qui ne contiennent pas d'élément nilpotent dans leur support. Soient $f$ une telle fonction et $\theta$ un quasicaractère sur $\mathfrak{h}(F)$. On a vu dans la section 9 que $J_N(\theta,f)$ et $J_{geom}(\theta,f)$ sont homogènes de même degré pour la transformation $\lambda\mapsto (\theta^\lambda,f^\lambda)$. On peut donc supposer que le support de $f$ est dans un bon voisinage $\omega$ de $0$ dans $\mathfrak{g}(F)$. Par l'exponentielle, on relève les fonctions $f$ et $\theta \mathbf{1}_\omega$ en des fonctions sur respectivement $G(F)$ et $H(F)$ que l'on note $\mathbf{f}$ et $\Theta_\omega$. On vérifie alors que $J_{geom,1,\omega}(\Theta_\omega,\mathbf{f} )=J_{geom}(\theta,f)$ et $J_{N,1,\omega}(\Theta_\omega,\mathbf{f})=J_N(\theta,f)$. Puisque $f$ ne contient pas d'élément nilpotent dans son support, $\mathbf{f}$ n'a pas d'élément unipotent dans son support et d'après le lemme précédent, on a donc $\lim\limits_{N\to \infty} J_{N,1,\omega}(\Theta_\omega,\mathbf{f})=J_{geom,1,\omega}(\Theta_\omega,\mathbf{f})$. \\

\subsection{Preuve de 5.4.1}

Soient $\theta$ un quasicaractère de $H(F)$ et $f\in C_c^\infty(G(F))$ une fonction très cuspidale. Soit $\omega$ un bon voisinage de $0$ dans $\mathfrak{g}$. Alors $f-f\mathbf{1}_{exp(\omega)}$ ne contient pas d'élément unipotent dans son support. D'après le lemme 10.0.1, il suffit donc d'établir que 
$\lim\limits_{N\to \infty} J_N(\theta,f\mathbf{1}_{exp(\omega)})=J_{geom}(\theta,f\mathbf{1}_{exp(\omega)})$. D'après les lemmes 6.1.1 et 6.2.1, on a les égalités $J_N(\theta,f\mathbf{1}_{exp(\omega)})=J_N(\theta_{1,\omega},f_{1,\omega})$ et $J_{geom}(\theta,f\mathbf{1}_{exp(\omega)})=J_{geom}(\theta_{1,\omega},f_{1,\omega})$. Le résultat sur l'algèbre de Lie que l'on vient d'établir permet alors de conclure.

\section{Décomposition de Cartan relative}

Dans cette section on se propose de montrer une sorte de "décomposition de Cartan relative" telle que celle qui est montrée dans [BO] et [DS] (pour les variétés symmétriques) ou dans [Sa] théorème 2.3.8 (pour des variétés sphériques associées à des groupes déployés). \\

Posons $\nu_0=h(v_0)$ et soit $w_0=v_0,w_1,\ldots,w_l$ une famille maximale de vecteurs deux à deux orthogonaux de $V_0$ vérifiant $h(w_i)=(-1)^i \nu_0$. \\
Posons $r_0=E(\frac{l+1}{2})$ et $r_1=E(\frac{l}{2})$. Introduisons $u_{\pm i}=w_{2i-2}\pm w_{2i-1}$, $i=1\ldots,r_0$ et $u'_{\pm i}=w_{2i-1}\pm w_{2i}$, $i=1,\ldots, r_1$. Les familles $(u_{\pm i})$ et $(u'_{\pm i})$ sont des familles hyperboliques maximales de $V_0$ et $W$ respectivement. Soient $V_{an,0}$ resp. $W_{an}$ l'orthogonal de $Eu_1\oplus Eu_{-1}\oplus\ldots\oplus Eu_{r_0}\oplus Eu_{-r_0}$ dans $V_0$ resp. de $Eu'_1\oplus Eu'_{-1}\oplus\ldots\oplus Eu'_{r_1}\oplus Eu'_{-r_1}$ dans $W$. Ce sont des sous-espaces hermitiens totalement anisotropes. Définissons $P_0$ et $P_H$ comme les sous-groupes paraboliques de $G_0$ et $H$ qui conservent les drapeaux de sous-espaces totalement isotropes

\begin{center}
$Eu_1\subset Eu_1\oplus Eu_2\subset \ldots \subset Eu_1\oplus\ldots\oplus Eu_{r_0}$ et $Eu'_1\subset \ldots\subset Eu'_1\oplus\ldots\oplus Eu'_{r_1}$ resp.
\end{center}
respectivement.

\begin{lem}
L'application $\overline{P}_0\times P_H\to G_0$, $(\overline{p}_0,p_H)\mapsto \overline{p_0}p_H$ est submersive à l'origine.
\end{lem}

\ul{Preuve}: Il s'agit de vérifier l'égalité $\overline{\mathfrak{p}}_0+\mathfrak{p}_H=\mathfrak{g}_0$. Démontrons cela par récurrence sur $l$. Pour $l=0$ c'est évident car $\overline{\mathfrak{p}}_0=\mathfrak{g}_0$. Soit $l\geqslant 1$ et supposons le résultat vérifié pour $l'=l-1$. Soit $W'$ l'orthogonal dans $V_0$ du sous-espace engendré par $u_{\pm 1}$, $H'$ son groupe unitaire et $P_{H'}$ le sous groupe parabolique de $H'$ qui conserve le drapeau

$$Eu_2\subset\ldots\subset Eu_2\oplus\ldots\oplus Eu_{r_0}$$

Alors par hypothèse de récurrence, on a $\overline{\mathfrak{p}_{H'}}+\mathfrak{p}_{H}=\mathfrak{h}$ et par conséquent $\mathfrak{h}\subset \overline{\mathfrak{p}_0}+\mathfrak{p}_H$. On en déduit que $(\overline{\mathfrak{p}_0}+\mathfrak{p}_H)^\perp\subset \overline{\mathfrak{p}_0}^\perp\cap \mathfrak{h}^\perp=\overline{\mathfrak{u}}_0\cap \mathfrak{h}^\perp$ où $\overline{U}_0$ est la radical unipotent de $\overline{P}_0$ et $\overline{\mathfrak{u}}_0$ est son algèbre de Lie. \\
Il suffit donc de vérifier que $\overline{\mathfrak{u}}_0\cap \mathfrak{h}^\perp=0$. On a $\mathfrak{h}^\perp=\{c(v_0,w), w\in W\}\oplus Fc(v_0,\eta v_0)$ où $\eta$ est un élément non nul de trace nulle dans $E$. Soit $N=c(v_0,w)+\lambda c(v_0,\eta v_0)\in \overline{\mathfrak{u}}_0\cap \mathfrak{h}^\perp$. On a alors $0=Nu_{-1}=\nu_0 w+h(w,w_1)v_0+\lambda\nu_0 (\eta-\overline{\eta})v_0$ d'où $w=0$ puis $\lambda=0$. $\blacksquare$

\vspace{4mm}

Soient $A_0$ et $A_H$ les tores déployés maximaux de $G_0$ et $H$ respectivement qui laissent stable les droites $Eu_{i}$ ($i=\pm1,\ldots,\pm r_0$) et $Eu'_{i}$ ($i=\pm1,\ldots,\pm r_1$) respectivement. On notre $A_0^{+}$ et $A_H^+$ les sous-ensembles de $A_0(F)$ et $A_H(F)$ qui contractent $P_0(F)$ et $P_H(F)$ respectivement. Posons

$$R_0=\mathcal{O}_Eu_1\oplus\mathcal{O}_Eu_{-1}\oplus\ldots\oplus \mathcal{O}_Eu_r\oplus \mathcal{O}_E u_{-r}\oplus R_{an,0}$$

\noindent où $R_{an,0}=\{v\in V_{an,0}; val_E(h(v))\geqslant val_E(\nu_0)-1\}$. Pour $\mathcal{V}$ un $E$-espace vectoriel, $\mathcal{R}$ un $\mathcal{O}_E$-réseau de $\mathcal{V}$ et $g\in GL_E(\mathcal{V})$ on note $val_{\mathcal{R}}(g)=inf\{val_{\mathcal{R}}(gv), v\in\mathcal{R}\}$ et $||g||_{\mathcal{R}}=|\pi_E|_E^{val_{\mathcal{R}}(g)}$.

\begin{lem}
Il existe un sous-ensemble compact $C_1$ de $G_0(F)$ vérifiant la propriété suivante
\begin{center}
Pour tous $v,v'\in V_0$ tels que $h(v)=h(v')=\nu_0$ et $val_{R_0}(v)=val_{R_0}(v')$, il existe $\gamma\in C_1$ tel que $\gamma v=v'$
\end{center}
\end{lem}

\ul{Preuve}: Si $l=0$, $G_0(F)$ est compact et il suffit de prendre $C_1=G_0(F)$. Supposons dorénavant $l\geqslant 1$. Soit $K_0=Stab_{G_0(F)}(R_0)$, c'est un sous-groupe ouvert-compact de $G_0(F)$. Pour tout $g\in G_0(F)$, on notera $||g||=||g||_{R_0}$. La fonction $g\mapsto ||g||$ est minorée par une constante strictement positive sur $G_0(F)$ et est biinvariante par $K_0$. Soit $v\in V_0$. Montrons

\vspace{2mm}

 (1) Il existe $k\in K_0$ tel que $kv\in Eu_1\oplus Eu_{-1}\oplus V_{0,an}$.

\vspace{2mm}

\noindent On a une décomposition $v=\displaystyle\sum_{i=\pm1,\ldots,\pm r} \lambda_i u_i+v_{an}$ où $v_{an}\in V_{0,an}$. Quitte à multiplier $v$ par un élément du groupe de Weyl de $A_0$ (identifié à un sous-groupe de $K_0$), on peut supposer que $val(\lambda_1)=inf_{i=\pm1,\ldots,\pm r} val(\lambda_i)$.\\
Soit $k_1\in K_0$ l'élément qui envoie  $u_1$ sur $u_1-\frac{\lambda_2}{\lambda_1}u_2-\ldots-\frac{\lambda_r}{\lambda_1} u_r$, $u_{-i}$ sur $u_{-i}+ \overline{(\frac{\lambda_i}{\lambda_1})}u_{-1}$ pour $i=2,\ldots,r$ et qui agit trivialement sur $V_{0,an}+Eu_{-1}+Eu_2+\ldots+Eu_r$. On a alors
$$k_1 v=\lambda_1u_1+\mu_1 u_{-1}+\lambda_{-2} u_{-2}+\ldots+\lambda_{-r}u_{-r}+v_{an}$$
Soit $k_2\in K_0$ l'élément qui envoie $u_1$ sur $u_1-\frac{\lambda_{-2}}{\lambda_1}u_{-2}-\ldots-\frac{\lambda_{-r}}{\lambda_1}u_{-r}$, $u_i$ sur $u_i+\overline{(\frac{\lambda_{-i}}{\lambda_1})}u_{-1}$ et qui agit trivialement sur $V_{0,an}+Eu_{-1}+\ldots+Eu_{-r}$. On a alors $k_2k_1 v\in Eu_1\oplus Eu_{-1}\oplus V_{0,an}$.

\vspace{2mm}

 (2) Pour tout $v_{an}\in V_{0,an}$, on a 
$$\displaystyle R_{an,0}\subset \frac{1}{2}\left((R_{an,0}\cap Ev_{an})\oplus (R_{an,0}\cap (v_{an})^\perp)\right)$$
où $(v_{an})^\perp$ est l'orthogonal de $v_{an}$ dans $V_{0,an}$.

\vspace{2mm}

\noindent En effet, soit $v\in R_{an,0}$ et $v=v_1+v_2$ sa décomposition suivant $V_{an,0}=Ev_{an}\oplus (v_{an})^{\perp}$. Supposons que $val_E(h(v_1)+h(v_2))> min\left(val_E(h(v_1)),val_E(h(v_2))\right)+val_E(4)$ alors $-\frac{h(v_1)}{h(v_2)}\in 1+4\mathfrak{p}_F$. Puisque $1+4\mathfrak{p}_F\subset \mathcal{O}_F^{\times,2}$, il existerait $\lambda\in\mathcal{O}_F^\times$ tel que $h(v_1+\lambda v_2)=0$, contredisant le fait que $V_{0,an}$ est anisotrope. Par conséquent, on a

$$min(val_E(h(2v_1)),val_E(h(2v_2)))\geqslant val_E(h(v_1+v_2))=val_E(h(v))\geqslant val_E(\nu_0)-1$$

\noindent d'où $v_1,v_2\in \frac{1}{2}R_{0,an}$.

\vspace{2mm}

 (3) Il existe $g\in G_0(F)$ tel que $gv=\lambda_1 u_1+ \lambda_{-1} u_{-1}+ v_{an}$ avec $v_{an}\in V_{0,an}$, \\

$val_E(\lambda_1)=inf\left(val_E(\lambda_1),val_E(\lambda_{-1}), val_{R_{an,0}}(v_{an})\right)$ et $||g||\leqslant |2|_E^{-1}$

\vspace{2mm}

\noindent D'après (1), on peut supposer que $v=\lambda_1 u_1+\lambda_{-1}u_{-1}+v_{an}$ et quitte à appliquer l'élément qui échange $u_1$ et $u_{-1}$ et laisse stable l'orthogonal du plan engendré par $u_1$ et $u_{-1}$, on peut aussi supposer que $val_E(\lambda_1)\leqslant val_E(\lambda_{-1})$. Si $val_E(\lambda_1)\leqslant val_{R_{an,0}}(v_{an})$, il n'y a rien à ajouter. Si au contraire $-d=val_{R_{an,0}}(v_{an})<val_F(\lambda_1)$, on considère l'élément $g\in G_0(F)$ qui envoie $u_1$ sur $u_1$, $v_{an}$ sur $v_{an}+\pi_E^{-d}u_1$, $u_{-1}$ sur $\displaystyle u_{-1}-\frac{2\overline{\pi_E}^{-d}\nu_0}{h(v_{an})}v_{an}-\frac{N(\pi_E)^{-d}\nu_0}{h(v_{an})}u_1$ et qui agit trivialement sur l'orthogonal de $Eu_1+Eu_{-1}+Ev_{an}$. Puisque $val_E(h(v_{an}))\leqslant -2d+val_E(\nu_0)$, on a 

$$g(R_0\cap (Eu_1+Eu_{-1}+Ev_{an}))=R_0\cap (Eu_1+Eu_{-1}+Ev_{an})$$

\noindent D'après (2) on a donc $||g||\leqslant |2|_E^{-1}$. Toujours en utilisant l'inégalité $val_E(h(v_{an}))\leqslant -2d+val_E(\nu_0)$, on vérifie aisément que $gv$ est bien de la forme désirée.

\vspace{2mm}

\begin{center}
(4) Il existe $g\in G_0(F)$ tel que $gv\in Eu_1+Eu_{-1}$ et $||g||\leqslant |\pi_E|_{E}^{-1}|2|_{E}^{-4}$.
\end{center}

D'après (3) on peut supposer que $v=\lambda_1 u_1+ \lambda_{-1} u_{-1}+ v_{an}$ avec $v_{an}\in V_{0,an}$, $val_E(\lambda_1)=inf(val_E(\lambda_1),val_E(\lambda_{-1}), val_{R_{an}}(v_{an}))$ et montrer le résultat avec un facteur 2 de moins. Soit $g\in G_0(F)$ l'élément qui agit ainsi sur $u_{\pm 1}$ et $v_{an}$:

$$\displaystyle u_{-1}\mapsto u_{-1},\; u_1\mapsto u_1-\frac{1}{\lambda_1}v_{an}-\frac{h(v_{an})}{4\nu_0N(\lambda_1)} u_{-1}$$

$$\displaystyle v_{an}\mapsto v_{an}+\frac{h(v_{an})}{2\overline{\lambda_1} \nu_0} u_{-1}$$
et agit trivialement sur l'orthogonal de $Eu_1+Eu_{-1}+Ev_{an}$. Alors $gv$ est bien de la forme désirée et d'après (2) on a

\[\begin{aligned}
val_{R_0}(g)\leqslant inf(0,val_{R_{an,0}}(v_{an})-val_E(\lambda_1),val_E(h(v_{an}))-val_E(4\nu_0N(\lambda_1)), \\
val_E(h(v_{an}))-val_E(2\lambda_1 \nu_0)-val_{R_{an,0}}(v_{an}))-val_E(2)
\end{aligned}\]

\noindent Puisque $val_E(h(v_{an}))\geqslant 2val_{R_{an,0}}(v_{an})+val_E(\nu_0)-1$, on a $val_{R_0}(g)\geqslant -1-3val_E(2)$.

\vspace{2mm}

Supposons maintenant que $h(v)=\nu_0$.

\begin{center}
(5) Il existe $g\in G_0(F)$ et $\lambda\in \mathcal{O}_E-\{0\}$ tels que $gv=\lambda u_1+\frac{\nu_0}{2}\overline{\lambda}^{-1} u_{-1}$ et $||g||\leqslant |\pi_E|_E^{-1}|2|_E^{-5}$.
\end{center}

D'après (4), on sait qu'il existe $g'\in G_0(F)$ de sorte que $g'v=\lambda u_1+\mu u_{-1}$ et $||g'||\leqslant |\pi_E|_E^{-1}|2|_E^{-4}$. Quitte à multiplier $g'$ à gauche par l'élément qui échange $u_1$ et $u_{-1}$ et agit trivialement sur l'orthogonal du plan engendré par ces deux derniers, on peut supposer que $val(\mu)\geqslant val(\lambda)$. Soit $\displaystyle u=\frac{\nu_0}{2N(\lambda)}-\frac{\mu}{\lambda}$. On a $\displaystyle \frac{\mu}{\lambda}\in\mathcal{O}_E$ et $\nu_0=h(g'v)=Tr(\lambda\overline{\mu})$ d'où $val_E(\nu_0)\geqslant val_E(\lambda)+val_E(\mu)\geqslant 2val_E(\lambda)=val_E(N(\lambda))$. Par conséquent $u\in \frac{1}{2}\mathcal{O}_E$. On vérifie facilement que $Tr(u)=0$. Soit $g''\in G_0(F)$ l'élément qui envoie $u_1$ sur $u_1+uu_{-1}$ et qui agit trivialement sur l'orthogonal de $u_{-1}$. C'est bien un élément de $G_0(F)$ car $Tr(u)=0$ et d'après ce que l'on vient de voir $||g''||\leqslant |2|_E^{-1}$. On a alors $g''g'v=\lambda u_{-1}+\overline{\lambda}^{-1}\frac{\nu_0}{2} u_1$ et $||g''g'||\leqslant |\pi_E|_E^{-1}|2|_E^{-5}$.

\vspace{2mm}

On peut maintenant terminer la preuve du lemme. Soient $v$ et $v'$ comme dans le lemme. Il existe $g_1,g_2\in G_0(F)$ et $\lambda_1,\lambda_2\in\mathcal{O}_E-\{0\}$ vérifiant

\begin{itemize}
\item $||g_1||,||g_2||\leqslant |\pi_E|_E^{-1}|2|_E^{-5}$
\item $g_1 v=\lambda_1 u_1+ \frac{\nu_0}{2}\overline{\lambda_1}^{-1} u_{-1}$ et $g_2 v'=\lambda_2 u_1+ \frac{\nu_0}{2} \overline{\lambda_2}^{-1}u_{-1}$
\end{itemize}

\noindent Il existe un entier $N_0$ tel que

$|val_{R_0}(gv_0)-val_{R_0}(v_0)|\leqslant -val_{R_0}(g)+N_0$

\noindent pour tout $g\in G_0(F)$ et pour tout $v_0\in V_0$. Puisque $val_{R_0}(g_1v)=-val_E(\lambda_1)-val_E(2)$ et $val_{R_0}(g_2v')=-val_E(\lambda_2)-val_E(2)$ et que $val_{R_0}(v)=val_{R_0}(v')$, on a 

$$|val_E(\lambda_1)-val_E(\lambda_2)|\leqslant 2+10val_E(2)+2N_0$$

\noindent Soit $a$ l'élément de $G_0(F)$ qui envoie $u_1$ sur $\displaystyle\frac{\lambda_2}{\lambda_1} u_1$, $u_{-1}$ sur $\displaystyle\overline{\frac{\lambda_1}{\lambda_2}} u_{-1}$ et qui agit trivialement sur l'orthogonal de $Eu_1+Eu_{-1}$. On a alors $ag_1v=g_2v'$. D'où $g_2^{-1}ag_1v=v'$ et $||g_2^{-1}ag_1||$ est borné par une constante. $\blacksquare$

\vspace{4mm}

Supposons $l\geqslant 1$. Pour tout $\lambda\in E^\times$, on note $a(\lambda)$ l'élément de $G_0(F)$ qui envoie $u_1$ sur $\lambda u_1$, $u_{-1}$ sur $\overline{\lambda}^{-1}u_{-1}$ et qui agit trivialement sur l'orthogonal de $Eu_1+Eu_{-1}$. Soit $R_{\sharp,H}$ un $\mathcal{O}_E$-réseau de l'orthogonal de $Ew_1\oplus\ldots\oplus Ew_l$ dans $W$. On pose alors $R_1=\mathcal{O}_Ew_1\oplus \ldots\oplus \mathcal{O}_Ew_l\oplus R_{\sharp,H}$ et $R_2=\mathcal{O}_Ew_2\oplus \ldots\oplus \mathcal{O}_Ew_l\oplus R_{\sharp,H}$.

\begin{lem}
Supposons $l\geqslant 1$. Il existe un sous-ensemble compact $C_2$ de $G_0(F)$ et une constante $c_0>0$ vérifiant la propriété suivante

\begin{center}
Pour tout $h\in H(F)$ et tout $\lambda\in \mathcal{O}_E-\{0\}$, il existe $h'\in h(a(\lambda)C_2a(\lambda)^{-1}\cap H(F))$ tel que
$$val_{R_1}(h')\geqslant val_{R_1}(h'w_1)-val_E(\lambda)-c_0$$
\end{center}
\end{lem}

\ul{Preuve}:Soit $e_1,\ldots,e_t$ une base orthogonale du $\mathcal{O}_E$-module $R_2$. Il existe un entier strictement positif $\alpha$ vérifiant les deux propriétés suivantes

\begin{itemize}
\item $\displaystyle\frac{h(e_i)}{h(w_1)}\in \mathfrak{p}_E^{-\alpha}$, pour $i=1,\ldots,t$;
\item pour tout entier $k\geqslant 0$ et pour tout $x\in (1+\mathfrak{p}_E^{k+\alpha})\cap\mathcal{O}_F$, il existe $y\in 1+\mathfrak{p}_E^k$ tel que $N(y)=x$.
\end{itemize}

Pour $i=1,\ldots,t$ et $\lambda\in\mathcal{O}_E-\{0\}$, on considère l'élément $\gamma_i(\lambda)\in H(F)$ qui agit trivialement sur $(Ew_1+Ee_i)^\perp$ et qui agit de la façon suivante sur $w_1$ et $e_i$

$$\displaystyle w_1\mapsto aw_1+\lambda \pi_E^{\alpha} e_i$$
$$\displaystyle e_i\mapsto \frac{h(e_i)}{h(w_1)}\overline{\lambda}\overline{\pi}_E^{\alpha}w_1-\overline{a} e_i$$

\noindent où $a\in 1+\lambda^2\mathcal{O}_E$ est tel que $\displaystyle N(a)=1-N(\lambda \pi_E^{\alpha})\frac{h(e_i)}{h(w_1)}$. Alors $\gamma_i(\lambda)$ laisse stable $\mathcal{O}_E w_1\oplus R_2$ donc reste dans un compact indépendant de $i$ et de $\lambda$. Vérifions que $a(\lambda)^{-1}\gamma_i(\lambda)a(\lambda)$ reste aussi dans un compact. Cet élément envoie $u_{-1}$ sur 

$$\overline{\lambda}^{-1} a(\lambda)^{-1} \gamma_i(\lambda)u_{-1}=\overline{\lambda}^{-1} a(\lambda)^{-1}(w_0-\gamma_i(\lambda)w_1)=u_{-1}+\overline{\lambda}^{-1} a(\lambda)^{-1}(w_1-\gamma_i(\lambda)w_1)$$

Mais $w_1-\gamma_i(\lambda)w_1$ est dans $\lambda^2 \mathcal{O}_E w_1\oplus R_2$. Ainsi $a(\lambda)^{-1}\gamma_i(\lambda)a(\lambda)u_{-1}\in \frac{1}{2}\mathcal{O}_Eu_{-1}+R_1$. On montre de la même façon que $a(\lambda)^{-1}\gamma_i(\lambda)a(\lambda)u_1$ reste borné. Enfin, puisque $a(\lambda)^{-1}\gamma_i(\lambda)a(\lambda)e_i=a(\lambda)^{-1}\gamma_i(\lambda)e_i$ et que $\gamma_i(\lambda)e_i\in \lambda \mathcal{O}_E w_1+R_2$, $a(\lambda)^{-1}\gamma_i(\lambda)a(\lambda)e_i$ reste aussi borné. \\

Soit $C_2$ un sous-ensemble compact de $G_0(F)$ qui contient tout les $a(\lambda)^{-1}\gamma_i(\lambda)a(\lambda)$ pour tout $i=1,\ldots,t$ et pour tout $\lambda\in\mathcal{O}_E-\{0\}$. Puisque les $\gamma_i(\lambda)$ sont bornés il existe une constante $c_1>0$ telle que pour tout $h\in H(F)$, pour tout $i=1,\ldots,t$ et pour tout $\lambda\in\mathcal{O}_E-\{0\}$ on ait

$$val_{R_1}(h\gamma_i(\lambda))\geqslant val_{R_1}(h)-c_1$$

On a $val_{R_1}(h\gamma_i(\lambda)w_1)=inf(val_{R_1}(hw_1),val_E(\lambda)+\alpha+val_{R_1}(he_i))$. Puisque $val_{R_1}(h)=inf(val_{R_1}(hw_1),val_{R_1}(he_1),\ldots,val_{R_1}(he_t))$, on a

$$inf_{i=1,\ldots,t}(val_{R_1}(h\gamma_i(\lambda)w_1))\leqslant val_E(\lambda)+\alpha+val_{R_1}(h)\leqslant val_E(\lambda)+\alpha+c_1+inf_{i=1,\ldots,t}val_{R_1}(h\gamma_i(\lambda))$$

 Pour obtenir le lemme, il suffit donc de prendre $c=\alpha+c_1$ et pour $h\in H(F)$ et $\lambda\in \mathcal{O}_E-\{0\}$ de choisir $h'=h\gamma_j(\lambda)$ où $j$ est tel que $inf_{i=1,\ldots,t}(val_{R_1}(h\gamma_i(\lambda)w_1))=val_{R_1}(h\gamma_j(\lambda)w_1)$ $\blacksquare$

\vspace{4mm}

\begin{prop}
Il existe des sous-ensembles compacts $C_0\subset G_0(F)$ et $C_H\subset H(F)$ tels que
$$G_0(F)=C_HA_H^{+}A_0^{+}C_0$$
\end{prop}

\ul{Preuve}: On démontre le résultat par récurrence sur $l$. Si $l=0$, c'est vrai car $G_0(F)$ est compact. On peut donc supposer que $l\geqslant 1$. Le lemme 11.0.2 nous fournit un compact $C_1\subset G_0(F)$. Le lemme 11.0.3 nous fournit un compact $C_2\subset G_0(F)$ ainsi qu'une constante $c_0$. Posons $C=C_2^{-1}.C_1$. Soit $g\in G_0(F)$. D'après le lemme 11.0.2, il existe $k_1\in C_1$ et $\lambda\in \mathcal{O}_E-\{0\}$ tels que

$$k_1g^{-1}v_0=a(\lambda)^{-1}v_0$$
 
\noindent c'est-à-dire $g\in H(F)a(\lambda)k_1$. En utilisant le lemme 11.0.3, on voit qu'il existe $k_2\in C_2$, $\lambda\in \mathcal{O}_E-\{0\}$ et $h\in H(F)$ tels que

\begin{center}
$g=ha(\lambda)k_2^{-1}k_1$ et 

$$\mbox{(1)} \;\;\; val_{R_1}(h)\geqslant val_{R_1}(hw_1)-val_E(\lambda)-c_0$$
\end{center}

Posons $k=k_2^{-1}k_1\in C$. Soit $H'$ le groupe unitaire de $W'=(Ew_0\oplus Ew_1)^\perp$, $P_{H'}$ le sous-groupe parabolique qui conserve le drapeau $Eu_2\subset\ldots\subset Eu_2\oplus\ldots\oplus Eu_r$, $A_{H'}$ le tore maximal déployé qui conserve les droites $Eu_i$ ($i=\pm2,\ldots,\pm r_0$) et $A_{H'}^{+}$ le sous-ensemble de $A_{H'}(F)$ des éléments positifs pour $P_{H'}$. Par hypothèse de récurrence, il existe des sous-ensembles compacts $C_{H}^\sharp\subset H(F)$ et $C_{H'}^\sharp\subset H'(F)$ tels que $H(F)=C_{H}^\sharp A_H^+A_{H'}^+C_{H'}^\sharp$. Suivant cette décomposition, on peut écrire $h=k_1^\sharp a_Ha_{H'}k_2^\sharp$. On a alors $g=k_1^\sharp a_Ha_{H'}a(\lambda)k_2^\sharp k$. On a bien $a_{H'}a(\lambda)\in A_0(F)$ mais rien n'assure qu'il soit dans la chambre positive pour $P_0$. On va utiliser l'inégalité (1) pour vérifier que $a_{H'}a(\lambda)$ est bien dans $A_0^+$ modulo un sous-ensemble fini.
 
\vspace{2mm}

Soit $a_{H',2}$ la valeur propre de $a_{H'}$ agissant sur $u_2$. Il s'agit de vérifier que $val_E(a_{H',2})-val_E(\lambda)$ est borné par une constante. Il existe des constantes positives $c_1$, $c_2$, $c_3$ et $c_4$ telles que

\begin{itemize}
\item $val_{R_2}(k^\sharp)\geqslant -c_1$ pour tout $k^\sharp\in C_{H'}^\sharp$
\item $val_{R_1}(h^{-1})\geqslant val_{R_1}(h)-c_2$ pour tout $h\in H(F)$
\item $val_{R_1}(hk^\sharp)\geqslant val_{R_1}(h)-c_3$ pour tout $h\in H(F)$ et $k^\sharp\in C_H^\sharp$
\item $val_{R_1}(hw_1)\geqslant -val_E(a_{H,1})-c_4$ si $h=k^\sharp a_Hh'$ avec $k^\sharp\in C_H^\sharp$, $a_H\in A_H(F)^+$ et $h'\in H'(F)$ et où on a noté $a_{H,1}$ la valeur propre de $a_H$ agissant sur $u'_1$.
\end{itemize}

On a alors

\[\begin{aligned}
-val_E(a_{H',2}) =val_{R_1}(a_{H'}^{-1}w_2)+val_E(2) & \geqslant val_{R_2}((k_2^\sharp)^{-1}a_{H'}^{-1}w_2)-c_1 \\
 & =val_{R_1}((k_2^\sharp)^{-1}a_{H'}^{-1}u'_1-w_1)-c_1 \\
 & =val_{R_1}(h^{-1}k_1^\sharp a_Hu'_1-w_1)-c_1 \\
 & \geqslant inf\left(0,val_{R_1}(h^{-1}k_1^\sharp a_Hu'_1)\right)-c_1
\end{aligned}\]

et

\[\begin{aligned}
val_{R_1}(h^{-1}k_1^\sharp a_Hu'_1)& =val_E(a_{H,1})+val_{R_1}(h^{-1}k_1^\sharp u'_1) \\
& \geqslant val_E(a_{H,1})+val_{R_1}(h)-c_2-c_3 \\
& \geqslant val_E(a_{H,1})+val_{R_1}(hw_1)-val_E(\lambda)-c_0-c_2-c_3 \\
& \geqslant -val_E(\lambda)-c_0-c_2-c_3-c_4
\end{aligned}\]

On en déduit que $-val_E(a_{H',2})\geqslant -val_E(\lambda)-c_0-c_1-c_2-c_3-c_4$ ce qu'il nous fallait. $\blacksquare$

\vspace{4mm}

On aura aussi besoin du lemme suivant

\begin{lem}
Soient $C_0\subset G_0(F)$ et $C_H\subset H(F)$ deux compacts-ouverts vérifiant

$$G_0(F)=C_HA_H^+A_0^+C_0$$

On a les majorations suivantes
\begin{enumerate}
\item $\Xi^{G_0}(a_H a_0)<<\Xi^{G_0}(a_H) \Xi^{G_0}(a_0)$ pour tous $a_0\in A_0^+,a_H\in A_H^+$
\item $mes(C_Ha_Ha_0C_0)<< \delta_{P_H}(a_H)^{-1}\delta_{P_0}(a_0)^{-1}$ pour tous $a_0\in A_0^+,a_H\in A_H^+$ 
\item $\sigma(a_0)<<\sigma(a_0a_H)$ pour tous $a_0\in A_0^+,a_H\in A_H^+$ 
\end{enumerate}
\end{lem}

\ul{Preuve}: 1)Soient $K_1\subset \overline{P}_0(F)$ et $K_2\subset P_H(F)$ deux sous-groupes ouverts vérifiant

\begin{itemize}
\item $K_1\subset a_0K_1a_0^{-1}$ pour tout $a_0\in A_0^+$;
\item $K_2\subset a_H^{-1} K_2a_H$ pour tout $a_H\in A_H^+$;
\item $\Xi^{G_0}$ est biinvariant par $K_1$ et $K_2$.
\end{itemize}

On a alors pour tous $a_0\in A_0^+,a_H\in A_H^+$ et pour tout $(k_1,k_2)\in K_1\times K_2$

$$\Xi^{G_0}(a_Ha_0)=\Xi^{G_0}(a_Hk_2k_1a_0)$$
D'après le lemme 11.0.1, $K_2K_1$ contient un sous-groupe ouvert $K_3$ de $G_0(F)$. On a alors
$$\Xi^{G_0}(a_Ha_0)=mes(K_3)^{-1}\displaystyle\int_{K_3} \Xi^{G_0}(a_Hk_3a_0)dk_3$$

Or, on a $\displaystyle\int_{K_3} \Xi^{G_0}(gk_3g')dk_3<<\Xi^{G_0}(g)\Xi^{G_0}(g')$ pour tous $g,g'\in G_0(F)$. On obtient ainsi la première partie du lemme. \\

2) Soient $K_1\subset G_0(F)$ et $K_2\subset H(F)$ deux sous-groupes compacts-ouverts tels que 
$$K_2K_1\subset a_H^{-1}C_Ha_Ha_0C_0a_0^{-1}$$
pour tous $a_0\in A_0^+,a_H\in A_H^+$ (de tels sous-groupes existent d'après le lemme 11.0.1). Soient $a_0\in A_0^+, a_H\in A_H^+$, on a alors $mes\left(C_Ha_Ha_0C_0\right)=mes\left((C_Ha_HK_2)(K_1a_0C_0)\right)$. L'ensemble $C_Ha_HK_2$ est union de $mes(C_Ha_HK_2)mes(K_2)^{-1}$ classes à gauche pour $K_2$ et l'ensemble $K_1a_0C_0$ est union de $mes(K_1a_0C_0)mes(K_1)^{-1}$ classes à droite pour $K_1$. On en déduit que

$$mes\left((C_Ha_HK_2)(K_1a_0C_0)\right)\leqslant mes(C_Ha_HK_2)mes(K_1a_0C_0)mes(K_1)^{-1}mes(K_2)^{-1}mes(K_2K_1)$$
Puisque $mes(C_Ha_HK_2)<<\delta_{P_W}(a_H)^{-1}$ et $mes(K_1a_0C_0)<<\delta_{P_0}(a_0)^{-1}$ pour tous $a_H\in A_H^+,a_0\in A_0^+$ on a bien 2) \\

3) Il existe une constante $c_1$ telle que $val_{R_0}(g)\leqslant val_{R_0}(g^{-1}v_0)+c_1$ pour tout $g\in G_0(F)$. On a donc $val_{R_0}(a_Ha_0)\leqslant val_{R_0}(a_0^{-1} v_0)+c_1$ pour tous $a_0\in A_0^+,a_H\in A_H^+$. Pour tout $a_0\in A_0(F)$, notons $a_{0,1}$ la valeur propre de $a_0$ agissant sur $u_1$. On a alors $val_{R_0}(a_0^{-1} v_0)=-val_E(a_{0,1})+val_E(2)$ et $val_{R_0}(a_0)=-val_E(a_{0,1})$ pour tout $a_0\in A_0^+$. Le résultat découle alors de ce que

$$sup\left(1,-val_{R_0}(g)\right)<<\sigma(g)<<sup\left(1,-val_{R_0}(g)\right)$$

\noindent pour tout $g\in G(F)$ $\blacksquare$

\vspace{4mm}

\section{Majorations d'intégrales cas r=0}

\begin{lem}
Il existe $\epsilon>0$ tel que $\Xi^{G_0}(h)<<exp(-\epsilon \sigma(h)) \Xi^H(h)$ pour tout $h\in H(F)$.
\end{lem}

\ul{Preuve}: Soient $(z_{\pm i})_{1\leqslant i\leqslant r_0}$ et $(z_{\pm i})_{1\leqslant i\leqslant r_1}$ des familles hyperboliques maximales de $V_0$ et $W$ respectivement. On a $r_0\geqslant r_1$. Soit $P_{G_0}$ le sous-groupe parabolique minimal de $G_0$ qui conserve le drapeau

$$Ez_1\subset\ldots\subset Ez_1\oplus\ldots\oplus Ez_{r_0}$$

\noindent Soit $A_{0,G_0}$ le tore des éléments qui stabilisent les droites $Ez_i$ pour $i=\pm 1,\ldots,\pm r_0$ et qui agissent trivialement sur l'orthogonal de l'espace engendré par ces droites. Soit $A_{0,G_0}^+$ l'ensemble des éléments de $A_{0,G_0}(F)$ qui contractent $P_{G_0}(F)$. Posons $A_{0,H}^+=A_{0,G_0}^+\cap H(F)$. D'après Bruhat-Tits, il existe un sous-ensemble compact $\Gamma\subset H(F)$ tel que $H(F)=\Gamma A_{0,H}^+\Gamma$. On peut donc se contenter de montrer le lemme pour $h=a\in A_{0,H}^+$. Pour un tel $a$, on note $a_i$ la valeur propre de $a$ agissant sur $z_i$ pour $i=\pm 1,\ldots,\pm r_1$. D'après le lemme II.1.1 de [W3] il existe un réel $d$ tel que \\
\noindent $\Xi^{G_0}(a)<<\delta_{P_{G_0}}(a)^{1/2}\sigma(a)^d$ pour tout $a\in A_{0,G_0}^+$. Un simple calcul montre que l'on a

$$\delta_{P_{G_0}}(a)=|a_1|_E^{d_V-1}|a_2|_E^{d_V-3}\ldots|a_{r_1}|_E^{d_V+1-2r_1}$$
et
$$\delta_{P_{0,H}}(a)=|a_1|_E^{d_W-1}|a_2|_E^{d_W-3}\ldots|a_{r_1}|_E^{d_W+1-2r_1}$$

\noindent On a donc $\delta_{P_{G_0}}(a)=\delta_{P_{0,H}}(a)(|a_1|_E\ldots|a_{r_1}|_E)$. On en déduit qu'il existe $\epsilon>0$ tel que

$$\Xi^{G_0}(a)^{1/2}<<\delta_{P_{0,H}}(a)^{1/2}exp(-\epsilon \sigma(a))$$

\noindent pour tout $a\in A_{0,H}^+$. D'après le lemme II.1.1 de [W3] on a aussi $\delta_{P_{0,H}}(a)^{1/2}<<\Xi^H(a)$ pour tout $a\in A_{0,H}^+$. Le lemme en découle $\blacksquare$

\vspace{4mm}

On reprend les notations de la section 11: $P_0,P_H,A_0,A_H,A_0^+,A_H^+$. Soient $C_0\subset G_0(F)$ et $C_H\subset H(F)$ des sous-ensembles compacts-ouverts tels que

$$G_0(F)=C_HA_H^+A_0^+C_0$$

\noindent (une telle décomposition existe d'après la proposition 11.0.1). Soient $A_0^1$ et $A_H^1$ les sous-groupes compacts maximaux de $A_0(F)$ et $A_H(F)$ respectivement. On peut toujours supposer que $A_0^1C_0=C_0$ et $C_HA_H^1=C_H$. Posons $\Lambda_0^+=A_0^+/A_0^1$ et $\Lambda_H^+=A_H^+/A_H^1$. On identifie $\Lambda_0^+$ et $\Lambda_H^+$ à des sous-ensembles de $A_0(F)$ et $A_H(F)$ via des sections. On a alors

$$G_0(F)=C_H\Lambda_H^+\Lambda_0^+C_0$$

\noindent Pour toute fonction $f$ mesurable à valeurs positives sur $G_0(F)$ on a

$$\mbox{(1)}\;\;\; \displaystyle\int_{G_0(F)} f(g)dg\leqslant \sum_{a_H\in\Lambda_H^+, a_0\in \Lambda_0^+} \int_{C_Ha_Ha_0C_0} f(g) dg$$

\begin{lem}
Soient $D$ un réel, $N\geqslant 1$ un entier et $g_0\in G_0(F)$ alors l'intégrale

$$\displaystyle\chi(g_0,N,D)=\int_{G_0(F)} \Xi^{G_0}(g_0^{-1}g) \Xi^{G_0}(g) \kappa_N(g) \sigma(g)^D dg$$

\noindent est convergente. De plus, à $D$ fixé, il existe un réel $R$ tel que

$$\displaystyle \chi(g_0,N,D)<<\Xi^{G_0}(g_0)\sigma(g_0)^RN^R$$
pour tous $g_0\in G_0(F)$ et $N\geqslant 1$.
\end{lem}

\ul{Preuve}:

On fixe le réel $D$. Introduisons un paramètre auxiliaire $\delta>0$ que l'on précisera plus tard. On a 

$$\chi(g_0,N,D)=\chi_{\leqslant \delta N}(g_0,N,D)+\chi_{>\delta N}(g_0,N,D)$$

\noindent où on a posé

$$\displaystyle \chi_{\leqslant \delta N}(g_0,N,D)=\int_{G_0(F)}\mathbf{1}_{\sigma\leqslant \delta N}(g) \Xi^{G_0}(g_0^{-1}g) \Xi^{G_0}(g) \kappa_N(g) \sigma(g)^D dg$$

et

$$\displaystyle \chi_{> \delta N}(g_0,N,D)=\int_{G_0(F)}\mathbf{1}_{\sigma> \delta N}(g) \Xi^{G_0}(g_0^{-1}g) \Xi^{G_0}(g) \kappa_N(g) \sigma(g)^D dg$$

Pour tout réel $D_1$ on a

$$\displaystyle\chi_{\leqslant \delta N}(g_0,N,D)\leqslant(\delta N)^{D_1} \int_{G_0(F)} \Xi^{G_0}(g_0^{-1}g)\Xi^{G_0}(g) \sigma(g)^{D-D_1} dg$$

Si $D_1$ est assez grand, l'intégrale précédente est convergente. Fixons un tel $D_1$. En introduisant une intégrale intérieure sur un sous-groupe compact-ouvert de $G_0(F)$, on voit qu'elle est essentiellement majorée par $\Xi^{G_0}(g_0)$ pour tout $g_0\in G_0(F)$. On a donc une majoration

$$\mbox{(2)}\;\;\; \chi_{\leqslant \delta N} (g_0,N,D)<<(\delta N)^{D_1} \Xi^{G_0}(g_0)$$

\noindent pour tous $g_0\in G_0(F)$, $N\geqslant 1$, $\delta>0$. Majorons à présent $\chi_{>\delta N}(g_0,N,D)$. Il existe un réel $\alpha>0$ tel que l'on ait la majoration $\Xi^{G_0}(g'^{-1}g)<<exp(\alpha \sigma(g')) \Xi^{G_0}(g)$ pour tous $g,g'\in G_0(F)$. On a par conséquent

$$\displaystyle\chi_{> \delta N}(g_0,N,D)<<exp(\alpha\sigma(g_0))\int_{G_0(F)} \mathbf{1}_{\sigma > \delta N}(g) \Xi^{G_0}(g)^2 \kappa_N(g) \sigma(g)^D dg$$

\noindent pour tous $g_0\in G_0(F)$, $N\geqslant 1$, $\delta>0$. Il existe un réel $c_1>0$ telle que pour tous $k_H\in C_H,k_0\in C_0,a_H\in A_H^+, a_0\in A_0^+$ on ait $\sigma(k_Ha_Ha_0k_0)\leqslant c_1+\sigma(a_H)+\sigma(a_0)$. Il existe aussi une constante $c_2>0$ telle que pour tous $k_H\in C_H,k_0\in C_0,a_H\in A_H^+, a_0\in A_0^+$ et pour tout $N\geqslant 1$ on ait $\kappa_N(k_Ha_Ha_0k_0)=1\Rightarrow val_E(a_{0,1})\leqslant c_2 N$ où $a_{0,1}$ désigne la valeur propre de $a_0$ agissant sur $u_1$. Puisque $\sigma(a_0)<<sup(1,val_E(a_{0,1}))$ pour tout $a_0\in A_0^+$, quitte à agrandir la constante $c_2$ on peut remplacer $val_E(a_{0,1})$ par $\sigma(a_0)$ dans l'inégalité précédente. D'après le lemme 11.0.4 et l'inégalité (1), on a donc

\[\begin{aligned}
\displaystyle \int_{G_0(F)} \mathbf{1}_{\sigma > \delta N}(g) \Xi^{G_0}(g)^2 \kappa_N(g) \sigma(g)^D dg << & \left (\sum_{a_0\in \Lambda_0^+, \sigma(a_0)\leqslant c_2N} \Xi^{G_0}(a_0)^2\delta_{P_0}(a_0)^{-1} \sigma(a_0)^D\right) \times\\
 & \left(\sum_{a_H\in \Lambda_H^+} \mathbf{1}_{\sigma> (\delta-c_2) N-c_1}(a_H)\Xi^{G_0}(a_H)^2 \delta_{P_H}(a_H)^{-1} \sigma(a_H)^D\right)
\end{aligned}\]

Il existe des réels $D_2$ et $\epsilon>0$ tels que l'on ait les majorations $\Xi^{G_0}(a_0)^2<<\delta_{P_0}(a_0) \sigma(a_0)^{D_2}$, $\Xi^{H}(a_H)^2<<\delta_{P_H}(a_H) \sigma(a_H)^{D_2}$ et $\Xi^{G_0}(h)<<exp(-\epsilon \sigma(h)) \Xi^H(h)$ pour tous $a_0\in A_0^+$, $a_H\in A_H^+$, $h\in H(F)$. Le terme $\chi_{> \delta N}(g_0,N,D)$ est donc essentiellement majoré par le produit de $exp(\alpha \sigma(g_0))$ et des deux sommes

$$\displaystyle\sum_{a_0\in\Lambda_0^+, \sigma(a_0)\leqslant c_2 N} \sigma(a_0)^{D+D_2}$$

et

$$\displaystyle\sum_{a_H\in\Lambda_H^+,\sigma(a_H)>(\delta-c_2)N-c_1} exp(-\epsilon\sigma(a_H)) \sigma(a_H)^{D+D_2}$$

La première somme est essentiellement majorée par une puissance de $N$ et la deuxième par $exp(-\epsilon'(\delta-c_2)N)$ pour un certain $\epsilon'>0$ qui ne dépend que de $\epsilon$. Rassemblant cette majoration avec la majoration (2), on obtient

$$\chi(g_0,N,D)<<(\delta N)^{D_1}\Xi^{G_0}(g_0)+exp(\alpha\sigma(g_0))N^{D_3}exp(-\epsilon'(\delta-c_2)N)$$

\noindent pour tous $g_0\in G_0(F)$, $N\geqslant 1$, $\delta>0$ (et pour un certain $D_3>0$). On vérifie qu'il existe un réel $\beta>0$ tel que $exp(-\beta \sigma(g))<<\Xi^{G_0}(g)$ pour tout $g\in G_0(F)$. Il suffit alors de prendre $\displaystyle\delta=\frac{(\alpha+\beta)\sigma(g_0)}{\epsilon'N}+c_2$ pour obtenir la majoration du lemme $\blacksquare$

\vspace{4mm}

\begin{lem}
Pour tout réel $D$, l'intégrale
$$\displaystyle\int_{H(F)} \Xi^H(h)\Xi^{G_0}(h) \sigma(h)^D dh$$
est convergente.
\end{lem}

\ul{Preuve}: D'après le lemme 12.0.5 pour tout $d>0$ on a une majoration
$$\Xi^{G_0}(h)<<\Xi^H(h)\sigma(h)^{-d}$$

\noindent pour tout $h\in H(F)$. Et d'après le lemme II.1.5 de [W3], pour $d$ assez grand l'intégrale $\displaystyle\int_{H(F)} \Xi^H(h)^2 \sigma(h)^{D-d} dh$ converge. $\blacksquare$

\vspace{4mm}

\begin{lem}
Pour tout réel $D$, l'intégrale
$$\displaystyle\int_{H(F)}\int_{H(F)} \Xi^{G_0}(h)\Xi^H(h'h)\Xi^{G_0}(h') \sigma(h)^D\sigma(h')^D dhdh'$$
est convergente
\end{lem}

\ul{Preuve}: Soit $K_H$ un sous-groupe compact-ouvert de $H(F)$. On a des majorations \\
\noindent $\displaystyle\int_{K_H} \Xi^H(h'kh)dk<<\Xi^H(h')\Xi^H(h)$, $\sigma(k_Hh)<<\sigma(h)$ et $\Xi^{G_0}(k_Hh)<<\Xi^{G_0}(h)$ pour tous $h,h'\in H(F)$ et pour tout $k_H\in K_H$. Par conséquent, en introduisant une intégrale intérieure sur $K_H$, l'intégrale de l'énoncé est majorée par une constante fois le carré de l'intégrale suivante
$$\displaystyle\int_{H(F)} \Xi^{G_0}(h)\Xi^H(h) \sigma(h)^D dh$$
qui est convergente d'après le lemme précédent. $\blacksquare$

\vspace{4mm}

Soient $D$ et $b>0$ deux réels et $N\geqslant 1$ un entier. Posons

$$I^1(N,D)=\displaystyle\int_{G_0(F)}\int_{H(F)} \Xi^H(h)\Xi^{G_0}(g)\Xi^{G_0}(h^{-1}g) \kappa_N(g) \sigma(h)^D\sigma(g)^D dh dg$$

et

$$I^1(N,D,b)=\displaystyle\int_{G_0(F)}\int_{H(F)} \mathbf{1}_{\sigma\geqslant b}(h)\Xi^H(h)\Xi^{G_0}(g)\Xi^{G_0}(h^{-1}g) \kappa_N(g) \sigma(h)^D\sigma(g)^D dh dg$$

\begin{lem}
Ces deux expressions sont convergentes. A $D$ fixé, il existe des réels $R$ et $\epsilon>0$ tels que
\begin{enumerate}
\item $I^1(N,D)<<N^R$ pour tout $N\geqslant 1$
\item $I^1(N,D,b)<<N^R exp(-\epsilon b)$ pour tout $N\geqslant 1$ et pour tout $b>0$
\end{enumerate}
\end{lem}

\ul{Preuve}: Fixons un réel $D$. D'après les lemmes 12.0.5 et 12.0.6, il existe un réel $R>0$ et $\epsilon_1>0$ tels que

$$I^1(N,D)<<N^R\displaystyle\int_{H(F)} \Xi^H(h)^2 exp(-\epsilon_1 \sigma(h))\sigma(h)^D dh$$
et
$$I^1(N,D,b)<<N^R\displaystyle\int_{H(F)} \mathbf{1}_{\sigma\geqslant b}(h)\Xi^H(h)^2exp(-\epsilon_1 \sigma(h))\sigma(h)^D dh$$

\noindent pour tout $N\geqslant 1$ et pour tout $b>0$. Pour tout $\epsilon'>0$, l'intégrale $\displaystyle\int_{H(F)} \Xi^H(h)^2 exp(-\epsilon' \sigma(h))\sigma(h)^D dh$ est convergente. On en déduit la première majoration de l'énoncé. Pour $h\in H(F)$ tel que $\sigma(h)\geqslant b$ on a $exp(-\epsilon_1 \sigma(h))\leqslant exp(-\epsilon_1b/2)exp(-\epsilon_1 \sigma(h)/2)$. On en déduit la deuxième majoration de l'énoncé avec $\epsilon=\epsilon_1/2$ $\blacksquare$

\section{Majorations d'intégrales cas général}

Pour $N\geqslant 1$ un entier et $D$ un réel, on pose
$$I(N,D)=\displaystyle\int_{G(F)} \Xi^G(g)^2 \kappa_N(g) \sigma(g)^Ddg$$

\begin{prop}
Cette intégrale est convergente et le réel $D$ étant fixé, il existe un réel $R$ tel que

$$I(N,D)<<N^R$$
pour tout $N\geqslant 1$.
\end{prop}

\ul{Preuve}: Fixons le réel $D$. Les fonctions $\Xi^G,\kappa_N$ et $\sigma$ étant invariantes à droite par $K$ et la fonction $\kappa_N$ étant invariante à gauche par $U(F)$, on a

$$I(N,D)=\displaystyle\int_{M(F)}\kappa_N(m) \int_{U(F)} \Xi^G(mu)^2\sigma(mu)^D du dm$$

D'après la proposition 2.1.1, il existe un réel $D'$ tel que $I(N,D)$ soit essentiellement majorée par

$$\displaystyle\int_{M(F)} \Xi^M(m)^2 \kappa_N(m)\sigma(m)^{D'} dm$$

\noindent pour tout $N\geqslant 1$. On a $M(F)=A(F)\times G_0(F)$ et

\begin{center}
 $\kappa_N(m)=\kappa_N(a)\kappa_N(g_0)$, $\sigma(m)<<\sigma(a)\sigma(g_0)$, $\Xi^M(m)=\Xi^{G_0}(g_0)$
\end{center} 
 
\noindent  pour tout $m=ag_0\in M(F)=A(F)G_0(F)$. L'intégrale précédente est donc essentiellement majorée par

$$\displaystyle\int_{A(F)}\kappa_N(a)\sigma(a)^{D'}da \int_{G_0(F)} \Xi^{G_0}(g_0)^2\kappa_N(g_0)\sigma(g_0)^{D'}dg_0$$

\noindent pour tout $N\geqslant 1$. D'après le lemme 12.0.6, l'intégrale sur $G_0(F)$ est essentiellement bornée par une puissance de $N$ et il n'est pas difficile de vérifier qu'il en va de même de l'intégrale sur $A(F)$. $\blacksquare$

\vspace{4mm}

Soit $c\geqslant 0$ un nombre réel. On définit $U(F)_c$ comme l'ensemble des éléments $u\in U(F)$ qui vérifient $val_F(Tr(h(uv_i,v_{-i-1})))\geqslant -c$ pour $i=0,\ldots,r-1$. Pour $c\geqslant 0$ un réel et $n\in \mathbb{N}$, on pose

$$A_{n,c}=\{u\in U(F)_c; \; q^{-n-1}<\Xi^M(m_{\overline{P}}(u))\delta_{\overline{P}}(m_{\overline{P}}(u))^{1/2}\leqslant q^{-n}\}$$

\begin{lem}
Il existe des réels $\epsilon>0$ et $\alpha>0$ tels que l'on ait la majoration
$$mes(A_{n,c})<<q^{n(1-\epsilon)+\alpha c}$$
pour tout $c\geqslant 0$ et tout $n\in \mathbb{N}$
\end{lem}

\ul{Preuve}: Si $V_0$ est anisotrope, c'est le corollaire 2.4.1 appliqué à $\overline{Q}=P$, à un certain sous-groupe parabolique minimal $P_{min}$ de $G$ tel que $P$ soit antistandard et à certaines formes linéaires $l_{\alpha}$, $\alpha\in \Delta$. Si $V_0$ n'est pas anisotrope, les résultats de la section 2.4 ne peuvent pas s'appliquer directement. On peut alors trouver deux vecteurs isotropes $e_1,e_{-1}\in V_0$ vérifiant $v_0=e_1+e_{-1}$ et $h(e_1,e_{-1})=\nu_0/2$. Fixons deux tels éléments $e_1,e_{-1}$. Pour $c\geqslant 0$ on définit $U(F)'_c$ comme l'ensemble des éléments $u\in U(F)$ qui vérifient $val_F(Tr(h(uv_i,v_{-i-1})))\geqslant -c$ pour $i=1,\ldots,r-1$ et $val_F(Tr(h(ue_1,v_{-1})))\geqslant -c$. Alors les résultats du paragraphe 2.4 s'appliquent à $U(F)'_c$. En particulier, il existe des nombres réels $\epsilon_1>0$ et $\alpha_1>0$ de sorte que l'on ait une majoration

$$mes\{u\in U(F)'_c; \; q^{-n-1}<\Xi^M(m_{\overline{P}}(u))\delta_{\overline{P}}(m_{\overline{P}}(u))^{1/2}\leqslant q^{-n}\}<< q^{n(1-\epsilon_1)+\alpha_1 c}$$

\noindent pour tout $c\geqslant 0$ et pour tout $n\in \mathbb{N}$. Soit $\delta>0$ un réel que l'on précisera plus tard. On a pour tout $n\in\mathbb{N}$ et tout $c\geqslant 0$

$$A_{n,c}=(A_{n,c}\cap U(F)'_{\delta n}) \sqcup (A_{n,c}\backslash U(F)'_{\delta n})$$

D'après ce qui précède, la mesure de $A_{n,c}\cap U(F)'_{\delta n}$ est essentiellement majorée par $q^{n(1-\epsilon_1+\alpha_1 \delta)}$. Soit $B_n$ le sous-ensemble de $U(F)$ des éléments $u$ qui vérifient $q^{-n-1}<\Xi^M(m_{\overline{P}}(u))\delta_{\overline{P}}(m_{\overline{P}}(u))^{1/2}\leqslant q^{-n}$. On a alors

$$\displaystyle mes(A_{n,c}\backslash U(F)'_{\delta n})=\int_{U(F)} \mathbf{1}_{B_n}(u) \mathbf{1}_{U(F)_c\backslash U(F)'_{\delta n}} (u) du$$

Pour $\lambda \in \mathcal{O}_F^\times$, soit $a(\lambda)$ l'élément de $G_0(F)$ qui envoie $e_1$ sur $\lambda e_1$, $e_{-1}$ sur $\lambda^{-1} e_{-1}$ et qui agit trivialement sur l'orthogonal de $Ee_1\oplus Ee_{-1}$ dans $V_0$. La conjugaison par $a(\lambda)$ pour $\lambda\in\mathcal{O}_F^\times$ laisse stable $U(F)$ et ne change pas la mesure. Par conséquent on a

$$\displaystyle mes(A_{n,c}\backslash U(F)'_{\delta n})=\int_{\mathcal{O}_F^\times}\int_{U(F)} \mathbf{1}_{B_n}(a(\lambda)ua(\lambda)^{-1}) \mathbf{1}_{U(F)_c\backslash U(F)'_{\delta n}} (a(\lambda)ua(\lambda)^{-1}) du d\lambda$$

\noindent où $d\lambda$ est la mesure de Haar sur $\mathcal{O}_F^\times$ qui donne une mesure totale de 1. On ne perd rien à supposer que les $a(\lambda)$ pour $\lambda\in\mathcal{O}_F^\times$ sont dans le sous-groupe compact spécial de $G(F)$ qui permet de définir les éléments $m_{\overline{P}}$ et qu'ils laissent stable $\Xi^M$ par translation à gauche. On a alors $\mathbf{1}_{B_n}(a(\lambda)ua(\lambda)^{-1})=\mathbf{1}_{B_n}(u)$ pour tous $u\in U(F),\lambda\in\mathcal{O}_F^\times$. On en déduit que

$$\displaystyle mes(A_{n,c}\backslash U(F)'_{\delta n})=\int_{U(F)} \mathbf{1}_{B_n}(u)\int_{\mathcal{O}_F^\times} \mathbf{1}_{U(F)_c\backslash U(F)'_{\delta n}} (a(\lambda)ua(\lambda)^{-1}) d\lambda du$$

Soit $u\in U(F)$ et posons $u_1=Tr(h(ue_1,v_{-1}))$ et $u_{-1}=Tr(h(ue_{-1},v_{-1}))$. Si l'intégrale intérieure est non nulle, on a $val_F(u_1)<-\delta n$. Supposons que ce soit le cas. On a alors l'égalité

$$\displaystyle\int_{\mathcal{O}_F^\times} \mathbf{1}_{U(F)_c\backslash U(F)'_{\delta n}} (a(\lambda)ua(\lambda)^{-1}) d\lambda=\int_{\mathcal{O}_F^\times} \mathbf{1}_{val_F\geqslant -c}(\lambda^{-1}u_1+\lambda u_{-1}) d\lambda$$

On vérifie aisément que l'intégrale précédente est essentiellement majorée par $q^{val_F(u_1)+c}$ pour tous $c\geqslant 0$, $u_1,u_{-1}\in F^\times$. On en déduit que $mes(A_{n,c}\backslash U(F)'_{\delta n})$ est essentiellement majorée par $mes(B_n) q^{-\delta n+c}$. D'après les lemmes II.3.4 et II.4.3 de [W3] alliés au fait que la fonction $\Xi^M$ est bornée par $1$, il existe un entier $d$ tel que l'on ait une majoration $mes(B_n)<<n^d q^n$ pour tout $n\in \mathbb{N}$. Pour $\delta$ qui vérifie $0<\delta<\epsilon_1/\alpha_1$ on obtient le résultat du lemme en sommant les majorations obtenues pour $mes(A_{n,c}\cap U(F)'_{\delta n})$ et $mes(A_{n,c}\backslash U(F)'_{\delta n})$ $\blacksquare$

\vspace{2mm}

\begin{lem}
Il existe un réel $\epsilon_1>0$ tel que pour tous réels $\epsilon_0,c,C$ vérifiant $c\geqslant 0,C>0$ et $\epsilon_0<\epsilon_1$, l'intégrale

$$\displaystyle\int_{U(F)_c} \mathbf{1}_{\sigma\geqslant C}(u) \Xi^M(m_{\overline{P}}(u)) \delta_{\overline{P}}(m_{\overline{P}}(u))^{1/2} exp(\epsilon_0\sigma(u)) du$$

\noindent est convergente. De plus, il existe un réel $\alpha >0$ tel que cette intégrale soit essentiellement majorée par $exp(-(\epsilon_1-\epsilon_0) C+\alpha c)$ pour tout $c\geqslant 0$ pour tout $C>0$ et pour tout $\epsilon_0<\epsilon_1$.
\end{lem}

\ul{Preuve}: D'après le lemme II.3.4 de [W3] et puisque la fonction $\Xi^M$ est bornée, il existe deux réels $c_1,c_2>0$ tels que $\Xi^M(m_{\overline{P}}(u)) \delta_{\overline{P}}(u)^{1/2}\leqslant c_1 exp(-c_2\sigma(u))$ pour tout $u\in U(F)$. La convergence et la majoration de l'énoncé sont alors des conséquences faciles du lemme précédent: il suffit de découper l'intégrale en un somme d'intégrales sur les $A_{n,c}$ pour $n\geqslant 0$, de remarquer que les termes de la somme pour $n\leqslant c_2 C-log(c_1)-1$ sont nuls et de majorer $\sigma(u)$ par $\displaystyle\frac{n+1+log(c_1)}{c_2}$ pour $u\in A_{n,c}$ $\blacksquare$

\vspace{5mm}

Pour $D,c>0$ deux réels et $m\in M(F)$. Posons

$$\displaystyle X(c,D,m)=\int_{U(F)_c} \Xi^G(um) \sigma(u)^Ddu$$

\begin{prop}
L'intégrale précédente est toujours convergente et à $D$ fixé, il existe un réel $R$ tel que l'on ait la majoration

$$X(c,D,m)<<c^R \Xi^M(m)\delta_P(m)^{1/2} \sigma(m)^R$$
pour tout $m\in M(F)$ et pour tout $c\geqslant 1$.
\end{prop}

\ul{Preuve}: Fixons un réel $D$. Introduisons un paramètre auxiliaire $b>0$ que l'on précisera plus tard. On a $X(c,D,m)=X_{\leqslant b}(c,D,m)+X_{>b}(c,D,m)$ où

$$\displaystyle X_{\leqslant b}(c,D,m)=\int_{U(F)_c} \mathbf{1}_{\sigma\leqslant b}(u)\Xi^G(um) \sigma(u)^Ddu$$

et

$$\displaystyle X_{>b}(c,D,m)=\int_{U(F)_c} \mathbf{1}_{\sigma> b}(u)\Xi^G(um) \sigma(u)^Ddu$$

Pour tout réel $D'$ on a

$$\displaystyle X_{\leqslant b}(c,D,m)\leqslant b^{D'} \int_{U(F)}\Xi^G(um) \sigma(u)^{D-D'}du$$

D'après la proposition II.4.5 de [W3], si $D'$ est assez grand, l'intégrale précédente est essentiellement majorée par $\delta_P(m)^{1/2} \Xi^M(m)$ pour tout $m\in M(F)$. Fixons un tel $D'$. Il existe un réel $\beta>0$ tel que l'on ait la majoration $\Xi^G(gg')<<\Xi^G(g) exp(\beta \sigma(g'))$ pour tous $g,g'\in G(F)$. On en déduit la majoration

$$\displaystyle X_{>b}(c,D,m)<<exp(\beta \sigma(m))\int_{U(F)_c} \mathbf{1}_{\sigma> b}(u)\Xi^G(u) \sigma(u)^Ddu$$

\noindent pour tout $m\in M(F)$ pour tout $c>0$ et pour tout $b>0$. D'après [W3], lemmes II.1.1 et II.3.2, il existe un réel $D''$ tel que

$$\Xi^G(g)<<\Xi^M(m_{\overline{P}}(g))\delta_{\overline{P}}(m_{\overline{P}}(g))^{1/2} \sigma(g)^{D''}$$

\noindent pour tout $g\in G(F)$. On en déduit que

$$\displaystyle X_{>b}(c,D,m)<<exp(\beta \sigma(m))\int_{U(F)_c}\mathbf{1}_{\sigma> b}(u)\Xi^M(m_{\overline{P}}(u))\delta_{\overline{P}}(m_{\overline{P}}(u))^{1/2} \sigma(u)^{D+D''} du$$
pour tout $m\in M(F)$ pour tout $c>0$ et pour tout $b>0$. D'après le lemme précédent, il existe deux réels $\epsilon>0$ et $\alpha>0$ tels que l'intégrale précédente soit essentiellement majorée par $exp(-\epsilon b+\alpha c)$ pour tous $b,c>0$. Il existe un réel $\gamma>0$ tel que $exp(-\gamma \sigma(m))<<\Xi^M(m)\delta_P(m)^{1/2}$ pour tout $m\in M(F)$. Il suffit de prendre $\displaystyle b=\frac{(\beta+\gamma)\sigma(m)+\alpha c}{\epsilon}$ pour obtenir la majoration de l'énoncé $\blacksquare$

\vspace{5mm}

\begin{cor}
Pour tous réels $D$ et $c\geqslant 0$, les intégrales
\begin{enumerate}
\item $$\displaystyle \int_{H(F)U(F)_c} \Xi^H(h)\Xi^G(hu)\sigma(hu)^D dudh$$
\item $$\displaystyle \int_{H(F)U(F)_c}\int_{H(F)U(F)_c} \Xi^G(hu)\Xi^G(h'h)\Xi^G(h'u')\sigma(hu)^D\sigma(h'u')^D du'dh'dudh$$
\end{enumerate}
sont convergentes.
\end{cor}

\ul{Preuve}: On peut effectuer dans les deux intégrales le changement de variable $u\mapsto h^{-1}uh$. Puisque la conjugaison par $H(F)$ préserve $U(F)_c$, on obtient les mêmes intégrales où $hu$ a été changé en $uh$. Le point 1. est alors une conséquence de la proposition 13.0.3 et du lemme 12.0.7 et le point 2. une conséquence de la proposition 13.0.3 et du lemme 12.0.8 $\blacksquare$

\vspace{5mm}

\begin{prop}
Soient $c>0$ et $\epsilon>0$ deux réels. Il existe un réel $\epsilon_1>0$ tel que pour tout $\epsilon_0<\epsilon_1$ et pour tout $h\in H(F)$, l'intégrale

$$\displaystyle\int_{U(F)_c} exp(\epsilon_0 \sigma(u)) \Xi^G(hu) du$$
soit convergente et que de plus elle soit essentiellement majorée par $exp(\epsilon \sigma(h))\Xi^{G_0}(h)$ pour tout $\epsilon_0<\epsilon_1$ et tout $h\in H(F)$.

\end{prop}

\ul{Preuve}: Fixons $c>0$ et $\epsilon>0$ deux réels. Introduisons un paramètre $b>0$ que l'on précisera plus tard. Pour tout réel $\epsilon_0$ on a l'égalité

\[\begin{aligned}
\displaystyle\int_{U(F)_c} exp(\epsilon_0 \sigma(u)) \Xi^G(hu) du = & \int_{U(F)_c} \mathbf{1}_{\sigma\leqslant b}(u) exp(\epsilon_0 \sigma(u)) \Xi^G(hu) du \\
 & +\int_{U(F)_c} \mathbf{1}_{\sigma>b}(u) exp(\epsilon_0 \sigma(u)) \Xi^G(hu) du
\end{aligned}\]

Notons $I_{\leqslant b}(\epsilon_0,h)$ et $I_{>b}(\epsilon_0,h)$ la première et la deuxième intégrale qui apparaîssent dans la somme ci-dessus. Pour tout réel $D'$ on a une majoration

$$\displaystyle I_{\leqslant b}(\epsilon_0,h)\leqslant b^{D'} exp(\epsilon_0 b) \int_{U(F)} \Xi^G(hu) \sigma(u)^{-D'} du$$

D'après la proposition II.4.5 de [W3], pour $D'$ assez grand l'intégrale précécente est essentiellement majorée par $\Xi^{G_0}(h)$ pour tout $h\in H(F)$. On fixe un tel $D'$. Il existe un réel $\beta>0$ tel que $\Xi^G(gg')<<exp(\beta \sigma(g)) \Xi^G(g')$ pour tous $g,g'\in G(F)$. On en déduit la majoration

$$\displaystyle I_{>b}(\epsilon_0,h)<<exp(\beta \sigma(h)) \int_{U(F)_c} \mathbf{1}_{\sigma>b}(u) exp(\epsilon_0\sigma(u)) \Xi^G(u) du$$
pour tout $\epsilon_0$, pour tout $b>0$ et pour tout $h\in H(F)$. D'après [W3], lemmes II.1.1 et II.3.2, il existe un réel $D''$ tel que

$$\Xi^G(g)<<\Xi^M(m_{\overline{P}}(g))\delta_{\overline{P}}(m_{\overline{P}}(g))^{1/2} \sigma(g)^{D''}$$

\noindent pour tout $g\in G(F)$. On en déduit la majoration

$$\displaystyle I_{>b}(\epsilon_0,h)<<exp(\beta \sigma(h)) \int_{U(F)_c} \mathbf{1}_{\sigma>b}(u) \Xi^M(m_{\overline{P}}(u))\delta_{\overline{P}}(m_{\overline{P}}(u))^{1/2} exp(\epsilon_0 \sigma(u)) \sigma(u)^{D''} du$$

\noindent pour tout $\epsilon_0$, pour tout $b>0$ et pour tout $h\in H(F)$. D'après le lemme 13.0.11, il existe un réel $\epsilon_1'>0$ tel que pour tout réel $\epsilon_0<\epsilon_1'$, l'intégrale précédente sans le $\sigma(u)^{D''}$ soit convergente et essentiellement majorée par $exp(-(\epsilon_1'-\epsilon_0)b)$ pour tout $b>0$ et pour tout $\epsilon_0<\epsilon_1'$. Soit $\epsilon_1'>\epsilon_2>0$ un réel, on peut majorer le terme $\sigma(u)^{D''}$ par une constante multipliée par $exp(\epsilon_2 \sigma(u))$. On en déduit que l'intégrale précédente est convergente pour $\epsilon_0<\epsilon_1'-\epsilon_2$ et qu'elle est essentiellement majorée par $exp(-(\epsilon_1'-\epsilon_2-\epsilon_0)b)$ pour tout $b>0$ et tout $\epsilon_0<\epsilon_1'-\epsilon_2$. En particulier on a la majoration suivante

$$I_{> b}(\epsilon_0,h)<<exp(\beta \sigma(h)) exp(-(\epsilon_1'-\epsilon_2)b/2)$$

\noindent pour tout $b>0$, pour tout $h\in H(F)$ et pour tout $\epsilon_0<(\epsilon_1'-\epsilon_2)/2$. Il existe un réel $\gamma>0$ tel que $exp(-\gamma \sigma(h))<<\Xi^{G_0}(h)$ pour tout $h\in H(F)$. Pour $\displaystyle b<\frac{\epsilon\sigma(h)}{\epsilon_0}$, on a la majoration de l'énoncé pour le terme $I_{\leqslant b}(\epsilon_0,h)$. Pour $\displaystyle b>\frac{2(\beta+\gamma)}{\epsilon_1'-\epsilon_2} \sigma(h)$ et $\epsilon_0<(\epsilon_1'-\epsilon_2)/2$ on obtient la majoration de l'énoncé pour le terme $I_{>b}(\epsilon_0,h)$. On peut trouver un $b$ tel que les trois inégalités précédentes soient vérifiées si

$$\epsilon_0<min\left((\epsilon_1'-\epsilon_2)/2,\epsilon(\epsilon_1'-\epsilon_2)/(2(\beta+\gamma))\right)$$

Il suffit donc de prendre $\epsilon_1=min((\epsilon_1'-\epsilon_2)/2,\epsilon(\epsilon_1'-\epsilon_2)/(2(\beta+\gamma)))$ pour obtenir le résultat de l'énoncé $\blacksquare$

\vspace{5mm}

Soient $D$ un réel, $c$ et $c'$ deux entiers tels que $c'\geqslant c>0$ et $m,m'\in M(F)$. Posons

$$X(c,D,m,m')=\displaystyle\int_{U(F)}\int_{U(F)_c} \Xi^G(um)\Xi^G(vum')\sigma(v)^D\sigma(u)^Ddvdu$$

$$X(c,c',D,m,m')=\displaystyle\int_{U(F)-U(F)_{c'}}\int_{U(F)_c}\Xi^G(um)\Xi^G(vum')\sigma(v)^D\sigma(u)^Ddvdu$$

\begin{lem}
Les intégrales ci-dessus sont convergentes et à $D$ fixé, il existe un réel $R$ et un réel $\epsilon>0$ tels qu'on ait les majorations

\begin{enumerate}
\item

$$X(c,D,m,m')<<c^R\sigma(m)^R\sigma(m')^R\delta_{P}(m)^{1/2}\delta_{P}(m')^{1/2}\Xi^M(m)\Xi^M(m')$$

\noindent pour tous $m,m'\in M(F)$ et tout $c\geqslant 1$.
\item

$$X(c,c',D,m,m')<<exp(-\epsilon c')\sigma(m)^R\sigma(m')^R\delta_{P}(m)^{1/2}\delta_{P}(m')^{1/2}\Xi^M(m)\Xi^M(m')$$

\noindent pour tous $m,m'\in M(F)$ et tout $c'\geqslant c\geqslant 1$.
\end{enumerate}

\end{lem}

\ul{Preuve}: Pour $c\in\mathbb{N}$ et $x\in F$ on pose $val_c(x)=min(0,val_F(\lambda)+c)$. Pour $u\in U(F)/U(F)_c$ on définit $val_c(u)=\displaystyle\sum_{i=0}^{r-1} val_c(Tr(h(uv_i,v_{-i-1})))$. Montrons que:

\vspace{2mm}

(1) Il existe un réel $D_1$ tel que l'on ait la majoration
$$\displaystyle\int_{U(F)_c}\Xi^G(vum)\sigma(vum)^Ddv<<(c-val_c(u))^{D_1}q_F^{val_c(u)}\sigma(m)^{D_1}\delta_{P}(m)^{1/2}\Xi^M(m)$$
pour tout $m\in M(F)$, tout $c\geqslant 1$ et tout $u\in U(F)$.

\vspace{2mm}

Pour $a\in A(F)\cap K$ on peut remplacer $\Xi^G(vum)$ et $\sigma(vum)$ par $\Xi^G(avuma^{-1})$ et $\sigma(avuma^{-1})$ et on peut alors intégrer sur $A(F)\cap K$. Puisque $a$ commute à $m$ et normalise $U(F)_c$, on obtient
$$\displaystyle\int_{U(F)_c}\Xi^G(vum)\sigma(vum)^Ddv<<\int_{A(F)\cap K}\int_{U(F)_c} \Xi^G(vaua^{-1}m)\sigma(vaua^{-1}m)^D dvda$$
L'application

$$A(F)\cap K\to U(F)/U(F)_c$$
$$a\mapsto aua^{-1}$$

\noindent a son image incluse dans $U(F)_{c-val_c(u)}/U(F)_c$ et son jacobien est borné par $q_F^{val_c(u)}$. On en déduit que

$$\displaystyle\int_{U(F)_c}\Xi^G(vum)\sigma(vum)^Ddv<<q_F^{val_c(u)}\int_{U(F)_{c-val_c(u)}} \Xi^G(vm)\sigma(vm)^Ddv$$

\noindent La proposition 13.0.3 permet alors d'obtenir la majoration (1). \\

Dans les expressions définissant $X(c,D,m,m')$ et $X(c',c,D,m,m')$ on peut écrire l'intégrale sur $U(F)$, respectivement sur $U(F)-U(F)_{c'}$, comme une double intégrale sur $U(F)_c$ et sur $U(F)/U(F)_c$, respectivement comme une double intégrale sur $U(F)_c$ et sur $(U(F)-U(F)_{c'})/U(F)_c$. D'après (1), on a les majorations

\[\begin{aligned}
X(c,D,m,m')<< & \sigma(m)^R\sigma(m')^R\delta_{P}(m)^{1/2}\delta_{P}(m')^{1/2}\Xi^M(m)\Xi^M(m') \\
 & \displaystyle\int_{U(F)/U(F)_c} q_F^{2val_c(u)}(c-val_c(u))^{D_1}du
\end{aligned}\]

et

\[\begin{aligned}
X(c',c,D,m,m') & <<\sigma(m)^R\sigma(m')^R\delta_{P}(m)^{1/2}\delta_{P}(m')^{1/2}\Xi^M(m)\Xi^M(m') \\
 & \displaystyle\int_{U(F)-U(F)_{c'}/U(F)_c} q_F^{2val_c(u)}(c-val_c(u))^{D_1}du
\end{aligned}\]

\noindent pour tous $m,m'\in M(F)$ et pour tous $c'\geqslant c\geqslant 1$. Les deux majorations de l'énoncé découlent alors des majorations suivantes

$$\displaystyle\int_{F/\mathfrak{p}_F^{-c}} q_F^{2val_c(x)}(c-val_c(x))^{D_1}dx<<c^{D_2}$$

et

$$\displaystyle\int_{(F-\mathfrak{p}_F^{-c'})/\mathfrak{p}_F^{-c}}q_F^{2val_c(x)}(c-val_c(x))^{D_1}dx<<exp(-\epsilon c')$$

\noindent qui sont vérifiées pour certains réels $D_2$ et $\epsilon>0$ $\blacksquare$

\vspace{5mm}

Soient $D,C>0,c>0$ trois réels et $N\geqslant 1$ un entier. Posons

$$\displaystyle\chi(c,C,N,D)=\int_{M(F)}\int_{U(F)_c} \mathbf{1}_{\sigma\geqslant C}(u) \Xi^M(m) \Xi^G(um) \kappa_N(m) \delta_P(m)^{-1/2} \sigma(u)^D \sigma(m)^D du dm$$

\begin{lem}
Cette expression est convergente. Le réel $D$ étant fixé, pour tout réel $R$ il existe un réel $\beta>0$ tel que

$$\chi(c,C,N,D)<<exp(-cR) N^{-R}$$
pour tout $c>0$, $N\geqslant 1$ et tout $C$ tel que $C\geqslant \beta(log(N)+c)$.
\end{lem}

\ul{Preuve}: Fixons le réel $D$. Il existe un réel $D'$ tel que

$$\Xi^G(g)=\Xi^G(g^{-1})<<\Xi^M(m_{\overline{P}}(g^{-1})) \delta_{\overline{P}}(m_{\overline{P}}(g^{-1}))^{1/2} \sigma(g)^{D'}$$
pour tout $g\in G(F)$. En particulier on a

$$\Xi^G(um)<<\Xi^M(m^{-1}m_{\overline{P}}(u^{-1})) \delta_P(m)^{1/2}\delta_{\overline{P}}(m_{\overline{P}}(u^{-1}))^{1/2} \sigma(u)^{D'} \sigma(m)^{D'}$$
pour tous $m\in M(F),u\in U(F)$. En effectuant le changement de variable $u\mapsto u^{-1}$ on en déduit que $\chi(c,C,N,D)$ est essentiellement majoré par

$$\displaystyle\int_{U(F)_c}\mathbf{1}_{\sigma\geqslant C}(u)\delta_{\overline{P}}(m_{\overline{P}}(u))^{1/2}\sigma(u)^{D+D'}\int_{M(F)}\Xi^M(m^{-1}m_{\overline{P}}(u)) \Xi^M(m)\sigma(m)^{D+D'} \kappa_N(m) dmdu$$
pour tout $c>0$, pour tout $C>0$ et pour tout $N\geqslant 1$. Réécrivons l'intégrale sur $M(F)$ en un intégrale sur $G_0(F)$ et $A(F)$. On a la majoration $\sigma(ag_0)<<\sigma(a)\sigma(g_0)$ pour tous $a\in A(F),g_0\in G_0(F)$. L'intégrale sur $M(F)$ est donc essentiellement majorée par le produit des intégrales

$$\displaystyle \int_{A(F)} \kappa_N(a) \sigma(a)^{D+D'} da$$

et

$$\displaystyle\int_{G_0(F)} \Xi^{G_0}(g_0) \Xi^{G_0}(g_0^{-1} g_0(m_{\overline{P}}(u))) \sigma(g_0)^{D+D'} \kappa_N(g_0) dg_0$$
où pour $m\in M(F)$ on définit $g_0(m)$ comme l'unique élément de $G_0(F)$ tel que $mg_0(m)^{-1}\in A(F)$. Comme on l'a déjà vu, la première intégrale est essentiellement majorée par $N^{R_1}$ pour un certain réel $R_1>0$. D'après le lemme 12.0.6, il existe un réel $R_2$ tel que la deuxième intégrale soit essentiellement majorée par $N^{R_2}\Xi^{G_0}(g_0(m_{\overline{P}}(u))) \sigma(g_0(m_{\overline{P}}(u)))^{R_2}$ pour tout $u\in U(F)$ et tout $N\geqslant 1$. On a évidemment $\Xi^{G_0}(g_0(m_{\overline{P}}(u)))=\Xi^M(m_{\overline{P}}(u))$ pour tout $u\in U(F)$. On a aussi une majoration $\sigma(g_0(m_{\overline{P}}(u)))<<\sigma(u)$ pour tout $u\in U(F)$. On en déduit que $\chi(c,C,N,D)$ est essentiellement majorée par

$$\displaystyle N^{R_1+R_2}\int_{U(F)_c} \mathbf{1}_{\sigma\geqslant C}(u)\delta_{\overline{P}}(m_{\overline{P}}(u))^{1/2}\Xi^M(m_{\overline{P}}(u)) \sigma(u)^{D+D'+R_2} du$$
pour tout $c>0$, pour tout $C>0$ et pour tout $N\geqslant 1$. Il ne reste plus qu'à évoquer le lemme 13.0.11 pour obtenir la majoration de l'énoncé $\blacksquare$

\vspace{2mm}

Soient $D$ et $C$ deux réels et $c,c',N$ trois entiers naturels tels que $C,c,c',N\geqslant 1$. Pour $m\in M(F),h\in H(F),u,u'\in U(F)$ on pose
$$\phi(m,h,u,u',D,N)=\Xi^H(h)\Xi^G(u'm)\Xi^G(uhu'm)\kappa_N(m)\sigma(u')^D\sigma(u)^D\sigma(h)^D\sigma(m)^D\delta_P(m)^{-1}$$

On pose aussi
$$I(c,N,D)=\displaystyle\int_{M(F)}\int_{H(F)U(F)_c}\int_{U(F)} \phi(m,h,u,u',D,N)du'dudhdm$$
$$I(c,c',N,D)=\displaystyle\int_{M(F)}\int_{H(F)U(F)_c}\int_{U(F)-U(F)_{c'}} \phi(m,h,u,u',D,N)du'dudhdm$$
$$I(c,c',N,C,D)=\displaystyle\int_{M(F)}\int_{H(F)U(F)_c}\int_{U(F)_{c'}} \mathbf{1}_{\sigma\geqslant C}(hu)\phi(m,h,u,u',D,N)du'dudhdm$$

\begin{prop}
Les expressions ci-dessus sont convergentes et à $c$ et $D$ fixés on a les majorations suivantes
\begin{enumerate}
\item Il existe un réel $R$ tel que
$$I(c,N,D)<<N^R$$
pour tout $N\geqslant 1$
\item Pour tout réel $R$, il existe $\alpha>0$ tel que
$$I(c,c',N,D)<<N^{-R}$$
pour tout $N\geqslant 2$ et tout $c'\geqslant \alpha log(N)$
\item Pour tout réel $R$, il existe $\alpha>0$ tel que
$$I(c,c',N,C,D)<<N^{-R}$$
pour tout $N\geqslant 1$, tout $c'\geqslant 1$ et tout $C\geqslant \alpha(log(N)+c')$
\end{enumerate}
\end{prop}

\ul{Preuve}: Effectuons dans la première intégrale le changement de variables $u'\mapsto h^{-1}u'h$ et majorons $\sigma(h^{-1}u'h)$ par $\sigma(h)^2\sigma(u')$. On reconnaît l'intégrale intérieure sur $U(F)_c$ et $U(F)$: c'est $X(c,D,m,hm)$. D'après le 1. du lemme 13.0.12, il existe un réel $R$ tel que $I(c,N,D)$ soit essentiellement majoré par

$$\displaystyle\int_{H(F)}\int_{M(F)} \Xi^H(h)\Xi^M(hm)\Xi^M(m)\sigma(m)^{D+2R}\sigma(h)^{D+R} \kappa_N(m) dm dh$$
pour tout $N\geqslant 1$. L'intégrale intérieure est essentiellement majorée par le produit de 

$$\displaystyle\int_{A(F)} \kappa_N(a) \sigma(a)^{D+2R} da$$

et

$$\displaystyle\int_{G_0(F)} \Xi^{G_0}(hg_0)\Xi^{G_0}(g_0) \sigma(g_0)^{D+2R} \kappa_N(g_0) dg_0$$

La première intégrale est essentiellement majorée par une puissance de $N$. D'après le lemme 12.0.6, il existe un réel $R'$ telle que la deuxième intégrale soit essentiellement majorée par $\Xi^{G_0}(h)\sigma(h)^{R'}N^{R'}$. D'après le lemme 12.0.7, l'intégrale

$$\displaystyle\int_{H(F)} \Xi^{G_0}(h) \Xi^H(h) \sigma(h)^{D+R+R'} dh$$
est convergente. On en déduit la première majoration de la proposition. En utilisant le 2. du lemme 13.0.12 pour majorer l'intégrale intérieure sur $U(F)_c$ et $U(F)$ dans $I(c,c',N,D)$ et en effectuant les même manipulations que précédemment on obtient la même majoration multipliée par $exp(-\epsilon c')$ pour un certain $\epsilon>0$. On en déduit la deuxième majoration de l'énoncé.

\vspace{2mm}

On a la majoration $I(c,c',N,C,D)\leqslant I(max(c,c'),max(c,c'),N,C,D)$. Pour établir 3., il nous suffit de majorer $I(c,c,N,C,D)$ à $D$ fixé. Introduisons un paramètre $b>0$ que l'on précisera plus tard. L'expression $I(c,c,N,C,D)$ est majorée par la somme de deux intégrales similaires $I_{\geqslant b}(c,c,N,C,D)$ et $I_{< b}(c,c,N,C,D)$ où on a échangé le terme $\mathbf{1}_{\sigma\geqslant C}(hu)$ par les termes $\mathbf{1}_{\sigma\geqslant b}(h)$ et $\mathbf{1}_{\sigma< b}(h)\mathbf{1}_{\sigma\geqslant C-b}(u)$ respectivement. D'après la proposition 13.0.3, il existe un réel $D'$ tel que $I_{\geqslant b}(c,c,N,C,D)$ soit essentiellement majoré par le produit d'une puissance de $c$ et de

$$\displaystyle\int_{H(F)}\int_{M(F)}\mathbf{1}_{\sigma\geqslant b}(h) \Xi^H(h)\Xi^M(hm)\Xi^M(m)\kappa_N(m)\sigma(h)^{D'}\sigma(m)^{D'}dmdh$$

\noindent pour tout $N\geqslant 1$ et pour tous $c,C,b>0$. On peut comme on l'a fait plusieurs fois décomposer l'intégrale sur $M(F)$ en un produit d'une intégrale sur $A(F)$ et d'une intégrale sur $G_0(F)$. Comme on l'a aussi vu plusieurs fois, l'intégrale sur $A(F)$ est essentiellement majorée par une puissance de $N$. On reconnaît l'intégrale sur $G_0(F)$: c'est $I^1(N,D',b)$. Cette intégrale est majorée par le lemme 12.0.9. Alliant cette majoration aux précédentes, on en déduit qu'il existe deux réel $R_1$ et $\epsilon_1>0$ tels que $I_{\geqslant b}(c,c,N,C,D)$ soit essentiellement majorée par $c^{R_1}N^{R_1}exp(-\epsilon_1 b)$ pour tout $N\geqslant 1$ et pour tous $c,C,b>0$. En particulier pour $\displaystyle b= \frac{(R+R_1)log(N)+R_1 log(c)}{\epsilon_1}$, l'expression $I_{\geqslant b}(c,c,N,C,D)$ est essentiellement majorée par $N^{-R}$. On fait dorénavant ce choix pour $b$.

\vspace{2mm}

 Il reste à majorer $I_{>b}(c,c,N,C,D)$, on effectue les changements de variable $u'\mapsto hu'h^{-1}$ et $u\mapsto uu'^{-1}$. On majore $\sigma(uu'^{-1})$ par $\sigma(u)\sigma(u')$ et $\mathbf{1}_{\geqslant C-b}(uu'^{-1})$ par $\mathbf{1}_{\geqslant (C-b)/2}(u)+\mathbf{1}_{\geqslant (C-b)/2}(u')$. On obtient que $I_{<b}(c,c,N,C,D)$ est essentiellement majoré par l'intégrale sur $H(F)$ du produit de $\mathbf{1}_{\sigma<b}(h)\Xi^H(h)\sigma(h)^D$ et de l'intégrale

$$\displaystyle\int_{M(F)}\int_{U(F)_c}\int_{U(F)_c} (\mathbf{1}_{\geqslant (C-b)/2}(u)+\mathbf{1}_{\geqslant (C-b)/2}(u'))\Xi^G(u'm) \Xi^G(uhm) \kappa_N(m)\sigma(u')^{2D}\sigma(u)^D \sigma(m)^D$$

$$\delta_P(m)^{-1} dudu'dm$$

\noindent pour tout $c>0$, pour tout $C>0$ et pour tout $N\geqslant 1$. Cette intégrale est somme de deux intégrales similaires. Il suffit donc de borner l'une des deux. Considérons l'intégrale où on a fait disparaître le termer $\mathbf{1}_{\geqslant (C-b)/2}(u)$, on reconnaît l'intégrale sur $u$: c'est $X(c,D,hm)$. D'après le lemme 13.0.12, elle est essentiellement majorée par $c^{R_2}\Xi^M(hm)\delta_P(m)^{1/2}\sigma(m)^{R_2}$ pour un certain réel $R_2$. On en déduit que $I_{>b}(c,c,N,C,D)$ est essentiellement majorée par l'intégrale sur $H(F)$ du produit de $c^{R_2}\mathbf{1}_{\sigma<b}(h)\Xi^H(h)\sigma(h)^D$ et de l'intégrale

$$\displaystyle\int_{M(F)}\int_{U(F)_c}\mathbf{1}_{\geqslant (C-b)/2}(u)\Xi^M(hm)\Xi^G(u'm)\kappa_N(m)\sigma(u')^{2D} \sigma(m)^{D+R_2} \delta_P(m)^{-1/2}du dm$$
pour tout $c$, pour tout $C>0$ et pour tout $N\geqslant 1$. Il existe un réel $\beta_1$ tel que $\Xi^M(mm')<<exp(\beta_1\sigma(m))\Xi^M(m')$ pour tous $m,m'\in M(F)$. L'intégrale précédente est donc essentiellement majorée par le produit de $exp(\beta_1 \sigma(h))$ et de $\chi(c,(C-b)/2,N,2D)$. D'après le lemme 13.0.13, pour tout réel $R'$, il existe un réel $\alpha_1$ tel que ce dernier terme soit essentiellement majoré par $exp(-cR')N^{-R'}$ pour tout $N\geqslant 1$ et pour tous $c>0,C>0$ vérifiant $C-b\geqslant \alpha_1(log(N)+c)$. L'intégrale sur $H(F)$ de $\mathbf{1}_{\sigma<b}(h)\Xi^H(h)\sigma(h)^Dexp(\beta_1\sigma(h))$ est majorée par le produit de $b^Dexp(\beta_1 b)$ et de
 
$$\displaystyle\int_{H(F)}\mathbf{1}_{\sigma<b}(h) dh$$

D'après l'inégalité 4.3(1) de [W2], cette intégrale est essentiellement majorée par $exp(R_3b)$ pour un certain réel $R_3$. Au final, $I_{<b}(c,c,N,C,D)$ est essentiellement majoré par

$$c^{R_2}exp(-cR')N^{-R'}b^Dexp((\beta_1+R_3) b)$$

\noindent pour tout $N\geqslant 1$ et pour tout $c>0,C>0$ vérifiant $C-b\geqslant \alpha_1(log(N)+c)$. D'après notre choix de $b$, le terme $b^Dexp((\beta_1+R_3) b)$ est essentiellement majorée par le produit d'une puissance de $c$ et de $N^{\frac{R+R_1}{2\epsilon_1}}$. La puissance de $c$ multipliée par $c^{R_2}exp(-cR')$ est essentiellement majorée par 1. Pour $R'=R+\frac{R+R_1}{2\epsilon_1}$, on obtient la majoration de l'énoncé $\blacksquare$

\section{Entrelacements tempérés}

Soient $\pi\in Temp(G)$, $\sigma\in Temp(H)$ et fixons des produits scalaires invariants sur $E_\pi$ et $E_\sigma$. Pour $e,e'\in E_\pi$, $\epsilon,\epsilon'\in E_\sigma$ et $c\in\mathbb{N}$, on pose

$$\mathcal{L}_{\pi,\sigma,c}(\epsilon'\otimes e',\epsilon\otimes e)=\displaystyle\int_{H(F)U(F)_c} (\sigma(h)\epsilon',\epsilon)(e',\pi(hu)e)\overline{\xi}(u) du dh$$

D'après le corollaire 13.0.1, $\mathcal{L}_{\pi,\sigma,c}(\epsilon'\otimes e',\epsilon\otimes e)$ est défini par une intégrale convergente. Pour tout réel $c'>0$ notons $\omega_A(c')$ l'ensemble des $a\in A(F)$ qui vérifient $val_F(a_i-1)\geqslant c'$ pour $i=1,\ldots,r$.

\begin{lem}
Soit $c'\geqslant 00$ un réel, il existe un réel $c_0\geqslant 0$ tel que pour tout $\pi\in Temp(G)$, pour tout $\sigma\in Temp(H)$, pour tous $e,e'\in E_\pi^{\omega_A(c')}$, pour tous $\epsilon,\epsilon'\in E_\sigma$ et pour tout $c\geqslant c_0$ on ait

$$\mathcal{L}_{\pi,\sigma,c}(\epsilon'\otimes e',\epsilon\otimes e)=\mathcal{L}_{\pi,\sigma,c_0}(\epsilon'\otimes e',\epsilon\otimes e)$$
\end{lem}

\ul{Preuve}: Fixons $c'>0$ et soient $e,e'\in E_\pi^{\omega_A(c')}$ et $\epsilon,\epsilon'\in E_\sigma$. On a alors, pour tout $c\geqslant 0$,

$$\displaystyle \mathcal{L}_{\pi,\sigma,c}(\epsilon'\otimes e',\epsilon\otimes e)=\int_{H(F)U(F)_c}(\sigma(h)\epsilon',\epsilon)(e',\pi(hu)e)\int_{\omega_A(c')}\overline{\xi}(aua^{-1})da du dh$$

\noindent Il existe un réel $c_0\geqslant 0$, tel que $\int_{\omega_A(c')}\overline{\xi}(aua^{-1})da=0$ pour tout $u\in U(F)-U(F)_{c_0}$, ce qui permet de conclure $\blacksquare$

\vspace{3mm}

Posons

$$\mathcal{L}_{\pi,\sigma}(\epsilon'\otimes e',\epsilon\otimes e)=lim_{c\to\infty}\mathcal{L}_{\pi,\sigma,c}(\epsilon'\otimes e',\epsilon\otimes e)$$

On a

$$\mathcal{L}_{\pi,\sigma}(\sigma(h')\epsilon'\otimes \pi(hu)e',\epsilon\otimes \pi(u')e)=\xi(u^{-1}u')\mathcal{L}_{\pi,\sigma}(\epsilon'\otimes e',\sigma(h^{-1})\epsilon\otimes\pi(h'^{-1}) e)$$

 Pour tous $h,h'\in H(F)$, $u,u'\in U(F)$ ,$e,e'\in E_\pi$, $\epsilon,\epsilon'\in E_\sigma$. Il est alors facile de vérifier que $\mathcal{L}_{\pi,\sigma}\neq0$ entraîne $Hom_{H,\xi}(\pi,\sigma)\neq 0$. L'objectif de cette section est d'établir la réciproque.

\vspace{4mm}

\subsection{Un lemme sur les entrelacements}

Fixons des données $(w_i)_{i=0,\ldots,l},P_0,P_H,A_0,A_H$ comme dans la section 11 et des compacts-ouverts $C_0\subset G_0(F)$ et $C_H\subset H(F)$ qui vérifient la conclusion de la proposition 11.0.1. On peut toujours supposer que $C_HA_H^1=C_H$ et $A_0^1C_0=C_0$, ce que l'on fait. On a alors $G_0(F)=C_HA_H^+A_0^+C_0$ d'où $M(F)=C_HA_H^+A_0^+A(F)C_0$. Posons $P_{0,G}=P_0AU$ et $A_{0,G}=AA_0$. Ce sont respectivement un sous-groupe parabolique minimal et un tore déployé maximal de $G$. On note $A_{0,G}^+$ l'ensemble des éléments de $A_{0,G}(F)$ qui contractent $P_{0,G}(F)$. On pose $\Lambda_H^+=A_H^+/A_H^1,\Lambda_0^+=A_0^+/A_0^1,\Lambda_G^+=A_{0,G}^+/A_{0,G}^1$ que l'on identifie à des sous-ensembles de $A_H^+,A_0^+$ et $A_{0,G}^+$ via le choix de sections.

\begin{lem}
Soient $\Gamma_1$ et $\Gamma_2$ des sous-groupes compacts-ouverts de $G(F)$ et $H(F)$ respectivement. Il existe un sous-ensemble compact-ouvert $C_M\subset M(F)$ et des sous-groupes ouverts $K_1\subset K$, $K_2\subset K_H$ tels que pour tous $\pi\in Irr(G),\sigma\in Irr(H)$ pour tout $l\in Hom_{H,\xi}(\pi,\sigma)$ on ait

\begin{enumerate}
\item Pour tout $e\in E_\pi^{\Gamma_1}$ et pour tout $m\in M(F)-C_HA_H^+A_{0,G}^+C_M$,

$$l(\pi(m)e)=0$$

\item Pour tous $e\in E_\pi^{\Gamma_1}$, $\epsilon^\vee\in E_{\sigma^\vee}^{\Gamma_2}$ et pour tous $h\in C_HA_H^+$, $m\in A_{0,G}^+C_M$, on a

$$<\epsilon^\vee,l(\pi(hm)e)>=<\sigma^\vee(e_{K_2})\sigma^\vee(h^{-1})\epsilon^\vee,l(\pi(e_{K_1})\pi(m)e)>$$
\end{enumerate}
\end{lem}

\ul{Preuve}:

1-Pour $m\in M(F)$ et une décomposition $m=k_Ha_Ha_0ak_0$ avec $k_H\in C_H,k_0\in C_0,a_H\in A_H^+,a_0\in A_0^+$ et $a\in A(F)$. On a alors $l(\pi(m)e)=\sigma(k_H a_H)l(\pi(a_0ak_0)e)$. Puisque $C_0$ est compact, il suffit de montrer qu'il existe un compact $C_{0,G}\subset A_{0,G}(F)$ ne dépendant que de $\Gamma_1$ tel que $l(\pi(a)e)=0$ pour tout $a\in A_0^+A(F)-A_{0,G}^+C_{0,G}$ et tout $e\in E_\pi^{\Gamma_1}$. Il existe certainement un compact $C_{0,G}\subset A_{0,G}(F)$ tel que pour tout $a\in A_0^+A(F)-A_{0,G}^+C_{0,G}$ on ait $a\Gamma_1a^{-1}\cap(U(F)-Ker\xi)\neq\emptyset$ (on utilise ici le fait que $v_0=(u_1+u_{-1})/2$ avec les notations de la section 11). Pour un tel $a$ et $u\in a\Gamma_1a^{-1}\cap(U(F)-Ker\xi)$ on a alors $l(\pi(a)e)=l(\pi(ua)e)=\xi(u)l(\pi(a)e)$ d'où $l(\pi(a)e)=0$.

\vspace{2mm}

2- Résulte de ce que l'application produit $U\times P_H\times \overline{P}_{0,G}\to G$ est submersive à l'origine (d'après le lemme 11.0.1) $\blacksquare$

\subsection{Entrelacements dans une famille d'induites}

Fixons un sous-groupe de Levi minimal de $H$ défini sur $F$ ainsi qu'un sous-groupe compact spécial $K_H$ de $H(F)$ en bonne position par rapport à ce Levi minimal. Soient $Q=LU_Q$ et $R=SU_S$ des sous-groupes paraboliques semistandards de $G$ et $H$ respectivement et $\rho$ et $\tau$ des représentations irréductibles de la série discrète de $L(F)$ et $S(F)$ respectivement. Pour tous $\lambda\in i\mathcal{A}_{L,F}^*$ et $\mu\in i\mathcal{A}_{S,F}^*$ on note $\pi_\lambda=i_Q^G(\rho_\lambda)$ et $\sigma_\mu=i^H_R(\tau_\mu)$. Les représentations $\pi_\lambda$ se réalisent dans un espace commun $\mathcal{K}_{Q,\rho}^G$ de fonctions sur $K$ et les représentations $\sigma_\mu$ se réalisent dans un espace commun $\mathcal{K}_{R,\tau}^H$ de fonctions sur $K_H$. On peut munir comme en 1.5 ces deux espaces sont munis de produits scalaires invariants. On utilise ces produits scalaires pour définir les formes sesquilinéaires $\mathcal{L}_{\pi_\lambda,\sigma_\mu}$ pour tout $(\lambda,\mu)\in i\mathcal{A}_{L,F}^*\times i\mathcal{A}_{S,F}^*$. On notera $\lambda\mapsto m(\rho_\lambda)$ et $\mu\mapsto m(\tau_\mu)$ les mesures de Plancherel.

\vspace{4mm}

\begin{prop} On conserve les notations précédentes. Alors \\
(i) Pour tous $e,e'\in\mathcal{K}_{Q,\rho}^G$ et $\epsilon,\epsilon'\in \mathcal{K}_{R,\tau}^H$, la fonction $(\lambda,\mu)\mapsto \mathcal{L}_{\pi_\lambda,\sigma_\mu}(\epsilon'\otimes e',\epsilon\otimes e)$ est analytique sur $i\mathcal{A}_{L,F}^*\times i\mathcal{A}_{S,F}^*$. \\
(ii) S'il existe $\lambda_0\in i\mathcal{A}_{L,F}^*$ et $\mu_0\in i\mathcal{A}_{S,F}^*$ tels que $\mathcal{L}_{\pi_{\lambda_0},\sigma_{\mu_0}}\neq 0$ alors pour tout $\lambda\in i\mathcal{A}_{L,F}^*$ et tout $\mu\in i\mathcal{A}_{S,F}^*$ on a $\mathcal{L}_{\pi_\lambda,\sigma_\mu }\neq 0$.
\end{prop}

\ul{Preuve}: (i) Fixons $e,e'\in\mathcal{K}_{Q,\rho}^G$ et $\epsilon,\epsilon'\in \mathcal{K}_{R,\tau}^H$. Pour tous $\lambda\in i\mathcal{A}_{L,F}^*$, $\mu\in i\mathcal{A}_{S,F}^*$ on a $\mathcal{L}_{\pi_\lambda,\sigma_\mu}(\epsilon'\otimes e',\epsilon\otimes e)=\mathcal{L}_{\pi_\lambda,\sigma_\mu,c}(\epsilon'\otimes e',\epsilon\otimes e)$ pour $c$ assez grand. On peut choisir un entier $c_0$ à partir duquel l'égalité précédente est vérifiée qui ne dépend que des stabilisateurs de $e$ et $e'$ dans $A(F)\cap K$. Par conséquent on peut trouver un entier $c$ tel que l'égalité soit vérifiée simultanément pour tous $(\lambda,\mu)$. On a alors

$$\mathcal{L}_{\pi_\lambda,\sigma_\mu}(\epsilon'\otimes e',\epsilon\otimes e)=\displaystyle\int_{H(F)U(F)_c} (\sigma_\mu(h)\epsilon',\epsilon)(e',\pi_\lambda(hu)e)\overline{\xi}(u) du dh$$
pour tous $(\lambda,\mu)\in i\mathcal{A}_{L,F}^*\times i\mathcal{A}_{S,F}^*$. \\

Il suffit de vérifier que cette intégrale est encore uniformément convergente pour $(\lambda,\mu)$ dans un voisinage de $i\mathcal{A}_L^*\times i\mathcal{A}_S^*\subset \mathcal{A}_{L,\mathbb{C}}^*\times \mathcal{A}_{S,\mathbb{C}}^*$. Soit $\epsilon>0$ pour $(\lambda,\mu)$ dans un voisinage assez petit de $i\mathcal{A}_L^*\times i\mathcal{A}_S^*$ on a des majorations uniformes

$$|(\sigma_\mu(h)\epsilon',\epsilon)|<<exp(\epsilon \sigma(h))\Xi^H(h)$$
$$|(e',\pi_\lambda(hu)e)|<<exp(\epsilon \sigma(hu)) \Xi^G(hu)$$
pour tout $h\in H(F)$ et pour tout $u\in U(F)$. Il suffit donc de vérifier que pour $\epsilon$ assez petit l'intégrale suivante est absolument convergente

\begin{center}
$\displaystyle\int_{H(F)U(F)_c} exp(\epsilon(\sigma(h)+\sigma(hu)))\Xi^H(h)\Xi^G(hu) dhdu$
\end{center}

C'est une conséquence de la proposition 13.0.4 et du lemme 12.0.5.

\vspace{2mm}

(ii) On s'inspire ici grandement de la démonstration de la proposition 6.4.1 de [SV]. On suppose que la fonction $(\lambda,\mu)\mapsto \mathcal{L}_{\pi_\lambda,\sigma_\mu}$ n'est pas identiquement nulle et on veut montrer qu'elle est non nulle en $(0,0)$. \\
 Soient $e\in \mathcal{K}_{Q,\rho}^G$ et $\epsilon\in \mathcal{K}_{R,\tau}^G$ tels que $(e,e)=(\epsilon,\epsilon)=1$ et tels que l'ordre d'annulation en $(0,0)$
 de $(\lambda,\mu)\mapsto \mathcal{L}_{\pi_\lambda,\sigma_\mu}(\epsilon\otimes e,\epsilon\otimes e)$ soit minimal. Reprenons les notations du paragraphe 1.7. Il y ait défini des ensembles $\mathcal{E}(\rho)$ et $\mathcal{E}(\tau)$. Soient $\phi$ et $\psi$ des fonctions $C^\infty$ sur $i\mathcal{A}_{L,F}^*$ et $i\mathcal{A}_{S,F}^*$ respectivement qui vérifient
 
\begin{center}
(1) $(w^{-1}Supp(\phi)+\mu)\cap Supp(\phi)=\emptyset$ pour tout $(w,\mu)\in \mathcal{E}(\rho)-\{(Id,0)\}$ et $(w^{-1}Supp(\psi)+\mu)\cap Supp(\psi)=\emptyset$ pour tout $(w,\mu)\in \mathcal{E}(\tau)-\{(Id,0)\}$.
\end{center}

\noindent On définit une fonction $\Phi$ sur $G(F)$ par la formule suivante

$$\Phi(g)=\displaystyle\int_{i\mathcal{A}_{L,F}^*\times i\mathcal{A}_{S,F}^*} \phi(\lambda)\psi(\mu)m(\rho_\lambda)m(\tau_\mu) \mathcal{L}_{\pi_\lambda,\sigma_\mu}(\epsilon\otimes\pi_\lambda(g)e,\epsilon\otimes e) d\mu d\lambda$$

\noindent A $g$ fixé, il existe un $c$ assez grand de sorte que

$$\mathcal{L}_{\pi_\lambda,\sigma_\mu}(\epsilon\otimes\pi_\lambda(g)e,\epsilon\otimes e)=\displaystyle\int_{H(F)U(F)_c} (\sigma_\mu(h)\epsilon,\epsilon)(\pi_\lambda(g)e,\pi_\lambda(hu)e)\overline{\xi}(u) du dh$$

\noindent pour tous $\lambda\in i\mathcal{A}_{L,F}^*$ et $\mu\in i\mathcal{A}_{S,F}^*$. Pour un tel $c$, on a alors

$$\Phi(g)=\displaystyle\int_{i\mathcal{A}_{L,F}^*\times i\mathcal{A}_{S,F}^*} \int_{H(F)U(F)_c} \phi(\lambda)\psi(\mu) m(\rho_\lambda)m(\tau_\mu) (\sigma_\mu(h)\epsilon,\epsilon)(\pi_\lambda(g)e,\pi_\lambda(hu)e)\overline{\xi}(u) du dh d\mu d\lambda$$

L'intégrale ci-dessus est absolument convergente, car il existe $C>0$ indépendant de $\lambda$ et $\mu$ tel que pour $h\in H(F)$ et $u\in U(F)$, $|(\sigma_\mu(h)\epsilon,\epsilon)(\pi_\lambda(g)e,\pi_\lambda(hu)e)|\leqslant C\Xi^H(h)\Xi^G(hu)$ et l'intégrale $\int_{H(F)U(F)_c} \Xi^H(h)\Xi^G(hu) dh du$ est absolument convergente. On peut donc changer l'ordre d'intégration et on obtient alors

\begin{center}
(2) $\Phi(g)=\displaystyle\int_{H(F)U(F)_c} f_{e,e,\phi}(u^{-1}h^{-1}g) f_{\epsilon,\epsilon,\psi}(h) \overline{\xi(u)} du dh$
\end{center}

\noindent (on repren ici les notations du paragraphe 1.7). Montrons

\begin{center}
(3) Il existe un entier $c_0$ tel que l'égalité (2) soit vérifiée pour tout $c\geqslant c_0$ et pour tout $g\in M(F)K$.
\end{center}

En effet, l'entier $c_0(g)$ à partir duquel la formule précédente est vérifiée ne dépend que d'un sous-groupe compact-ouvert de $A(F)\cap K$ qui laisse stable $e$ et $\pi_\lambda(g)e$. Comme $A(F)$ commute à $M(F)$ on peut trouver un sous-groupe compact-ouvert de $A(F)$ qui laisse stable ces éléments pour tout $g\in M(F)K$.

\begin{center}
(4) La fonction $(m,k)\mapsto |\Phi(mk)|^2 \delta_P(m)^{-1}$ est intégrable sur $M(F)\times K$.
\end{center}

En effet, les fonctions $f_{e,e,\phi}$ et $f_{\epsilon,\epsilon,\psi}$ sont dans $\mathcal{S}(G(F))$ et $\mathcal{S}(H(F))$ respectivement. Par conséquent pour tout réel $D'>0$ on a une majoration

$$|\Phi(mk)|<<\displaystyle\int_{H(F)U(F)} \Xi^G(uh^{-1}m)\Xi^H(h) \sigma(uh^{-1}m)^{-D'}\sigma(h)^{-D'}dh du$$
pour tout $m\in M(F)$ et pour tout $k\in K$. D'après le lemme II.4.5 de [W3], pour tout $D>0$, il existe $D'>0$ tel que le terme précédent soit essentiellement majoré par

$$\delta_{P}(m)^{1/2}\displaystyle\int_{H(F)} \Xi^M(h^{-1}m)\Xi^H(h) \sigma(h^{-1}m)^{-D} \sigma(h)^{-D}dh$$
On a $\delta_{P}(m)=\delta_P(a)$, $\Xi^M(h^{-1}m)=\Xi^{G_0}(h^{-1}g_0)$ et $\sigma(a)^{1/2}<<\sigma(h^{-1}m)$ pour tout $h\in H(F)$ et pour tout $m=ag_0\in M(F)=A(F)G_0(F)$. Il suffit donc de montrer que les deux intégrales suivantes sont absolument convergentes dès que $D$ est assez grand

$$\displaystyle\int_{G_0(F)} \int_{H(F)\times H(F)} \Xi^{G_0}(h^{-1}g_0)\Xi^H(h)\Xi^{G_0}(h'^{-1}g_0)\Xi^H(h') \sigma(h)^{-D}\sigma(h')^{-D}dhdh'dg_0$$

\noindent et

$$\displaystyle\int_{A(F)} \sigma(a)^{-D} da$$

\noindent la deuxième intégrale ne pose pas de problème et d'après le lemme 12.0.6, la première est essentiellement majorée par

$$\displaystyle\int_{H(F)\times H(F)} \Xi^H(h^{-1}h')\Xi^H(h)\Xi^H(h') \sigma(h)^{-D}\sigma(h')^{-D}dhdh'$$
qui est une intégrale convergente pour $D$ assez grand.

\vspace{4mm}

Dans l'égalité (2) on peut faire tendre le réel $c$ vers l'infini. Comme on l'a vu lors de la preuve de (4) le terme sous l'intégrale est intégrable sur $H(F)U(F)$. On en déduit que l'on a aussi

\begin{center}
(5) $\displaystyle \Phi(g)=\int_{H(F)U(F)} f_{e,e,\phi}(u^{-1}h^{-1}g)f_{\epsilon,\epsilon,\psi}(h)\overline{\xi(u)} dudh$
\end{center}
pour tout $g\in G(F)$. Posons $S=\int_{M(F)\times K}|\Phi(mk)|^2 \delta_P(m)^{-1} dk dm$. Explicitons cette intégrale en remplaçant $\Phi$ par son expression (5) et $\overline{\Phi}$ par son expression (2). On obtient

\[\begin{aligned}
S & =\displaystyle\int_{M(F)\times K}\int_{H(F)U(F)}\int_{H(F)U(F)_c} f_{e,e,\phi}(u'h'^{-1}mk)\overline{f_{e,e,\phi}(uh^{-1}mk)} f_{\epsilon,\epsilon,\psi}(h')\overline{f_{\epsilon,\epsilon,\psi}(h)} \\ 
 & \delta_P(m)^{-1}\xi(u'u^{-1})dudhdu'dh'dkdm \\
 & =\displaystyle\int_{G(F)}\int_{H(F)}\int_{H(F)U(F)_c} f_{e,e,\phi}(h'^{-1}g)\overline{f_{e,e,\phi}(uh^{-1}g)} f_{\epsilon,\epsilon,\psi}(h')\overline{f_{\epsilon,\epsilon,\psi}(h)} \overline{\xi(u)} dudhdh'dg
\end{aligned}\]

\noindent pour tout $c\geqslant c_0$. D'après les majorations déjà effectuées, cette triple intégrale est absolument convergente. L'hypothèse (1) permet d'appliquer l'égalité 1.7(3), on a donc

$$\displaystyle\int_{G(F)} f_{e,e,\phi}(h'^{-1}g)\overline{f_{e,e,\phi}(uh^{-1}g)} dg=\int_{i\mathcal{A}_{L,F}^*} m(\rho_\lambda) |\phi(\lambda)|^2(e,\pi_\lambda(uh^{-1}h')e) d\lambda$$

On obtient, après le changement de variable $h'\mapsto hh'$

$$S=\displaystyle\int_{H(F)}\int_{H(F)U(F)_c}\int_{i\mathcal{A}_{L,F}^*} m(\rho_\lambda) |\phi(\lambda)|^2 (e,\pi_\lambda(uh')e) f_{\epsilon,\epsilon,\psi}(hh') \overline{f_{\epsilon,\epsilon,\psi}(h)} \overline{\xi(u)} d\lambda du dh dh'$$

L'intégrale ci-dessus est encore absolument convergente. Appliquons à nouveau l'égalité 1.7(3), on a

$$\displaystyle\int_{H(F)} f_{\epsilon,\epsilon,\psi}(hh')\overline{f_{\epsilon,\epsilon,\psi}(h)} dh=\int_{i\mathcal{A}_{S,F}^*} m(\tau_\mu)|\psi(\mu)|^2 (\sigma_\mu(h')\epsilon,\epsilon) d\mu$$

\noindent On en déduit que

\[\begin{aligned}
\mbox{(6)} \;\;\; S & =\displaystyle\int_{H(F)U(F)_c}\int_{i\mathcal{A}_{L,F}^*\times i\mathcal{A}_{S,F}^*} m(\rho_\lambda)m(\tau_\mu)|\phi(\lambda)|^2 |\psi(\mu)|^2 (\sigma_\mu(h)\epsilon,\epsilon)(e,\pi_\lambda(uh)e) \overline{\xi(u)} d\lambda d\mu dudh \\
 & =\displaystyle\int_{i\mathcal{A}_{L,F}^*\times i\mathcal{A}_{S,F}^*}m(\rho_\lambda)m(\tau_\mu)|\phi(\lambda)|^2 |\psi(\mu)|^2 \mathcal{L}_{\pi_\lambda,\sigma_\mu}(\epsilon\otimes e,\epsilon \otimes e) d\lambda d\mu
\end{aligned}\]

\vspace{4mm}

Soient $\Gamma_1\subset K$ et $\Gamma_2\subset K_H$ des sous-groupes ouverts distingués stabilisant $e$ et $\epsilon$. Soient $C_M\subset M(F)$ un compact et $K_1\subset K$, $K_2\subset K_H$ des sous-groupes ouverts qui vérifient les conclusions du lemme 14.1.1. On peut toujours supposer que $A_{0,G}^1C_M=C_M$. Puisque l'on a $\mathcal{L}_{\pi_\lambda,\sigma_\mu}\in Hom_{H,\xi}(E_{\pi_\lambda},E_{\sigma_\mu})\otimes \overline{Hom_{H,\xi}(E_{\pi_\lambda},E_{\sigma_\mu})}$ pour tout $(\lambda,\mu)\in i\mathcal{A}_{L,F}^*\times i\mathcal{A}_{S,F}^*$, le 1. du lemme 14.1.1 implique que le support de $(m,k)\mapsto \Phi(mk)$ est contenu dans $(C_H\Lambda_H^+\Lambda_G^+ C_M)\times K$.
Soient $\mathcal{B}^{K_1}$ et $\mathcal{B}^{K_2}$ des bases orthonormales de $(\mathcal{K}_{Q,\rho}^G)^{K_1}$ et $(\mathcal{K}_{R,\tau}^H)^{K_2}$ respectivement. \\
D'après le 2. du lemme 14.1.1, pour tout $m=k_Ha_Ha_Gk_M\in C_H\Lambda_H^+\Lambda_G^+ C_M$ et tout $k\in K$, on a

$$\Phi(mk)=\displaystyle\sum_{e'\in \mathcal{B}^{K_1},\epsilon'\in \mathcal{B}^{K_2}} F_{e',\epsilon'}(a_Gk_Mk,k_Ha_H)$$

\noindent où on a posé pour $e'\in \mathcal{B}^{K_1}$ et $\epsilon'\in \mathcal{B}^{K_2}$

$$F_{e',\epsilon'}(g,h)=\displaystyle\int_{i\mathcal{A}_{L,F}^*\times i\mathcal{A}_{S,F}^*} \phi(\lambda)\psi(\mu) m(\rho_\lambda)m(\tau_\mu)\mathcal{L}_{\pi_\lambda,\sigma_\mu}(\epsilon\otimes e',\epsilon'\otimes e) (\pi_\lambda(g)e,e')(\sigma_\mu(h)\epsilon',\epsilon) d\lambda d\mu$$

\noindent Soit $F=\sum_{e'\in \mathcal{B}^{K_1},\epsilon'\in\mathcal{B}^{K_2}} F_{e',\epsilon'}$. D'après ce qui précède, le 3- du lemme 11.0.4 et le fait que $\delta_P\delta_{P_0}=\delta_{P_{0,G}}$, on a une majoration

$$S<<\displaystyle\sum_{a_G\in \Lambda_G^+,a_H\in\Lambda_H^+} \delta_{P_{0,G}}(a_G)^{-1}\delta_{P_H}(a_H)^{-1} \int_{C_H\times C_M\times K} |F(a_Gk_Mk,k_Ha_H)|^2 dkdk_Mdk_H$$

\noindent où la constante implicite ne dépend pas de $\phi$ et $\psi$. Soit $C_G$ un sous-ensemble compact de $G(F)$ qui contient $C_MK$. On a alors

\[\begin{aligned}
\displaystyle\int_{C_H\times C_M\times K} |F(a_Gk_Mk,k_Ha_H)|^2 dkd & k_Mdk_H  <<\int_{C_H\times C_G} |F(a_Gk_G,k_Ha_H)|^2 dk_Gdk_H\\
 & << \int_{(C_H)^2\times (C_G)^2} |F(k_G^1a_Gk_G^2,k_H^1a_Hk_H^2)|^2 dk_G^1dk_G^2dk_H^1dk_H^2
\end{aligned}\]

Pour toute fonction $f$ intégrable et mesurable sur $G(F)$ on a

$$\displaystyle\int_{KA_{0,G}^+K} f(g)dg=\sum_{a_G\in\Lambda_G^+} mes(Ka_GK)\int_{K\times K} f(k_1a_Gk_2) dk_1dk_2$$

D'après [W3] I.(5), on a $\delta_{P_{0,G}}(a_G)^{-1}<<mes(Ka_GK)$ pour tout $a_G\in\Lambda_{0,G}^+$. Les mêmes résultats sont bien sûrs valables pour $H(F)$. On en déduit la majoration

\[\begin{aligned}
\displaystyle\sum_{a_G\in \Lambda_G^+,a_H\in\Lambda_H^+} \delta_{P_{0,G}}(a_G)^{-1}\delta_{P_H}(a_H)^{-1}\int_{(C_H)^2\times (C_G)^2} |f(k_G^1a_Gk_G^2,k_H^1a_Hk_H^2)|^2 dk_G^1dk_G^2dk_H^1dk_H^2 \\
<<\int_{G(F)\times H(F)} |f(g,h)|^2 dhdg
\end{aligned}\]

\noindent pour toute fonction $f$ mesurable et intégrable sur $G(F)\times H(F)$. En particulier, on a

$$\mbox{(7)} \!\;\;\; S<<\int_{G(F)\times H(F)} |F(g,h)|^2 dhdg$$

\noindent la constante implicite ne dépendant pas de $\phi$ et $\psi$. Pour $(\lambda,\mu)\in i\mathcal{A}_{L,F}^*\times i\mathcal{A}_{S,F}^*$, considérons la représentation $\pi_\lambda\boxtimes \sigma_\mu=i^{G\times H}_{Q\times R}(\rho_\lambda\boxtimes \sigma_\mu)$. Avec les notations de la section 1.7, on a alors $F_{e',\epsilon'}=f_{e'\otimes \epsilon,e\otimes \epsilon',\varphi_{e',\epsilon'}}$, où

$$\varphi_{e',\epsilon'}(\lambda,\mu)=\phi(\lambda)\psi(\mu)\mathcal{L}_{\pi_\lambda,\sigma_\mu}(\epsilon\otimes e',\epsilon'\otimes e)$$

\noindent pour tout $(\lambda,\mu)\in i\mathcal{A}_{L,F}^*\times i\mathcal{A}_{S,F}^*$. D'après l'hypothèse (1) faite sur les supports de $\phi$ et $\psi$ et 1.7(3), les fonctions $(F_{e',\epsilon'})_{e'\in \mathcal{B}^{K_1},\epsilon'\in \mathcal{B}^{K_2}}$ forment une famille orthogonale pour le produit scalaire $L^2$ et la norme $L^2$ de $F_{e',\epsilon'}$ est égale à

$$\displaystyle\int_{i\mathcal{A}_{L,F}^*\times i\mathcal{A}_{S,F}^*} \vert \phi(\lambda)\vert ^2 \vert \psi(\mu)\vert ^2 \vert \mathcal{L}_{\pi_\lambda,\sigma_\mu}(\epsilon\otimes e',\epsilon'\otimes e)\vert ^2 d\mu d\lambda$$

\noindent On en déduit que

\[\begin{aligned}
\displaystyle\int_{G(F)\times H(F)} |F(g,h)|^2 dg dh=\int_{i\mathcal{A}_{L,F}^*\times i\mathcal{A}_{S,F}^*} & |\phi(\lambda)|^2|\psi(\mu)|^2m(\rho_\lambda) m(\tau_\mu) \\
 & \left(\sum_{e'\in\mathcal{B}^{K_1},\epsilon'\in\mathcal{B}^{K_2}} |\mathcal{L}_{\pi_\lambda,\sigma_\mu}(\epsilon\otimes e',\epsilon'\otimes e) | ^2 \right)d\lambda d\mu
\end{aligned}\]

D'après l'égalité (6) et l'inégalité (7), il existe une constante $C$, telle que pour tous $\phi$ et $\psi$ qui vérifient l'hypothèse (1) on ait

\[\begin{aligned}
\displaystyle\int_{i\mathcal{A}_{L,F}^*\times i\mathcal{A}_{S,F}^*}m(\rho_\lambda)m(\tau_\mu)|\phi(\lambda)|^2 |\psi(\mu)|^2 \mathcal{L}_{\pi_\lambda,\sigma_\mu}(\epsilon\otimes e,\epsilon \otimes e) d\lambda d\mu \\
 \leqslant C\int_{i\mathcal{A}_{L,F}^*\times i\mathcal{A}_{S,F}^*} |\phi(\lambda)|^2|\psi(\mu)|^2m(\rho_\lambda) m(\tau_\mu)\sum_{e'\in\mathcal{B}^{K_1},\epsilon'\in\mathcal{B}^{K_2}} |\mathcal{L}_{\pi_\lambda,\sigma_\mu}(\epsilon & \otimes e',\epsilon'\otimes e)|^2 d\lambda d\mu
\end{aligned}\]

Pour presque tout $\lambda\in i\mathcal{A}_{L,F}^*$ il existe un voisinage $\omega_\lambda$ de $\lambda$ tel que $(w^{-1}\omega_\lambda+\mu)\cap\omega_\lambda=\emptyset$ pour tout $(w,\mu)\in \mathcal{E}(\rho)-\{(Id,0)\}$. On a le même résultat sur $i\mathcal{A}_{S,F}^*$. On a donc l'inégalité

$$\mathcal{L}_{\pi_\lambda,\sigma_\mu}(\epsilon\otimes e,\epsilon \otimes e)\leqslant C\displaystyle\sum_{e'\in\mathcal{B}^{K_1},\epsilon'\in\mathcal{B}^{K_2}} |\mathcal{L}_{\pi_\lambda,\sigma_\mu}(\epsilon\otimes e',\epsilon'\otimes e)|^2$$
pour tout $(\lambda,\mu)\in i\mathcal{A}_{L,F}^*\times i\mathcal{A}_{S,F}^*$. Puisque $\epsilon$ et $e$ ont été choisis tels que l'ordre d'annulation en $(0,0)$ de $(\lambda,\mu)\mapsto \mathcal{L}_{\pi_\lambda,\sigma_\mu}(\epsilon\otimes e,\epsilon \otimes e)$ soit minimal, on déduit de l'inégalité précédente que $\mathcal{L}_{\pi_0,\sigma_0}(\epsilon\otimes e,\epsilon \otimes e)\neq 0$ $\blacksquare$

\vspace{4mm}

\begin{prop}

i) Soient $\lambda_0\in i\mathcal{A}_{L,F}^*$ et $\mu_0\in i\mathcal{A}_{S,F}^*$ tels que $\mathcal{L}_{\pi_{\lambda_0},\sigma_{\mu_0}}\neq 0$. Il existe alors un unique couple $(\pi',\sigma')$ constitué de sous-représentations irréductibles de $\pi_{\lambda_0}$ et $\sigma_{\mu_0}$ respectivement tel que $\mathcal{L}_{\pi',\sigma'}\neq 0$. \\

ii) Soit $\sigma\in Temp(H)$ et supposons que $\mathcal{L}_{\pi_\lambda,\sigma}$ pour $\lambda\in i\mathcal{A}_{L,F}^*$ ne soit pas identiquement nulle, alors

\vspace{2mm}

  a) Pour tout $\lambda\in i\mathcal{A}_{L,F}^*$, $\mathcal{L}_{\pi_\lambda,\sigma}$ est non nul; \\ 
   
  b) Il existe des familles finies $(\epsilon_i)_{i=1,\ldots,n}$ et $(\epsilon'_i)_{i=1,\ldots,n}$ d'éléments de $E_\sigma$, des familles finies $(e_i)_{i=1,\ldots,n}$ et $(e'_i)_{i=1,\ldots,n}$ d'éléments de $\mathcal{K}^G_{Q,\tau}$ et une famille finie $(\varphi_i)_{i=1,\ldots,n}$ de fonctions holomorphes sur un voisinage de $i\mathcal{A}_{L,F}^*$ dans $\mathcal{A}_{L,\mathbb{C}}^*/i\mathcal{A}_{L,F}^\vee$, telles que pour tout $\lambda\in i\mathcal{A}_{L,F}^*$ on ait

$$\displaystyle\sum_{i=1}^n \varphi_i(\lambda) \mathcal{L}_{\pi_\lambda,\sigma}(\epsilon'_i\otimes e'_i,\epsilon_i\otimes e_i)=1$$
\end{prop}

\ul{Preuve}: i) On ne perd rien à supposer que $(\lambda_0,\mu_0)=(0,0)$.
\begin{center}
(1) Pour tous $e',e\in\mathcal{K}^G_{Q,\rho}$ et pour tous $\epsilon',\epsilon\in \mathcal{K}^H_{R,\tau}$ on a
$$|\mathcal{L}_{\pi_0,\sigma_0}(e'\otimes \epsilon',e\otimes \epsilon)|^2=|\mathcal{L}_{\pi_0,\sigma_0}(e'\otimes \epsilon,e'\otimes \epsilon)||\mathcal{L}_{\pi_0,\sigma_0}(e\otimes \epsilon',e\otimes \epsilon')|$$
\end{center}

\noindent En effet, d'après le 1) de la proposition précédente, il suffit de vérifier la même égalité en remplaçant $\pi_0$ et $\sigma_0$ par $\pi_\lambda$ et $\sigma_\mu$ pour $(\lambda,\mu)$ dans un ouvert dense de $i\mathcal{A}_{L,F}^*\times i\mathcal{A}_{S,F}^*$. On peut donc supposer $\pi_\lambda$ et $\sigma_\mu$ irréductibles. Mais alors $dim Hom_{H,\xi}(\pi_\lambda,\sigma_\mu)\leqslant 1$, donc il existe $l\in Hom_{H,\xi}(\pi_\lambda,\sigma_\mu)$ et $c\in\mathbb{C}$ tels que

$$\mathcal{L}_{\pi_\lambda,\sigma_\mu}(e'\otimes \epsilon',e\otimes \epsilon)=c(l(e'),\epsilon)(\epsilon',l(e))$$

\noindent Pour tous $e',e\in\mathcal{K}^G_{Q,\rho}$ et pour tous $\epsilon',\epsilon\in \mathcal{K}^H_{R,\tau}$. L'égalité (1) est alors facile à vérifier. L'existence d'un couple comme dans l'énoncé est évidente. Supposons donc qu'il existe deux tels couples $(\pi',\sigma')$ et $(\pi'',\sigma'')$. D'après (1), on peut trouver $e_1\in E_{\pi'}\subset \mathcal{K}^G_{Q,\rho}$, $e_2\in E_{\pi''}\subset \mathcal{K}^G_{Q,\rho}$, $\epsilon_1\in E_{\sigma'}\subset \mathcal{K}^H_{R,\tau}$ et $\epsilon_2\in E_{\sigma''}\subset \mathcal{K}^H_{R,\tau}$ tels que

$$\mathcal{L}_{\pi_0,\sigma_0}(e_1\otimes \epsilon_1,e_1\otimes \epsilon_1)\neq 0$$

et

$$\mathcal{L}_{\pi_0,\sigma_0}(e_2\otimes \epsilon_2,e_2\otimes \epsilon_2)\neq 0$$

D'après (1) on a aussi $\mathcal{L}_{\pi_0,\sigma_0}(e_1\otimes \epsilon_2,e_2\otimes \epsilon_1)\neq 0$. On a $\pi'\neq \pi''$ ou $\sigma'\neq\sigma''$. Dans le premier cas $E_{\pi'}$ et $E_{\pi''}$ sont orthogonaux pour le produit scalaire et dans le deuxième cas $E_{\sigma'}$ et $E_{\sigma''}$ sont orthogonaux pour le produit scalaire. Dans les deux cas on a $\mathcal{L}_{\pi_0,\sigma_0}(e_1\otimes \epsilon_2,e_2\otimes \epsilon_1)= 0$. On aboutit donc à une contradiction. \\

ii) On peut supposer que $\sigma$ est une sous-représentation de $\sigma_0$. La fonction $\lambda\mapsto \mathcal{L}_{\pi_\lambda,\sigma}$ est non nulle sur un ouvert non vide de $i\mathcal{A}_{L,F}^*$. D'après le i), cela implique que pour toute autre sous-représentation $\sigma'$ de $\sigma_0$ la fonction $\lambda\mapsto \mathcal{L}_{\pi_\lambda,\sigma'}$ s'annule sur un ouvert non vide. Comme elle est analytique d'après la proposition précédente, elle est donc identiquement nulle. Par conséquent d'après le ii) de la proposition précédente la fonction $\lambda\mapsto \mathcal{L}_{\pi_\lambda,\sigma}$ n'est jamais nulle. Cela prouve a). \\

 On peut donc certainement trouver des familles finies $(\epsilon_i)_{i=1,\ldots,n}$ et $(\epsilon'_i)_{i=1,\ldots,n}$ d'éléments de $E_\sigma$, et des familles finies $(e_i)_{i=1,\ldots,n}$ et $(e'_i)_{i=1,\ldots,n}$ d'éléments de $\mathcal{K}^G_{Q,\tau}$ telles que pour tout $\lambda\in i\mathcal{A}_{L,F}^*$ la famille

$$\left(\mathcal{L}_{\pi_\lambda,\sigma}(\epsilon'_i\otimes e'_i,\epsilon_i\otimes e_i)\right)_{i=1,\ldots,n}$$

\noindent contienne un élément non nul. Pour $i=1,\ldots,n$ les fonctions $\psi_i: \lambda\mapsto \mathcal{L}_{\pi_\lambda,\sigma}(\epsilon'_i\otimes e'_i,\epsilon_i\otimes e_i)$ admettent un prolongement holomorphe dans un voisinage de $i\mathcal{A}_{L,F}^*$ dans $\mathcal{A}_{L,\mathbb{C}}^*/i\mathcal{A}_{L,F}^\vee$ que l'on note aussi $\psi_i$. Les fonctions $\lambda\mapsto \overline{\psi_i(\lambda)}=\overline{\psi_i(-\overline{\lambda})}$ pour $\lambda\in i\mathcal{A}_L^*$ admettent donc aussi un prolongement holomorphe dans un voisinage de $i\mathcal{A}_{L,F}^*$. Pour obtenir le b), il suffit de prendre pour $i=1,\ldots,n$

$$\displaystyle \phi_i(\lambda)=\overline{\psi_i(\lambda)} . \left(\sum_{i=1}^n |\psi_i(\lambda)|^2\right)^{-1} $$

$\blacksquare$

\subsection{Tout entrelacement est tempéré}

\begin{theo}
Soient $\pi\in Temp(G)$ et $\sigma\in Temp(H)$ alors $\mathcal{L}_{\pi,\sigma}\neq0$ si et seulement si $Hom_{H,\xi}(\pi,\sigma)\neq 0$.
\end{theo}

\ul{Preuve}: D'après une remarque déjà faite, il suffit d'établir que si $Hom_{H,\xi}(\pi,\sigma)\neq 0$ alors $\mathcal{L}_{\pi,\sigma}\neq0$. Soit donc $l\in Hom_{H,\xi}(\pi,\sigma)$ non nul. Il existe des sous-groupes paraboliques semistandards $Q=LU_Q$ et $R=SU_S$ de $G$ et $H$ respectivement et $\rho$ et $\tau$ des représentations irréductibles de la série discrète de $L(F)$ et $S(F)$ tel qu'avec les notations du paragraphe 14.2, $\pi$ et $\sigma$ soient des sous-représentations de $\pi_0$ et $\sigma_0$. Soient $\phi \in C^\infty(i\mathcal{A}_{L,F}^*)$ et $e,e'\in \mathcal{K}^G_{Q,\rho}$. On pose $f=f_{e,e',\phi}$. Soient $\epsilon\in E_\sigma$ et $e_0\in E_\pi$.

\begin{center}
(1) l'intégrale $I(\epsilon,e_0,f)=\displaystyle\int_{G(F)} (\epsilon,l(\pi(g)e_0)) f(g)dg$ converge absolument.
\end{center}

En effet, on décompose l'intégrale en

$$\displaystyle\int_{U(F)\times M(F)\times K} |(\epsilon,l(\pi(umk)e_0))| |f(umk)| \delta_{P}(m)^{-1} dkdmdu$$

Puisque $f$ est invariante à droite par un sous-groupe compact-ouvert et que le stabilisateur de $e_0$ est ouvert, on peut oublier l'intégrale sur $K$. On a $|(\epsilon,l(\pi(ug)e_0))|=|(\epsilon,l(\pi(g)e_0))|$ pour tout $u\in U(F)$ et $g\in G(F)$. La fonction $f$ appartient à $\mathcal{S}(G(F))$, d'après la proposition II.4.5 de [W3], on a donc pour tout $d>0$

$$\displaystyle\int_{U(F)} |f(um)| du<< \delta_{P}(m)^{1/2} \Xi^M(m) \sigma(m)^{-d}$$
pour tout $m\in M(F)$. Il suffit donc de montrer que pour $d$ assez grand l'intégrale suivante converge

$$\displaystyle\int_{M(F)} |(\epsilon,l(\pi(m)e_0))|\Xi^M(m)\delta_P(m)^{-1/2} \sigma(m)^{-d} dm$$
Appliquons le lemme 14.1.1 à des sous-groupes ouverts $\Gamma_1\subset K$ et $\Gamma_2\subset K_H$ qui laissent stables $e_0$ et $\epsilon$. On en déduit l'existence d'un compact $C_M\subset M(F)$ et de sous-groupes ouverts $K_1\subset K$, $K_2\subset K_H$ qui vérifient les conclusions de ce lemme. On peut toujours supposer que $A_{0,G}^1C_M=C_M$. Alors l'intégrale précédente est à support dans $C_H\Lambda_H^+\Lambda_G^+C_M$. D'après le 2- du lemme 14.1.1 on a une majoration $|(\epsilon,l(\pi(m)e))|<< \Xi^G(a_G)\Xi^H(a_H)$ pour tout $m=k_Ha_Ha_Gk_M\in C_H\Lambda_H^+\Lambda_G^+C_M$. Pour $a_G\in A_{0,G}(F)$, on note $(a_G)_0$ l'unique élément de $A_0(F)$ tel que $a_G(a_G)_0^{-1}\in A(F)$. On a alors la majoration $\sigma(ha_G)>>\sigma(a_G(a_G)_0^{-1})^{1/2}\sigma(h(a_G)_0)^{1/2}$ pour tous $a_G\in A_{0,G}(F)$, $h\in H(F)$. D'après le 3- du lemme 11.0.4, on a donc la majoration

$$\sigma(k_Ha_Ha_Gk_M)>>\sigma((a_G)_0)^{1/2}\sigma(a_G(a_G)_0^{-1})^{1/2}>>\sigma(a_G)^{1/2}$$
pour tous $a_G\in A_{0,G}^+$, $a_H\in A_H^+$, $k_H\in C_H$ et $k_M\in C_M$. D'après le lemme 11.0.4 on a aussi

\[\begin{aligned}
mes(C_Ha_Ha_GC_M)=mes(C_Ha_H(a_G)_0C_M)<<\delta_{P_H}(a_H)^{-1} & \delta_{P_0}((a_G)_0)^{-1} \\
 & =\delta_{P_H}(a_H)^{-1}\delta_{P_{0,G}\cap M}(a_G)^{-1}
\end{aligned}\]

\noindent pour tous $a_G\in A_{0,G}^+$, $a_H\in A_H^+$. Par conséquent, on a

\[\begin{aligned}
\displaystyle\int_{M(F)} & |(\epsilon,l(\pi(m)e_0))| \Xi^M(m)\delta_P(m)^{-1/2} \sigma(m)^{-d} dm  \\
 & << \displaystyle\sum_{a_G\in \Lambda_G^+,a_H\in \Lambda_H^+}\delta_{P_{0,G}\cap M}(a_G)^{-1} \delta_{P_H}(a_H)^{-1} \Xi^G(a_G) \Xi^H(a_H)\Xi^M(a_Ha_G) \delta_P(a_G)^{-1/2}\sigma(a_G)^{-d/2}
\end{aligned}\]

D'après le lemme 11.0.4, on a

$$\Xi^M(a_Ha_G)=\Xi^{G_0}(a_H(a_G)_0)<<\Xi^{G_0}(a_H)\Xi^{G_0}((a_G)_0)=\Xi^{G_0}(a_H)\Xi^M(a_G)$$
pour tout $a_H\in \Lambda_H^+$ et pour tout $a_G\in \Lambda_G^+$. La somme ci-dessus est donc essentiellement majorée par la produit de

$$\displaystyle\sum_{a_H\in \Lambda_H^+} \delta_{P_H}(a_H)^{-1}\Xi^{G_0}(a_H)\Xi^H(a_H)$$
et de

$$\displaystyle\sum_{a_G\in \Lambda_G^+} \delta_{P_{0,G}\cap M}(a_G)^{-1}\delta_P(a_G)^{-1/2} \Xi^G(a_G)\Xi^M(a_G)\sigma(a_G)^{-d/2}$$

D'après le lemme 12.0.5 et la majoration $\Xi^H(a_H)<<\delta_{P_H}(a_H)^{1/2}\sigma(a_H)^{d'}$ pour un certain réel $d'$, la première somme converge absolument. D'après le lemme II.1.1 de [W3] il existe des réels $d_1,d_2$ tels qu'on ait des majorations $\Xi^G(a_G)<<\delta_{P_{0,G}}(a_G)^{1/2}\sigma(a_G)^{d_1}$ et $\Xi^M(a_G)<<\delta_{P_{0,G}\cap M}(a_G)^{1/2} \sigma(a_G)^{d_2}$ pour tout $a_G\in\Lambda_G^+$. Puisque $\delta_{P_{0,G}}=\delta_{P_{0,G}\cap M}\delta_P$, la deuxième somme est essentiellement majorée par

$$\displaystyle\sum_{a_G\in\Lambda_G^+} \sigma(a_G)^{d_1+d_2-d/2}$$
qui est une série absolument convergente pour $d$ assez grand.

\vspace{4mm}

On peut calculer $I(\epsilon,e_0,f)$ de deux façons. La première consiste à choisir une suite exhaustive de sous-ensembles compact-ouverts $K$-biinvariants $(\Omega_n)$ de $G(F)$ et de la calculer comme la limite de l'intégrale restreinte à $\Omega_n$ lorsque $n$ tend vers l'infini. On trouve alors que

\begin{center}
(2) $I(\epsilon,e_0,f)=(\epsilon,l(\pi(f)e_0))$
\end{center}

On peut aussi écrire

$$I(\epsilon,e_0,f)=\displaystyle\int_{U(F)}\int_{G_0(F)}\int_K\int_{A(F)} (\epsilon,l(\pi(ag_0k)e_0)) f(uag_0k) \delta_P(a)^{-1} \xi(u)dadkdg_0du$$

Puisque $f$ est biinvariante par un sous-groupe compact-ouvert de $A(F)$, il existe un entier $c$ tel que l'intégrale intérieure soit nulle pour $u\in U(F)-U(F)_c$ (cf preuve du lemme 14.0.14). On a alors

\[\begin{aligned}
I(\epsilon,e_0,f) & = \\
 & \displaystyle\int_K\int_{A(F)}\int_{H(F)\backslash G_0(F)}\int_{H(F)U(F)_c} (\epsilon,l(\pi(hag_0k)e_0)) f(uhag_0k) \delta_P(a)^{-1} \xi(u)dudhdg_0dadk
\end{aligned}\]

\noindent et

\[\begin{aligned}
\int_{H(F)U(F)_c} (\epsilon,l(\pi(hag_0k)e_0)) f(uhag_0k) \xi(u)dudh & \\
=\int_{H(F)U(F)_c}\int_{i\mathcal{A}_{L,F}^*} \phi(\lambda)m(\rho_\lambda) (\pi_\lambda(uhag_0k)e',e) & (\epsilon,\sigma(h)l(\pi(ag_0k)e_0))\xi(u)d\lambda dudh
\end{aligned}\]

Comme on l'a déjà vérifié dans la preuve précédente, ce genre d'intégrale est absolument convergente. En permutant les deux intégrales, on voit apparaître le terme \\
\noindent $\mathcal{L}_{\pi_\lambda,\sigma,c}(\epsilon\otimes \pi_\lambda(ag_0k)e',l(\pi(ag_0k)e_0)\otimes e)$. On peut prendre $c$ aussi grand que l'on veut et comme déjà expliqué dans la preuve de la proposition 14.2.1, on peut choisir $c$ tel que pour tous $a\in A(F),g_0\in G_0(F)$, $k\in K$ et pour tout $\lambda\in i\mathcal{A}_{L,F}^*$, on ait
$$\mathcal{L}_{\pi_\lambda,\sigma,c}(\epsilon\otimes \pi_\lambda(ag_0k)e',l(\pi(ag_0k)e_0)\otimes e)=\mathcal{L}_{\pi_\lambda,\sigma}(\epsilon\otimes \pi_\lambda(ag_0k)e',l(\pi(ag_0k)e_0)\otimes e)$$
(Car $A$ commute à $AG_0=M$). On en déduit que

\[\begin{aligned}
\mbox{(3)} \;\; I(\epsilon,e_0,f)=\displaystyle\int_K\int_{A(F)}\int_{H(F)\backslash G_0(F)}\int_{i\mathcal{A}_{L,F}^*} &  \phi(\lambda)m(\rho_\lambda)\mathcal{L}_{\pi_\lambda,\sigma}(\epsilon\otimes \pi_\lambda(ag_0k)e',l(\pi(ag_0k)e_0)\otimes e) \\
 & \delta_P(a)^{-1}dudhdg_0dadk
\end{aligned}\]

Choisissons $\epsilon \in E_\sigma$ et $e_0\in E_\pi$ tels que $(l(e_0),\epsilon)\neq 0$. On applique les calculs précédents au cas où $e=e'=e_0$ et $\phi$ à support dans un voisinage $\omega$ assez petit de $0$ tel que $\phi(0)\neq 0$. Alors $\pi(f)e_0$ est un multiple non nul de $e_0$. D'après (2) $I(\epsilon,e_0,f)$ n'est pas nulle. Alors (3) implique l'existence de $\lambda$ tel que $\mathcal{L}_{\pi_\lambda,\sigma}$ ne soit pas nul. D'après la proposition 14.2.1, $\mathcal{L}_{\pi_0,\sigma_0}$ n'est pas nul non plus. Si $\pi_0$ et $\sigma_0$ sont irréductibles alors $\pi=\pi_0$, $\sigma=\sigma_0$ et on a ce que l'on voulait. Sinon d'après la proposition 14.2.2 i), il existe un unique couple $(\pi',\sigma')$ de sous-représentations irréductibles de $\pi_0$ et $\sigma_0$ tel que $\mathcal{L}_{\pi',\sigma'}$ soit non nul. Soient alors $e_1\in E_{\pi'}\subset \mathcal{K}^G_{Q,\rho}$ et $\epsilon_1\in E_{\sigma'}\subset \mathcal{K}^H_{R,\tau}$ tels que $\mathcal{L}_{\pi_0,\sigma_0}(\epsilon_1\otimes e_1,\epsilon_1\otimes e_1)\neq 0$. Quitte à restreindre $\omega$, on peut supposer que pour tout $\lambda\in \omega$, $\mathcal{L}_{\pi_\lambda,\sigma_0}(\epsilon_1\otimes e_1,\epsilon_1\otimes e_1)$ est non nul. Soit $\phi'$ la fonction à support dans $\omega$ définie par

$$\phi'(\lambda)=\phi(\lambda) \mathcal{L}_{\pi_\lambda,\sigma_0}(\epsilon\otimes e_1,\epsilon_1\otimes e_0)\mathcal{L}_{\pi_\lambda,\sigma_0}(\epsilon_1\otimes e_1,\epsilon_1\otimes e_1)^{-1}$$

Posons $f'=f_{e_1,e_0,\phi'}$. On a alors l'égalité

\begin{center}
(4) $I(\epsilon,e_0,f)=I(\epsilon_1,e_0,f')$
\end{center}

En effet, d'après (3) il suffit de vérifier que pour tout $\lambda\in i\mathcal{A}_{L,F}^*$ et $g\in G(F)$, on a l'égalité

\[\begin{aligned}
\mathcal{L}_{\pi_\lambda,\sigma_0}(\epsilon\otimes e_1,\epsilon_1\otimes e_0)\mathcal{L}_{\pi_\lambda,\sigma_0}(\epsilon_1\otimes \pi_\lambda(g)e_0,l(\pi(g)e_0)\otimes e_0) \\
=\mathcal{L}_{\pi_\lambda,\sigma_0}(\epsilon_1\otimes e_1,\epsilon_1\otimes e_1)\mathcal{L}_{\pi_\lambda,\sigma_0}(\epsilon\otimes \pi_\lambda(g)e_0,l(\pi(g)e_0)\otimes e_0)
\end{aligned}\]

Il suffit même de vérifier l'égalité précédente en remplaçant $\sigma_0$ par $\sigma_\mu$ pour $\lambda$ et $\mu$ en position générale. On peut alors supposer $\pi_\lambda$ et $\sigma_\mu$ irréductibles. Alors il existe $\underline{l}\in Hom_{H,\xi}(\pi_\lambda,\sigma_\mu)$ et $c\in \mathbb{C}^\times$ tels que pour tout $e,e'\in E_{\pi_\lambda}$ et $\epsilon,\epsilon'\in E_{\sigma_\mu}$

$$\mathcal{L}_{\pi_\lambda,\sigma_\mu}(\epsilon'\otimes e',\epsilon\otimes e)=c(\underline{l}(e'),\epsilon)(\epsilon',\underline{l}(e))$$

Les deux membres de l'égalité que l'on cherche à prouver valent alors

$$c^2(\underline{l}(e_1),\epsilon_1)(\epsilon,\underline{l}(e_0))(\epsilon_1,\underline{l}(e_1))(\underline{l}(\pi_\lambda(g)e_0),l(\pi(g)e_0))$$

Cela prouve (4). En particulier on a $I(\epsilon_1,e_0,f')\neq 0$. D'après (2) cela implique $\pi(f')e_0\neq 0$ donc notamment $\phi'(0)\neq 0$. On a par conséquent $\mathcal{L}_{\pi_0,\sigma_0}(\epsilon\otimes e_1,\epsilon_1\otimes e_0)\neq 0$. Mais si $\pi'\neq \pi$ ou $\sigma'\neq \sigma$, il est facile de vérifier que $\mathcal{L}_{\pi_0,\sigma_0}(\epsilon\otimes e_1,\epsilon_1\otimes e_0)=0$. On a donc $\pi'=\pi$ et $\sigma'=\sigma$ et le résultat est démontré. $\blacksquare$

\section{Induction et multiplicités}

Dans cette section on se propose de montrer qu'en un certain sens la multiplicité $m(\pi,\sigma)$ est compatible à l'induction. Pour cela il est commode d'étendre un peu la définition de la multiplicité $m(\pi,\sigma)$. Soient $(V,h)$ et $(V',h')$ deux espaces hermitiens de dimension $d$ et $d'$, on dira qu'ils sont compatibles si $d$ et $d'$ sont de parités différentes et si le plus petit des deux peut s'injecter dans le plus gros. Supposons que $(V,h)$ et $(V',h')$ soient compatibles. Notons $G$ resp. $G'$ les groupes unitaires de $V$ rep. $V'$. Fixons une injection de $(V',h')$ dans $(V,h)$ ou de $(V,h)$ dans $(V',h')$ (suivant que $d>d'$ ou $d'>d$). On peut alors toujours trouver une décomposition $V=V'\oplus^\perp D\oplus^\perp (Z_+\oplus Z_-)$ ou $V'=V\oplus^\perp D\oplus^\perp (Z_+\oplus Z_-)$ où $D$ est une droite et $Z_+$, $Z_-$ sont des sous-espaces totalement isotropes. Soit $P=MU$ le sous-groupe parabolique de $G$ ou $G'$ (suivant que $d>d'$ ou $d'>d$) qui fixe un drapeau complet de sous-espaces de $Z_+$. On construit comme dans la section 4 un caractère $\xi$ de $U(F)$. Soient $\pi$ et $\pi'$ des représentations irréductibles lisses de $G(F)$ et $G'(F)$ respectivement. On définit la multiplicité $m'(\pi,\pi')$ comme étant $m(\pi,\pi')$ si $d>d'$ et $m(\pi',\pi)$ si $d'>d$, cette multiplicité ne dépend pas des divers choix effectués. Dorénavant on notera aussi $m(\pi,\sigma)$ cette multiplicité généralisée. \\

\begin{lem}
Soient $(V,h)$ et $(V',h')$ deux espaces hermitiens compatibles de groupes unitaires $G$ et $G'$ respectivement. Pour $\pi\in Irr(G)$ et $\pi'\in Irr(G')$ on a
$$m(\pi^\vee,\pi'^\vee)=m(\pi,\pi')$$
\end{lem}

\ul{Preuve}: On ne perd rien à supposer que $V'\subset V$. Soit $\delta$ un automorphisme $F$-linéaire de $V$ tel que $h(\delta u,\delta v)=h(v,u)$ pour tout $u,v\in V$. On peut clairement trouver un tel élément tel que $\delta(V')=V'$. D'après [MVW] $\pi^\vee$ est isomorphe à $\pi^\delta$ où $\pi^\delta(x)=\pi(\delta x\delta^{-1})$ et de la même façon $\pi'$ est isomorphe à $\pi'^{\delta'}$ où $\delta'$ est la restriction de $\delta$ à $V'$. Modulo ces isomorphismes on a les égalités $Hom_{H,\xi}(\pi^\vee,\pi'^\vee)=Hom_{H,\xi}(\pi^\delta,\pi'^{\delta'})=Hom_{H,\xi}(\pi,\pi')$ et donc l'égalité $m(\pi^\vee,\pi'^\vee)=m(\pi,\pi')$ $\blacksquare$

\vspace{2mm}

Soient $(V,h)$ un espace hermitien et $W$ un sous-espace hermitien de $V$ de sorte que $V$ et $W$ soient compatibles. Fixons une décomposition orthogonale $V=(Z_+\oplus Z_-)\oplus^\perp D\oplus^\perp W$ où $Z_+,Z_-$ sont des sous-espaces totalement isotropes et $D$ est une droite. Fixons aussi un générateur $v_0$ de $D$ et posons $r=dim(Z_+)=dim(Z_-)$. Cette situation sera conservée jusqu'en 15.3 inclus. Soit $k\geqslant 1$ un entier et introduisons quelques notations relatives à des sous-groupes de $GL_k$. Pour $i=1,\ldots,k$, on note $Q_{k-i,i}$ le sous-groupe parabolique standard de $GL_k$ de composante de Levi standard $GL_{k-i}\times GL_1\times\ldots\times GL_1$ et $U_{k-i,i}$ son radical unipotent. On définit $P_{k-i,i}$ comme le sous-groupe des éléments de $Q_{k-i,i}$ dont les projections sur les $i$ facteurs $GL_1$ sont triviales. En particulier $P_{k-1,1}$ est le sous-groupe mirabolique de $GL_k$. On définit un caractère $\psi_i$ de $P_{k-i,i}(E)$ par

$$\displaystyle \psi_i(p)=\psi_E(\sum_{j=k-i+1}^{k-1} p_{j,j+1})$$

\noindent pour tout $p\in P_{k-i,i}(E)$ où $p_{j,l}$ désigne les coefficients de $p$.

\subsection{Induction de $\pi$ et multiplicité I}

Soit $Y_+$ un sous-espace totalement isotrope de dimension $k$ de $V$ et notons $Q$ son stabilisateur dans $G$. C'est un sous-groupe parabolique de $G$, on note $N$ son radical unipotent. Fixons des sous-espaces $\tilde{V}$ et $Y_-$ de sorte que $V=(Y_+\oplus Y_-)\oplus^\perp \tilde{V}$ et notons $L$ la composante de Levi de $Q$ qui stabilise $\tilde{V}$ et $Y_-$. On a alors un isomorphisme $L\simeq GL(Y_+)\times \tilde{G}$. Soient $\tilde{\pi}\in Temp(\tilde{G})$, $\pi_+\in Temp(GL(Y_+))$, $\sigma\in Temp(H)$ et posons $\pi=i^G_Q(\pi_+\otimes \tilde{\pi})$.

\begin{prop}
Supposons $r=0$ alors
$$m(\pi,\sigma)=m(\tilde{\pi},\sigma)$$
\end{prop}

\ul{Preuve}: Pour tout $\lambda\in i\mathbb{R}/(\frac{2i\pi}{log(q)}\mathbb{Z})$, posons $\pi_\lambda=i^G_Q((\pi_+ |det|^\lambda)\otimes \tilde{\pi})$. La représentation $\pi_\lambda$ se réalise par translation à droite sur l'espace $V_{\pi_\lambda}$ des fonctions $\varphi: G(F)\to V_{\pi_+}\otimes V_{\tilde{\pi}}$ lisses vérifiant

$$\varphi(g_+\tilde{g}ng)=|det(g_+)|_E^\lambda \delta_Q(g_+)^{1/2}(\pi_+(g_+)\otimes \tilde{\pi}(\tilde{g}))\varphi(g)$$

\noindent pour tous $g_+\in GL(Z_+), \tilde{g}\in\tilde{G}(F), n\in N(F), g\in G(F)$. On a $m(\pi_\lambda,\sigma)=m(\pi,\sigma)$ pour tout $\lambda$ d'après la proposition 14.2.2 ii) a) et le théorème 14.3.1. L'espace quotient $Q(F)\backslash G(F)$ paramétrise les sous-espaces totalement isotropes de dimension $k$ de $V$. L'action de $H(F)$ sur $Q(F)\backslash G(F)$ admet deux orbites: l'une ouverte $\mathcal{U}$ correspond à l'ensemble des sous-espaces totalement isotropes $Y'$ de dimension $k$ tels que $dim(Y'\cap W)=k-1$, l'autre fermée $\mathcal{Y}$ correspond à l'ensemble des sous-espaces totalement isotropes $Y'$ de dimension $k$ inclus dans $W$. Définissons $V_{\pi_\lambda,\mathcal{U}}$ comme le sous-espace de $V_{\pi_\lambda}$ des fonctions à support dans $\mathcal{U}$ et notons $V_{\pi_\lambda,\mathcal{Y}}=V_{\pi_\lambda}/V_{\pi_\lambda,\mathcal{U}}$. On a alors la suite exacte de $H(F)$-représentations
$$0\to V_{\pi_\lambda,\mathcal{U}}\to V_{\pi_\lambda}\to V_{\pi_\lambda,\mathcal{Y}}\to 0$$

\begin{lem}
On a $Hom_H(V_{\pi_\lambda,\mathcal{Y}},\sigma)=0$ et $Ext^1(V_{\pi_\lambda,\mathcal{Y}},\sigma)=0$.
\end{lem}

\ul{Preuve}: On peut toujours supposer que $Y_+\subset W$. Alors $Q\cap H=Q_H=L_HN_H$ est un sous-groupe parabolique de $H$ et on a un isomorphisme de $H(F)$-représentations

$$V_{\pi_\lambda,\mathcal{Y}}\simeq i^H_{Q_H}(\delta_Q^{1/2}\delta_{Q_H}^{-1/2} (\pi_+|det|^\lambda\otimes\tilde{\pi})_{|Q_H})$$

\noindent D'après le second théorème d'adjonction de Bernstein, on a

$$Hom_H(V_{\pi_\lambda,\mathcal{Y}},\sigma)=Hom_{L_H}(\delta_Q^{1/2}\delta_{Q_H}^{-1/2} (\pi_+|det|^\lambda\otimes\tilde{\pi}),(\sigma)_{\overline{Q}_H})$$

\noindent et

$$Ext^1_H(V_{\pi_\lambda,\mathcal{Y}},\sigma)=Ext^1_{L_H}(\delta_Q^{1/2}\delta_{Q_H}^{-1/2} (\pi_+|det|^\lambda\otimes\tilde{\pi}),(\sigma)_{\overline{Q}_H})$$

\noindent La représentation $(\sigma)_{\overline{Q}_H}$ est de longueur finie et puisque $\sigma$ est tempérée, les caractères centraux de tous ses sous-quotients ont leur partie réelle dans un certain cône. La représentation $\delta_Q^{1/2}\delta_{Q_H}^{-1/2} (\pi_+|det|^\lambda\otimes\tilde{\pi})$ est irréductible et son caractère central n'est pas dans ce cône. On en déduit la nullité des deux espaces précédents $\blacksquare$

D'après le lemme on a $Hom_H(\pi_\lambda,\sigma)=Hom_H(V_{\pi_\lambda,\mathcal{U}},\sigma)$. Quitte à conjuguer, on peut toujours supposer que $dim(Y_+\cap W)=k-1$ et donc $\mathcal{U}=Q(F)\backslash Q(F)H(F)$. La restriction à $H(F)$ de $V_{\pi_\lambda,\mathcal{U}}$ est alors une induite compacte à partir du sous-groupe $H(F)\cap Q(F)$. Décrivons plus en détail ce sous-groupe. On peut fixer une décomposition

$$W=(Y'_+\oplus Y'_-)\oplus^\perp D'\oplus^\perp \tilde{V}$$

\noindent où $Y'_+=Y_+\cap W$ et $D'$ est une droite. Soit $\tilde{G}'$ le groupe unitaire de $D'\oplus \tilde{V}$ et $Q'$ le sous-groupe parabolique de $H$ des éléments qui stabilisent $Y'_+$. On a alors une décomposition $Q'=L'N'$ où $L'$ est le Levi qui stabilise $D'\oplus \tilde{V}$ et $Y'_-$, ce Levi s'identifie donc naturellement à $GL(Y_+')\times \tilde{G}'$. On a alors $Q\cap H=(GL(Y'_+)\times \tilde{G})N'$. La représentation $\delta_Q^{1/2}(\pi_+|det|^\lambda\otimes\tilde{\pi})_{|Q(F)\cap H(F)}$ n'est pas triviale sur $N'(F)$ mais elle est triviale sur $N'_\sharp(F)$ où $N'_\sharp$ est le sous-groupe des éléments de $N'$ qui agissent comme l'identité sur $D'$. On a $N'_\sharp=N\cap H$ et c'est un sous-groupe distingué de $Q\cap H$. Le quotient $(Q(F)\cap H(F))/N'_\sharp(F)$ est naturellement un sous-groupe de $L(F)=GL(Y_+)\times \tilde{G}(F)$. On peut fixer une base de $Y_+$ de sorte que l'isomorphisme déduit $GL(Y_+)\simeq GL_k(E)$ identifie $(Q(F)\cap H(F))/N'_\sharp(F)$ avec $P_{k-1,1}(E)\times \tilde{G}(F)$. On a alors un isomorphisme de $H(F)$-représentations

$$V_{\pi_\lambda,\mathcal{U}}\simeq c-ind_{(P_{k-1,1}(E)\times \tilde{G}(F))N'_\sharp(F)}^{H(F)}((\delta_Q^{1/2}\pi_+|det|^\lambda)_{|P_{k-1,1}(E)}\otimes \tilde{\pi})$$

\noindent D'après [BZ] 3.5, la représentation $(\delta_Q^{1/2}\pi_+|det|^\lambda)_{|P_{k-1,1}(E)}$ possède une filtration

$$\{0\}=\tau_{k+1,\lambda}\subset \tau_{k,\lambda}\subset\ldots\subset \tau_{1,\lambda}=(\delta_Q^{1/2}\pi_+|det|^\lambda)_{|P_{k-1,1}(E)}$$

\noindent Les quotients de cette filtration vérifiant

$$\tau_{i,\lambda}/\tau_{i+1,\lambda}\simeq c-ind_{P_{k-i,i}(E)}^{P_{k-1,1}(E)}(\Delta^i(\delta_Q^{1/2}\pi_+|det|^\lambda)\otimes \psi_i)$$

\noindent L'induction à support compact étant un foncteur exact, on en déduit une filtration

$$\{0\}=\mu_{k+1,\lambda}\subset \mu_{k,\lambda}\subset \ldots\subset \mu_{1,\lambda}=V_{\pi_\lambda,\mathcal{U}}$$

\noindent où pour $i=1,\ldots,k$ on a

$$\mu_{i,\lambda}/\mu_{i+1,\lambda}\simeq c-ind_{(P_{k-i,i}(E)\times \tilde{G}(F))N'_\sharp(F)}^{H(F)}(\Delta^i(\delta_Q^{1/2}\pi_+|det|^\lambda)\otimes \psi_i\otimes \tilde{\pi})$$

\noindent Pour $i=1,\ldots,k$, soit $Q_{k-i}$ le sous-groupe parabolique de $H$ qui fixe l'espace engendré par les $k-i$ premiers vecteurs de la base fixée de $Y_+$, il admet un Levi naturellement isomorphe à $Res_{E/F}(GL_{k-i})\times \tilde{G}_{k-i}$ où $\tilde{G}_{k-i}$ est un groupe unitaire contenant $\tilde{G}$. La représentation $\Delta^i(\delta_Q^{1/2}\pi_+|det|^\lambda)\otimes \psi_i\otimes \tilde{\pi}$ est triviale sur le radical unipotent de $Q_{k-i}$. Posons

$$\mu'_i=c-ind_{(U_{k-i,i}(E)\tilde{G}(F)N'_\sharp(F))\cap \tilde{G}_{k-i}}^{\tilde{G}_{k-i}}(\psi_i\otimes \tilde{\pi})$$

\noindent On a alors

$$\mu_{i,\lambda}/\mu_{i+1,\lambda}\simeq i_{Q_{k-i}}^H(\Delta^i(\delta_Q^{1/2}\pi_+|det|^\lambda)\otimes \mu'_i)$$

\begin{lem}
Pour $\lambda$ générique, on a pour $i=1,\ldots,k-1$

$$Hom_H(\mu_{i,\lambda}/\mu_{i+1,\lambda},\sigma)=0$$

et

$$Ext^1_H(\mu_{i,\lambda}/\mu_{i+1,\lambda},\sigma)=0$$
\end{lem}

\ul{Preuve}:
On a $\Delta^i(\delta_Q^{1/2}\pi_+|det|^\lambda)=\Delta^i(\delta_Q^{1/2}\pi_+)|det|^\lambda$ et $\Delta^i(\delta_Q^{1/2}\pi_+)$ est de longueur finie. Par conséquent pour $\lambda$ générique tous les sous-quotients de $\mu_{i,\lambda}/\mu_{i+1,\lambda}$ ont des supports cuspidaux différents de celui de $\sigma$ $\blacksquare$

\vspace{3mm}

D'après le lemme on a donc pour $\lambda$ générique $Hom_H(\pi,\sigma)=Hom_H(\mu_{k,\lambda},\sigma)$. D'après [BZ], on a $\Delta^k(\delta_Q^{1/2}\pi_+|det|^\lambda)=1$ (car $\pi_+$ est tempérée donc générique) et $\tilde{G}_{k-i}=Q_{k-i}=H$. Par conséquent $\mu_{k,\lambda}=c-ind_{\tilde{G}P_{0,k}(E)N'_\sharp(F)}^H(\psi_k\otimes \tilde{\pi})$. La contragrédiente de l'induite à supports compacts d'une représentation admissible étant l'induite ordinaire de la contragrédiente, on a par réciprocité de Frobenius

$$Hom_H(\pi,\sigma)=Hom_{\tilde{G}P_{0,k}(E)N'_\sharp(F)}(\sigma^\vee,\psi_k^{-1}\otimes \tilde{\pi}^\vee)$$

Remarquons que $P_{0,k}(E)N'_\sharp(F)$ n'est autre que le radical unipotent d'un sous-groupe parabolique de $H$ de composante de Lévi $Res_{E/F}(GL_1)\times \tilde{G}$ et que $\psi_k^{-1}$ a alors une définition analogue à celle de $\xi$ pour une normalisation convenable. La multiplicité $m(\tilde{\pi},\sigma)=m(\tilde{\pi}^\vee,\sigma^\vee)$ est donc la dimension de l'espace $Hom_{\tilde{G}P_{0,k}(E)N'_\sharp(F)}(\sigma^\vee,\psi_k^{-1}\otimes \tilde{\pi}^\vee)$. D'où le résultat $\blacksquare$

\subsection{Induction de $\sigma$ et multiplicité}

Soit $W=(Y_{+,H}\oplus Y_{-,H})\oplus^\perp \tilde{W}$ une décomposition orthogonale où $Y_{+,H}$ et $Y_{-,H}$ sont totalement isotropes. Notons $Q_H$ le sous-groupe parabolique de $H$ des éléments qui stabilisent $Y_{+,H}$, $N_H$ son radical unipotent et $L_H$ la composante de Levi qui stabilise $Y_{-,H}$. Notons $\tilde{H}$ le groupe unitaire de $\tilde{W}$, on a alors $L_H=GL(Y_{+,H})\times \tilde{H}$. Soient $\pi\in Temp(G)$, $\sigma_+\in Temp(GL(Y_{+,H}))$, $\tilde{\sigma}\in Temp(\tilde{H})$ et posons $\sigma=i^H_{Q_H}(\sigma_+\otimes \tilde{\sigma})$.

\begin{prop}
On a $m(\pi,\sigma)=m(\pi,\tilde{\sigma})$
\end{prop}

\ul{Preuve}: Soit $(D',h_{D'})$ un espace hermitien de dimension 1 engendré par un vecteur $v_0'$ tel que $h_{D'}(v_0')=-h(v_0)$. Posons $V'$ comme étant la somme orthogonale de $V$ et $D'$ et notons $G'$ son groupe unitaire. Soient $Z'_+=Z_+\oplus E(v_0+v_0')$ et $Z'_-=Z_-\oplus E(v_0-v_0')$, on a alors $V'=(Z'_+\oplus Z'_-)\oplus^\perp W$ et $Z'_+,Z'_-$ sont des sous-espaces totalement isotropes de $V'$. Notons $P'$ le sous-groupe parabolique de $G'$ des éléments qui stabilisent $Z'_+$ et $M'$ sa composante de Levi qui stabilise $Z'_-$. On a alors $M'=GL(Z'_+)\times H$. Fixons une représentation tempérée arbitraire $\pi'_+$ de $GL(Z'_+)$. Posons $\sigma'=i_{P'}^{G'}(\pi'_+\otimes \sigma)$. Alors, d'après la proposition 15.1.1 $m(\pi,\sigma)=m(\pi,\sigma')$. Soit $Q'$ le sous-groupe parabolique de $G'$ des éléments qui stabilisent $Z'_+\oplus Y_{+,H}$ et $L'$ sa composante de Levi qui stabilise $Z'_-\oplus Y_{-,H}$. Alors $L'=GL(Z'_+\oplus Y_{+,H})\times \tilde{H}$ et on a un isomorphisme $\sigma'\simeq i^{G'}_{Q'}((\pi'_+\times \sigma_+)\otimes \tilde{\sigma})$. Appliquant une nouvelle fois la proposition 15.1.1 on a $m(\pi,\sigma')=m(\pi,\tilde{\sigma})$ $\blacksquare$

\subsection{Induction de $\pi$ et multiplicité II}

On reprend ici les notations du paragraphe 15.1: on a une décomposition $V=(Y_+\oplus Y_-)\oplus^\perp \tilde{V}$, $\pi_+\in Temp(GL(Y_+))$, $\tilde{\pi}\in Temp(\tilde{G})$, $\pi=i_Q^G(\pi_+\otimes \tilde{\pi})$ et $\sigma\in Temp(H)$.

\begin{prop}
On a $m(\pi,\sigma)=m(\tilde{\pi},\sigma)$.
\end{prop}

\ul{Preuve}: On démontre le résultat par récurrence sur $dim(V)$. Si $dim(V)=1$ il n'y a rien à dire. Supposons donc le résultat établi pour tout les couples $(V',W')$ avec $dim(V')<dim(V)$. Posons $V'=Z_+\oplus Z_-\oplus W$, $G'$ le groupe unitaire de $V'$, $P'$ le sous-groupe parabolique de $G'$ qui préserve $Z_+$ et $M'$ la composante de Levi qui préserve $Z_-$. On a alors un isomorphisme naturel $M'\simeq GL(Z_+)\times H$. Soit $\sigma_+\in Temp(GL(Z_+))$ et posons $\pi'=i^{G'}_{P'}(\sigma_+\otimes \sigma)$. D'après la proposition 15.2.1 on a $m(\pi,\sigma)=m(\pi,\pi')$. D'après la proposition 15.1.1 on a $m(\pi,\pi')=m(\pi',\tilde{\pi})$. Enfin d'après l'hypothèse de récurrence appliquée au couplus $(V',\tilde{V})$, on a $m(\pi',\tilde{\pi})=m(\tilde{\pi},\sigma)$ $\blacksquare$

\section{Le développement spectral}

On reprend les notations de la section 4 :$V,W,P,M,U,\xi,\ldots$. Soit $\sigma\in Temp(H)$ et $f\in C_c^\infty(G(F))$ une fonction très cuspidale. Le but de cette section est de donner une expression spectrale de la limite $\lim\limits_{N\to \infty} J_N(\theta_\sigma,f)$.

\subsection{La formule}

Soit $L$ un Lévi de $G$ et $\mathcal{O}\in \Pi_{ell}(L)$ une $i\mathcal{A}_{L,F}^*$ orbite de représentations irréductibles elliptiques de $L(F)$. On a une décomposition $L\simeq R_{E/F} GL_{n_1}\times \ldots\times R_{E/F} GL_{n_k}\times \tilde{G}$ où $\tilde{G}$ est le groupe unitaire d'un sous-espace hermitien $\tilde{V}$ de $V$. Soit $\pi\in \mathcal{O}$, on a une décomposition analogue de $\pi$ en produit tensoriel $\pi\simeq \pi_1\otimes\ldots\otimes \pi_k\otimes\tilde{\pi}$ où pour $j=1,\ldots,k$, $\pi_j$ est une représentation irréductible de la série discrète de $GL_{n_j}(E)$ et $\tilde{\pi}$ est une représentation irréductible elliptique de $\tilde{G}(F)$. Alors $\tilde{\pi}$ ne dépend pas du choix de $\pi$ et on peut poser $t(\mathcal{O})=t(\tilde{\pi})$ avec la notation du paragraphe 3.1. Soit $\sigma\in Temp(H)$ une représentation irréductible tempérée de $H(F)$. Les espaces hermitiens $W$ et $\tilde{V}$ sont compatibles, on peut donc poser $m(\mathcal{O},\sigma)=m(\tilde{\pi},\sigma)$.

\begin{duf}
Soient $\sigma\in Temp(H)$ une représentation tempérée irréductible de $H(F)$ et $f\in C_c^\infty(G(F))$ une fonction très cuspidale. On définit alors la quantité $J_{spec}(\sigma,f)$ par la formule

\[\begin{aligned}
J_{spec}(\sigma,f) & =\displaystyle\sum_{L\in\mathcal{L}(M_{min})} |W^L||W^G|^{-1}(-1)^{a_L}\sum_{\mathcal{O}\in\{\Pi_{ell}(L)\}; m(\mathcal{O},\sigma)=1} \\
 & [i\mathcal{A}_\mathcal{O}^\vee:i\mathcal{A}_{L,F}^\vee]^{-1}t(\mathcal{O})^{-1}\int_{i\mathcal{A}_{L,F}^*} J^G_L(\pi_\lambda,f) d\lambda
\end{aligned}\]

\noindent où pour chaque orbite $\mathcal{O}$ on a fixé un point base $\pi\in\mathcal{O}$.
\end{duf}

Le but de cette section est de prouver le résultat suivant

\begin{theo}
Soient $\sigma$ et $f$ comme précédemment alors on a
$$\lim\limits_{N\to\infty} J_N(\theta_\sigma,f)=J_{spec}(\sigma,f)$$
\end{theo}

\subsection{Utilisation de la formule de Plancherel}
 On peut exprimer $f$ grâce à la formule de Plancherel-Harish-Chandra. On a l'égalité

$$f(g)=\displaystyle\sum_{L\in \mathcal{L}(M_{min})} |W^L||W^G|^{-1} \sum_{\mathcal{O}\in \Pi_2(L)} f_\mathcal{O}(g)$$
où pour $\mathcal{O}\in \Pi_2(L)$ on a posé $f_\mathcal{O}(g)=[i\mathcal{A}_\mathcal{O}^\vee:i\mathcal{A}_{L,F}^\vee]^{-1} \displaystyle\int_{i\mathcal{A}_{L,F}^*} m(\tau_\lambda) Tr(i_Q^G(\tau_\lambda,g^{-1})i_Q^G(\tau_\lambda,f))d\lambda$ avec $Q\in \mathcal{P}(L)$ et $\tau \in \mathcal{O}$ quelconques. La fonction $f_\mathcal{O}$ appartient à l'espace de Schwarz-Harish-Chandra et elle est nulle pour presque tout $\mathcal{O}$. Cela permet d'obtenir l'expression suivante de ${}^g f^\xi$
 
$$\mbox{(1)} \;\;\; {}^g f^\xi(h)=\displaystyle\sum_{L\in \mathcal{L}(M_{min})} |W^L||W^G|^{-1}\sum_{\mathcal{O}\in \Pi_2(L)} \int_{U(F)} f_\mathcal{O}(g^{-1}hug) \xi(u) du$$

On a interverti l'intégrale et la double somme. Cela est possible car l'intégrale 

$$\displaystyle\int_{U(F)} f_\mathcal{O}(g^{-1}hug) \xi(u) du$$
est absolument convergente d'après la proposition II.4.5 de [W3]. Fixons un produit scalaire invariant sur $E_\sigma$ et une base orthonormale $\mathcal{B}_\sigma$ pour ce produit scalaire. On a alors l'égalité

$$J(\theta_\sigma,f,g)=\displaystyle\sum_{\epsilon\in \mathcal{B}_\sigma} \int_{H(F)} (\epsilon,\sigma(h)\epsilon) {}^g f^\xi(h)dh$$
Presque tout les termes de cette somme sont nuls. Substituant l'égalité (1) on obtient

$$J(\theta_\sigma,f,g)=\displaystyle\sum_{\epsilon\in\mathcal{B}_\sigma} \int_{H(F)} (\epsilon,\sigma(h)\epsilon) \sum_{L\in \mathcal{L}(M_{min})} |W^L||W^G|^{-1} $$
$$\displaystyle\sum_{\mathcal{O}\in\Pi_2(L)}\int_{U(F)} f_\mathcal{O}(g^{-1}hug)\xi(u) du dh$$

Pour $\epsilon\in E_\sigma$, $L\in \mathcal{L}(M_{min})$, $\mathcal{O}\in \Pi_2(L)$ et $g\in G(F)$, posons

\begin{center}
(2) $J_{L,\mathcal{O}}(\epsilon,f,g)=\displaystyle\int_{H(F)\times U(F)} (\epsilon,\sigma(h)\epsilon) f_\mathcal{O}(g^{-1}hug)\xi(u)dudh$
\end{center}

Alors le lemme 12.0.7 et la proposition II.4.5 de [W3] montrent que cette expression est absolument convergente. On a donc l'égalité

\begin{center}
(3) $J(\theta_\sigma,f,g)=\displaystyle\sum_{\epsilon\in\mathcal{B}_\sigma}\sum_{L\in \mathcal{L}(M_{min})}\sum_{\mathcal{O}\in \Pi_2(L)} |W^L||W^G|^{-1} J_{L,\mathcal{O}}(\epsilon,f,g)$
\end{center}

Fixons provisoirement $L\in \mathcal{L}(M_{min})$ et $\mathcal{O}\in \Pi_2(L)$. On a le fait suivant dont la démonstration est tout à fait analogue à celle du point 14.2(3) (le point crucial étant la centralité de $A(F)$ dans $M(F)$)

\begin{center}
(4) Il existe un entier naturel $c_0$ tel que pour tout $c\geqslant c_0$, $g\in M(F)K$ et $h\in H(F)$, on ait l'égalité
$$\displaystyle\int_{U(F)} f_\mathcal{O}(g^{-1}hug)\xi(u)du=\int_{U(F)_c} f_\mathcal{O}(g^{-1}hug)\xi(u)du$$
\end{center}

Fixons $\tau\in\mathcal{O}$ un point base et un parabolique $Q\in\mathcal{P}(L)$. On notera $\pi_\lambda=i_Q^G(\tau_\lambda)$ pour $\lambda\in i\mathcal{A}_{L,F}^*$ et on réalise toutes ces représentations sur l'espace commun $\mathcal{K}_{Q,\tau}^G$. On munit $E_\tau$ d'un produit scalaire invariant, alors $\mathcal{K}_{Q,\tau}^G$ hérite d'un produit scalaire induit. Soit $K_f\subset K$ un sous-groupe ouvert tel que $f$ soit biinvariante par $K_f$ et soit $\mathcal{B}_\mathcal{O}^{K_f}$ une base orthonormée de $(\mathcal{K}_{Q,\tau}^G)^{K_f}$. Pour tout $g\in G(F)$ on a

$$f_\mathcal{O}(g)=[i\mathcal{A}_\mathcal{O}^\vee:i\mathcal{A}_{L,F}^\vee]^{-1}\displaystyle\sum_{e\in \mathcal{B}_\mathcal{O}^{K_f}}\int_{i\mathcal{A}_{L,F}^*} m(\tau_\lambda) (e,\pi_\lambda(g^{-1})\pi_\lambda(f)e)d\lambda$$

Soit $c_0$ un entier tel que (4) soit vérifié. Pour $c\geqslant c_0$, $g\in M(F)K$ et $\epsilon\in E_\sigma$, remplaçons dans (2) $f_\mathcal{O}$ par son expression précédente. Après les changements de variables $h\mapsto h^{-1}$ et $u\mapsto u^{-1}$, on obtient

\[\begin{aligned}
J_{L,\mathcal{O}}(\epsilon,f,g)=[i\mathcal{A}_\mathcal{O}^\vee:i\mathcal{A}_{L,F}^\vee]^{-1}\displaystyle\int_{H(F)\times U(F)_c} (\sigma(h)\epsilon,\epsilon) \\
\sum_{e\in \mathcal{B}_\mathcal{O}^{K_f}}\int_{i\mathcal{A}_{L,F}^*} m(\tau_\lambda) (\pi_\lambda(g)e,\pi_\lambda(hug)\pi_\lambda(f)e)\overline{\xi(u)}d\lambda dudh
\end{aligned}\]

A $g$ fixé le coefficient $(\pi_\lambda(g)e,\pi_\lambda(hug)\pi_\lambda(f)e)$ est essentiellement majoré par $\Xi^G(hu)$ indépendamment de $\lambda$. On en déduit que l'expression précédente est absolument convergente ce qui permet de permuter les intégrales et on reconnaît alors l'intégrale intérieure: c'est $\mathcal{L}_{\pi_\lambda,\sigma,c}(\epsilon\otimes \pi_\lambda(g)e,\epsilon\otimes \pi_\lambda(g)\pi_\lambda(f)e)$. Quitte à accroître $c_0$ c'est aussi $\mathcal{L}_{\pi_\lambda,\sigma}(\epsilon\otimes \pi_\lambda(g)e,\epsilon\otimes \pi_\lambda(g)\pi_\lambda(f)e)$. On a alors

\[\begin{aligned}
J_{L,\mathcal{O}}(\epsilon,f,g)=[i\mathcal{A}_\mathcal{O}^\vee:i\mathcal{A}_{L,F}^\vee]^{-1}\displaystyle \sum_{e\in \mathcal{B}_\mathcal{O}^{K_f}}\int_{i\mathcal{A}_{L,F}^*} m(\tau_\lambda)\\
 \mathcal{L}_{\pi_\lambda,\sigma}(\epsilon\otimes \pi_\lambda(g)e,\epsilon\otimes \pi_\lambda(g)\pi_\lambda(f)e)d\lambda
\end{aligned}\]

En particulier si $m(\mathcal{O},\sigma)=0$ alors $J_{L,\mathcal{O}}(\epsilon,f,g)=0$ pour tous $\epsilon\in E_\sigma$ et $g\in M(F)K$. Supposons que $m(\mathcal{O},\sigma)=1$ et fixons des familles $(\epsilon'_j)_{j=1,\ldots,n}$, $(\epsilon_j)_{j=1,\ldots,n}$, $(e'_j)_{j=1,\ldots,n}$, $(e_j)_{j=1,\ldots,n}$, $(\varphi_j)_{j=1,\ldots,n}$ vérifiant le ii)b) de la proposition 14.2.2. Pour $\lambda\in i\mathcal{A}_{L,F}^*$, $g\in M(F)K$ et $e\in \mathcal{K}^G_{Q,\tau}$ considérons la somme

$$X_\lambda(e,g)=\displaystyle\sum_{\epsilon\in \mathcal{B}_\sigma} \mathcal{L}_{\pi_\lambda,\sigma}(\epsilon\otimes \pi_\lambda(g)e, \epsilon\otimes \pi_\lambda(g)\pi_\lambda(f)e)$$

Pour $j=1,\ldots,n$ et $c,c'\in\mathbb{N}$ posons

$$X_{\lambda,j,c,c'}(e,g)=\displaystyle\int_{H(F)U(F)_c}\int_{H(F)U(F)_{c'}} (\sigma(h)\epsilon'_j,\epsilon_j) (\pi_\lambda(h'u'g)e,\pi_\lambda(hu)e_i)$$
$$(e'_j,\pi_\lambda(h'u'g)\pi_\lambda(f)e)\overline{\xi}(u)du'dh'dudh$$

C'est une intégrale absolument convergente. Alors la même discussion que dans [W2] p.91 montre que pour $c$ et $c'$ assez grand on a

$$X_\lambda(e,g)=\displaystyle\sum_{j=1}^n \varphi_j(\lambda) X_{\lambda,j,c,c'}(e,g)$$

On pose $J_{L,\mathcal{O}}(\theta_\sigma,f,g)=\displaystyle\sum_{\epsilon\in\mathcal{B}_\sigma} J_{L,\mathcal{O}}(\epsilon,f,g)$. Ce que l'on vient de dire implique alors l'égalité suivante

$$J_{L,\mathcal{O}}(\theta_\sigma,f,g)=[i\mathcal{A}_\mathcal{O}^\vee:i\mathcal{A}_{L,F}^\vee]^{-1}\displaystyle\sum_{e\in \mathcal{B}_\mathcal{O}^{K_f}}\sum_{j=1}^n \int_{i\mathcal{A}_{L,F}^*} m(\tau_\lambda) \varphi_j(\lambda) X_{\lambda,j,c,c'}(e,g) d\lambda$$

D'après (3) on a aussi

$$J(\theta_\sigma,f,g)=\displaystyle\sum_{L\in\mathcal{L}(M_{min})}\sum_{\mathcal{O}\in \Pi_2(L)} |W^G||W^L|^{-1} J_{L,\mathcal{O}}(\theta_\sigma,f,g)$$

$J_N(\theta_\sigma,f)$ est l'intégrale de $J(\theta_\sigma,f,mk)\kappa_N(mk)\delta_P(m)^{-1}$ sur $m\in H(F)\backslash M(F)$ et $k\in K$. Pour $L\in \mathcal{L}(M_{min})$, $\mathcal{O}\in \Pi_2(L)$ et $N,C\in\mathbb{N}$ posons

$$J_{L,\mathcal{O},N,C}(\theta_\sigma,f)=[i\mathcal{A}_\mathcal{O}^\vee:i\mathcal{A}_{L,F}^\vee]^{-1}\displaystyle\sum_{e\in \mathcal{B}_\mathcal{O}^{K_f}}\sum_{j=1}^n \int_{i\mathcal{A}_{L,F}^*} m(\tau_\lambda) \varphi_j(\lambda)$$
$$\displaystyle\int_{H(F)U(F)_c} \mathbf{1}_{\sigma <Clog(N)}(hu)(\sigma(h)\epsilon'_j,\epsilon_j)\int_{G(F)}(\pi_\lambda(g)e,\pi_\lambda(hu)e_j)$$
$$(e'_j,\pi_\lambda(g)\pi_\lambda(f)e)\kappa_N(g) dgdudhd\lambda$$

\begin{lem}
(i) Cette expression est absolument convergente.\\
(ii) Il existe $C$ tel qu'on ait la majoration

$$|J_N(\theta_\sigma,f)-\displaystyle\sum_{L\in\mathcal{L}(M_{min})}\sum_{\mathcal{O}\in \Pi_2(L)} |W^G||W^L|^{-1} J_{L,\mathcal{O},N,C}(\theta_\sigma,f)|<<N^{-1}$$
pour tout $N\geqslant 2$.
\end{lem}

\ul{Preuve}: C'est exactement la même que celle du lemme 6.4 de [W2] où on remplace les majorations 4.3(6) (7) et (8) de [W2] par les majorations de la proposition 13.0.5 $\blacksquare$

On fixe dorénavant un $C$ qui vérifie le (ii) du lemme précédent.

\subsection{Changement de fonction de troncature}

On fixe jusqu'au paragraphe 16.5 des données $L\in\mathcal{L}(M_{min})$, $Q\in \mathcal{P}(L)$, $\mathcal{O}\in\{\Pi_2(L)\}$ et $\tau\in\mathcal{O}$. Soit $Y\in \mathcal{A}_{P_{min}}^+$. On en déduit une $(G,M_{min})$-famille orthogonale positive $\mathcal{Y}=(Y_P)_{P\in \mathcal{P}(M_{min})}$. Soit $g\mapsto u(g,\mathcal{Y})$ la fonction caractéristique de l'ensemble des $g\in G(F)$ qui s'écrivent $g=k_1mk_2$ avec $k_1,k_2\in K$ et $m\in M_{min}(F)$ qui vérifie $\sigma^G_{M_{min}}(H_{M_{min}}(m),\mathcal{Y})=1$. Soient $e',e''\in \mathcal{K}^G_{Q,\tau}$ et $\varphi$ une fonction holomorphe sur un voisinage de $i\mathcal{A}_{L,F}^*$ dans $\mathcal{A}_{L,\mathbb{C}}^*/i\mathcal{A}_{L,F}^\vee$. Pour $e\in \mathcal{K}^G_{Q,\tau}$, $g,g'\in G(F)$ et $\lambda\in i\mathcal{A}_{L,F}^*$ posons

$$\Phi(e,g,g',\lambda)=(\pi_\lambda(g)e,\pi_\lambda(g')e')(e'',\pi_\lambda(g)\pi_\lambda(f)e)$$

On définit alors

$$\Phi_N(g')=\displaystyle\sum_{e\in \mathcal{B}_{\mathcal{O}}^{K_f}}\int_{i\mathcal{A}_{L,F}^*}\phi(\lambda) m(\tau_\lambda)\int_{G(F)} \Phi(e,g,g',\lambda) \kappa_N(g) dg d\lambda$$

$$\Phi_Y(g')=\displaystyle\sum_{e\in \mathcal{B}_{\mathcal{O}}^{K_f}}\int_{i\mathcal{A}_{L,F}^*}\phi(\lambda) m(\tau_\lambda)\int_{G(F)} \Phi(e,g,g',\lambda) u(g,\mathcal{Y}) dg d\lambda$$

\begin{prop}
Les deux expressions ci-dessus sont absolument convergentes. Soit $R$ un réel, il existe deux réels $c_1,c_2>0$ tels que

$$|\Phi_N(g')-\Phi_Y(g')|<<N^{-R}$$
pour tout $N\geqslant 2$ et pour tous $g'\in G(F)$, $Y\in \mathcal{A}_{P_{min}}^+$ vérifiant $\sigma(g')\leqslant Clog(N)$ et $c_1log(N)\leqslant \alpha(Y)\leqslant c_2 N$ pour tout $\alpha\in \Delta_{min}$
\end{prop}

\ul{Preuve}: Il suffit de reprendre la preuve de la proposition 6.6 de [W2] et d'y apporter trois modifications:

\begin{itemize}
\item Il faut remplacer la majoration 4.3(2) de [W2] par la majoration de la proposition 13.0.2.

\item Pour établir l'analogue du point (5) de la preuve dans [W2], on n'a besoin que de la propriété suivante sur la fonction de troncature $\kappa_N$:
\begin{center}
Il existe une constante $c_3>0$ tel que pour tout $g\in G(F)$ vérifiant $\sigma(g)\leqslant c_3 N$ on a $\kappa_N(g)=1$
\end{center}

\item Il y a une erreur dans la preuve de la proposition 6.6 de [W2] qui a été corrigée dans [W5]. L'erreur se trouve dans l'utilisation de l'inégalité (11) pour majorer la fonction $f_5$. Si la fonction $f_5$ est seulement $C^\infty$ cela ne suffit pas, car ses coefficients de Fourier ne seront qu'à décroissance rapide. Puisque l'on a supposé ici que la fonction $\phi$ admet un prolongement holomorphe sur un voisinage de $i\mathcal{A}_{L,F}^*$, la fonction $f_4$ de la référence admettra elle aussi un prolongement holomorphe dans un tel voisinage. Ainsi les coefficients de Fourier de $f_4$ seront à décroissance exponentielle. Il suffit alors de remplacer l'inégalité (11) de la référérence par la suivante

\vspace{4mm}

\textbf{(11')} Pour tout $c>0$ on a la minoration

$$|\zeta(xm',ym)|>c log(N)$$
pour tout $x,y\in K_1$, pourvu que $c_1$ soit assez grand.

\vspace{4mm}

On peut alors majorer $f_5$ par une puissance négative de $N$ aussi grande que l'on veut, pourvu que $c_1$ soit assez grand. Avec cette correction mineure la preuve de [W2] devient correcte.
\end{itemize}

$\blacksquare$

\subsection{Utilisation des calculs spectraux d'Arthur}

Pour tout $\epsilon>0$ on note $\mathcal{D}(\epsilon)$ l'ensemble des $Y\in \mathcal{A}_{P_{min}}^+$ tels que

$$inf\{\alpha(Y);\alpha\in\Delta\}>\epsilon sup\{\alpha(Y); \alpha\in\Delta\}$$

Fixons une norme $|.|$ sur $\mathcal{A}_{M_{min}}$. Pour $L'\in \mathcal{L}(L)$, on note $\Lambda^{L'}_{\mathcal{O},ell}$ l'ensemble des $\lambda\in i\mathcal{A}_L^*$ tels que $R^{L'}(\tau_\lambda)\cap W^{L'}(L)_{reg}\neq \emptyset$. Cet ensemble est stable par translation par $i\mathcal{A}_{L,F}^\vee+i\mathcal{A}_{L'}^*$. Pour un tel $L'$ et un tel $\lambda$ on dispose d'une décomposition

$$i^{L'}_{L'\cap Q}(\tau_\lambda)=\displaystyle\bigoplus_{\zeta\in R^{L'}(\tau_\lambda)^\vee} i^{L'}_{L'\cap Q}(\tau_\lambda,\zeta)$$

Fixons $S'=L'U'\in\mathcal{P}(L')$ et posons $Q(S')=(L'\cap Q)U'$. On en déduit par induction une décomposition analogue de $E^G_{Q(S'),\tau_\lambda}$ que l'on peut aussi écrire

$$\mathcal{K}^G_{Q(S'),\tau}=\displaystyle\bigoplus_{\zeta\in R^{L'}(\tau_\lambda)^\vee} \mathcal{K}^G_{Q(S'),\tau,\zeta}$$
On note $proj_{\lambda,\zeta}$ la projection de $\mathcal{K}^G_{Q(S'),\tau}$ sur $\mathcal{K}^G_{Q(S'),\tau,\zeta}$ par rapport aux autres facteurs. Pour tout $g'\in G(F)$ on pose

$$\Phi(g')=\displaystyle\sum_{L'\in\mathcal{L}(L)} (-1)^{a_{L'}}\sum_{\lambda\in \Lambda^{L'}_{\mathcal{O},ell}/(i\mathcal{A}_{L,F}^\vee+i\mathcal{A}_{L'}^*)} |R^{L'}(\tau_\lambda)|2^{a_{L'}-a_L}\sum_{\zeta\in R^{L'}(\tau_\lambda)^\vee}$$

$$\displaystyle\int_{i\mathcal{A}_{L',F}^*} \left(proj_{\lambda,\zeta}\circ R_{Q(S')|Q}(\tau_{\lambda+\mu})e'',proj_{\lambda,\zeta}\circ R_{Q(S')|Q}(\tau_{\lambda+\mu})i^G_Q(\tau_{\lambda+\mu},g')e'\right)$$

$$J^G_{L'}(i^{L'}_{L'\cap Q}(\tau_{\lambda+\mu},\zeta),f)\varphi(\lambda+\mu) d\mu$$

\begin{prop}
Soit $\epsilon>0$ et $R\geqslant 1$. On a

$$|\Phi_Y(g')-\Phi(g')|<<\sigma(g')^R\Xi^G(g')|Y|^{-R}$$
pour tout $g'\in G(F)$ et tout $Y\in \mathcal{D}(\epsilon)\cap \mathcal{A}_{M_{min},F}$.
\end{prop}

\ul{Preuve}: C'est la même que celle des sections 6.7 et 6.8 de [W2]. La proposition 6.7 de [W2] reprend les calculs spectraux d'Arthur pour la formule des traces locale (dans [A3] p.69 à 88) qui donnent une approximation de notre terme $\Phi_Y(1)$ dans le cas où $\varphi=1$. Waldspurger montre qu'on peut alors glisser une fonction $\varphi$ tout le long des calculs d'Arthur et que ceci donne une approximation du terme $\Phi_Y(g')$ même lorsque $g'\neq 1$. Le lemme 6.8 de [W2] démontre que cette approximation est alors exactement notre terme $\Phi(g')$. Ces deux démonstrations sont tout à fait générales et s'appliquent à des groupes $G$ réductifs connexes généraux sauf le dernier paragraphe de la preuve du lemme 6.8 où il est fait usage de propriétés particulières des $R$-groupes dans le cas des groupes spéciaux orthogonaux. Ces propriétés sont aussi vérifiées dans le cas unitaire cf section 3.2 $\blacksquare$

\subsection{Evaluation d'une limite}

\begin{lem}
On a l'égalité

$$\lim\limits_{N\to \infty} J_{L,\mathcal{O},N,C}(\theta_\sigma,f)=[i\mathcal{A}_\mathcal{O}^\vee: i\mathcal{A}_{L,F}^\vee]^{-1}\displaystyle\sum_{L'\in\mathcal{L}(L)} (-1)^{a_{L'}}\sum_{\lambda\in \Lambda^{L'}_{\mathcal{O},ell}/(i\mathcal{A}_{L,F}^\vee+i\mathcal{A}_{L'}^*)}$$
$$|R^{L'}(\tau_\lambda)|2^{a_{L'}-a_L}\displaystyle\sum_{\zeta\in R^{L'}(\tau_\lambda)^\vee \\ m(i^{L'}_{L'\cap Q}(\tau_\lambda,\zeta),\sigma)=1} \int_{i\mathcal{A}_{L',F}^*} J_{L'}^G(i^{L'}_{L'\cap Q}(\tau_{\lambda+\mu},\zeta),f)d\mu$$
\end{lem}

\ul{Preuve}: Encore une fois c'est exactement la même que celle du lemme 6.9 de [W2]. Il suffit d'utiliser la majoration du corollaire 13.0.1(1) à la place de la majoration 4.3(4) de [W2], et de remarquer que le théorème 14.3.1 et la proposition 15.3.1 entraînent les analogues des lemmes 5.3(ii) et 5.4 de [W2] $\blacksquare$

\subsection{Preuve du théorème}

Il suffit de reprendre le paragraphe 6.10 de [W2]. Fixons $L'\in\mathcal{L}(M_{min})$ et $\mathcal{O}'\in \{\Pi_{ell}(L')\}$. D'après les lemmes 16.2.1 et 16.5.1, il s'agit essentiellement de compter les quadruplets $(L,\mathcal{O},\lambda,\zeta)$ où $L\in\mathcal{L}^{L'}(M_{min})$, $\mathcal{O}\in \{\Pi_{ell}(L)\}$, $\lambda\in \Lambda^{L'}_{\mathcal{O},ell}/(i\mathcal{A}_{L,F}^\vee+i\mathcal{A}_{L'}^*)$ et $\zeta\in R^{L'}(\tau_\lambda)^\vee$ ($\tau$ est un point base de $\mathcal{O}$), tels que

$$\{i^{L'}_{L'\cap Q}(\tau_\lambda,\zeta)_\mu; \; \mu\in i\mathcal{A}_{L'}^*\}=\mathcal{O}'$$
(Q est n'importe quel élément de $\mathcal{P}(L)$) $\blacksquare$

\section{Une formule pour la multiplicité}

\subsection{Le théorème}

Soient $(V,h_V)$ et $(W,h_W)$ deux espaces hermitiens compatibles de groupes unitaires respectifs $G$ et $H$. Soient $\theta$ et $\theta'$ des quasicaractères sur $G(F)$ et $H(F)$ respectivement. Supposons que $d_V>d_W$. On reprend alors les notations de la section 5. Il y est défini un ensemble $\mathcal{T}$ de tores de $H$. On pose alors

$$m_{geom}(\theta,\theta')=\displaystyle\sum_{T\in\mathcal{T}} |W(H,T)|^{-1} \lim\limits_{s\to 0^+} \int_{T(F)} c_{\theta'}(t)c_\theta(t)D^H(t)^{1/2}D^G(t)^{1/2}\Delta(t)^{s-1/2} dt$$

Cette expression a un sens d'après les lemmes 5.2.2 et 5.3.1. Si $d_W>d_V$ on pose $m_{geom}(\theta,\theta')=m_{geom}(\theta',\theta)$. Pour $\pi\in Temp(G)$ et $\sigma\in Temp(H)$ on notera $c_\pi=c_{\theta_\pi}$, $c_\sigma=c_{\theta_\sigma}$, $m_{geom}(\theta,\sigma)=m_{geom}(\theta,\theta_{\sigma^\vee})$ et $m_{geom}(\pi,\sigma)=m_{geom}(\theta_\pi,\theta_{\sigma^\vee})$.

\vspace{3mm}

\begin{theo}
Pour tout $\pi\in Temp(G)$ et $\sigma\in Temp(H)$ on a l'égalité

$$m(\pi,\sigma)=m_{geom}(\pi,\sigma)$$
\end{theo}

\vspace{3mm}

Soit $f\in C_c^\infty(G(F))$ une fonction cuspidale et $\sigma\in Temp(H)$. Posons

$$m_{spec}(f,\sigma)=\displaystyle\sum_{\pi\in \Pi_{ell}(G); m(\pi,\sigma^\vee)=1} t(\pi)^{-1} \theta_\pi(f)$$

Le théorème précédent découlera alors du théorème suivant

\vspace{3mm}

\begin{theo}
Soient $f$ et $\sigma$ comme ci-dessus. On a

$$m_{spec}(f,\sigma)=m_{geom}(I\theta_f,\sigma)$$
\end{theo}

\vspace{3mm}

Ces deux théorèmes seront démontrés en 16.4 et 16.5.

\subsection{Induction pour la multiplicité géométrique}

Soit $k\geqslant 1$ un entier. Pour $\theta$ un quasicaractère de $GL_k(E)$ on pose $m_{geom}(\theta)=c_{\theta,\mathcal{O}_k}(1)$ où $\mathcal{O}_k$ est l'unique orbite nilpotente régulière de $\mathfrak{gl}_k(E)$. Remarquons que ce terme ne dépend que de la restriction de $\theta$ au sous-ensemble $\{H_{R_{E/F}GL_k}=0\}$ de $GL_k(E)$, car ce sous-ensemble est un voisinage de $1$. \\

Fixons deux espaces hermitiens compatibles $(V,h_V)$ et $(W,h_W)$ de groupes unitaires respectifs $G$ et $H$ et soit $L$ un Lévi de $G$. On a alors une décomposition

$$L\simeq R_{E/F}GL_{n_1}\times\ldots\times R_{E/F} GL_{n_k}\times \tilde{G}$$

\noindent où $\tilde{G}$ est le groupe unitaire d'un sous-espace hermitien $\tilde{V}$ de $V$. Pour $j=1,\ldots,k$ soit $\theta_j$ un quasicaractère de $R_{E/F} GL_{n_j}$ et $\tilde{\theta}$ un quasicaractère de $\tilde{G}$. On note $\theta^L=\theta_1\otimes \ldots\otimes \theta_k\otimes \tilde{\theta}$ le quasicaractère de $L(F)$ donné par $\theta^L(g_1,\ldots,g_k,\tilde{g})=\theta_1(g_1)\ldots\theta_k(g_k)\tilde{\theta}(\tilde{g})$. Enfin, on pose $m_{geom}(\theta^L,\rho)=m_{geom}(\tilde{\theta},\rho)m_{geom}(\theta_1)\ldots m_{geom}(\theta_k)$. Soit $\sigma\in Temp(H)$.

\begin{lem}
Supposons que $\theta=Ind^G_L(\theta^L)$ alors on a

$$m_{geom}(\theta,\sigma)=m_{geom}(\theta^L,\sigma)$$
\end{lem}

\ul{Preuve}: Posons $n=n_1+\ldots+n_k$. Fixons une injection $(W,h_W)$ dans $(V,h_V)$. On distingue trois cas

\vspace{2mm}

\ul{Cas $d_V>d_{\tilde{V}}>d_W$}: Quitte à conjuguer $L$, on peut supposer que $W\subset \tilde{V}$. Les ensembles de tores à partir desquels sont définis $m_{geom}(\theta,\sigma)$ et $m_{geom}(\tilde{\theta},\sigma)$ sont alors les mêmes : c'est l'ensemble noté $\mathcal{T}$ dans la section 5. On a par définition

$$m_{geom}(\theta,\sigma)=\displaystyle\sum_{T\in\mathcal{T}} |W(H,T)|^{-1} \lim\limits_{s\to 0^+} \int_{T(F)} c_{\sigma^\vee}(t)c_\theta(t)D^H(t)^{1/2}D^G(t)^{1/2}\Delta(t)^{s-1/2} dt$$

\noindent et

\[\begin{aligned}
\displaystyle m_{geom}(\theta^L,\sigma)=c_{\theta_1,\mathcal{O}_{n_1}}(1)\ldots c_{\theta_k,\mathcal{O}_{n_k}}(1) & \sum_{T\in\mathcal{T}} |W(H,T)|^{-1} \\
 & \lim\limits_{s\to 0^+} \int_{T(F)} c_{\sigma^\vee}(t)c_{\tilde{\theta}}(t)D^H(t)^{1/2}D^{\tilde{G}}(t)^{1/2}\Delta(t)^{s-1/2} dt
\end{aligned}\]

Soit $T\in\mathcal{T}$. On a alors une décomposition orthogonale $W=W'\oplus W''$ tel que $T$ soit un tore maximal anisotrope de $H'$ le groupe unitaire de $W'$. On note $H''$, $G''$ et $\tilde{G}''$ les groupes unitaires des supplémentaires orthogonaux de $W'$ dans $W$, $V$ et $\tilde{V}$ respectivement. Rappelons que $T_\natural$ est l'ensemble des $t\in T$ dont toutes les valeurs propres sur $W'$ sont différentes de $1$ et de multiplicité $1$. On a alors $Z_G(t)=G''\times T$. Il suffit de montrer que pour tout $t\in T_\natural(F)$ on a 

$$\mbox{(1)}\;\;\; D^G(t)^{1/2}c_\theta(t)=c_{\tilde{\theta}}(t)c_{\theta_1,\mathcal{O}_{n_1}}(1)\ldots c_{\theta_k,\mathcal{O}_{n_k}}(1) D^{\tilde{G}}(t)^{1/2}$$

\noindent Rappelons que l'on a

$$c_\theta(t)=\displaystyle\frac{1}{|Nil(\mathfrak{g}''(F))_{reg}|}\sum_{\mathcal{O}\in Nil(\mathfrak{g}''(F))_{reg}} c_{\theta,\mathcal{O}}(t)$$

Soit $\mathcal{O}\in Nil(\mathfrak{g}''(F))_{reg}$. Le lemme 2.3 de [W2] permet d'exprimer $c_{\theta,\mathcal{O}}(t)$ à partir de $\theta^L$. Explicitons la formule dans notre cas particulier. Il y figure une somme sur $\mathcal{X}^L(t)$ un ensemble de représentants des classes de conjugaison par $L(F)$ des éléments de $L(F)$ conjugués à $t$ par $G(F)$.

\begin{center}
(2) On peut prendre $\mathcal{X}^L(t)=\{t\}$
\end{center}

En effet, soit $g\in G(F)$ tel que $gtg^{-1}\in L(F)$. Alors $A_L$ commute à $gtg^{-1}$ donc $A_L\subset gZ_G(t)g^{-1}=g(G''\times T)g^{-1}$. Le tore $T$ étant anisotrope on en déduit $A_L\subset gG''g^{-1}$. L'intersection des noyaux des $a-1$ dans $V$ pour $a\in A_L(F)$ est $\tilde{V}$ d'où $gW'\subset \tilde{V}$. Comme on a aussi $W'\subset\tilde{V}$, d'après le théorème de Witt quitte à multiplier $g$ par un élément de $\tilde{G}(F)\subset L(F)$ on peut supposer que $gW'=W'$. Alors $g$ stabilise aussi $V''$ et si on note $g'$ la restriction de $g$ à $W'$ on a $gtg^{-1}=g'tg'^{-1}$ et $g'\in H'(F)\subset L(F)$ d'où (2). \\

La deuxième somme du lemme 2.3 de [W2] porte alors sur $\Gamma_t/G_t(F)$ où $\Gamma_t=Z_{G(F)}(t)$ et on a vu que ce centralisateur est connexe. La deuxième somme est donc elle aussi triviale. Enfin la dernière somme porte sur les orbites $\mathcal{O}'\in Nil(\mathfrak{l}_t)$ telles que $[\mathcal{O}:\mathcal{O}']=1$ c'est-à-dire telles que $\mathcal{O}$ soit dans l'induite de $\mathcal{O}'$. On a

$$L_t=R_{E/F} GL_{n_1}\times\ldots\times R_{E/F} GL_{n_k}\times \tilde{G}''\times T$$

\noindent et

$$Nil(\mathfrak{l}_t)=Nil(\mathfrak{gl}_{n_1}(E))\times\ldots\times Nil(\mathfrak{gl}_{n_k}(E))\times Nil(\tilde{\mathfrak{g}}'')$$

 De plus si $[\mathcal{O}:\mathcal{O}']=1$, alors l'orbite nilpotente $\mathcal{O}'$ est régulière. Comme $Nil(\mathfrak{gl}_{n_j}(E))_{reg}=\{\mathcal{O}_{n_j}\}$ pour $j=1,\ldots,k$, cette dernière somme porte en fait sur les $\mathcal{O}'\in Nil(\tilde{\mathfrak{g}}'')_{reg}$ telles que $[\mathcal{O}:\mathcal{O}_{n_1}\times\ldots\times\mathcal{O}_{n_k}\times  \mathcal{O}']=1$. On vérifie aisément à partir de la description en 3.3 des orbites nilpotentes régulières d'un groupe unitaire qu'il existe une unique telle orbite $\mathcal{O}'$ et que l'application $\mathcal{O}\mapsto \mathcal{O}'$ est une bijection entre $Nil(\mathfrak{g}'')_{reg}$ et $Nil(\tilde{\mathfrak{g}}'')_{reg}$. Enfin le terme $[Z_L(t)(F):L_t(F)]$ apparaissant dans le lemme 2.3 de [W2] est trivial car $Z_L(t)=R_{E/F}GL_{n_1}\times\ldots\times R_{E/F} GL_{n_k}\times \tilde{G}''\times T=L_t$. Finalement on obtient

$$\displaystyle\sum_{\mathcal{O}\in Nil(\mathfrak{g}''(F))_{reg}} c_{\theta,\mathcal{O}}(t)=D^G(t)^{-1/2}D^L(t)^{1/2}c_{\theta_1,\mathcal{O}_{n_1}}(1)\ldots c_{\theta_k,\mathcal{O}_{n_k}}(1)\sum_{\mathcal{O}\in Nil(\tilde{\mathfrak{g}}''(F))_{reg}} c_{\tilde{\theta},\mathcal{O}}(t)$$

Puisque $D^L(t)=D^{\tilde{G}}(t)$ et $|Nil(\mathfrak{g}''(F))_{reg}|=|Nil(\tilde{\mathfrak{g}}''(F))_{reg}|$, c'est l'égalité (1), ce qu'il nous fallait.

\vspace{2mm}

\ul{Cas $d_V>d_W>d_{\tilde{V}}$}: A nouveau, quitte à conjuguer $L$, on peut supposer que $\tilde{V}\subset W$. Les ensembles de tores définissant $m_{geom}(\theta,\sigma)$ et $m_{geom}(\tilde{\theta},\sigma)$ sont cette fois différents. On les note $\mathcal{T}$ et $\tilde{\mathcal{T}}$ respectivement: ce sont des classes de représentants dans des ensembles de tores $\underline{\mathcal{T}}$ et $\underline{\tilde{\mathcal{T}}}$ pour l'action par conjugaison de $H(F)$ et $\tilde{G}(F)$ respectivement. On a alors

$$m_{geom}(\theta,\sigma)=\displaystyle\sum_{T\in\mathcal{T}} |W(H,T)|^{-1} \lim\limits_{s\to 0^+} \int_{T(F)} c_{\sigma^\vee}(t)c_\theta(t)D^H(t)^{1/2}D^G(t)^{1/2}\Delta(t)^{s-1/2} dt$$

\noindent et

\[\begin{aligned}
\displaystyle m_{geom}(\theta^L,\sigma)=c_{\theta_1,\mathcal{O}_{n_1}}(1)\ldots c_{\theta_k,\mathcal{O}_{n_k}}(1) & \sum_{T\in\tilde{\mathcal{T}}} |W(\tilde{G},T)|^{-1} \\
 & \lim\limits_{s\to 0^+} \int_{T(F)} c_{\sigma^\vee}(t)c_{\tilde{\theta}}(t)D^{\tilde{G}}(t)^{1/2}D^H(t)^{1/2}\Delta(t)^{s-1/2} dt
\end{aligned}\]

Notons $m_T(\theta,\sigma)$ le terme correspondant à $T\in\mathcal{T}$ dans la première somme.

\begin{center}
(3) Si $T\in\mathcal{T}$ n'est conjugué à aucun élément de $\tilde{\mathcal{T}}$ par $H(F)$ alors $m_T(\theta,\sigma)=0$
\end{center}

Soient $W=W'_T\oplus W''_T$ la décomposition orthogonale associée à $T$ et $t\in T(F)$ en position générique. Il suffit de montrer que $t$ n'est conjugué à aucun élément de $L(F)$ par $G(F)$ car alors $t$ n'est pas dans le support de $\theta=Ind^G_L(\theta^L)$. Or si $g\in G(F)$ est tel que $gtg^{-1}\in L(F)$ alors par le même raisonnement que pour (2) on a $gW'_T\subset \tilde{V}\subset W$. D'après le théorème de Witt, il existe $h\in H(F)$ tel que $hW'_T=gW'_T\subset \tilde{V}$. Mais alors $hTh^{-1}\in\underline{\tilde{\mathcal{T}}}$ ce qui contredit l'hypothèse.

\begin{center}
(4) Si $T\in \mathcal{T}$ est conjugué par $H(F)$ à un élément $T'\in\tilde{\mathcal{T}}$, alors cet élément est unique et on a $|W(H,T)|=|W(\tilde{G},T')|$.
\end{center}

Soient $T_1,T_2\in \tilde{\mathcal{T}}$ et $h\in H(F)$ tel que $hT_1h^{-1}=T_2$. On a des décompositions orthogonales $\tilde{V}=\tilde{V}_1'\oplus \tilde{V}_1''$ et $\tilde{V}=\tilde{V}_2'\oplus \tilde{V}_2''$ associées à $T_1$ et $T_2$ respectivement. On a alors $h\tilde{V}'_1=\tilde{V}'_2$. D'après le théorème de Witt, quitte à multiplier $h$ par un élément de $\tilde{G}(F)$, on peut supposer que $\tilde{V}'_1=\tilde{V}'_2$. Soit $h'$ la restriction de $h$ à $\tilde{V}'_1$ on a alors $h'T_1h'^{-1}=T_2$ et $h'\in \tilde{G}(F)$ ce qui prouve le premier point. Soit $T'\in \tilde{\mathcal{T}}$ et $\tilde{V}=\tilde{V}'\oplus \tilde{V}''$ la décomposition orthogonale associée. On note $W''$ le supplémentaire orthogonal de $\tilde{V}'$ dans $W$ et $H'$, $\tilde{G}''$ $H''$ les groupes unitaires de $\tilde{V}'$, $\tilde{V}''$ et $W''$ respectivement. Pour établir le deuxième point il suffit de remarquer que $Norm_{H(F)}(T')=H''\times Norm_{H'(F)}(T)$ et $Norm_{\tilde{G}(F)}(T')=\tilde{G}''\times Norm_{H'(F)}(T)$ d'où des isomorphismes $W(H,T')\simeq W(H',T)\simeq W(\tilde{G},T')$. \\

D'après les points (3) et (4) on peut faire porter la somme définissant $m_{geom}(\theta,\sigma)$ sur $T\in\tilde{\mathcal{T}}$. Soit $T\in\tilde{\mathcal{T}}$ un tel tore. Soit $\tilde{V}=\tilde{V}'\oplus\tilde{V}''$ la décomposition orthogonale associée à $T$, $V''$ l'orthogonal de $\tilde{V}'$ dans $V$ et $\tilde{G}''$, $G''$ les groupes unitaires de $\tilde{V}''$ et $V''$ respectivement. D'après (4), pour établir l'égalité voulue, il suffit de voir que pour $t\in T_\natural(F)$, on a 

$$\mbox{(5)}\;\;\; c_\theta(t)D^G(t)^{1/2}=c_{\tilde{\theta}}(t)c_{\theta_1,\mathcal{O}_{n_1}}(1)\ldots c_{\theta_k,\mathcal{O}_{n_k}}(1)D^{\tilde{G}}(t)^{1/2}$$

\noindent On a toujours

$$c_\theta(t)=\displaystyle\frac{1}{|Nil(\mathfrak{g}''(F))_{reg}|}\sum_{\mathcal{O}\in Nil(\mathfrak{g}''(F))_{reg}} c_{\theta,\mathcal{O}}(t)$$

Soit $\mathcal{O}\in Nil(\mathfrak{g}''(F))_{reg}$ et exprimons $c_{\theta,\mathcal{O}}(t)$ à partir du lemme 2.3 de [W2]. Le même raisonnement que dans le premier cas montre qu'on a $\Gamma_t=G_t(F)$ et qu'on peut prendre $\mathcal{X}^L(t)=\{t\}$. De plus, si $\tilde{V}\neq 0$, il existe une unique orbite $\mathcal{O}'\in Nil(\tilde{\mathfrak{g}}'')_{reg}$ telle que $[\mathcal{O}:\mathcal{O}_{n_1}\times\ldots\times\mathcal{O}_{n_k}\times  \mathcal{O}']=1$ et l'application $\mathcal{O}\mapsto \mathcal{O}'$ est une bijection de $Nil(\mathfrak{g}'')_{reg}$ sur $Nil(\tilde{\mathfrak{g}}'')_{reg}$. Si $\tilde{V}=0$, alors $Nil(\tilde{\mathfrak{g}}'')_{reg}$ ne contient qu'une orbite $\mathcal{O}'$ et pour toute $\mathcal{O}\in Nil(\mathfrak{g}'')_{reg}$ on a $[\mathcal{O}:\mathcal{O}_{n_1}\times\ldots\times\mathcal{O}_{n_k}\times  \mathcal{O}']=1$. Dans tout les cas, on en déduit que

$$c_\theta(t)=c_{\tilde{\theta}}(t)c_{\theta_1,\mathcal{O}_{n_1}}(1)\ldots c_{\theta_k,\mathcal{O}_{n_k}}(1)D^G(t)^{-1/2}D^L(t)^{1/2}$$

Puisque $D^L(t)=D^{\tilde{G}}(t)$ c'est l'égalité (5).

\vspace{2mm}

\ul{Cas $d_W>d_V>d_{\tilde{V}}$}: On peut alors trouver une décomposition orthogonale $W=(Z_+\oplus Z_-)\oplus D\oplus V$ où $D$ est une droite et $Z_+,Z_-$ sont deux sous-espaces totalement isotropes. Soit $D'$ une droite munie d'une forme hermitienne opposée à $h_D$. Considérons l'espace hermitien $V'$ somme orthogonale de $W$ et $D'$ et notons $G'$ son groupe unitaire. Alors $V'$ est somme orthogonale de $Z'_+\oplus Z'_-$ et de $V$ où $Z'_+$ et $Z'_-$ sont totalement isotropes. Posons $L'=GL(Z'_+)\times G$, c'est naturellement un sous-groupe de Lévi de $G'$. Soit $\mathcal{O}_{GL}$ l'unique orbite nilpotente régulière de $\mathfrak{gl}(Z'_+)$ et $\theta_{GL}$ un quasicaractère de $GL(Z'_+)$ tel que $c_{\theta_{GL},\mathcal{O}_{GL}}(1)=1$ (il en existe). Posons $\theta'=Ind^{G'}_{L'}(\theta_{GL}\otimes \theta)$. D'après le deuxième cas on a $m_{geom}(\theta,\sigma)=m_{geom}(\theta',\sigma)$. Soit $L''=GL(Z'_+)\times L$, c'est naturellement un sous-groupe de Lévi de $G'$. Puisque $\theta=Ind^G_L(\theta^L)$, on a aussi $\theta'=Ind^{G'}_{L''}(\theta_{GL}\otimes \theta_1\otimes\ldots\otimes \theta_k\otimes\tilde{\theta})$. Appliquant à nouveau le deuxième cas, on obtient l'égalité voulue $\blacksquare$

\subsection{Pseudo-coefficients}

Soit $(V,h_V)$ un espace hermitien de groupe unitaire $G$. Soit $L$ un Levi de $G$ et $\tau$ une représentation irréductible de la série discrète de $L(F)$ tel que $R(\tau)\cap W(L)_{reg}\neq\emptyset$. On a vu dans le paragraphe 3.1 qu'alors $R(\tau)\cap W(L)_{reg}$ ne contient qu'un seul élément que l'on avait noté $t$. Pour tout $\zeta\in R(\tau)^\vee$, on pose $\pi(\zeta)=i_Q^G(\tau,\zeta)$ où $Q\in\mathcal{P}(L)$. Les considérations des paragraphes 7.3 et 7.4 de [W2] valent toujours car elles n'utilisaient que le fait que $A_G=\{1\}$. En particulier on a

\begin{lem}
Il existe une fonction cuspidale $f\in C^\infty_c(G(F))$ vérifiant

\begin{enumerate}
\item $I\theta_f=\displaystyle\sum_{\zeta\in R(\tau)^\vee}\zeta(t) \theta_{\pi(\zeta)}$
\item $\theta_{\pi(\zeta)^\vee}(f)=\zeta(t) t(\pi(\zeta)^\vee)$ pour tout $\zeta\in R(\tau)^\vee$
\item $\theta_{\sigma^\vee}(f)=0$ pour tout $\sigma\in Temp(G)$ qui n'est pas l'une des représentations $\pi(\zeta)$.
\end{enumerate}

\end{lem}

On définit $T_{ell}(G)$ comme l'ensemble des représentations virtuelles tempérées de $G(F)$ de la forme $\displaystyle\sum_{\zeta\in R(\tau)^\vee} \zeta(t)\pi(\zeta)$, où $\tau$ est comme précédemment.
 
\subsection{Le cas du groupe linéaire}

Soient $k\geqslant 1$ un entier et $G=R_{E/F} GL_k$. Soit $f\in C_c^\infty(G(F))$ une fonction cuspidale. Posons

$$\displaystyle m_{spec}(f)=\sum_{\mathcal{O}\in\{\Pi_{ell}(G)\}} [i\mathcal{A}_{\mathcal{O}}^\vee:i\tilde{\mathcal{A}}_{G,F}^\vee]^{-1} \theta_\pi(f\mathbf{1}_{H_G=0})$$

où on a fixé un point base $\pi$ dans chaque orbite.

\begin{lem}
Pour toute fonction cuspidale $f\in C_c^\infty(G(F))$, on a l'égalité

$$m_{geom}(I\theta_f)=m_{spec}(f)$$
\end{lem}

\ul{Preuve}: C'est le lemme 7.6 de [W2] $\blacksquare$

\subsection{Le théorème 17.1.2 implique le théorème 17.1.1}

Supposons le théorème 17.1.2 vérifié et montrons que le théorème 17.1.1 l'est aussi. Le résultat est établi par récurrence. On veut montrer que $m_{geom}(\pi,\sigma)=m(\pi,\sigma)$ pour tout $(\pi,\sigma)\in Temp(G)\times Temp(H)$. Fixons $\sigma$. On peut étendre la définition des multiplicités $m(.,\sigma)$ et $m_{geom}(.,\sigma)$ par linéarité à l'ensemble des représentations virtuelles tempérées. Il suffit de montrer l'égalité sur une famille génératrice de cet espace. Une telle famille est formée par l'union de $T_{ell}(G)$ et de l'ensemble des représentations tempérées qui sont des induites propres paraboliques de représentations tempérées irréductibles. \\

 Traitons d'abord le cas où $\pi$ est une induite propre. On peut alors trouver un Levi $L$ de la forme $R_{E/F}GL_k\times \tilde{G}$ où $\tilde{G}$ est le groupe unitaire d'un sous-espace hermitien $\tilde{V}$ de $V$, un parabolique $Q\in\mathcal{P}(L)$ et une représentation $\mu\otimes \tilde{\pi}\in Temp(L)=Temp(GL_k(E))\times Temp(\tilde{G})$ de sorte que $\pi\simeq i_Q^G(\mu\otimes\tilde{\pi})$. D'après la proposition 15.3.1 on a alors $m(\pi,\sigma)=m(\tilde{\pi},\sigma)$. D'après le lemme 17.2.1 on a $m_{geom}(\pi,\sigma)=m_{geom}(\tilde{\pi},\sigma)c_{\theta_\mu,\mathcal{O}_k}(1)$. D'après un résultat de Rodier ([Ro]) le coefficient $c_{\theta_\mu,\mathcal{O}_k}(1)$ vaut $1$ si $\mu$ admet un modèle de Whittaker et $0$ sinon. Puisque $\mu$ est tempérée, $\mu$ admet un modèle de Whittaker et ce coefficient vaut 1. D'après l'hypothèse de récurrence, on a $m(\tilde{\pi},\sigma)=m_{geom}(\tilde{\pi},\sigma)$. D'où l'égalité $m(\pi,\sigma)=m_{geom}(\pi,\sigma)$. \\

 Il reste à établir l'égalité pour $\pi\in T_{ell}(G)$. On peut alors trouver $L$ et $\tau$ comme dans le paragraphe 17.3 tels qu'avec les notations de ce paragraphe on ait $\pi=\displaystyle\sum_{\zeta\in R(\tau)^\vee} \zeta(t)\pi(\zeta)$. Soit $f\in C_c^\infty(G(F))$ une fonction cuspidale vérifiant les conclusions du lemme 17.3.1. On a alors

\[\begin{aligned}
m_{spec}(f,\sigma) & =\displaystyle\sum_{\pi'\in \Pi_{ell}(G); m(\pi',\sigma^\vee)=1} t(\pi')^{-1}\theta_{\pi'}(f) \\
 & =\sum_{\zeta\in R(\tau)^\vee} m(\pi(\zeta)^\vee,\sigma^\vee)\zeta(t) \\
 & =\sum_{\zeta\in R(\tau)^\vee} m(\pi(\zeta),\sigma) \zeta(t) \\
 & =m(\pi,\sigma)
\end{aligned}\]
où à la troisième ligne on a utilisé le lemme 15.0.1. Enfin, d'après le premier point du lemme 17.3.1 on a aussi $I\theta_f=\theta_\pi$ d'où $m_{geom}(I\theta_f,\sigma)=m_{geom}(\pi,\sigma)$. Le théorème 17.1.2 appliqué à $f$ et $\sigma$ donne alors l'égalité voulue.

\subsection{Preuve du théorème 17.1.2}

 La démonstration de ce théorème procède aussi par récurrence. Les multiplicités $m_{geom}(I\theta_f,\sigma)$ et $m_{spec}(f,\sigma)$ vues comme distributions en $f$ sont invariantes. D'après le lemme 1.8.1, on peut donc supposer $f$ très cuspidale. Les théorèmes 5.4.1 et 16.1.1 calculent de deux façons la limite lorsque $N$ tend vers l'infini de $J_N(\theta_{\sigma^\vee},f)$. Ces deux expressions sont donc égales, autrement dit on a $J_{geom}(\theta_{\sigma^\vee},f)=J_{spec}(\sigma^\vee,f)$. On a $J_{geom}(\theta_{\sigma^\vee},f)=m_{geom}(\theta_f,\sigma)$. D'après la définition 16.1.1 on a

$$J_{spec}(\sigma^\vee,f)=\displaystyle\sum_{L\in\mathcal{L}(M_{min})} |W^L||W^G|^{-1}(-1)^{a_L} J_{spec,L}(\sigma^\vee,f)$$
où on a posé

$$J_{spec,L}(\sigma^\vee,f)=\displaystyle\sum_{\mathcal{O}\in\{\Pi_{ell}(L)\}; m(\mathcal{O},\sigma^\vee)=1} [i\mathcal{A}_\mathcal{O}^\vee:i\mathcal{A}_{L,F}^\vee]^{-1}t(\mathcal{O})^{-1} \int_{i\mathcal{A}_{L,F}^*} J^G_L(\pi_{\lambda},f) d\lambda$$

On reconnaît $J_{spec,G}(\sigma^\vee,f)$: c'est $m_{spec}(f,\sigma)$. Plus généralement, soit $L\in\mathcal{L}(M_{min})$. Reprenons les notations et les résultats d'Arthur rappelés dans la section 1.8. On pose $f_L=\phi_L(f)\mathbf{1}_{H_L=0}$. C'est une fonction cuspidale de $L(F)$. Puisqu'on a fixé les mesures de sorte que $i\tilde{\mathcal{A}}_{L,F}^*$ soit de mesure $1$, l'égalité 1.8(1) nous donne

$$J_{spec,L}(\sigma^\vee,f)=\displaystyle\sum_{\mathcal{O}\in\{\Pi_{ell}(L)\}; m(\mathcal{O},\sigma^\vee)=1} [i\mathcal{A}_\mathcal{O}^\vee:i\tilde{\mathcal{A}}_{L,F}^\vee]^{-1}t(\mathcal{O})^{-1}\theta_\pi(f_L)$$

Le Levi $L$ admet une décomposition $L\simeq R_{E/F} GL_{n_1}\times\ldots\times R_{E/F}GL_{n_k}\times \tilde{G}$ où $\tilde{G}$ est le groupe unitaire d'une sous-espace hermitien $\tilde{V}$ de $V$. L'espace des fonctions cuspidales sur $L(F)$ est produit tensoriel des espaces des fonctions cuspidales sur les $GL_{n_j}(E)$ et de l'espace de fonctions cuspidales sur $\tilde{G}(F)$. On définit sur cet espace la forme linéaire $m_{spec}(.,\sigma)$ comme le produit tensoriel des formes linéaires $m_{spec}$ sur les $GL_{n_j}(E)$ et de la forme linéaire $m_{spec}(.,\sigma)$ sur $\tilde{G}(F)$. On vérifie alors aisément l'égalité

$$J_{spec,L}(\sigma^\vee,f)=m_{spec}(f_L,\sigma)$$

Si $L\neq G$, on a aussi $m_{spec}(f_L,\sigma)=m_{geom}(I\theta_{f_L},\sigma)$. En effet, il suffit de vérifier l'égalité sur les tenseurs purs et c'est alors une conséquence du lemme 17.4.1 et de l'hypothèse de récurrence sur les espaces compatibles $\tilde{V}$ et $W$. On a donc l'égalité

$$\mbox{(1)} \;\;\; \displaystyle m_{spec}(f,\sigma)=m_{geom}(\theta_f,\sigma)-\sum_{L\in\mathcal{L}(M_{min}); L\neq G} |W^L||W^G|^{-1} (-1)^{a_L} m_{geom}(I\theta_{f_L},\sigma)$$

D'après le lemme 1.8.2 on a l'égalité

$$m_{geom}(\theta_f,\sigma)=\displaystyle\sum_{L\in\mathcal{L}(M_{min})} |W^L||W^G|^{-1}(-1)^{a_L}m_{geom}(Ind_L^G(I\theta_{\phi_L(f)}),\sigma)$$

D'après le lemme 17.2.1, pour tout $L\in \mathcal{L}(M_{min})$, on a 

$$m_{geom}(Ind_L^G(I\theta_{\phi_L(f)}),\sigma)=m_{geom}(I\theta_{\phi_L(f)},\sigma)$$

 Le quasicaractère $I\theta_{\phi_L(f)}$ est la somme des quasi caractères $I\theta_{\phi_L(f) \mathbf{1}_{H_L}=Z}$ pour $Z\in\mathcal{A}_{L,F}$. D'après une remarque faite au début du paragraphe 17.2, seul le terme pour $Z=0$ intervient de façon non nul dans $m_{geom}(I\theta_{\phi_L(f)},\sigma)$. On a donc $m_{geom}(I\theta_{\phi_L(f)},\sigma)=m_{geom}(I\theta_{f_L},\sigma)$. La plupart des termes de l'égalité (1) se simplifient et on obtient

$$m_{spec}(f,\sigma)=m_{geom}(I\theta_f,\sigma)$$

\section{Une application à la conjecture de Gross-Prasad}

\subsection{Hypothèses sur les $L$-paquets}

On reprend la situation de la section 4. On note $d$ la dimension de $V$ et $d_W$ celle de $W$. On suppose désormais $G$ et $H$ quasidéployés et on affecte les objets que l'on avait définis d'un indice $i$: $G_i$, $H_i$, $V_i$, $W_i$... Il existe à isomorphisme près un unique espace hermitien $V_a$ non isomorphe à $V_i$ et de même dimension que $V_i$. Si $d_W\geqslant 1$, il existe aussi un unique espace hermitien $W_a$ non isomorphe à $W_i$ et de même dimension que $W_i$. Si $d_W=0$, on pose simplement $W_a=0$. Notons $G_a$ et $H_a$ les groupes unitaires de $V_a$ et $W_a$ respectivement. Dans le cas où $d_W\geqslant 1$, les espaces hermitiens $V_a$ et $Z\oplus^\perp D\oplus^\perp W_a$ sont isomorphes et on fixe un tel isomorphisme. Pour tous $\pi\in Temp(G_a)$, $\sigma\in Temp(H_a)$, on dispose donc de la multiplicité $m(\pi,\sigma)$ (indépendante de l'isomorphisme choisi). Si $d_W=0$, il n'existe alors pas d'isomorphisme entre $V_a$ et $Z\oplus D\oplus W_a$ et pour tous $\pi\in Temp(G_a)$, $\sigma\in Temp(H_a)$ on pose $m(\pi,\sigma)=0$. On affectera d'un indice $a$ les objets définis à partir du couple $(V_a,W_a)$.\\

D'après des conjectures d'Arthur et Langlands il existe des décompositions de $Temp(G_i)$ et $Temp(G_a)$ en sous-ensembles finis disjoints appelés $L$-paquets vérifiant les conditions suivantes

\begin{enumerate}
\item Si $\Pi$ est un $L$-paquet de $Temp(G_i)$ ou $Temp(G_a)$ alors $\theta_\Pi=\displaystyle\sum_{\pi\in \Pi} \theta_\pi$ est une distribution stable.
\item Il existe une application injective des $L$-paquets de $Temp(G_a)$ dans les $L$-paquets de $Temp(G_i)$ qui vérifie les deux conditions suivantes:

\begin{enumerate}
\item Si le $L$-paquet $\Pi_i$ est l'image du $L$-paquet $\Pi_a$ alors $(-1)^{d+1}\theta_{\Pi_a}$ est le transfert à $G_a(F)$ de $\theta_{\Pi_i}$.
\item Si le $L$-paquet $\Pi_i$ n'est l'image d'aucun $L$-paquet de $Temp(G_a)$ alors le transfert de la distribution $\theta_{\Pi_i}$ à $G_a(F)$ est nul.
\end{enumerate}

\item Pour tout $L$-paquet $\Pi_i$ de $Temp(G_i)$ et toute orbite nilpotente régulière $\mathcal{O}$ de $\mathfrak{g}_i(F)$, il existe un et un seul élément de $\Pi$ qui admet un modèle de Whittaker relativement à $\mathcal{O}$.
\end{enumerate}

On suppose que ces conjectures sont vérifiées pour le couple $(G_i,G_a)$ ainsi que pour le couple $(H_i,H_a)$ si $d_W\geqslant 1$.

\subsection{Un calcul de fonction $\hat{j}$}

\begin{lem}
Soient $B$ un sous-groupe de Borel de $G_i$ défini sur $F$, $T_{qd}$ un tore maximal de $B$ défini sur $F$ et $X_{qd}\in\mathfrak{t}_{qd}(F)\cap\mathfrak{g}_{i,reg}(F)$. Posons $W^{G_i}=W(G_i,T)$. Pour tout $\mathcal{O}\in Nil_{reg}(\mathfrak{g}_i(F))$, on a
$$\hat{j}(\mathcal{O},X_{qd})=|Nil_{reg}(\mathfrak{g}_i(F))|^{-1} |W^{G_i}| D^{G_i}(X_{qd})^{-1/2}$$
\end{lem}

\ul{Preuve}: Pour tout $\lambda\in F^{\times, 2}$ on a $\hat{j}(\mathcal{O},\lambda X_{qd})=|\lambda|_{F}^{(dim(T)-dim(G_i))/2}\hat{j}(\mathcal{O},X_{qd})$ et $D^{G_i}(\lambda X_{qd})=|\lambda|_{F}^{(dim(T)-dim(G_i))/2}D^{G_i}(X_{qd})$. Il suffit donc d'établir le résultat pour $X_{qd}$ dans un voisinage de $0$. Pour $X_{qd}\in \mathfrak{t}_{qd}(F)\cap\mathfrak{g}_{i,reg}(F)$ assez proche de 0, on a

$$\hat{j}(X_{qd},X_{qd})=\displaystyle\sum_{\mathcal{O}\in Nil(\mathfrak{g}_i(F))} \Gamma_\mathcal{O}(X_{qd}) \hat{j}(\mathcal{O},X_{qd})$$

Appliquons l'égalité 1.4(1) à $M=T_{qd}$ et $X=Y=X_{qd}$. On obtient

$$\hat{j}(X_{qd},X_{qd})=|W^{G_i}| D^{G_i}(X_{qd})^{-1/2}$$

D'après le lemme 9.3.1, on a donc

$$\displaystyle\sum_{\mathcal{O}\in Nil(\mathfrak{g}_i(F))_{reg}}\hat{j}(\mathcal{O},X_{qd})=|W^{G_i}| D^{G_i}(X_{qd})^{-1/2}$$

Si $dim(V_i)$ est impaire, alors il n'y a qu'une seule orbite nilpotente régulière dans $\mathfrak{g}_i(F)$ et on a le résultat voulu. Si $dim(V_i)$ est paire alors il y a deux orbites nilpotentes régulières que l'on note $\mathcal{O}_+$ et $\mathcal{O}_-$. Soit $g$ un élément non central du groupe des similitudes de $V_i$ qui commute à $T_{qd}$. Alors la conjugaison par $g$ échange $\mathcal{O}_+$ et $\mathcal{O}_-$. On a alors pour tout $X_{qd}\in \mathfrak{t}_{qd}(F)\cap\mathfrak{g}_{i,reg}(F)$

$$\hat{j}(\mathcal{O}_-,X_{qd})=\hat{j}(g\mathcal{O}_+g^{-1},gX_{qd}g^{-1})=\hat{j}(\mathcal{O}_+,X_{qd})$$

$\blacksquare$

\vspace{3mm}

\subsection{Classes de conjugaison stable de tores}

Dans la suite, on notera $\underline{\mathcal{T}}_\flat=\underline{\mathcal{T}}(W_\flat,D)$ où $\flat=i$ ou $a$. Fixons un $E\otimes_F \overline{F}$-isomorphisme $\Phi: W_a\otimes_F \overline{F}\to W_i\otimes_F \overline{F}$ tel que $h_{W_i}(\Phi(w),\Phi(w'))=h_{W_a}(w,w')$ pour tous $w,w'\in W_a\otimes_F \overline{F}$. Pour $h\in H_a$, on notera $\phi(h)=\Phi\circ h\circ \Phi^{-1}\in H_i$. Soient $T,T'\in \underline{\mathcal{T}}_a\cup \underline{\mathcal{T}}_i$, on dit que $T$ et $T'$ sont stablement conjugués si et seulement si l'une des deux conditions suivantes est vérifiée:

\begin{enumerate}
\item $T,T'\in\underline{\mathcal{T}}_\flat$ pour un certain indice $\flat$ et il existe $h\in H_\flat$ tel que $hTh^{-1}=T'$ et l'isomorphisme $t\mapsto hth^{-1}$ de $T$ sur $T'$ est défini sur $F$.
\item Quitte à échanger $T$ et $T'$, on a $T\in\underline{\mathcal{T}}_a$ et $T'\in\underline{\mathcal{T}}_i$ et il existe $h\in H_i$ tel que $h\phi(T)h^{-1}=T'$ et l'isomorphisme $t\mapsto h\phi(t)h^{-1}$ de $T$ sur $T'$ est défini sur $F$.
\end{enumerate}

Soit $\flat=i$ ou $a$. Pour $T\in \underline{\mathcal{T}}_\flat$, on introduit le groupe de cohomologie $H^1(T)=H^1(F,T)$. On définit alors

$$h(T)=\left\{
    \begin{array}{ll}
        |H^1(T)|/2 & \mbox{si } T\neq \{1\}\\
        1 & \mbox{sinon.}
    \end{array}
\right.$$

Soit $\overline{W}(H_\flat,T)=Norm_{H_\flat}(T)/Z_{H_\flat}(T)$, c'est un groupe algébrique défini sur $F$. On note $\overline{W}_F(H_\flat,T)$ l'ensemble de ses points sur $F$.

\begin{lem}
(i) Soient $T,T'\in\underline{\mathcal{T}}_a\cup \underline{\mathcal{T}}_i$ tels que $T\sim_{stab} T'$. Alors $W''_T$ est isomorphe à $W''_{T'}$.\\

De plus, si on est dans le cas 1) alors on peut trouver $h$ qui vérifie 1) et dont la restriction à $W''_T$ soit un isomorphisme défini sur $F$ avec $W''_{T'}$, si on est dans le cas 2) alors on peut trouver $h$ qui vérifie 2) et tel que la restriction à $W''_T$ de $h\circ\Phi$ soit un isomorphisme défini sur $F$ avec $W''_{T'}$. \\

(ii) Soit $T\in\mathcal{T}_\flat$ où $\flat=i$ ou $a$. On a alors

$$\displaystyle\sum_{T'\in\mathcal{T}_\flat, T'\sim_{stab} T} |W(H_\flat,T')|^{-1}=h(T)|\overline{W}_F(H_\flat,T)|^{-1}$$
et les termes $h(T)$ et $|\overline{W}_F(H_\flat,T)|$ ne dépendent que de la classe de conjugaison stable de $T$. \\

(iii) Toutes les classes de conjugaison stable dans $\mathcal{T}_a\cup\mathcal{T}_i$ coupent $\mathcal{T}_a$ et $\mathcal{T}_i$ sauf la classe du tore $\{1\}\in\mathcal{T}_i$ qui ne coupe pas $\mathcal{T}_a$.
\end{lem}

\ul{Preuve}: \\
(i) Dans les deux cas les espaces hermitiens $W''_T$ et $W''_{T'}$ sont de même dimension et leurs groupes unitaires sont quasi-déployés sur $F$. Si $dim(W''_T)$ est pair cela suffit à prouver que $W''_T$ et $W''_{T'}$ sont isomorphes sur $F$. Si $dim(W''_T)$ est impair, on utilise le fait que $W''_T\oplus D$ et $W''_{T'}\oplus D$ ont aussi des groupes unitaires quasi-déployés sur $F$ donc sont isomorphes et d'après le théorème de Witt, c'est aussi le cas de $W''_T$ et $W''_{T'}$. \\
Si on est dans le cas (1), il suffit de changer la restriction de $h$ à $W''_T\otimes \overline{F}$, qui est un isomorphisme entre $W''_T\otimes \overline{F}$ et $W''_{T'}\otimes \overline{F}$, par un isomorphisme entre $W''_T$ et $W''_{T'}$ défini sur $F$: le nouvel élément $h$ vérifie alors (1) et le (i) du lemme. \\
Si on est dans le cas (2), on remplace la restriction de $h$ à $\Phi(W''_T\otimes \overline{F})$ par $h''\circ \Phi^{-1}_{|W''_T\otimes \overline{F}}$ où $h''$ est un isomorphisme défini sur $F$ de $W''_T$ sur $W''_{T'}$.

\vspace{2mm}

(ii) Soit $T\in \underline{\mathcal{T}}_\flat$, alors d'après (i) pour tout $T'\in\underline{\mathcal{T}}_\flat$ qui est stablement conjugué à $T$, on peut trouver $h\in H_\flat$ qui vérifie à la fois (1) et (i). Pour un tel $h$, on a alors $\sigma(h)^{-1}h\in T$ pour tout $\sigma\in \Gamma_F$. Réciproquement, pour $h\in H_\flat$ vérifiant $\sigma(h)^{-1}h\in T$ pour tout $\sigma\in \Gamma_F$, le tore $hTh^{-1}$ est défini sur $F$, appartient à $\underline{\mathcal{T}}_\flat$ et est stablement conjugué à $T$. Soit $\mathcal{H}_T=\{h\in H_\flat : \forall\sigma\in \Gamma_F, \sigma(h)^{-1}h\in T\}$, on a alors une surjection

$$Z:H_\flat(F)\backslash \mathcal{H}_T/T\to \{T'\in \mathcal{T}_\flat: T'\sim_{stab} T\}$$

\begin{center}
$H_\flat(F)hT\mapsto$ l'unique élément de $\mathcal{T}_\flat$ qui est conjugué à $hTh^{-1}$ par un élément de $H_\flat(F)$.
\end{center}

On vérifie facilement que l'application $H_\flat(F)\backslash \mathcal{H}_T/T \to H^1(F,T)$ qui à $H_\flat(F)hT$ associe la classe du 1-cocyle $\sigma\mapsto \sigma(h)^{-1}h$ induit une bijection entre $H_\flat(F)\backslash \mathcal{H}_T/T$ et $Ker(H^1(F,T)\to H^1(F,H_\flat))$ dont le cardinal vaut $h(T)$. \\
D'autre part, le cardinal de $H_\flat(F)\backslash \mathcal{H}_T/T$ est aussi la somme des cardinaux des fibres de l'application $Z$. Soit $T'\in\mathcal{T}_\flat$ vérifiant $T'\sim_{stab} T$ et $h\in \mathcal{H}_T$ tel que $hTh^{-1}=T'$. La fibre de $Z$ au dessus de $T'$ est alors l'image de l'application

$$\mathcal{H}_{T'}\cap Norm_{H_\flat}(T')\to H_\flat(F)\backslash \mathcal{H}_T/T$$

$$n\mapsto H_\flat(F) nh T$$

\noindent qui se quotiente en une bijection de la fibre avec $(Norm_{H_\flat(F)}(T')T')\backslash\mathcal{H}_{T'}\cap Norm_{H_\flat}(T')$. Puisque l'on a une bijection naturelle

$$(H''_{T'}(F) T')\backslash Norm_{H_\flat(F)}(T')T' \simeq Z_{H_\flat(F)}(T')\backslash Norm_{H_\flat(F)}(T')=W(H_\flat,T')$$

\noindent ,le nombre d'éléments dans la fibre est $|W(H_\flat,T')|^{-1}|H''_{T'}(F)T'\backslash\mathcal{H}_{T'}\cap Norm_{H_\flat}(T')|$. On a une application naturelle

$$H''_{T'}(F)T'\backslash\mathcal{H}_{T'}\cap Norm_{H_\flat}(T')\to \overline{W}_F(H_\flat,T')$$

On vérifie aisément qu'elle est injective. Montrons qu'elle est aussi surjective. Soit $w\in\overline{W}_F(H_\flat,T')$ et $n\in Norm_{H_\flat}(T')$ qui relève $w$, pour tout $\sigma\in \Gamma_F$ on a alors $\sigma(n)^{-1}n\in Z_{H_\flat}(T')=H''_{T'} T'$. Soit $n'$ la restriction de $n$ à $W'_{T'}\otimes \overline{F}$: c'est un élément de $H'_{T'}$ que l'on considère comme un élément de $H_\flat$ en le laissant agir comme l'identité sur $W''_{T'}\otimes \overline{F}$ et qui relève alors encore $w$. On a encore $\sigma(n')^{-1}n'\in T'$ pour tout $\sigma\in \Gamma_F$ donc $n'\in \mathcal{H}_{T'}\cap Norm_{H_\flat}(T')$. L'application est donc bijective et on a par conséquent

$$|H''_{T'}(F)T'\backslash\mathcal{H}_{T'}\cap Norm_{H_\flat}(T')|= |\overline{W}_F(H_\flat,T')|$$

D'où la formule suivante

$$\displaystyle\sum_{T'\in\mathcal{T}_\flat, T'\sim_{stab} T} |W(H_\flat,T')|^{-1} |\overline{W}_F(H_\flat,T')|=h(T)$$
Pour obtenir le résultat annoncé, il suffit donc de montrer que $|\overline{W}_F(H_\flat,T)|$ ne dépend que de la classe de conjugaison stable de $T$. Or si on est dans le cas (1), alors on a une bijection $\overline{W}_F(H_\flat,T)\to \overline{W}_F(H_\flat,T')$ donnée par $w\mapsto hwh^{-1}$. On a une bijection analogue dans le cas (2). Enfin, $h(T)$ ne dépend aussi que de la conjugaison stable de $T$ car si $T\sim_{stab} T'$ alors $T$ et $T'$ sont isomorphes sur $F$. Cela établit (ii).

\vspace{2mm}

(iii) Soit $T\in \mathcal{T}_a$ alors puisqu'un des deux groupes $G_a$ et $H_a$ n'est pas quasi-déployé, on a $W'_T\neq 0$ donc la classe de conjugaison stable de $\{1\}\in\mathcal{T}_i$ ne coupe pas $\mathcal{T}_a$. De plus, il existe un unique espace hermitien $W'_i$ de même dimension que $W'_T$ mais non isomorphe à $W'_T$ et $W_i$ est alors isomorphe à $W'_i\oplus W''_T$. On peut alors transférer $T$ en un tore maximal $T'$ du groupe unitaire de $W'_i$ et on a alors $T'\in \underline{\mathcal{T}}_i$ donc la classe de conjugaison stable de $T$ coupe $\mathcal{T}_i$. \\

De la même façon soit $T\in \mathcal{T}_i$ tel que $W'_T\neq 0$. Alors il existe un unique espace hermitien $W'_a$ de même dimension que $W'_T$ mais qui ne lui est pas isomorphe et $W_a$ est alors isomorphe à $W'_a\oplus W''_T$. On peut alors transférer $T$ en un tore maximal $T'$ du groupe unitaire de $W'_a$ et on a alors $T'\in \underline{\mathcal{T}}_a$ donc la classe de conjugaison stable de $T$ coupe $\mathcal{T}_a$. $\blacksquare$

\vspace{4mm}

\subsection{Un résultat dans le sens de la conjecture de Gross-Prasad}

\begin{theo}
Soit $\Pi_i$ un $L$-paquet de $Temp(G_i)$ et $\Sigma_i$ un $L$-paquet de $Temp(H_i)$. Si $\Pi_i$ est l'image d'un $L$-paquet de $Temp(G_a)$ alors on note $\Pi_a$ ce $L$-paquet, sinon on pose $\Pi_a=\emptyset$. On définit de même $\Sigma_a$. Il existe alors un unique couple $(\sigma,\pi)\in (\Pi_i\times \Sigma_i)\cup (\Pi_a\times \Sigma_a)$ tel que $m(\pi,\sigma^\vee)=1$.
\end{theo}

\ul{Preuve}: Le cas $d_W=0$ découle de notre hypothèse 18.1.3 sur les $L$-paquets de $Temp(G_i)$. On suppose dorénavant que $d_W\geqslant 1$. Soit $\flat=i$ ou $a$. On pose $m_{\Pi_\flat,\Sigma_\flat}=\displaystyle\sum_{\pi\in\Sigma_\flat,\sigma\in\Sigma_\flat} m(\pi,\sigma^\vee)$. Pour $T\in\mathcal{T}_\flat$ on définit les fonctions $c_{\Pi_\flat}$ et $c_{\Sigma_\flat}$ sur $T_\natural(F)$ par

$$c_{\Pi_\flat}(t)=\displaystyle\sum_{\pi\in\Pi_\flat} c_\pi(t)$$

\noindent et

$$c_{\Sigma_\flat}(t)=\displaystyle\sum_{\sigma\in\Sigma_\flat} c_\sigma(t)$$

\noindent D'après le théorème 17.1.1, on a

$$\mbox{(1)}\;\;\; m_{\Pi_\flat,\Sigma_\flat}= \displaystyle\sum_{T\in\mathcal{T}_\flat} |W(H_\flat,T)|^{-1} \nu(T) \lim\limits_{s\to 0^+}\int_{T(F)} c_{\Pi_\flat}(t) c_{\Sigma_\flat}(t) D^{H_\flat}(t)^{1/2}D^{G_\flat}(t)^{1/2}\Delta(t)^{s-1/2} dt$$

Soient $T\in\mathcal{T}_\flat$ et $t\in T_\natural(F)$. Soit $B\subset G''_{\flat,T}$ un sous-groupe de Borel défini sur $F$, $T_B\subset B$ un tore maximal défini sur $F$ et $X_{qd}\in \mathfrak{t}_B(F)\cap \mathfrak{g}''_{\flat,T,reg}(F)$. On fixe de la même façon un élément $X_{qd}^H\in \mathfrak{h}''_{\flat,T,reg}(F)$. Pour $\lambda\in F^{\times, 2}$ assez petit, on a

$$\theta_{\Pi_\flat}(texp(\lambda X_{qd}))=\displaystyle\sum_{\mathcal{O}\in Nil(\mathfrak{g}''_{\flat,T})} c_{\theta_{\Pi_\flat},\mathcal{O}}(t) |\lambda|_F^{-dim(\mathcal{O})/2} \hat{j}(\mathcal{O},X_{qd})$$

\noindent D'après le lemme 18.2.1, on a

$$c_{\Pi_\flat}(t)=\lim\limits_{\substack{\lambda\in F^{\times} \\ \lambda\to 0}} |W^{G''_{\flat,T}}|^{-1}|\lambda|_F^{\delta(G''_{\flat,T})/2} \theta_{\Pi_\flat}(texp(\lambda X_{qd}))$$

\noindent et de la même façon

$$c_{\Sigma_\flat}(t)=\lim\limits_{\substack{\lambda\in F^\times \\ \lambda\to 0}} |W^{H''_{\flat,T}}|^{-1}|\lambda|_F^{\delta(H''_{\flat,T})/2} \theta_{\Sigma_\flat}(texp(\lambda X_{qd}^H))$$

Soient $T'\in\mathcal{T}_\flat$ stablement conjugué à $T$, $h\in H_\flat$ qui vérifie 18.3.1(i) et posons $t'=hth^{-1}\in T_{\natural}'(F)$, $X'_{qd}=hX_{qd}h^{-1}$ et $X_{qd}'^{H}=hX_{qd}^H h^{-1}$. Alors $X'_{qd}$ et $X_{qd}'^{H}$ vérifient les mêmes hypothèses par rapport à $T'$ que $X_{qd}$ et $X_{qd}^H$ par rapport à $T$. On a donc

$$c_{\Pi_\flat}(t')=\lim\limits_{\substack{\lambda\in F^\times \\ \lambda\to 0}} |W^{G''_{\flat,T'}}|^{-1}|\lambda|_F^{\delta(G''_{\flat,T'})/2} \theta_{\Pi_\flat}(t'exp(\lambda X'_{qd}))$$

\noindent et

$$c_{\Sigma_\flat}(t')=\lim\limits_{\substack{\lambda\in F^\times \\ \lambda\to 0}} |W^{H''_{\flat,T'}}|^{-1}|\lambda|_F^{\delta(H''_{\flat,T'})/2} \theta_{\Sigma_\flat}(t'exp(\lambda X_{qd}'^{H}))$$

Puisque $\theta_{\Pi_\flat}$ et $\theta_{\Sigma_\flat}$ sont des distributions stables que $t'exp(\lambda X'_{qd})=htexp(X_{qd})h^{-1}$ et $|W^{G''_{\flat,T}}|=|W^{G''_{\flat,T'}}|$ (car $G''_{\flat,T}$ et $G''_{\flat,T'}$ sont isomorphes sur $F$), on a $c_{\Pi_\flat}(t)=c_{\Pi_\flat}(t')$ et de la même façon $c_{\Sigma_\flat}(t)=c_{\Sigma_\flat}(t')$. Comme de plus, $D^{H_\flat}(t')=D^{H_\flat}(t)$, $D^{G_\flat}(t')=D^{G_\flat}(t)$, $\Delta(t')=\Delta(t)$ et que l'isomorphisme $T(F)\to T'(F)$, $t\mapsto hth^{-1}$ envoie la mesure $\nu(T)dt$ sur la mesure $\nu(T')dt'$, les termes indexés par $T$ et $T'$ dans (1) sont les mêmes. Soit $\mathcal{T}_{\flat,stab}$ un ensemble de représentants des classes de conjugaisons stables dans $\mathcal{T}_\flat$. D'après 18.3.1(ii), la formule (1) peut se réécrire

$$m_{\Pi_\flat,\Sigma_\flat}=\displaystyle\sum_{T\in\mathcal{T}_{\flat,stab}} h(T) |\overline{W}_F(H_\flat,T)|^{-1} \lim\limits_{s\to 0^+} \int_{T(F)} c_{\Pi_\flat}(t) c_{\Sigma_\flat}(t) D^{H_\flat}(t)^{1/2}D^{G_\flat}(t)^{1/2} \Delta(t)^{s-1/2} dt $$

Soient $T'\in\mathcal{T}_i$ et $T\in\mathcal{T}_a$ stablement conjugués et $h\in H_i$ qui vérifie 18.3.1(i). Soit $t\in T_\natural(F)$ et posons $t'=h\phi(t)h^{-1}$, $X_{qd}'=h\phi(X_{qd})h^{-1}$ et $X_{qd}'^H=h\phi(X_{qd}')h^{-1}$. Puisque le transfert de $\theta_{\Pi_a}$ à $G_i$ est $(-1)^d\theta_{\Pi_i}$, on montre comme précédemment que $c_{\Pi_i}(t')=(-1)^dc_{\Pi_a}(t)$. De la même façon, on a $c_{\Sigma_i}(t')=(-1)^{d_W}c_{\Sigma_a}(t)$. On a aussi $D^{H_i}(t')=D^{H_a}(t)$, $D^{G_i}(t')=D^{G_a}(t)$ et $\Delta(t')=\Delta(t)$. L'isomorphisme $T(F)\to T'(F)$, $t\mapsto hth^{-1}$ envoie la mesure $\nu(T)dt$ sur la mesure $\nu(T') dt'$. Puisque $(-1)^{d+d_W}=-1$, on a par conséquent

\[\begin{aligned}
\displaystyle \int_{T(F)} c_{\Pi_a}(t) & c_{\Sigma_a}(t) D^{H_a}(t)^{1/2}D^{G_a}(t)^{1/2} \Delta(t)^{s-1/2} dt= \\
 & -\int_{T'(F)} c_{\Pi_i}(t) c_{\Sigma_i}(t) D^{H_i}(t)^{1/2}D^{G_i}(t)^{1/2} \Delta(t)^{s-1/2}  dt
\end{aligned}\]

\noindent pour tout $s\in \mathbb{C}$ vérifiant $Re(s)>0$. D'après 18.3.1(ii), on a aussi $ h(T) |\overline{W}_F(H_a,T)|^{-1}= h(T') |\overline{W}_F(H_i,T')|^{-1}$. Par conséquent, la contribution de la classe de conjugaison stable de $T$ dans $m_{\Pi_a,\Sigma_a}$ est l'opposé de la contribution de la classe de conjugaison stable de $T'$ dans $m_{\Pi_i,\Sigma_i}$. Puisque la seule classe de conjugaison stable qui ne coupe pas $\mathcal{T}_i$ et $\mathcal{T}_a$ est celle de $\{1\}$, seule cette classe contribue dans la somme $m_{\Pi_i,\Sigma_i}+m_{\Pi_a,\Sigma_a}$. On a donc

$$m_{\Pi_i,\Sigma_i}+m_{\Pi_a,\Sigma_a}=c_{\Pi_i}(1)c_{\Sigma_i}(1)$$

D'après un résultat de Rodier ([Ro]), pour une représentation $\pi\in Irr(G)$ et $\mathcal{O}\in Nil_{reg}(\mathfrak{g}(F))$, le coefficient $c_{\theta_\pi,\mathcal{O}}(1)$ vaut $1$ si et seulement si 
$\pi$ admet un modèle de Whittaker par rapport à $\mathcal{O}$ et vaut $0$ sinon. D'après les hypothèses qu'on a fait sur les $L$-paquets, on a donc

$$m_{\Pi_i,\Sigma_i}+m_{\Pi_a,\Sigma_a}=1$$

Il y a par conséquent un seul terme qui vaut $1$ dans la somme 

$$m_{\Pi_i,\Sigma_i}+m_{\Pi_a,\Sigma_a}=\displaystyle\sum_{(\sigma,\pi)\in (\Sigma_i\times\Pi_i)\cup(\Sigma_a\times\Pi_a)} m(\pi,\sigma^\vee)$$

$\blacksquare$

\bigskip

{\bf Bibliographie}

\bigskip

[AGRS] A. Aizenbud, D. Gourevitch, S. Rallis, G. Schiffmann: {\it Multiplicity one theorems}, Ann. of Math. (2) 172 (2010), no. 2, 1407-1434

[A1] J. Arthur: {\it The trace formula in invariant form}, Annals of Math. 114 (1981), p.1-74

[A2] -----------: {\it The invariant trace formula I. Local theory}, J. AMS 1 (1988), p. 323-383

[A3] -----------: {\it A local trace formula}, Publ. IHES 73 (1991), p.5-96

[A4] -----------: {\it Intertwining operators and residues I. Weighted characters}, J. Funct. Analysis 84 (1989), p. 19-84

[A5] -----------: {\it On elliptic tempered characters}, Acta Math. 171 (1993), p.73-138

[A6] -----------: {\it On a family of distributions obtained from Eisenstein series II: explicit formulas}, Amer. J. Math. 104 (1982), p.1289-1336

[BDKV] J.-N. Bernstein: {\it Le "`centre"' de Bernstein}, Ed. by P. Deligne. Travaux en Cours, Représentations des groupes réductifs sur un corps local, 1-32, Hermann, Paris, 1984

[BO] Y. Benoist, H. Oh: {\it Polar decomposition for p-adic symmetric spaces}, Int. Math. Res. Not. IMRN 2007, no. 24, Art. ID rnm121, 20 pp.

[BZ] I. Bernstein, A. Zelevinsky: {\it Induced representations of reductive p-adic groups I}, Ann. Sc. ENS 10 (1977), p.441-472

[CS] W. Casselman, J. Shalika: {\it The unramified principal series of p-adic groups II. The Whittaker function}, Comp. Math. 41, n°2 (1980), p.207-231

[DS] P. Delorme, V. Sécherre: {\it An analogue of the Cartan decomposition for p-adic reductive symmetric spaces of split p-adic reductive groups}, Pacific J. Math. 251 (2011), no. 1, 1-21

[GGP]  W. T. Gan, B. Gross, D. Prasad:{\it Symplectic local root numbers, central critical $L$-values and restriction problems in the representation theory of classical groups}, in Astérisque 346 (2012)
 
[GP] B. Gross, D. Prasad: {\it On irreducible representations of $SO_{2n+1}\times SO_{2m}$}, Can. J. Math. 46 (1994), p.930-950

[HCDeBS]: Harish-Chandra, S. DeBacker, P. Sally: {\it Admissible invariant distributions on reductive $p$-adic groups}, AMS Univ. lecture series 16 (1999)

[II] A. Ichino, T. Ikeda: {\it On the periods of automorphic forms on special orthogonal groups and the Gross-Prasad conjecture} , Geom. Funct. Anal. 19 (2010), no. 5, 1378-1425

[MVW] C. Moeglin, M.-F. Vignéras, J.-L. Waldspurger: {\it Correspondances de Howe sur un corps p-adique}, Springer LN 1291 (1987)

[Ro] F. Rodier: {\it Modèle de Whittaker et caractères de représentations}, in Non commutative harmonic analysis, J. Carmona, J. Dixmier, M. Vergne éd. Springer LN 466 (1981), p.151-171

[Sa] Y. Sakellaridis: {\it Spherical varieties and integral representations of L-functions}, prépublication 2009

[SV] Y. Sakellaridis, A. Venkatesh: {\it Periods and harmonic analysis on spherical varieties}, prépublication 2012

[W1] J.-L. Waldspurger: {\it Une formule intégrale reliée à la conjecture de Gross-Prasad}, Compos. Math. 146 (2010), no. 5, p. 1180-1290

[W2] -----------------------: {\it Une formule intégrale reliée à la conjecture locale de Gross-Prasad, $2^{eme}$ partie: extension aux représentations tempérées}, in Astérisque 346 (2012)

[W3] -----------------------: {\it La formule de Plancherel pour les groupes $p$-adiques, d'après Harish-Chandra}, J. of the Inst. of Math. Jussieu 2 (2003), p.235-333

[W4] -----------------------: {\it Intégrales orbitales nilpotentes et endoscopie pour les groupes classiques non ramifiés}, Astérisque 269 (2001)

\bigskip

Institut de mathématiques de Jussieu 2 place Jussieu 75005 Paris \\
 e-mail: rbeuzart@math.jussieu.fr
\end{document}